\documentstyle[12pt,leqno,twoside,personal,bezier,eufrak,rotating]{article} 
\renewcommand{\arraystretch}{1.5}

\begin{document}
\title{A two-box-shift morphism between Specht modules}
\author{Matthias K\"unzer}
\maketitle

\vspace{1cm}

\begin{small}
\begin{quote}
\begin{center}{\bf Resum\'e}\end{center}\vspace*{2mm}

Soit $n\geq 1$, soit $\lambda$ une partition de $n$, soit $\mu$ une partition obtenue \`a partir de $\lambda$ par d\'ecalage \`a gauche de deux cases situ\'es 
en bas d'une colonne. Nous donnons une formule pour un morphisme $\Z\Sl_n$-lin\'eaire d'ordre $m$ entre les modules de Specht sur $\Z/(m)$ correspondants, 
o\`u $m$ d\'esigne la longueur du d\'ecalage des cases (divis\'ee par deux en certains cas, specifi\'es combinatoirement).
En d'autres termes, nous obtenons une extension des modules de Specht sur $\Z$ correspondants, d'ordre $m$ en tant qu'\'el\'ement de $\Ext^1$. 
\end{quote}
\end{small}

\vspace{2cm}

\begin{small}
\begin{quote}
\begin{center}{\bf Abstract}\end{center}\vspace*{2mm}

Let $n\geq 1$, let $\lambda$ be a partition of $n$, let $\mu$ be a partition arising from $\lambda$ by a downwards shift of two boxes situated at the bottom of a column. We give a formula for
a $\Z\Sl_n$-linear morphism of order $m$ between the corresponding Specht modules over $\Z/(m)$, where $m$ is the box shift length (divided by two in certain combinatorially specified cases).
Reformulated, this yields an extension of the corresponding Specht modules over $\Z$ of order $m$ in $\Ext^1$. 
\end{quote}
\end{small}

\setcounter{section}{-1}

\section{Introduction}
\label{SecIntr}

\subsection{Known results}
\label{SubsecKno}

The irreducible modules of the complex group algebra of the symmetric group on $n$ letters possess the Specht lattices $S^\lambda$ as combinatorially described integral models, indexed by 
partitions $\lambda$ of $n$. Let $d\geq 1$ be an integer. Suppose given a partition $\lambda$ of $n$ and integers $1\leq g\leq k\leq\lambda_1 - 1$ such that the shift of $d$ 
boxes situated at the bottom of column $k+1$ of $\lambda$ downwards to column $g$ yields a partition again, which we denote by $\mu$. Let $\lambda'$ be the partition transposed to $\lambda$. The number 
\[
m_0 := \lambda'_g - \lambda'_{k+1} + (k+1) - g + d
\]
of steps by which we shift leftwards and downwards is called the {\it box shift length.} 

For instance, we may shift $d = 2$ boxes of $\lambda = (5,4,4,2,1,1)$ from column $k+1 = 4$ downwards to column $g = 2$ to obtain the partition $\mu = (5,3,3,2,2,2)$. We have
$\lambda' = (6,4,3,3,1)$, hence the box shift length equals $m_0 = 5$.

{\bf Theorem.} Based on [CL 74], {\sc Carter} and {\sc Payne} [CP 80] have shown that
\[
\Hom_{K\Sl_n}(K\ts_{\sZ} S^\lambda,K\ts_{\sZ} S^\mu) \neq 0,
\]
$K$ being an infinite field of characteristic $p$ such that $d < p^{v_p(m_0)}$, i.e.\ such that the number of shifted boxes is smaller than the $p$-part of the box shift length.
In case $d = 1$, this condition on $p$ translates to $p\; |\; m_0$. In case $d = 2$, it translates for $p\geq 3$ to $p \; |\; m_0$, and for $p = 2$ to $4 \;|\; m_0$. 

{\bf Theorem.} Suppose $d = 1$. In [K 99, 4.3.31, cf.\ 0.7.1], the existence of an element of order $m_0$ in
\[
\Hom_{(\sZ/(m_0))\Sl_n}(\Z/(m_0)\ts_{\sZ} S^\lambda, \Z/(m_0)\ts_{\sZ} S^\mu)
\]
has been shown via construction of such a morphism. 

\subsection{Results}
\label{SubsecRes}

We maintain the situation from Section \ref{SubsecKno}.

Suppose $d = 2$. Let $m$ be either $m_0$ or $m_0/2$, the latter value being taken in case $m_0$ is even, a similar combinatorial integer is even and the diagram of $\lambda$ takes a certain 
shape near the column $g$. 

{\bf Theorem} (\ref{ThC16}). There exists a morphism of order $m$ in
\[
\Hom_{(\sZ/(m))\Sl_n}(\Z/(m)\ts_{\sZ} S^\lambda, \Z/(m)\ts_{\sZ} S^\mu).
\]
As a consequence, we recover the according part of the result of {\sc Carter} and {\sc Payne} [CP 80]. An account of the construction of this morphism, including the definition of the modulus $m$, 
is given at the end of Section \ref{SecSpecht}. In Section \ref{SubsecCon}, we attempt to give an informal description of that construction. Some examples may be found in Section \ref{SubSecEx}. 

The proof by means of our techniques requires consideration of (essentially) $9^2$ cases - the base $9$ resulting from the structure of the one-step Garnir relations, the
exponent $2$ being the number of shifted boxes (cf.\ \ref{LemM2}). Moreover, the connection between the modulus and the box shift length results from the calculation only. 

A morphism modulo the box shift length $m_0$ itself need not exist. For instance, we obtain 
\[
\Hom_{\sZ\Sl_5}(S^{(2,2,1)},\Z/(5!)\ts_{\sZ} S^{(1,1,1,1,1)}) \llaiso \Hom_{\sZ\Sl_5}(S^{(2,2,1)},\Z/(2)\ts_{\sZ} S^{(1,1,1,1,1)}) \iso \Z/2.
\]

Let $d$ be arbitrary, but suppose $g = k$. This is to say, we shift $d$ boxes from column $g+1$ to column $g$.

{\bf Theorem} (\ref{ThLW6}). Let $\ilog_p d := \max\{ i\in\Z_{\geq 0}\; |\; p^i\leq d\}$. There exists a morphism of order 
$m := m_0\cdot\prodd{p \mbox{\rm\scr\ prime, } p | m_0} p^{-\min(v_p(m_0),\silog_p(d))}$ in 
\[
\Hom_{(\sZ/(m))\Sl_n}(\Z/(m)\ts_{\sZ} S^\lambda, \Z/(m)\ts_{\sZ} S^\mu).
\]
Again, the according part of the result of {\sc Carter} and {\sc Payne} [CP 80] is recovered. Roughly speaking, this morphism is given by picking $d$ entries from column $g + 1$ of a given 
$\lambda$-polytabloid, by placing them at the bottom of column $g$, and by mapping this $\lambda$-polytabloid to an alternating sum of the resulting $\mu$-polytabloids over the various 
possibilities to pick those entries.

{\bf Remark.} For an abelian group $A$ and an integer $m\geq 2$, we write $A[m] := \Kern(A\lraa{m} A)$. There is an isomorphism
\[
\Hom_{(\sZ/(m))\Sl_n}(\Z/(m)\ts_{\sZ} S^\lambda, \Z/(m)\ts_{\sZ} S^\mu) \lraiso \Ext^1_{\sZ\Sl_n}(S^\lambda, S^\mu)[m]
\]
which translates both results on the existence of a morphism of order $m$ into assertions on the existence of an extension of order $m$ in $\Ext^1$. 
Cf.\ Section \ref{SubsecMot}.

{\bf Remark}. Let $m\geq 2$ be an integer. For arbitrary partitions $\lambda$, $\mu$ of $n$, we dispose of the {\it transposition isomorphism} 
\[
\Hom_{(\sZ/(m))\Sl_n}(\Z/(m)\ts_{\sZ} S^\lambda, \Z/(m)\ts_{\sZ} S^\mu) \lraiso \Hom_{(\sZ/(m))\Sl_n}(\Z/(m)\ts_{\sZ} S^{\mu'}, \Z/(m)\ts_{\sZ} S^{\lambda'}),
\]
given by dualization, followed by alternation and isomorphic substitution (\ref{PropSST11}). Based on the simple assertion (\ref{PropSST2}), we give in (\ref{CorSST17})
an explicit, but nevertheless not quite satisfactory formula for the transpose of the two-box-shift morphism in (\ref{ThC16}) in terms of {\it tabloids.} Whereas the image of a polytabloid under
that transpose is known to be an element of the target Specht module, hence known to be expressible as a linear combination of {\it polytabloids,} a formula for such an expression is 
lacking. Some examples are given in Section \ref{SubsecHorEx}, cf.\ also (\ref{ExLW10}).

\subsection{Construction (informal description)}
\label{SubsecCon}

The following description might serve as an intuitive thread through our formalism concerning the two-box-shift morphism in (\ref{ThC16}). We maintain the situation of the beginning of
Section \ref{SubsecRes} and refer to the usual way to depict partitions and tableaux [J 78].

\begin{footnotesize}
We write the Specht lattice $S^\lambda$ as a quotient of the free $\Z\Sl_n$-module $F^\lambda$ on one generator, displayed as having as 
$\Z$-linear basis the set of $\lambda$-tableaux equipped with the action of the $\Sl_n$ on the tableau entries. The kernel of the canonical epimorphism $F^\lambda\lra S^\lambda$ 
is generated, over $\Z\Sl_n$, by signed column transpositions and one-step Garnir relations. Given a tableau, a signed column transposition is the sum of the tableau and the tableau having two 
entries of a column interchanged. A one-step Garnir relation depends on subsets of two subsequent columns of a tableau $[a]$ the sizes of which add up to (the length of the left column)$+1$. 
After factoring out the signed column transpositions, such a Garnir relation can be written, up to a scalar, as the alternating sum of the tableaux arising from $[a]$ by permutation inside the 
union of the chosen subsets of the columns, the sign being that of the afforded permutation. Our morphism is constructed as a factorization of a morphism $(F^\lambda \lra S^\mu \lra S^\mu/mS^\mu)$ 
over the canonical epimorphism $(F^\lambda\lra S^\lambda\lra S^\lambda/mS^\lambda)$. To a $\lambda$-tableau $[a]$, we may apply certain place operations, indexed by 
{\it double paths,} which produce $\mu$-polytabloids. A double path is a pair of sequences of positions in the union of the diagrams of $\lambda$ and $\mu$ that do not intersect, that run 
strictly from right to left, that start arbitrarily in column $k+1$ and that end at the the position of the shifted boxes in column $g$. (So at least $4$ and at most $4 + 2(k - g)$ positions 
are occupied.)
Given such a double path, the place operation it gives rise to is obtained by inserting entries of column $k+1$ into the double path from the right and by subsequently pushing entries through 
along each path separately. In particular, the last two entries forced to move are installed at the positions occupied by the shifted boxes in the diagram of $\mu$. Cf.\ (\ref{ExM-1}). Each 
double path has a weight attached, i.e.\ a tuple recording the number of positions it occupies in each column, and a sign, depending on the positions it  
occupies in column $k+1$. Sending $[a]$ to the accordingly signed sum of the resulting $\mu$-polytabloids indexed by double paths of a {\it fixed weight} yields a $\Z\Sl_n$-linear map 
from $F^\lambda$ to $S^\mu$ that annihilates the signed column transpositions. Forming a linear combination of these maps, indexed {\it over the weights,} equipped with certain coefficients 
which are polynomial in the combinatorial data, and dividing by an essentially polynomial factor of redundancy yields a $\Z\Sl_n$-linear map $F^\lambda\lraa{f''} S^\mu$ that annihilates
all one-step Garnir relations except for those involving column $g$ and column $g+1$, which vanish only modulo $mS^\mu$. Thus $(F^\lambda \lraa{f''} S^\mu \lra S^\mu/mS^\mu)$
factors over a $\Z\Sl_n$-linear map $S^\lambda/mS^\lambda\lraa{f} S^\mu/mS^\mu$. Since the
image of $F^\lambda\lraa{f''} S^\mu$ has not been contained in a strict multiple $t S^\mu \tme S^\mu$, $t\geq 2$, the resulting morphism $S^\lambda/mS^\lambda\lraa{f} S^\mu/mS^\mu$ 
is of order $m$ in its $\Hom$-group.
\end{footnotesize}

\subsection{Motivation}
\label{SubsecMot}

We consider the integral group ring $\Z\Sl_n$ as a subring of a product of integral matrix 
rings via {\sc Wedderburn}'s embedding, sending a group element $g\in\Sl_n$ to the tuple $(\rho^\lambda(g))_\lambda$ of its operating matrices on the Specht lattices with respect to a choice of 
integral bases. Let $m\geq 2$ be an integer.
A $\Z$-linear map $S^\lambda\lraa{f}S^\mu$ that becomes $\Z\Sl_n$-linear modulo $m\geq 2$, i.e.\ that satisfies $f\c\rho^\lambda(g) \con_m \rho^\mu(g)\c f$ for all $g\in\Sl_n$, imposes this  
congruence as a necessary condition on an element of our product of integral matrix rings to lie in the image of that embedding. More generally, unscrewing the regular lattice into simple 
lattices via a binary tree of short exact sequences, we obtain a necessary and sufficient system of modular morphisms by means of the correspondence
\[
\Hom_{\sZ\Sl_n}(X,Y/mY) \lraiso \Ext^1_{\sZ\Sl_n}(X,Y)[m],
\]
where $X$ and $Y$ are rationally disjoint $\Z\Sl_n$-lattices. This correspondence attaches to a morphism $X\lra Y/mY$ the extension obtained as pullback of 
\[
0\lra Y \lraa{m} Y \lra Y/mY \lra 0.
\]
Conversely, a retraction up to $m$ to the inclusion of an extension of $X$ by $Y$ yields the corresponding modular morphism as being induced on the
cokernels. Since the $\Ext^1$-order of a short exact sequence that occurs in that binary tree is not necessarily square-free, modular means modulo prime powers. See [K 99] for an 
elaboration on this theme.

\subsection{Acknowledgements}
\label{SubsecAck}

I'd like to thank Prof.\ {\sc Rentschler} for kind hospitality during the time of the preparation of this article. Large parts of the work have been supported by the 
EU TMR-network `Algebraic Lie Representations', grant no.\ ERB FMRX-CT97-0100.  I'd like to thank my fellow postdocs in Paris for help of various kinds. I'd like to thank Prof.\ 
{\sc Steffen K\"onig} for encouragement, advice and support during and after the time of the preparation of my thesis, upon the methods of which the present result relies.

\subsection{Leitfaden}
\label{SubsecLeit}

\begin{center}
\begin{picture}(1800,1170)
\put( 500,1115){\scr\ref{SecGarnir}\ \ul{Garnir}}
\put( 500,1015){\scr\ref{SecSpecht}\ Specht lattices}
\put( 490, 990){\line(-1,-1){130}}
\put( 540, 990){\line(0,-1){130}}
\put( 760, 990){\line(1,-1){130}}
\put(   0, 915){\scr\ref{SecTwoBox}\ \ul{A two-box-shift morphism}}
\put( 200, 815){\scr\ref{SubSecMS}\ Double paths}
\put( 290, 790){\line(-1,-1){130}}
\put( 500, 815){\scr\ref{SecGar}\ A Garnir formula}
\put( 490, 790){\line(-1,-1){130}}
\put( 540, 790){\line(0,-1){730}}
\put( 900, 915){\scr\ref{SecSemi}\ \ul{Semistandard morphisms}}
\put( 900, 815){\scr\ref{SubsecCorr}\ Correspondoids}
\put( 930, 790){\line(0,-1){530}}
\put( 960, 790){\line(1,-1){130}}
\put(   0, 755){\scr\ref{AppTwoLem}.\ Two lemmata}
\put( 140, 730){\line(0,-1){70}}
\put(   0, 615){\scr\ref{SubSecLinComb}\ Morphism, unreduced version}
\put( 140, 590){\line(0,-1){130}}
\put( 970, 615){\scr\ref{SubsecRemTr2}\ Remarks on transposition in characteristic $2$}
\put(1300, 590){\line(0,-1){530}}
\put( 900, 215){\scr\ref{SubsecComp}\ Paths as correspondoids}
\put( 890, 190){\line(-1,-1){130}}
\put( 960, 190){\line(1,-1){130}}
\put(   0, 415){\scr\ref{SubSecRedund}\ Morphism, reduced version}
\put( 140, 390){\line(0,-1){130}}
\bezier{300}(360,390)(440,370)(520,350)
\bezier{300}(560,340)(720,300)(880,260)
\put( 500,  15){\scr\ref{SecTCSB}\ \ul{Two columns, several boxes}}
\put(1050,  15){\scr\ref{SubsecHorEx}\ Horizontal examples}
\put(   0, 215){\scr\ref{SubSecEx}\ Vertical examples}
\bezier{300}(250,190)(385,167)(520,145)
\bezier{300}(560,138)(665,121)(770,103)
\bezier{300}(830,93)(930,77)(1030,60)
\end{picture}
\end{center}
\section{Garnir}
\label{SecGarnir}

\subsection{Specht lattices}
\label{SecSpecht}

For $a,b\in\Z$, we let $[a,b] := \{i\in\Z\; |\; a\leq i\mbox{ and } i\leq b\}$.

Let $n\geq 1$. The symmetric group on the set $[1,n]$ is denoted by $\Sl_n := \mbox{Aut}_{\mbox{\scr Sets}} [1,n]$. 
The sign of a permutation $\sigma\in\Sl_n$ is denoted by $\eps_\sigma$.
Maps are written in various manners, on the right, on the left, using indices etc. The symmetric group acts on the right. Occasionally, we shall allow ourselves to 
treat tuples consisting of pairwise different entries and single elements as sets. Intervals are to be read as subsets of $\Z$. Conjugation in a group is also written as $h^g := g^{-1}hg$.
Given a module $X$ over a commutative ring $R$ and an element $r\in R$, we shall usually write $X/r := X/rX$ (cf.\ e.g.\ \ref{CorSST8}). If $X$ and $Y$ are modules over an $R$-algebra
$A$, we identify $\Hom_A(X/r,Y/r) = \Hom_A(X,Y/r)$ via composition with the residue class map $X\lra X/r$.

Let 
\[
\begin{array}{rcl}
\N & \lraa{\lambda} & \N_0 \\
i  & \lra           & \lambda_i \\
\end{array}
\] 
be a {\it partition of $n$,} i.e.\ assume $\sum_i \lambda_i = n$ and $\lambda_i \geq \lambda_{i+1}$ for $i\in \N$. Usually, we write a partition as a tuple, i.e.\ 
$\lambda = (\lambda_1,\lambda_2,\dots)$. Sometimes, we abbreviate repeated entries by exponents that indicate their multiplicity, such as e.g.\ $(9,2,2,2,1,1) = (9,2^3,1^2) \neq (9,8,1)$.
Let $[\lambda] := \{ i\ti j \in\N\ti\N \; |\; \lambda_i\geq j \}$ be the {\it diagram} of $\lambda$. A {\it $\lambda$-tableau} is a bijection 
\[
\begin{array}{rcl}
[\lambda] & \lraisoa{[a]} & [1,n] \\
i\ti j    & \lra          & a_{j,i} \\
\end{array}.
\]
The {\it transposed partition} of $\lambda$ is denoted by $\lambda'$, i.e.\ $i\ti j \in [\lambda] \equ j\ti i\in [\lambda']$. For a $\lambda$-tableau $[a]$, we denote by $[a']$
the $\lambda'$-tableau obtained by composition with $i\ti j\lra j\ti i$, mapping $[\lambda']$ onto $[\lambda]$. 
A permutation $\sigma\in\Sl_n$ acts on the set of $\lambda$-tableaux $T^\lambda$ via composition $[a] \lraa{\sigma} [a]\sigma$. Let $F^\lambda$ be the free $\Z$-module on $T^\lambda$, endowed 
with the induced action of the $\Sl_n$. Let 
\[
\begin{array}{rclcrcl}
[\lambda]  & \lraa{\pi_R^\lambda} & \N & , \hspace*{2cm} & [\lambda] & \lraa{\pi_C^\lambda} & \N \\
i\ti j     & \lra                 & i  &                 & i\ti j    & \lra                 & j  \\
\end{array}
\]
denote the projections. To a $\lambda$-tableau $[a]$ we attach a {\it $\lambda$-tabloid}
\[
\{ a\} := \{ [a]\}:= [a]^{-1}\pi^\lambda_R \in \N^{[1,n]}.
\]
In general, $X$ being a set, an element $\sigma\in\Sl_n$ acts on $f\in X^{[1,n]}$ via $(i)(f\sigma) := (i\sigma^{-1})f$. In particular, the free $\Z$-module on the set of 
$\lambda$-tabloids, denoted by $M^\lambda$, carries a structure as a $\Z\Sl_n$-lattice, viz.\ $\{ a\}\sigma = \{ [a]\sigma\}$. Let 
\[
\begin{array}{rcl}
C_{[a]} & := & \{ \kappa\in\Sl_n \; |\; [a]^{-1}\pi^\lambda_C = ([a]\kappa)^{-1}\pi^\lambda_C \} \\
R_{[a]} & := & \{ \rho\in\Sl_n \; |\; [a]^{-1}\pi^\lambda_R = ([a]\rho)^{-1}\pi^\lambda_R \} \\
\end{array}
\]
be the {\it column stabilizer} and the {\it row stabilizer} of $[a]$, respectively. Note that $C_{[a]\sigma} = (C_{[a]})^\sigma$ and $R_{[a]\sigma} = (R_{[a]})^\sigma$. Let the {\it Specht 
lattice} $S^\lambda$ be the $\Z\Sl_n$-sublattice of $M^\lambda$ generated $\Z$-linearly by the {\it $\lambda$-polytabloids}
\[
\spi{a} := \spi{[a]} := \sum_{\kappa\in C_{[a]}} \{ a\}\cdot\kappa\eps_\kappa \in M^\lambda.
\]
Let $\zeta\tm [1,n]$. We denote by $\Sl_\zeta := C_{\Sl_n}([1,n]\ohne\zeta)$ the subgroup of $\Sl_n$ consisting of permutations that act merely on $\zeta$. For an element $x$ of a $\Z\Sl_n$-module
$X$ we denote 
\[
x\c\zeta := \sum_{\sigma\in S_\zeta} x\cdot\sigma\eps_\sigma.
\] 
Given $p\geq 1$ and $1\leq j\leq k\leq\lambda'_p$, we write $a_{p,[j,k]} := (a_{p,i})_{i\in [j,k]}$ and abbreviate $a_p := a_{p,[1,\lambda'_p]}$.
Let $\xi\tm a_p$, $\eta\tm a_{p+1}$ be such that $\#\xi + \#\eta = \lambda'_p + 1$. Letting $\Sl_\xi\ti\Sl_\eta\<\Sl_{\xi\cup\eta}$ denote a chosen set of representatives of 
right cosets, the expression
\[
G''_{[a],\xi,\eta} := \sumd{\sigma\in \Sl_\xi\ti\Sl_\eta\<\Sl_{\xi\cup\eta}} [a] \sigma\eps_\sigma
\]
is called a {\it one-step Garnir relation}. For $p\geq 1$, $u,v\in a_p$, $u\neq v$, the expression 
\[
[a] + [a](u,v)
\]
is called a {\it signed column transposition.} From [J 78, 7.2] and from the proof of [J 78, 8.4] we take that the (finite) set
\[
\begin{array}{rl}
     & \{ [a] + [a](u,v) \;|\; [a]\in T^\lambda,\; p\geq 1,\; u,v\in a_p,\; u\neq v \} \\
\cup & \{ G''_{[a],\xi,\eta}\; |\; [a] \in T^\lambda,\; p\geq 1,\; s,t\geq 1,\; s + t = \lambda'_p + 1, \; \xi = a_{p,[\lambda'_p - s + 1,\lambda'_p]},\; \eta = a_{p+1,[1,t]} \} \\
\end{array}
\]
$\Z\Sl_n$-linearly generates the kernel of the epimorphism
\[
\begin{array}{rcl}
F^\lambda & \lraa{\nu^\lambda} & S^\lambda \\
{[a]}     & \lra               & \spi{a}.
\end{array}
\]
Given a $\Z\Sl_n$-module $X$, its {\it alternated module} is written $X^- := X\ts_{\sZ}\Z^-$, $\Z^-$ being the abelian group $\Z$, equipped with $1\cdot\sigma := \eps_\sigma$, $\sigma\in\Sl_n$. So
$\Z^- \iso S^{(1^n)}$. For each partition $\lambda$ of $n$, we fix a $\lambda$-tableau $[a_\lambda]$ and let
\[
\{[a_\lambda']\sigma\}^- := (\{a_\lambda'\} \ts 1)\sigma \in M^{\lambda',-}
\]
for $\sigma\in\Sl_n$. Note that $\{ a'\}^-\sigma = \{[a']\sigma\}^-$ and that $\{ a'\}^-\kappa = \{ a'\}\eps_\kappa$ for $[a]\in T^\lambda$, $\sigma\in\Sl_n$, $\kappa\in C_{[a]}$. We obtain a 
factorization into epimorphisms
\[ 
\left(
\begin{array}{rcl}
F^\lambda & \lraa{\nu^\lambda} & S^\lambda \\
{[a]}     & \lra               & \spi{a} \\
\end{array}
\right)
= 
\left(
\begin{array}{rclcl}
F^\lambda & \lraa{\nu^\lambda_M} & M^{\lambda',-} & \lraa{\nu^\lambda_S} & S^\lambda \\
{[a]}     & \lra                 & \{ a'\}^-      & \lra                 & \spi{a}   \\
\end{array}
\right)
\]
by first factoring out the signed column transpositions. Using this factorization, we may write
\[
G'_{[a],\xi,\eta} := (G''_{[a],\xi,\eta})\nu^\lambda_M = \fracd{1}{\#\xi!\#\eta!}\cdot \{ a'\}^-\c (\xi\cup\eta).
\]
A $\lambda$-tableau $[a]$ is called {\it standard} if $a_{j,i}\leq a_{j',i'}$ for $i\ti j,\; i'\ti j'\in [\lambda]$ such that $i\leq i'$ and $j\leq j'$. The
tuple $(\spi{a} \; |\; [a]\mbox{ is a standard $\lambda$-tableau})$ is a $\Z$-linear basis of $S^\lambda$ [J, 8.4]. Its elements are called {\it standard $\lambda$-polytabloids}. Let 
$[\ck a_\lambda]$ be the standard $\lambda$-tableau determined by $\ck a_{\lambda,j,i+1} = \ck a_{\lambda,j,i} + 1$ whenever $i\ti j,\; (i+1)\ti j\in [\lambda]$. (We do not require 
a priori that $[a_\lambda] = [\ck a_\lambda]$, a requirement which is convenient for practical purposes, however.)

\begin{footnotesize}
{\bf Account of the result (\ref{ThC16}).} The notation deviates from our working notation further down. In particular, the following definitions are valid {\it only}
in the remainder of this subsection (and will be repeated further down in case they coincide nonetheless).

Let $1\leq g < k\leq \lambda_1 - 1$ such that
\[
\mu'_i := 
\left\{
\begin{array}{ll}
\lambda'_i + 2 & \mbox{for } i = g \\
\lambda'_i - 2 & \mbox{for } i = k+1 \\
\lambda'_i     & \mbox{else} \\
\end{array}
\right.
\] 
defines a partition $\mu$. (For the case of $g = k$, see (\ref{PropM14})).

A {\it double path} $\gamma$ is a pair of integers $l(\delta)\geq 1$ together with a pair of maps
\[
\begin{array}{rcl}
\gamma_\delta \; : \; {[0,l(\delta)]} & \lra & [\lambda]\cup [\mu]  \\
                      i               & \lra & \alpha(\delta,i)\ti\beta(\delta,i), \\
\end{array}
\]
where $\delta\in [1,2]$, subject to the following conditions.
\begin{itemize}
\item[(i)] $i < i'$ implies $\beta(\delta,i) < \beta(\delta,i')$ for $\delta\in [1,2]$ and $i,i'\in [0,l(\delta)]$.
\item[(ii)] $\gamma_1([0,l(1)])\cap \gamma_2([0,l(2)]) = \leer$.
\item[(iii)] $\alpha(\delta,0)\ti\beta(\delta,0) = (\lambda'_g + \delta)\ti g$ for $\delta\in [1,2]$.
\item[(iv)] $\beta(\delta,l(\delta)) = k+1$ for $\delta\in [1,2]$.
\item[(v)] $\alpha(1,l(1)) < \alpha(2,l(2))$.
\end{itemize}
The set of double paths is denoted by $\Gamma$. Suppose given $\gamma\in\Gamma$. Let $\eps_\gamma := (-1)^{\alpha(1,l(1)) + \alpha(2,l(2))}$. For $j\in [g+1,k]$, we let 
\[
e(\gamma,j) := \sumd{\delta\in [1,2]}\#\gamma_\delta^{-1}([1,\lambda'_j]\ti \{ j\}).
\]
Let
\[
\begin{array}{rcl}
[0,\lambda'_{k+1}]       & \lra        & [0,\mu'_{k+1}] \\
i                        & \lraa{\phi} & \# \Big( [1,i]\ohne\{\alpha(1,l(1)),\alpha(2,l(2))\}\Big) \\
\min(\phi^{-1}(\{ j\}))  & \llaa{\psi} & j. \\
\end{array}
\]
Given a $\lambda$-tableau $[a]$, we let the $\mu$-tableau $[a^\gamma]$ be defined by 
\[
\begin{array}{rcll}
a^\gamma_{j,i}                              & := & a_{j,i}                                  & \mbox{ for } i\ti j\in [\mu]\ohne(\gamma_1([0,l(1)-1])\cup \gamma_2([0,l(2)-1])\cup\N\ti\{ k+1\}) \\
a^\gamma_{\beta(\delta,i),\alpha(\delta,i)} & := & a_{\beta(\delta,i+1),\alpha(\delta,i+1)} & \mbox{ for } \delta\in [1,2],\; i\in [0,l(\delta) - 1] \\
a^\gamma_{k+1,i}                            & := & a_{k+1,\psi(i)}                          & \mbox{ for } i\in [1,\mu'_{k+1}]. \\
\end{array}
\]
For $j\in [g,k+1]$, we let
\[
X_j := (\lambda'_j - j) - (\lambda'_{k+1} - (k+1)).
\]
For $i\in [0,1]$, we denote
\[
L(i) := \{ j\in [g+1,k-1] \; |\; \lambda'_{j+1} = \lambda'_j - i \}.
\]
We write $[g+1,k] = \Cup_{\kappa\in [1,K]} [p(\kappa),q(\kappa)]$ such that $p(1) = g+1$, such that $p(\kappa)\leq q(\kappa)$, $[p(\kappa),q(\kappa) - 1]\tm L(0)\cup L(1)$ and 
$q(\kappa)\not\in L(0)\cup L(1)$ for $\kappa\in [1,K]$, such that $q(\kappa) + 1 = p(\kappa + 1)$ for $\kappa\in [1,K-1]$, and such that $q(K) = k$. For $Z\in\Z$ and $i\geq 0$, we write 
$Z^{(i)} := (Z+i-1)!/(Z-1)!$. We let
\[
R := \Big(\prodd{i\in [0,1]}\;\;\prodd{j\in L(i)} X_j^{(2-i)}\Big)\Big(\prodd{\auf{\scm\kappa\in [1,K],}{[p(\kappa),q(\kappa)-1]\tm L(1)}} \gcd(2,X_{q(\kappa)})\Big)
\]
and
\[
m := \left\{\begin{array}{ll}
     X_g + 2                       & \mbox{ in case } [p(1),q(1)-1]\not\tm L(1) \\
     (X_g + 2)/\gcd(2,X_g,X_{g+1}) & \mbox{ in case } [p(1),q(1)-1]\tm L(1) \\
     \end{array}\right. .
\]
The $\Z\Sl_n$-linear map
\[
\begin{array}{rcl}
F^\lambda & \lraa{f''} & S^\mu \\
{[a]}     & \lra       & \sumd{\gamma\in\Gamma} \Big(\prodd{j\in [g+1,k]} X_j^{(2 - e(\gamma,j))}\Big) \spi{a^\gamma}\eps_\gamma 
\end{array}
\]
is divisible by $R$, and its quotient by $R$ factors $\Z\Sl_n$-linearly as 
\[
\left(
\begin{array}{rclcl}
F^\lambda & \lraa{f''/R} & S^\mu   & \lra & S^\mu/m \\
          &              & \spi{b} & \lra & \spi{b} + m S^\mu \\ 
\end{array}
\right)
= 
\left(
\begin{array}{rclcl}
F^\lambda & \lraa{\nu^\lambda} & S^\lambda & \lraa{f} & S^\mu/m \\
{[a]}     & \lra               & \spi{a}   &          &         \\
\end{array}
\right).
\]
The resulting morphism $S^\lambda \lraa{f} S^\mu/m$ is of order $m$ as an element of $\Hom_{\ssZ\Sl_n}(S^\lambda,S^\mu/m)$.
\end{footnotesize}          
\subsection{A Garnir formula}
\label{SecGar}

Let $\lambda$ be a partition of $n$, let $[a]$ be a $\lambda$-tableau, let $1\leq p < q\leq \lambda_1$, let $d\geq 1$.
Let $\xi\tm a_p$, $s := \#\xi$, $\b\xi := a_p\ohne\xi$ and $\eta\tm a_q$, $t := \#\eta$, $\b\eta := a_q\ohne\eta$, $u := \#\b\eta$, be such that $s + t = \lambda'_p + 1 - d$. We choose a 
disjoint decomposition $\b\xi = \b\xi_0 \cup \b\xi_1$ such that $\#\b\xi_0 = d - 1$ and $\#\b\xi_1 = t$.

\begin{Lemma}[{cf.\ [K 99, 4.3.4]}]
\label{LemG1}
Assume given $x\in\phi\tm a_p$, $y\in\psi\tm a_q$. Then
\[
\begin{array}{rcl}
\spi{a}\c (\phi\cup\psi)
& = & \#\phi\cdot\spi{a}\c (\phi\cup\psi\ohne x) - \#\psi\cdot\spi{a}\cdot (x,y)\c (\phi\cup\psi\ohne x) \\
& = & \#\psi\cdot\spi{a}\c (\phi\cup\psi\ohne y) - \#\phi\cdot\spi{a}\cdot (x,y)\c (\phi\cup\psi\ohne y). \\
\end{array}
\] 

\rm
In fact,
\[
\begin{array}{rcl}
\spi{a}\c (\phi\cup\psi)
& = & \sumd{z\in\phi\cup\psi\;\;}\sumd{\sigma\in\Sl_{\phi\cup\psi},\; z\sigma = x} \spi{a}\cdot\sigma\eps_\sigma \\
& = &  \Big(\sumd{z\in\phi\ohne x\;\;}\sumd{\sigma\in\Sl_{\phi\cup\psi},\; z\sigma = x} \spi{a}\cdot\sigma\eps_\sigma\Big) \\
& + &  \Big(\spi{a}\c(\phi\cup\psi\ohne x)\Big) \\
& + &  \Big(\sumd{z\in\psi\ohne y\;\;}\sumd{\sigma\in\Sl_{\phi\cup\psi},\; z\sigma = x} \spi{a}\cdot\sigma\eps_\sigma\Big) \\
& + &  \Big(\sumd{\sigma\in\Sl_{\phi\cup\psi},\; y\sigma = x} \spi{a}\cdot\sigma\eps_\sigma\Big)  \\
& = &  \Big((\#\phi - 1)\cdot\spi{a}\c (\phi\cup\psi\ohne x)\Big) \\
& + &  \Big(\spi{a}\c(\phi\cup\psi\ohne x)\Big) \\
& - &  \Big((\#\psi - 1)\cdot\spi{a}\cdot(x,y)\c (\phi\cup\psi\ohne x)\Big) \\
& - &  \Big(\spi{a}\cdot(x,y)\c (\phi\cup\psi\ohne x)\Big). \\
\end{array}
\]
The second equality follows by symmetry.
\end{Lemma}

Suppose given disjoint subsets $\phi,\psi \tm [1,n]$ and a bijection $\phi\lraisoa{\alpha}\psi$. We denote
\[
(\phi,\psi)_\alpha := \prodd{x\in\phi} (x,(x)\alpha) \in\Sl_n.
\]
In case $\phi\tm a_p$, $\psi\tm a_q$, $\#\phi = \#\psi$, the element $\spi{a}\cdot(\phi,\psi) := \spi{a}\cdot(\phi,\psi)_\alpha$ is independent of the choice of a bijection 
$\phi\lraisoa{\alpha}\psi$.

\begin{Lemma}[{[K 99, 4.3.5]}]
\label{LemG3}
In case $d = 1$ we have
\[
\spi{a}\c (\xi\cup\eta) = s!\, t!\cdot\spi{a}\cdot(\bar\xi,\eta).
\]

\rm
We calculate
\[
\begin{array}{rcl}
\spi{a}\c (\xi\cup\eta) 
& \auf{\mbox{\scr (\ref{LemG1})}}{=} & \frac{t}{s+1}\cdot\spi{a}\cdot(x,y)\c (\xi\cup\eta) \\
& \auf{\mbox{\scr (\ref{LemG1})}}{=} & \frac{t(t-1)}{(s+1)(s+2)}\cdot\spi{a}\cdot(x,y)\cdot (x',y')\c (\xi\cup\eta) \\
& \auf{\mbox{\scr (\ref{LemG1})}}{=} & \cdots \\
\end{array}
\]
\[
\begin{array}{rcl}
& \auf{\mbox{\scr (\ref{LemG1})}}{=} & \frac{s!\, t!}{(s+t)!}\cdot\spi{a}\cdot(\b\xi,\eta)\c (\xi\cup\eta) \\
& = & s!\, t!\cdot\spi{a}\cdot(\b\xi,\eta),\\
\end{array}
\]
where $\b\xi = \{ x,x',\dots\}$, $\eta = \{ y,y',\dots\}$.
\end{Lemma}

\begin{Lemma}
\label{LemG2}
Let $\phi\tm a_p$ and $\psi\tm a_q$ be given such that $\#\phi = \#\psi$ and let $\sigma\in\Sl_{a_p}$. We obtain
\[
\spi{a}\cdot(\phi,\psi)\cdot\sigma = \spi{a}\cdot\Big((\phi)\sigma,\psi\Big)\eps_\sigma.
\]

\rm
Factorization reduces us to the consideration of a transposition $\sigma = (x,y)$, where, moreover, we may assume $x\in\phi$ and $y\in a_p\ohne\phi$. Choosing $z\in\psi$, we conclude
\[
\begin{array}{rcl}
\spi{a}\cdot (\phi,\psi)\cdot (x,y) 
& = & \spi{a}\Big((\phi\ohne x),(\psi\ohne z)\Big)(x,z)(x,y) \\
& = & - \spi{a}\Big((\phi\ohne x),(\psi\ohne z)\Big)(y,z) \\
& = & - \spi{a}\Big((\phi)(x,y),\psi\Big).
\end{array}
\]
\end{Lemma}

\begin{Lemma}
\label{LemG4}
The formula
\[
\spi{a}\c (\xi\cup\eta) = \frac{s!}{(d-1)!}\cdot\spi{a}\cdot (\b\xi_1,\eta)\c \b\xi
\]
holds.

\rm
The claimed equality may be reformulated to
\[
\spi{a}\c (\xi\cup\eta)\c \b\xi = \frac{s!(t+d-1)!}{(s+t)!(d-1)!} \cdot\spi{a}\cdot (\b\xi_1,\eta)\c (\xi\cup\eta)\c\b\xi.
\]
We perform an induction over $d$, the case $d = 1$ being treated in (\ref{LemG3}). We assume $d\geq 2$ and perform an induction over $t$, the case $t = 0$ being trivial. We assume $t\geq 1$,
choose $x_0\in\b\xi_0$, $x_1\in\b\xi_1$, $y\in\eta$ and denote $\b\xi_2 := (\b\xi_1)(x_0,x_1)$. Evaluating in two ways, we obtain
\[
\begin{array}{cl}
  & \spi{a}\c (\xi\cup x_1\cup\eta)\c (\b\xi\ohne x_1) \\
\auf{\mbox{\scr 1., case $d-1$, $t$}}{=} & \frac{(s+1)!(t+d-2)!}{(s+t+1)!(d-2)!}\cdot\spi{a}\cdot (\b\xi_2,\eta)\c (\xi\cup x_1\cup\eta)\c (\b\xi\ohne x_1) \vspace*{1mm}\\
= & \frac{(s+1)!(t+d-2)!}{(s+t)!(d-2)!} \cdot\spi{a}\cdot (\b\xi_2,\eta)\c (\xi\cup\eta)\c (\b\xi\ohne x_1) \\
\auf{\mbox{\scr 2., (\ref{LemG1})}}{=} & (s+1)\cdot\spi{a}\c (\xi\cup\eta)\c (\b\xi\ohne x_1) \\
- & t \Big[\spi{a}\cdot (x_1,y)\c \Big( (\xi\cup y)\cup (\eta\ohne y)\Big)\c (\b\xi\ohne x_1)\Big] \\  
\auf{\mbox{\scr case $d$, $t-1$}}{=} & \frac{s+1}{t+d-1}\cdot\spi{a}\c (\xi\cup\eta)\c \b\xi \vspace*{1mm}\\
- & t \Big[\frac{(s+1)!(t+d-2)!}{(s+t)!(d-1)!}\cdot\spi{a}\cdot (x_1,y)\cdot (\b\xi_1\ohne x_1,\eta\ohne y)\c \Big( (\xi\cup y)\cup (\eta\ohne y)\Big)\c (\b\xi\ohne x_1)\Big]. \\  
\end{array}
\]
Therefore,
\[
\begin{array}{rcl}
\spi{a}\c (\xi\cup\eta)\c\b\xi 
& = & \frac{s!(t+d-1)!}{(s+t)!(d-2)!}\cdot\spi{a}\cdot(\b\xi_2,\eta)\c (\xi\cup\eta)\c (\b\xi\ohne x_1) \vspace*{1mm}\\
& + & t\frac{s!(t+d-1)!}{(s+t)!(d-1)!}\cdot\spi{a}\cdot(\b\xi_1,\eta)\c (\xi\cup\eta)\c (\b\xi\ohne x_1) \\
& \auf{\mbox{\scr (\ref{LemG2})}}{=} & -(d-1)\frac{s!(t+d-1)!}{(s+t)!(d-1)!}\cdot\spi{a}\cdot(\b\xi_1,\eta)\c (\xi\cup\eta)\cdot(x_0,x_1)\c (\b\xi_0\cup\b\xi_1\ohne x_1) \vspace*{1mm}\\
& + & t\frac{s!(t+d-1)!}{(s+t)!(d-1)!}\cdot\spi{a}\cdot(\b\xi_1,\eta)\c (\xi\cup\eta)\c (\b\xi_0\cup\b\xi_1\ohne x_1) \vspace*{1mm}\\
& \auf{\mbox{\scr (\ref{LemG1})}}{=} & \frac{s!(t+d-1)!}{(s+t)!(d-1)!}\cdot\spi{a}\cdot (\b\xi_1,\eta)\c (\xi\cup\eta)\c\b\xi. \\ 
\end{array}
\]
\end{Lemma}

\begin{Proposition}[Garnir formula]
\label{PropG6}
Assume given a subset $\phi\tm\Cup_{j\in [1,p-1]} a_j$. We have
\[
\spi{a}\c (\phi\cup\xi\cup\eta) = \frac{s!}{(s+t)!(d-1)!}\cdot\spi{a}\cdot(\b\xi_1,\eta)\c\b\xi\c (\phi\cup\xi\cup\eta).
\]

\rm
Let $A$ be the set of injections of $\phi$ into $\phi\cup\xi\cup\eta$. We write
\[
\begin{array}{rcl}
\spi{a}\c (\phi\cup\xi\cup\eta)                           & = & \sumd{\alpha\in A}\;\;\sumd{\sigma\in\Sl_{\phi\cup\xi\cup\eta},\;\sigma|_\phi = \alpha} \spi{a}\cdot \sigma\eps_\sigma \\
\spi{a}\cdot(\b\xi_1,\eta)\c\b\xi\c (\phi\cup\xi\cup\eta) 
& = & \sumd{\alpha\in A}\;\;\sumd{\sigma\in\Sl_{\phi\cup\xi\cup\eta},\;\sigma|_\phi = \alpha} \spi{a}\cdot(\b\xi_1,\eta)\c\b\xi\cdot \sigma\eps_\sigma \\ 
\end{array}
\]
and fix an embedding $\alpha\in A$ in order to compare the summands. We choose a permutation $\rho\in\Sl_{\phi\cup\xi\cup\eta}$ that restricts to $\rho|_{\phi} = \alpha$.
On the one hand, we obtain
\[
\begin{array}{rcl}
\sumd{\sigma\in\Sl_{\phi\cup\xi\cup\eta},\;\sigma|_\phi = \alpha} \spi{a}\cdot \sigma\eps_\sigma 
& \auf{\sigma\; =\;\sigma'\rho}{=}   & \sumd{\sigma'\in\Sl_{\xi\cup\eta}} \spi{a}\cdot\sigma'\eps_{\sigma'}\cdot\rho\eps_\rho \\ 
& =                                  & \spi{a}\c (\xi\cup\eta)\cdot\rho\eps_\rho \\ 
& \auf{\mbox{\scr (\ref{LemG4})}}{=} & \fracd{s!}{(d-1)!}\cdot\spi{a}\cdot (\b\xi_1,\eta)\c\b\xi\cdot\rho\eps_\rho, \\
\end{array}
\]
on the other hand, we have
\[
\begin{array}{rcl}
\sumd{\sigma\in\Sl_{\phi\cup\xi\cup\eta},\;\sigma|_\phi = \alpha} \spi{a}\cdot(\b\xi_1,\eta)\c\b\xi\cdot \sigma\eps_\sigma 
& \auf{\sigma\; =\;\sigma'\rho}{=} & \sumd{\sigma'\in\Sl_{\xi\cup\eta}} \spi{a}\cdot(\b\xi_1,\eta)\c\b\xi\cdot\sigma'\eps_{\sigma'}\cdot\rho\eps_\rho \\
& =                               & (s+t)!\cdot\spi{a}\cdot(\b\xi_1,\eta)\c\b\xi\cdot\rho\eps_\rho, \\
\end{array}
\]
whence the result follows by comparison. Independence of the right hand side of the choice of $\b\xi_1$ can also be seen by application of (\ref{LemG2}) to an appropriate permutation
$\tau\in\Sl_{\b\xi}$, for we have $(-)\cdot\tau\c\b\xi = (-)\c\b\xi\eps_\tau$. 
\end{Proposition}

\begin{Lemma}
\label{LemG7}
For $i_0\in [0,s]$, we write
\[
B_{[a],\xi,\eta}(i_0) := \sumd{s_0\in [i_0,s]} (-1)^{s_0} \sumd{\auf{\scm\xi_0\tm\xi,}{\#\xi_0 \; =\; s_0}}\;\;\sumd{\auf{\scm\eta_0\tm\eta,}{\#\eta_0 \; =\; s_0 - i_0}}\;\;
                         \sumd{\auf{\scm\phi_0\tm\b\eta,}{\#\phi_0 \; =\; i_0}} \spi{a}\cdot (\xi_0,\eta_0\cup\phi_0).
\]
Let $i\in [1,u]$. We define and obtain 
\[
\begin{array}{rcl}
C_{[a],\xi,\eta}(i) 
& := & \sumd{\auf{\scm\phi\tm\b\eta,}{\#\phi \; =\; i}} \spi{a}\c (\xi\cup\eta\cup\phi) \\
&  = & s!(t+i)! \sumd{i_0\in [0,\min(i,s)]} \smatze{u - i_0}{i - i_0} B_{[a],\xi,\eta}(i_0).
\end{array}
\]

\rm
We fix a subset $\phi\tm\b\eta$ of cardinality $\#\phi = i$.  Indexing by subsets that become interchanged by a certain element of $\Sl_{\xi\cup\eta\cup\phi}$, we obtain
\[
\begin{array}{l}
\spi{a}\c (\xi\cup\eta\cup\phi) \\
= s!(t+i)!\sumd{i_0\in [0,\min(i,s)]}\;\;\sumd{s_0\in [i_0,s]}\;\;\sumd{\auf{\scm\xi_0\tm\xi,}{\#\xi_0 \; =\; s_0}}\;\;\sumd{\auf{\scm\eta_0\tm\eta,}{\#\eta_0 \; =\; s_0 - i_0}}\;\;
                  \sumd{\auf{\scm\phi_0\tm\phi,}{\#\phi_0 \; =\; i_0}} \spi{a}\cdot (\xi_0,\eta_0\cup\phi_0)\eps_{(\xi_0,\eta_0\cup\phi_0)}. \\
\end{array}
\]
The number of subsets $\phi$ of cardinality $\#\phi = i$ that lie in between $\phi_0\tm\phi\tm\b\eta$, where $\phi_0$ is a given subset of $\b\eta$ of cardinality $\#\phi_0 = i_0$, amounts 
to $\smatze{u - i_0}{i - i_0}$. Whence the result is obtained by collecting
\[
\sumd{\auf{\scm\phi\tm\b\eta,}{\#\phi \; =\; i}}\;\;\sumd{\auf{\scm\phi_0\tm\phi,}{\#\phi_0 \; =\; i_0}} U(\phi_0) 
= \smatze{u-i_0}{i-i_0} \sumd{\auf{\scm\phi_0\tm\b\eta,}{\#\phi_0 \; =\; i_0}} U(\phi_0),
\]
$U(\phi_0)$ being some expression independent of $\phi$.
\end{Lemma}

\begin{Lemma}
\label{LemG8}
Suppose given in addition a subset $\psi\tm\xi$, of cardinality $\#\psi =: v$. For $i_0\in [0,s-v]$, we write
\[
B'_{[a],\xi,\psi,\eta}(i_0) := \frac{(-1)^{i_0}}{v!t!} \sumd{\auf{\scm\xi_0\tm\xi\ohne\psi,}{\#\xi_0 \; =\; i_0}}\;\;\sumd{\auf{\scm\phi_0\tm\eta,}{\#\phi_0 \; =\; i_0}}\;\;
                               \spi{a}\cdot (\xi_0,\phi_0)\c (\psi\cup\eta).
\]
Let $i\in [1,u]$. We define and obtain 
\[
\begin{array}{rcl}
C'_{[a],\xi,\psi,\eta}(i) 
& := & \sumd{\auf{\scm\phi\tm\eta,}{\#\phi \; =\; i}} \spi{a}\c (\xi\cup\phi) \c (\psi\cup\eta) \\
&  = & t!s!(i+v)!\cdot \sumd{i_0\in [0,\min(i,s-v)]}\smatze{t - i_0}{i - i_0} B'_{[a],\xi,\psi,\eta}(i_0). \\
\end{array}
\]

\rm
We fix a subset $\phi\tm\eta$ of cardinality $\#\phi = i$ and obtain
\[
\begin{array}{cl}
  & \frac{1}{s!i!} \spi{a}\c (\xi\cup\phi) \c (\psi\cup\eta) \\
= & \frac{1}{(s-v)!i!} \smatze{i+v}{v}\spi{a}\c ((\xi\ohne\psi)\cup\phi) \c (\psi\cup\eta) \\
= & \smatze{i+v}{v}\sumd{i_0\in [0,\min(i,s-v)]}\;\;
    \sumd{\auf{\scm\xi_0\tm\xi\ohne\psi,}{\#\xi_0 \; =\; i_0}}\;\;\sumd{\auf{\scm\phi_0\tm\phi,}{\#\phi_0 \; =\; i_0}} (-1)^{i_0}\spi{a}\cdot (\xi_0,\phi_0)\c (\psi\cup\eta), \\
\end{array}
\]
whence we conclude as in (\ref{LemG7}) that
\[
\begin{array}{cl}
  & \sumd{\auf{\scm\phi\tm\eta,}{\#\phi \; =\; i}} \spi{a}\c (\xi\cup\phi) \c (\psi\cup\eta) \\
= & \sumd{\auf{\scm\phi\tm\eta,}{\#\phi \; =\; i}} \frac{s!(i+v)!}{v!} \sumd{i_0\in [0,\min(i,s-v)]}\;\; 
    \sumd{\auf{\scm\xi_0\tm\xi\ohne\psi,}{\#\xi_0 \; =\; i_0}}\;\;\sumd{\auf{\scm\phi_0\tm\phi,}{\#\phi_0 \; =\; i_0}} (-1)^{i_0}\spi{a}\cdot (\xi_0,\phi_0)\c (\psi\cup\eta) \\
= & \frac{s!(i+v)!}{v!} \sumd{i_0\in [0,\min(i,s-v)]}\;\;\sumd{\auf{\scm\xi_0\tm\xi\ohne\psi,}{\#\xi_0 \; =\; i_0}}\;\;
    \smatze{t - i_0}{i - i_0}\sumd{\auf{\scm\phi_0\tm\eta,}{\#\phi_0 \; =\; i_0}} (-1)^{i_0}\spi{a}\cdot (\xi_0,\phi_0)\c (\psi\cup\eta). \\

\end{array}
\]
\end{Lemma}          
\section{A two-box-shift morphism}
\label{SecTwoBox}

\subsection{Double paths}
\label{SubSecMS}

Let $\lambda$ be a partition of $n$. We assume given integers $g,k$ such that $1\leq g\leq k\leq \lambda_1 - 1$ and such that
\[
\mu'_i := 
\left\{
\begin{array}{ll}
\lambda'_i + 2 & \mbox{for } i = g \\
\lambda'_i - 2 & \mbox{for } i = k+1 \\
\lambda'_i     & \mbox{else} \\
\end{array}
\right.
\] 
defines a partition $\mu$. A {\it weight} $e$ is a map
\[
\begin{array}{rcl}
[1,\lambda_1] & \lraa{e} & [0,2] \\
j             & \lra     & e_j \\
\end{array}
\]
that maps $g$ and $k+1$ to $e_g = e_{k+1} = 2$ and that maps each $j\in [1,\lambda_1]\ohne [g,k+1]$ to $e_j = 0$. The set of weights is denoted by $E$. A {\it pattern} $\Xi$ of weight $e$ is a 
subset $\Xi\tm [1,2]\ti [g,k+1]$ such that 
\[
\#(\Xi\cap ([1,2]\ti\{ j\})) = e_j
\]
for $j\in [g,k+1]$. A {\it double path} $\gamma$ of weight $e$ is an injection from a pattern $\Xi$ of weight $e$ to 
$[\lambda]\cup [\mu]$ of the form
\[
\begin{array}{rcl}
\Xi    & \lraa{\gamma} & [\lambda]\cup [\mu] \\
i\ti j & \lra          & \b\gamma(j,i)\ti j    \\
\end{array}
\]  
such that
\[
\begin{array}{ll}
\b\gamma(g,i) = \lambda'_g + i     & \mbox{for } i\in [1,2] \\
\b\gamma(k+1,1) < \b\gamma(k+1,2). & \\
\end{array}
\]
Sometimes, we denote its pattern by $\Xi_\gamma := \Xi$. Its {\it sign} is given by $\eps_\gamma := (-1)^{\b\gamma(k+1,1) + \b\gamma(k+1,2)}$, not to be confused with the sign of a permutation. 
The set of double paths of weight $e$ is denoted by $\dGamma(e)$.

A double path $\gamma$ gives rise to a place operation in the following manner. Let
\[
\begin{array}{rcl}
[0,\lambda'_{k+1}]       & \lra        & [0,\mu'_{k+1}] \\
i                        & \lraa{\phi} & \#\Big( [1,i]\ohne\{\b\gamma(k+1,1),\b\gamma(k+1,2)\}\Big) \\
\min(\phi^{-1}(\{ j\}))  & \llaa{\psi} & j. \\
\end{array}
\]
Given a weight $e\in E$, we define an operation of sets
\[
\begin{array}{rclcl}
\dGamma(e) & \ti & T^\lambda & \lra & T^\mu \\
\gamma     & \ti & [a]       & \lra & [a^\gamma] \\
\end{array}
\]
by 
\[
\begin{array}{rcll}
a^\gamma_{j,i}              & := & a_{j,i} \mbox{ for } i\ti j\in [\mu]\ohne (\Xi)\gamma \mbox{ and } j\neq k+1, \\ 
a^\gamma_{j,\b\gamma(j,i)}  & := & a_{j',\b\gamma(j',i)} \mbox{ for $i\in [1,2]$ and for $j\in [g,k]$,} \\
                            &    & \mbox{where $j'$ is minimal in $[j+1,k+1]$ with $i\ti j'\in\Xi$}, \\
a^\gamma_{k+1,i}            & := & a_{k+1,\psi(i)} \mbox{ for } i\in [1,\mu'_{k+1}]. \\
\end{array}
\]
We define a $\Z\Sl_n$-linear map
\[
\begin{array}{rcl}
F^\lambda & \lraa{f''_e}    & S^\mu \\
{[a]}     & \lra            & \sumd{\gamma\in\dGamma(e)} \spi{a^\gamma}\eps_\gamma. \\
\end{array}
\]

\begin{Example}
\label{ExM-1}
\rm
For example, consider the case $\lambda = (7,7,7,4)$, $\mu = (7,6,6,4,1,1)$, $g = 1$, $k = 6$. Let $e = (2,1,1,1,2,0,2)$, and let $\gamma$ be the (non-ordered) double path of weight $e$ given 
pictorially by
\begin{center}
\begin{picture}(400,300)
\put(   0, 250){x}
\put(  50, 250){x}
\put( 100, 250){v}
\put( 150, 250){x}
\put( 200, 250){x}
\put( 250, 250){x}
\put( 300, 250){u}

\put(   0, 200){x}
\put(  50, 200){u}
\put( 100, 200){x}
\put( 150, 200){x}
\put( 200, 200){v}
\put( 250, 200){x}
\put( 300, 200){x}

\put(   0, 150){x}
\put(  50, 150){x}
\put( 100, 150){x}
\put( 150, 150){u}
\put( 200, 150){u}
\put( 250, 150){x}
\put( 300, 150){v}

\put(   0, 100){x}
\put(  50, 100){x}
\put( 100, 100){x}
\put( 150, 100){x}

\put(   0,  50){u}

\put(   0,   0){v,}
\end{picture}
\end{center}
a `u' indicating the $\gamma$-image of some $1\ti j$, a `v' indicating the $\gamma$-image of some $2\ti j$. The operation of $\gamma$ yields e.g.\
\begin{center}
\begin{picture}(1000,250)
\put(-150, 170){$\gamma\;\;\;\;\ti$}

\put(   0, 250){\scr 1}
\put(  50, 250){\scr 5}
\put( 100, 250){\scr 9}
\put( 150, 250){\scr 13}
\put( 200, 250){\scr 17}
\put( 250, 250){\scr 20}
\put( 300, 250){\scr 23}

\put(   0, 200){\scr 2}
\put(  50, 200){\scr 6}
\put( 100, 200){\scr 10}
\put( 150, 200){\scr 14}
\put( 200, 200){\scr 18}
\put( 250, 200){\scr 21}
\put( 300, 200){\scr 24}

\put(   0, 150){\scr 3}
\put(  50, 150){\scr 7}
\put( 100, 150){\scr 11}
\put( 150, 150){\scr 15}
\put( 200, 150){\scr 19}
\put( 250, 150){\scr 22}
\put( 300, 150){\scr 25}

\put(   0, 100){\scr 4}
\put(  50, 100){\scr 8}
\put( 100, 100){\scr 12}
\put( 150, 100){\scr 16}

\put( 380, 180){\vector(1,0){150}}

\put( 600, 250){\scr 1}
\put( 650, 250){\scr 5}
\put( 700, 250){\scr 18}
\put( 750, 250){\scr 13}
\put( 800, 250){\scr 17}
\put( 850, 250){\scr 20}
\put( 900, 250){\scr 24}

\put( 600, 200){\scr 2}
\put( 650, 200){\scr 15}
\put( 700, 200){\scr 10}
\put( 750, 200){\scr 14}
\put( 800, 200){\scr 25}
\put( 850, 200){\scr 21}

\put( 600, 150){\scr 3}
\put( 650, 150){\scr 7}
\put( 700, 150){\scr 11}
\put( 750, 150){\scr 19}
\put( 800, 150){\scr 23}
\put( 850, 150){\scr 22}

\put( 600, 100){\scr 4}
\put( 650, 100){\scr 8}
\put( 700, 100){\scr 12}
\put( 750, 100){\scr 16}

\put( 600,  50){\scr 6}

\put( 600,   0){\scr 9}
\end{picture}
\end{center}
I.e.\ the double path $\gamma$ subsequently pushes $\scm 23\ra 19\ra 15\ra 6\ra$ and $\scm 25\ra 18\ra 9\ra$. 
\end{Example}

\begin{Lemma}[path switch]
\label{LemM0}
For $q\in [g+1,k+1]$, we dispose of an involution
\[
\begin{array}{rclcl}
{[1,2]} & \ti & [g,k+1] & \lraa{\iota_q} & \hspace*{10mm} [1,2] \;\;\ti\;\; [g,k+1] \\
i       & \ti & j       & \lra           & \left\{\begin{array}{rcll}
                                         i\cdot (1,2) & \ti & j & \mbox{for } j\in [g+1,q-1] \\ 
                                         i            & \ti & j & \mbox{else} \\         
                                         \end{array}.\right. \\
\end{array} 
\]
The composition $\iota_p\gamma$ of this involution, restricted to $\iota_p^{-1}(\Xi_\gamma)$, with a double path $\gamma$ furnishes a double path which we
denote by $\iota_p\gamma$, being a slight abuse of notation. So $\Xi_{\iota_p\gamma} = \iota_p^{-1}(\Xi_\gamma)$.

Let $[a]\in T^\lambda$, $e\in E$, let $p\in [g+1,k]$ such that $e_p = 2$, let $\gamma\in\dGamma(e)$, $u := a_{p,\b\gamma(p,1)}$ and $v := a_{p,\b\gamma(p,2)}$. We assert the following 
equalities.
\[
\begin{array}{rrcl}
\mbox{\rm (i)}   & \spi{a^\gamma}\cdot (u,v)   & = & -\spi{a^{\iota_p\gamma}} \\
\mbox{\rm (ii)}  & \spi{a^{\iota_{p+1}\gamma}} & = & \spi{a^\gamma}           \\
\mbox{\rm (iii)} & ([a]f''_e)(u,v)             & = & -[a]f''_e                \\
\end{array}
\]

\rm
The operation of the double path $\iota_p\gamma$ on $[a]$ can be performed using the operation of $\gamma$, but with entries $u$ and $v$ interchanged, and, moreover, with entries placed by 
$\gamma$ in column $g$ interchanged, yielding the sign in (i). Analoguously (ii). Since composition with $\iota_p$ is an involution on $\dGamma(e)$, multiplication by $(u,v)$ attaches a sign 
to $[a]f''_e$ by (i), whence (iii).
\end{Lemma}

\begin{Proposition}
\label{PropM1}
Let $e\in E$. There is a factorization 
\[
(F^\lambda\lraa{f''_e}S^\mu) = (F^\lambda\lraa{\nu^\lambda_M} M^{\lambda',-}\lraa{f'_e}S^\mu).
\]

\rm
We need to show that a signed column transposition $[a] + [a](u,v)$, $[a]$ a $\lambda$-tableau, $p\geq 1$, $u,v\in a_p$, $u\neq v$, vanishes under $f''_e$, for which we may assume $p$ to lie in 
$[g+1,k+1]$. We abbreviate the double path image of $\gamma$ in $[a]$ by $I_\gamma := ((\Xi)\gamma)[a]$.

{\bf Case $p\in [g+1,k]$.}  

{\it Subcase $u, v\in I_\gamma$.} See (\ref{LemM0}, iii).

{\it Subcase $u =: a_{p,\b\gamma(p,i)}\in I_\gamma$, $v =: a_{p,l}\not\in I_\gamma$.} Abbreviating by $\tau_\gamma$ the transposition that interchanges the elements 
$\b\gamma(p,i)\ti p$ and $l\ti p$ of $[\lambda]\cup [\mu]$, we conclude from
\[
\spi{a^\gamma}\cdot (u,v) = -\spi{a^{\gamma\tau_\gamma}}
\]
that the sum over these summands, multiplied by $(u,v)$, yields minus the sum over the summands of the subcase $u\not\in I_\gamma, v\in I_\gamma$. In fact,  
$\gamma\lra\gamma\tau_\gamma$ is a bijection between the sets of double paths occurring in these subcases which preserves the sign. 

{\it Subcase $u, v\not\in I_\gamma$.} $\spi{a^\gamma}\cdot (u,v) = -\spi{a^{\gamma}}$ shows these summands to yield zero.

{\bf Case $p = k+1$.} The arguments of the case $p\in [g+1,k]$ work in the {\it subcases} $u, v\in I_\gamma$ and $u, v\not\in I_\gamma$ as well.  

{\it Subcase $u =: a_{k+1,\b\gamma(k+1,i)}\in I_\gamma$, $v =: a_{k+1,l}\not\in I_\gamma$.} Writing the transposition $\tau_\gamma$ as above, we conclude from 
\[
\spi{a^\gamma}\cdot (u,v) = (-1)^{\b\gamma(p,i) - l + 1}\spi{a^{\gamma\tau_\gamma}}
\]
that the sum over these summands yields minus the result of the sum over the summands of the subcase $u \not\in I_\gamma, v\in I_\gamma$. In fact, $\gamma\lra\gamma\tau_\gamma$ 
is a bijection between the sets of double paths occurring in these subcases which changes the sign by $\eps_{\gamma\tau_\gamma}/\eps_\gamma = (-1)^{l - \b\gamma(p,i)}$. 
\end{Proposition}

\subsection{Morphism, unreduced version}
\label{SubSecLinComb}

\begin{footnotesize}
We shall show the vanishing, modulo our modulus, of the required Garnir relations. Rather than giving short versions of all calculations, we give
full versions, but only of the two `large' cases (\ref{LemM2}, \ref{LemM8}), and only in Appendix \ref{AppTwoLem}. 
\end{footnotesize}

For $i\in [g,k+1]$, we let 
\[
X_i := (\lambda'_i - i) - (\lambda'_{k+1} - (k+1))
\]
and denote $X_i^{(j)} := (X_i + j - 1)!/(X_i - 1)!$, $j\in [0,2]$. Let $M^{\lambda',-}\lraa{f'}S^\mu$ be the $\Z\Sl_n$-linear map defined by
\[
f' := \sumd{e\in E} \left(\prodd{i\in [g+1,k]} X_i^{(2-e_i)}\right) f'_e.
\]
Let $[a]$ be a $\lambda$-tableau, $p\in [g,k]$, $\xi\tm a_p$, $s := \#\xi$, $\b\xi := a_p\ohne\xi$ and $\eta\tm a_{p+1}$, $t := \#\eta$, 
$\b\eta := a_{p+1}\ohne\eta$ such that $s + t = \lambda'_p + 1$. Denote $u := \#\b\eta = \lambda'_{p+1} - t$. 

\begin{Lemma}[cf.\ (\ref{LemM2})]
\label{LemM2short}
Suppose $g < p < k$ and $s,t\geq 2$. The map $f'$ annihilates $G'_{[a],\xi,\eta}$.

\rm
We fix a map
\[
\begin{array}{rcl}
[1,\lambda_1]\ohne\{p,p+1\} & \lraa{\w e} & [0,2] \\
j                           & \lra        & \w e_j \\
\end{array}
\]
that maps $g$ and $k+1$ to $e_g = e_{k+1} = 2$, and that maps $j\in [1,\lambda_1]\ohne [g,k+1]$ to $e_j = 0$. For $\alpha,\beta\in [0,2]$, we denote by $\w e\alpha\beta$ be the prolongation 
of $\w e$ to $[1,\lambda_1]$ defined by $\w e\alpha\beta|_{[1,\lambda_1]\ohne\{p,p+1\}} := \w e$, $(\w e\alpha\beta)_p := \alpha$ and $(\w e\alpha\beta)_{p+1} := \beta$. We contend that
\[
\sumd{\alpha,\beta\in [0,2]} X_p^{(2-\alpha)} X_{p+1}^{(2-\beta)} G'_{[a],\xi,\eta} f'_{\w e\alpha\beta} = 0,
\]
from which the lemma ensues.

There exist elements $x,y\in\xi$, $x\neq y$, $x', y'\in\eta$, $x'\neq y'$, $z\in\b\xi$, which we choose and fix. For $\gamma\in\dGamma(\w e00)$, we let $x_\gamma := a_{j,\b\gamma(j,1)}$, where
$j\in [p+2,k+1]$ is minimal with $1\ti j\in\Xi_\gamma$, and $y_\gamma := a_{j,\b\gamma(j,2)}$, where $j\in [p+2,k+1]$ is minimal with $2\ti j\in\Xi_\gamma$. I.e.\ we pick the entries 
$x_\gamma, y_\gamma$ that `cross the columns' $p$ and $p+1$ under the operation of $\gamma$. We write 
\[
\begin{array}{rcl}
U_\gamma     & := & (s+t-2)!^{-1} \cdot\spi{a^\gamma}\eps_\gamma\cdot (\b\xi\ohne z,\eta\ohne\{ x',y'\})\cdot (x_\gamma,x,x')\cdot (y_\gamma,y,y') \c (\xi\cup\eta) \\
V_{1,\gamma} & := & (s+t-1)!^{-1} \cdot\spi{a^\gamma}\eps_\gamma\cdot (\b\xi,\eta\ohne x')\cdot (x_\gamma,x,x')\c (\xi\cup\eta) \\
V_{2,\gamma} & := & (s+t-1)!^{-1} \cdot\spi{a^\gamma}\eps_\gamma\cdot (\b\xi,\eta\ohne y')\cdot (y_\gamma,y,y')\c (\xi\cup\eta), \\
\end{array}
\]
and let
\[
\begin{array}{lcl}
A   & := & \sumd{\gamma\in \dGamma(\w e00)} U_\gamma\cdot\sumd{w'\in\b\eta} (w',z) \\
B   & := & \sumd{\gamma\in \dGamma(\w e00)} U_\gamma\cdot(1 - \sumd{w\in\b\xi\ohne z} (w,z)) \\
C_1 & := & \sumd{\gamma\in \dGamma(\w e00)} U_\gamma\cdot (z,x_\gamma) \\
C_2 & := & \sumd{\gamma\in \dGamma(\w e00)} U_\gamma\cdot (z,y_\gamma) \\
D   & := & \sumd{\gamma\in \dGamma(\w e00)} V_{1,\gamma} \cdot \sumd{w\in\b\xi} (w,y_\gamma) \\
    &  = & \sumd{\gamma\in \dGamma(\w e00)} V_{2,\gamma} \cdot \sumd{w\in\b\xi} (w,x_\gamma) \\
H   & := & \sumd{\gamma\in \dGamma(\w e00)} V_{1,\gamma} \cdot \sumd{w'\in\b\eta} (w',y_\gamma) \\
    &  = & \sumd{\gamma\in \dGamma(\w e00)} V_{2,\gamma} \cdot \sumd{w'\in\b\eta} (w',x_\gamma) \\
F_1 & := & \sumd{\gamma\in \dGamma(\w e00)} V_{1,\gamma} \\
F_2 & := & \sumd{\gamma\in \dGamma(\w e00)} V_{2,\gamma}. \\
\end{array}
\]
Calculations carried out in (\ref{LemM2}) yield the following table.
\[
\begin{array}{rcl}
G'_{[a],\xi,\eta}f'_{\w e22} & = & 2(s - u)A + (s - u)(s - u - 1)B + 2(s - u)(s - u + 1)D \\
G'_{[a],\xi,\eta}f'_{\w e12} & = & - 2A - 2(s - u)B - 2(s - u + 1)D + 2(s - u + 1)H \\
G'_{[a],\xi,\eta}f'_{\w e02} & = & B - 2H \\
G'_{[a],\xi,\eta}f'_{\w e21} & = & 2A + 2(s - u - 1)B - (s - u - 1)(C_1 + C_2) + 4(s - u)D \\
G'_{[a],\xi,\eta}f'_{\w e11} & = & - 2B + (C_1 + C_2) - 2D + 2H + (s - u)(F_1 + F_2) \\
G'_{[a],\xi,\eta}f'_{\w e01} & = & -(F_1 + F_2) \\
G'_{[a],\xi,\eta}f'_{\w e20} & = & B - (C_1 + C_2) + 2D \\
G'_{[a],\xi,\eta}f'_{\w e10} & = & (F_1 + F_2) \\
G'_{[a],\xi,\eta}f'_{\w e00} & = & 0 \\
\end{array}
\]
We note that $s - u = X_p - X_{p+1}$ and evaluate the linear combination
\[
\begin{array}{l}
\sumd{\alpha,\beta\in [0,2]} X_p^{(2-\alpha)} X_{p+1}^{(2-\beta)} G'_{[a],\xi,\eta} f'_{\w e\alpha\beta} \\
\begin{array}{clll}
= & 1            & \cdot 1                    & \cdot\Big( 2(X_p - X_{p+1})A + (X_p - X_{p+1})(X_p - X_{p+1} - 1)B \\
  &              &                            & \;\; + 2(X_p - X_{p+1})(X_p - X_{p+1} + 1)D\Big) \\
+ & X_p          & \cdot 1                    & \cdot\Big( - 2A - 2(X_p - X_{p+1})B - 2(X_p - X_{p+1} + 1)D \\
  &              &                            & \;\; + 2(X_p - X_{p+1} + 1)H\Big) \\
+ & X_p(X_p + 1) & \cdot 1                    & \cdot\Big( B - 2H\Big) \\
+ & 1            & \cdot X_{p+1}              & \cdot\Big( 2A + 2(X_p - X_{p+1} - 1)B  \\
  &              &                            & \;\; - (X_p - X_{p+1} - 1)(C_1 + C_2) + 4(X_p - X_{p+1})D\Big) \\
+ & X_p          & \cdot X_{p+1}              & \cdot\Big( - 2B + (C_1 + C_2) - 2D + 2H \\
  &              &                            & \;\; + (X_p - X_{p+1})(F_1 + F_2)\Big)\\
+ & X_p(X_p + 1) & \cdot X_{p+1}              & \cdot\Big( -(F_1 + F_2)\Big)\\
+ & 1            & \cdot X_{p+1}(X_{p+1} + 1) & \cdot\Big( B - (C_1 + C_2) + 2D\Big) \\
+ & X_p          & \cdot X_{p+1}(X_{p+1} + 1) & \cdot\Big( (F_1 + F_2)\Big)\\
+ & X_p(X_p + 1) & \cdot X_{p+1}(X_{p+1} + 1) & \cdot\Big( 0 \Big) \\
= & 0.\\
\end{array} \\
\end{array}
\]
\end{Lemma}

\begin{Lemma}
\label{LemM3}
Suppose $g < p < k$ and $s = \lambda'_p$, $t = 1$. The map $f'$ annihilates $G'_{[a],\xi,\eta}$.

\rm
Concerning the map $\w e$ and its prolongations, as well as concerning the integers $x_\gamma$, $y_\gamma$, we continue to use the notation of (\ref{LemM2short}). We fix such a map $\w e$ and
need to show that
\[
\sumd{\alpha,\beta\in [0,2]} X_p^{(2-\alpha)} X_{p+1}^{(2-\beta)} G'_{[a],\xi,\eta} f'_{\w e\alpha\beta} = 0.
\]
There exist elements $x, y\in\xi$, $x\neq y$, which we choose and fix.  For $\gamma\in\dGamma(\w e00)$, we write 
\[
\begin{array}{rcl}
V_{1,\gamma} & := & \frac{1}{(s-1)!}  \cdot\spi{a^\gamma}\eps_\gamma\cdot (y_\gamma,y,\eta) \c (\xi\cup\eta) \\
V_{2,\gamma} & := & \frac{1}{(s-1)!}  \cdot\spi{a^\gamma}\eps_\gamma\cdot (x_\gamma,x,\eta) \c (\xi\cup\eta) \\
\end{array}
\]
and let
\[
\begin{array}{lcl}
A   & := & \sumd{\gamma\in \dGamma(\w e00)} V_{1,\gamma}\cdot\sumd{w'\in\b\eta} (x_\gamma,x,w') \\ 
    &  = & \sumd{\gamma\in \dGamma(\w e00)} V_{2,\gamma}\cdot\sumd{w'\in\b\eta} (y_\gamma,y,w') \\ 
C_1 & := & - \sumd{\gamma\in \dGamma(\w e00)} V_{1,\gamma}\cdot (x_\gamma,x) \\
C_2 & := & - \sumd{\gamma\in \dGamma(\w e00)} V_{2,\gamma}\cdot (y_\gamma,y) \\
\end{array}
\]
\[
\begin{array}{lcl}
H   & := & s^{-1}\cdot\sumd{\gamma\in \dGamma(\w e00)} V_{1,\gamma}\cdot \sumd{w'\in\b\eta} (x_\gamma,w') \\
    &  = & s^{-1}\cdot\sumd{\gamma\in \dGamma(\w e00)} V_{2,\gamma}\cdot \sumd{w'\in\b\eta} (y_\gamma,w'). \\
F_1 & := & s^{-1}\cdot\sumd{\gamma\in \dGamma(\w e00)} V_{1,\gamma} \\
F_2 & := & s^{-1}\cdot\sumd{\gamma\in \dGamma(\w e00)} V_{2,\gamma} \\
\end{array}
\] 
Calculations similar to those of (\ref{LemM2}) yield the following table, to be compared to the table of (\ref{LemM2short}).
\[
\begin{array}{rcl}
G'_{[a],\xi,\eta}f'_{\w e22} & = & 2(s - u)A \\
G'_{[a],\xi,\eta}f'_{\w e12} & = & - 2A + 2(s - u + 1)H \\
G'_{[a],\xi,\eta}f'_{\w e02} & = & - 2H \\
G'_{[a],\xi,\eta}f'_{\w e21} & = & 2A - (s - u - 1)(C_1 + C_2) \\
G'_{[a],\xi,\eta}f'_{\w e11} & = & (C_1 + C_2) + 2H + (s - u)(F_1 + F_2) \\
G'_{[a],\xi,\eta}f'_{\w e01} & = & -(F_1 + F_2) \\
G'_{[a],\xi,\eta}f'_{\w e20} & = & -(C_1 + C_2) \\
G'_{[a],\xi,\eta}f'_{\w e10} & = & (F_1 + F_2) \\
G'_{[a],\xi,\eta}f'_{\w e00} & = & 0 \\
\end{array}
\]
\end{Lemma}

\begin{Lemma}
\label{LemM4}
Suppose $g < p < k$ and $s = 1$, $t = \lambda'_p = \lambda'_{p+1}$. The map $f'$ annihilates $G'_{[a],\xi,\eta}$.

\rm
Concerning $\w e$, $x_\gamma$, $y_\gamma$, we continue to use the notation of (\ref{LemM2short}). We fix such a map $\w e$ and need to show that
\[
\sumd{\alpha,\beta\in [0,2]} X_p^{(2-\alpha)} X_{p+1}^{(2-\beta)} G'_{[a],\xi,\eta} f'_{\w e\alpha\beta} = 0.
\]
There exist elements $x', y'\in\eta$, $x'\neq y'$, $z\in\b\xi$, which we choose and fix. For $\gamma\in\dGamma(\w e00)$, we write
\[
\begin{array}{rcl}
U_{1,\gamma} & := & - \frac{1}{(t-1)!} \cdot\spi{a^\gamma}\eps_\gamma\cdot (\b\xi\ohne z,\eta\ohne\{ x',y'\})\cdot (x_\gamma,\xi,x')\cdot(y_\gamma,y')\c (\xi\cup\eta) \\
U_{2,\gamma} & := & - \frac{1}{(t-1)!} \cdot\spi{a^\gamma}\eps_\gamma\cdot (\b\xi\ohne z,\eta\ohne\{ x',y'\})\cdot (y_\gamma,\xi,y')\cdot(x_\gamma,x')\c (\xi\cup\eta) \\
V_{1,\gamma} & := & \frac{1}{t!} \cdot\spi{a^\gamma}\eps_\gamma\cdot (\b\xi,\eta\ohne x')\cdot (x_\gamma,\xi,x')\c (\xi\cup\eta) \\ 
V_{2,\gamma} & := & \frac{1}{t!} \cdot\spi{a^\gamma}\eps_\gamma\cdot (\b\xi,\eta\ohne y')\cdot (y_\gamma,\xi,y')\c (\xi\cup\eta) \\ 
\end{array}
\]
and let
\[
\begin{array}{lcl}
B   & := & \sumd{\gamma\in \dGamma(\w e00)} U_{1,\gamma}\cdot(1 - \sumd{w\in\b\xi\ohne z} (w,z)) \\ 
    &  = & \sumd{\gamma\in \dGamma(\w e00)} U_{2,\gamma}\cdot(1 - \sumd{w\in\b\xi\ohne z} (w,z)) \\ 
C_1 & := & \sumd{\gamma\in \dGamma(\w e00)} U_{1,\gamma}\cdot (z,x_\gamma) \\ 
\end{array}
\]
\[
\begin{array}{lcl}
C_2 & := & \sumd{\gamma\in \dGamma(\w e00)} U_{2,\gamma}\cdot (z,y_\gamma) \\ 
D   & := & \sumd{\gamma\in \dGamma(\w e00)} V_{1,\gamma}\cdot \sumd{w\in\b\xi} (w,y_\gamma) \\ 
    &  = & \sumd{\gamma\in \dGamma(\w e00)} V_{2,\gamma}\cdot \sumd{w\in\b\xi} (w,x_\gamma) \\ 
F_1 & := & \sumd{\gamma\in \dGamma(\w e00)} V_{1,\gamma} \\
F_2 & := & \sumd{\gamma\in \dGamma(\w e00)} V_{2,\gamma} \\
\end{array}
\]
Calculations similar to those of (\ref{LemM2}), using the Garnir relation $B = C_1 + C_2$, yield the following table, to be compared to the table of (\ref{LemM2short}).
\[
\begin{array}{rcl}
G'_{[a],\xi,\eta}f'_{\w e22} & = & 4D \\
G'_{[a],\xi,\eta}f'_{\w e12} & = & - 2B - 4D \\
G'_{[a],\xi,\eta}f'_{\w e02} & = & B \\
G'_{[a],\xi,\eta}f'_{\w e21} & = & 4D \\
G'_{[a],\xi,\eta}f'_{\w e11} & = & - B - 2D + (F_1 + F_2)\\
G'_{[a],\xi,\eta}f'_{\w e01} & = & - (F_1 + F_2)\\
G'_{[a],\xi,\eta}f'_{\w e20} & = & 2D \\
G'_{[a],\xi,\eta}f'_{\w e10} & = & (F_1 + F_2)\\
G'_{[a],\xi,\eta}f'_{\w e00} & = & 0 \\
\end{array}
\]
Note that $s - u = X_p - X_{p+1} = 1$.
\end{Lemma}

\begin{Lemma}
\label{LemM5}
Suppose $g < p = k$ and $s,t\geq 2$, $\eta = a_{k+1,[1,t]}$. The map $f'$ annihilates $G'_{[a],\xi,\eta}$.

\rm
We fix a map
\[
\begin{array}{rcl}
[1,\lambda_1]\ohne \{ k\} & \lraa{\w e} & [0,2] \\
j                         & \lra        & \w e_j \\
\end{array}
\]
that maps $g$ and $k+1$ to $e_g = e_{k+1} = 2$, and that maps $j\in [1,\lambda_1]\ohne [g,k+1]$ to $e_j = 0$. For $\alpha\in [0,2]$, $\w e\alpha$ denotes the map which prolongs $\w e$ to 
$[1,\lambda_1]$ via $(\w e\alpha)_k = \alpha$. We need to show that
\[
\sumd{\alpha\in [0,2]} X_k^{(2-\alpha)} G'_{[a],\xi,\eta} f'_{\w e\alpha} = 0.
\]
There exist elements $x, y\in\xi$, $x\neq y$, $z\in\b\xi$, which we choose and fix. Let $x' := a_{g+1,1}$ and $y' := a_{g+1,2}$, so that $x',y'\in\eta$. For $v,w\in a_{k+1}$, $v\neq w$, we denote 
\[
\dGamma(\w e,a,v,w) := \{\gamma\in\dGamma(\w e0)\; |\; a_{k+1,\b\gamma(k+1,1)} = v,\; a_{k+1,\b\gamma(k+1,2)} = w\}
\]
and let
\[
\begin{array}{lcl}
A   & := & \frac{1}{(s+t-2)!} \sumd{\gamma\in \dGamma(\w e,a,x',y')} \spi{a^\gamma}\eps_\gamma\cdot (\b\xi\ohne z,\eta\ohne\{ x',y'\})\cdot (x,x')\cdot (y,y') \c (\xi\cup\eta)
           \cdot\sumd{w'\in\b\eta} (w',z) \\ 
B   & := & \frac{1}{(s+t-2)!} \sumd{\gamma\in \dGamma(\w e,a,x',y')} \spi{a^\gamma}\eps_\gamma\cdot \\
    &    & \cdot(\b\xi\ohne z,\eta\ohne\{ x',y'\})\cdot (x,x')\cdot (y,y') \c (\xi\cup\eta)
           \cdot(1 - \sumd{w\in\b\xi\ohne z} (w,z)) \\ 
D   & := & \frac{1}{(s+t-1)!} \sumd{\gamma\in\dGamma(\w e,a,x',y')}\;\;\sumd{w\in\b\xi}\spi{a^\gamma}\eps_\gamma\cdot (\b\xi\ohne w,\eta\ohne\{ x',y'\})
           \cdot (x,x')\cdot (w,y') \c (\xi\cup\eta) \\ 
    &  = & \frac{1}{(s+t-1)!} \sumd{\gamma\in\dGamma(\w e,a,x',y')}\;\;\sumd{w\in\b\xi}\spi{a^\gamma}\eps_\gamma\cdot (\b\xi\ohne w,\eta\ohne\{ x',y'\})
           \cdot (y,y')\cdot (w,x') \c (\xi\cup\eta) \\ 
H   & := & \frac{1}{(s+t-1)!}\sumd{w'\in\b\eta}\;\; \sumd{\gamma\in\dGamma(\w e,a,x',w')} \spi{a^\gamma}\eps_\gamma 
           \cdot (\b\xi,\eta\ohne x')\cdot (x,x')\c (\xi\cup\eta) \\  
    &  = & - \frac{1}{(s+t-1)!}\sumd{w'\in\b\eta}\;\;\sumd{\gamma\in\dGamma(\w e,a,y',w')}\spi{a^\gamma}\eps_\gamma 
           \cdot (\b\xi,\eta\ohne y')\cdot (y,y')\c (\xi\cup\eta). \\  
\end{array}
\]
Calculations similar to those of (\ref{LemM2}) yield the following table.
\[
\begin{array}{rcl}
G'_{[a],\xi,\eta}f'_{\w e2} & = & \frac{1}{2}\Big( 2(s - u)A + (s - u)(s - u - 1)B + 2(s-u)(s-u+1)D \Big) \\
G'_{[a],\xi,\eta}f'_{\w e1} & = & \frac{1}{2}\Big( - 2A - 2(s - u)B - 2(s - u + 1)D + 2(s - u + 1)H \Big) \\
G'_{[a],\xi,\eta}f'_{\w e0} & = & \frac{1}{2}\Big( B - 2H \Big) \\
\end{array}
\]
Comparison with the table in (\ref{LemM2short}) is possible since $X_{k+1} = 0$.
\end{Lemma}

\begin{Lemma}
\label{LemM6}
Suppose $g < p = k$ and $s = \lambda'_k$, $t = 1$, $\eta = a_{k+1,1}$. The map $f'$ annihilates $G'_{[a],\xi,\eta}$.

\rm
Concerning $\w e$ and $\dGamma$, we continue to use the notation of (\ref{LemM5}). We fix such a map $\w e$. There exist elements $x, y\in\xi$, $x\neq y$, which we choose and fix. We let
\[
\begin{array}{lcl}
A   & := & \frac{1}{(s-1)!} \sumd{w'\in\b\eta}\;\;\sumd{\gamma\in \dGamma(\w e,a,\eta,w')} \spi{a^\gamma}\eps_\gamma
           \cdot (x,\eta)\cdot (y,w') \c (\xi\cup\eta) \\
H   & := & \frac{1}{s!} \sumd{w'\in\b\eta}\;\;\sumd{\gamma\in \dGamma(\w e,a,\eta,w')} \spi{a^\gamma}\eps_\gamma
           \cdot (x,\eta) \c (\xi\cup\eta)
\end{array}
\]
Calculations similar to those of (\ref{LemM2}) yield the following table, to be compared to the table of (\ref{LemM5}).
\[
\begin{array}{rcl}
G'_{[a],\xi,\eta}f'_{\w e2} & = & \frac{1}{2}\Big( 2(s - u)A \Big) \\
G'_{[a],\xi,\eta}f'_{\w e1} & = & \frac{1}{2}\Big( - 2A + 2(s - u + 1)H \Big) \\
G'_{[a],\xi,\eta}f'_{\w e0} & = & \frac{1}{2}\Big( - 2H \Big) \\
\end{array}
\]
\end{Lemma}

\begin{Lemma}
\label{LemM7}
Suppose $g < p = k$ and $s = 1$, $t = \lambda'_k = \lambda'_{k+1}$. The map $f'$ annihilates $G'_{[a],\xi,\eta}$.

\rm
Concerning $\w e$ and $\dGamma$, we continue to use the notation of (\ref{LemM5}). We fix such a map $\w e$. There exist elements $x'\neq y'\in\eta$, $x'\neq y'$, $z\in\b\xi$, which we choose and 
fix. We let
\[
\begin{array}{lcl}
B   & := & \frac{1}{(t-1)!} \sumd{\gamma\in \dGamma(\w e,a,x',y')} \spi{a^\gamma}\eps_\gamma\cdot (\b\xi\ohne z,\eta\ohne\{ x',y'\})\c (\xi\cup\eta)
           \cdot(1 - \sumd{w\in\b\xi\ohne z} (w,z)) \\ 
D   & := & \frac{1}{t!} \sumd{\gamma\in\dGamma(\w e,a,x',y')}\;\;\sumd{w\in\b\xi}\spi{a^\gamma}\eps_\gamma\cdot (\b\xi\ohne w,\eta\ohne\{ x',y'\})
           \cdot (\xi,x')\cdot (w,y') \c (\xi\cup\eta) \\ 
    &  = & \frac{1}{t!} \sumd{\gamma\in\dGamma(\w e,a,x',y')}\;\;\sumd{w\in\b\xi}\spi{a^\gamma}\eps_\gamma\cdot (\b\xi\ohne w,\eta\ohne\{ x',y'\})
           \cdot (\xi,y')\cdot (w,x') \c (\xi\cup\eta). \\ 
\end{array}
\]
Calculations similar to those of (\ref{LemM2}) yield the following table.
\[
\begin{array}{rcl}
G'_{[a],\xi,\eta}f'_{\w e2} & = & \frac{1}{2}\Big( 4D \Big) \\
G'_{[a],\xi,\eta}f'_{\w e1} & = & \frac{1}{2}\Big( - 2B - 4D \Big) \\
G'_{[a],\xi,\eta}f'_{\w e0} & = & \frac{1}{2}\Big( B \Big) \\
\end{array}
\]
In comparison with the table of (\ref{LemM5}), we note that $s - u = 1$.
\end{Lemma}

\begin{Lemma}[cf.\ (\ref{LemM8})]
\label{LemM8short}
Suppose $g = p < k$ and $s,t\geq 2$. There exist elements $x', y'\in\eta$, $x'\neq y'$, $z\in\b\xi$, which we choose and fix. 
For $\gamma\in\dGamma(\w e 0)$, we let $x_\gamma := a_{j,\b\gamma(j,1)}$, where $j\in [g+2,k+1]$ is minimal with $1\ti j\in\Xi_\gamma$, and 
$y_\gamma := a_{j,\b\gamma(j,2)}$, where $j\in [g+2,k+1]$ is minimal with $2\ti j\in\Xi_\gamma$. I.e.\ we pick the entries $x_\gamma, y_\gamma$ that `cross the column' $g+1$ under the 
operation of $\gamma$. The set of maps
\[
\begin{array}{rcl}
[1,\lambda_1]\ohne\{ g+1\} & \lraa{\w e} & [0,2] \\
j                          & \lra        & \w e_j \\ 
\end{array}
\]
that send $g$ and $k+1$ to $\w e_g = \w e_{k+1} = 2$, and that map $j\in [1,\lambda_1]\ohne [g,k+1]$ to $e_j = 0$, is denoted by $\w E$. For $\w e\in\w E$ and $\beta\in [0,2]$, we denote by 
$\w e\beta$ the prolongation of $\w e$ to $[g,k+1]$ by $(\w e\beta)_{g+1} = \beta$. For $\gamma\in\dGamma(\w e0)$, we write 
\[
U_\gamma := \spi{a^\gamma}\eps_\gamma\cdot(\b\xi\ohne z,\eta\ohne\{ x',y'\})\cdot(x_\gamma,x')\cdot(y_\gamma,y') \\
\]
and let
\[
\begin{array}{lcl}
A_{[a],\xi,\eta,\w e}   & := & \sumd{\gamma\in\dGamma(\w e0)} U_\gamma\cdot\sumd{w'\in\b\eta} (w',z) \\
B_{[a],\xi,\eta,\w e}   & := & \sumd{\gamma\in\dGamma(\w e0)} U_\gamma\cdot(1 - \sumd{w\in\b\xi\ohne z} (w,z)) \\
C_{1,[a],\xi,\eta,\w e} & := & \sumd{\gamma\in\dGamma(\w e0)} U_\gamma\cdot (z,x_\gamma) \\
C_{2,[a],\xi,\eta,\w e} & := & \sumd{\gamma\in\dGamma(\w e0)} U_\gamma\cdot (z,y_\gamma). \\
\end{array}
\]
We obtain
\[
\begin{array}{r}
G'_{[a],\xi,\eta} f' = (X_g + 2)\sumd{\w e\in\w E}\Big(\prodd{j\in [g+2,k]} X_j^{(2-\w e_j)}\Big) \Big( 2 A_{[a],\xi,\eta,\w e} + (X_g + 1)B_{[a],\xi,\eta,\w e} \\
- X_{g+1}(C_{1,[a],\xi,\eta,\w e} + C_{2,[a],\xi,\eta,\w e})\Big).\\
\end{array}
\]

\rm
Calculations carried out in (\ref{LemM8}) yield the following table.
\[
\begin{array}{rcl}
G'_{[a],\xi,\eta} f'_{\w e2} & = & 2(s - u + 2)A_{[a],\xi,\eta,\w e} + (s - u + 2)(s - u + 1)B_{[a],\xi,\eta,\w e} \\
G'_{[a],\xi,\eta} f'_{\w e1} & = & 2A_{[a],\xi,\eta,\w e} + 2(s - u + 1)B_{[a],\xi,\eta,\w e} - (s - u + 1)(C_{1,[a],\xi,\eta,\w e} + C_{2,[a],\xi,\eta,\w e}) \\
G'_{[a],\xi,\eta} f'_{\w e0} & = & B_{[a],\xi,\eta,\w e} - (C_{1,[a],\xi,\eta,\w e} + C_{2,[a],\xi,\eta,\w e}) \\
\end{array}
\]
We note that $s - u = X_g - X_{g+1}$ and evaluate the linear combination
\[
\begin{array}{l}
\sumd{\beta\in [0,2]} X_{g+1}^{(2 - \beta)} G'_{[a],\xi,\eta} f'_{\w e\beta} \\
\begin{array}{cll}
= & 1                    & \cdot\Big( 2(X_g - X_{g+1} + 2)A_{[a],\xi,\eta,\w e} \\
  &                      & + (X_g - X_{g+1} + 2)(X_g - X_{g+1} + 1)B_{[a],\xi,\eta,\w e} \Big)  \\
+ & X_{g+1}              & \cdot\Big( 2A_{[a],\xi,\eta,\w e} + 2(X_g - X_{g+1} + 1)B_{[a],\xi,\eta,\w e} \\
  &                      & - (X_g - X_{g+1} + 1)(C_{1,[a],\xi,\eta,\w e} + C_{2,[a],\xi,\eta,\w e}) \Big) \\
+ & X_{g+1}(X_{g+1} + 1) & \cdot\Big( B_{[a],\xi,\eta,\w e} - (C_{1,[a],\xi,\eta,\w e} + C_{2,[a],\xi,\eta,\w e}) \Big)  \\
\end{array} \\
= 2(X_g + 2)A_{[a],\xi,\eta,\w e} + (X_g + 2)(X_g + 1)B_{[a],\xi,\eta,\w e} - X_{g+1}(X_g + 2)(C_{1,[a],\xi,\eta,\w e} + C_{2,[a],\xi,\eta,\w e}). \\
\end{array}
\]
\end{Lemma}

\begin{Lemma}
\label{LemM9}
Suppose $g = p < k$ and $s = \lambda'_g$, $t = 1$. Concerning $\w e$, $x_\gamma$, $y_\gamma$, we let the notation be as in (\ref{LemM8short}). For $\gamma\in\dGamma(\w e0)$, we write
\[
\begin{array}{rcl}
V_{1,\gamma} & := & \spi{a^\gamma}\eps_\gamma\cdot(y_\gamma,\eta) \\
V_{2,\gamma} & := & \spi{a^\gamma}\eps_\gamma\cdot(x_\gamma,\eta) \\
\end{array}
\]
and let
\[
\begin{array}{lcl}
A_{[a],\xi,\eta,\w e}   & := &   \sumd{\gamma\in\dGamma(\w e0)} V_{1,\gamma}\cdot \sumd{w'\in\b\eta} (x_\gamma,w') \\ 
           &  = &   \sumd{\gamma\in\dGamma(\w e0)} V_{2,\gamma}\cdot \sumd{w'\in\b\eta} (y_\gamma,w') \\ 
B_{[a],\xi,\eta,\w e}   & := &   0 \\ 
C_{1,[a],\xi,\eta,\w e} & := & - \sumd{\gamma\in\dGamma(\w e0)} V_{1,\gamma} \\ 
C_{2,[a],\xi,\eta,\w e} & := & - \sumd{\gamma\in\dGamma(\w e0)} V_{2,\gamma}. \\ 
\end{array}
\]
We obtain
\[
\begin{array}{r}
G'_{[a],\xi,\eta} f' = (X_g + 2)\sumd{\w e\in\w E}\Big(\prodd{j\in [g+2,k]} X_j^{(2-\w e_j)}\Big) \Big(2A_{[a],\xi,\eta,\w e} + (X_g + 1)B_{[a],\xi,\eta,\w e} \\
- X_{g+1}(C_{1,[a],\xi,\eta,\w e} + C_{2,[a],\xi,\eta,\w e})\Big). \\
\end{array}
\]

\rm
Calculations similar to those of (\ref{LemM8}) yield the following table, to be compared to the table of (\ref{LemM8short}).
\[
\begin{array}{rcl}
G'_{[a],\xi,\eta} f'_{\w e2} & = & 2(s - u + 2)A_{[a],\xi,\eta,\w e}  \\
G'_{[a],\xi,\eta} f'_{\w e1} & = & 2A_{[a],\xi,\eta,\w e} - (s - u + 1)(C_{1,[a],\xi,\eta,\w e} + C_{2,[a],\xi,\eta,\w e}) \\
G'_{[a],\xi,\eta} f'_{\w e0} & = & - (C_{1,[a],\xi,\eta,\w e} + C_{2,[a],\xi,\eta,\w e}) \\
\end{array}
\]
\end{Lemma}

\begin{Lemma}
\label{LemM10}
Suppose $g = p < k$ and $s = 1$, $t = \lambda'_g = \lambda'_{g+1}$. Concerning $\w e$, $x_\gamma$, $y_\gamma$, we let the notation be as in (\ref{LemM8short}). There exist elements 
$x',y'\in\eta$, $x'\neq y'$, $z\in\b\xi$, which we choose and fix. For $\gamma\in\dGamma(\w e0)$, we write
\[
U_{\gamma} := \spi{a^\gamma}\eps_\gamma\cdot(\b\xi\ohne z,\eta\ohne\{ x',y'\})\cdot(x_\gamma,x')\cdot(y_\gamma,y') 
\]
and let
\[
\begin{array}{lcl}
A_{[a],\xi,\eta,\w e}   & := & 0 \\
B_{[a],\xi,\eta,\w e}   & := & \sumd{\gamma\in\dGamma(\w e0)} U_\gamma\cdot (1 - \sumd{w\in\b\xi\ohne z} (w,z)) \\
C_{1,[a],\xi,\eta,\w e} & := & \sumd{\gamma\in\dGamma(\w e0)} U_\gamma\cdot (x_\gamma,z) \\
C_{2,[a],\xi,\eta,\w e} & := & \sumd{\gamma\in\dGamma(\w e0)} U_\gamma\cdot (y_\gamma,z). \\
\end{array}
\]
We obtain
\[
\begin{array}{r}
G'_{[a],\xi,\eta} f' = (X_g + 2)\sumd{\w e\in\w E}\Big(\prodd{j\in [g+2,k]} X_j^{(2-\w e_j)}\Big) \Big(2A_{[a],\xi,\eta,\w e} + (X_g + 1)B_{[a],\xi,\eta,\w e} \\
- X_{g+1}(C_{1,[a],\xi,\eta,\w e} + C_{2,[a],\xi,\eta,\w e})\Big). \\
\end{array}
\]

\rm
Calculations similar to those of (\ref{LemM8}) yield the following table.
\[
\begin{array}{rcl}
G'_{[a],\xi,\eta}f'_{\w e2} & = & 6B_{[a],\xi,\eta,\w e}  \\
G'_{[a],\xi,\eta}f'_{\w e1} & = & 4B_{[a],\xi,\eta,\w e} - 2(C_{1,[a],\xi,\eta,\w e} + C_{2,[a],\xi,\eta,\w e}))  \\
G'_{[a],\xi,\eta}f'_{\w e0} & = & B_{[a],\xi,\eta,\w e} - (C_{1,[a],\xi,\eta,\w e} + C_{2,[a],\xi,\eta,\w e}) \\
\end{array}
\]
In comparison to the table of (\ref{LemM8short}), we note that $X_g - X_{g+1} = 1$.
\end{Lemma}

\begin{Lemma}
\label{LemM11}
Suppose $g = p = k$ and $s\geq 1$, $t\geq 2$, $\eta = a_{g+1,[1,t]}$. Let $\gamma^1$ be the double path defined by $\gamma^1_{g+1,1} := 1$, $\gamma^1_{g+1,2} := 2$. Let $x' := a_{g+1,1}$ and 
$y' := a_{g+1,2}$, so that $x', y'\in\eta$. There exists an element $z\in\b\xi$, which we choose and fix. We write
\[
U := \spi{a^{\gamma^1}}\eps_{\gamma^1}\cdot (\bar\xi\ohne z,\eta\ohne\{ x, y\}) 
\]
and denote
\[
\begin{array}{rcl}
A & := & U \cdot \sumd{w'\in\b\eta} (w',z) \\
B & := & U \cdot (1 - \sumd{w\in\b\xi\ohne z} (w,z)) \\
\end{array}
\]
We obtain 
\[
G'_{[a],\xi,\eta} f' = \frac{1}{2}\Big((X_g + 1)(X_g + 2)B + 2(X_g + 2)A\Big),
\]
whence the map $f'$ annihilates $G'_{[a],\xi,\eta}$ modulo $X_g + 2$ in case $X_g + 2$ is odd, modulo $(X_g + 2)/2$ in case $X_g + 2$ is even.

\rm
Calculations similar to those of (\ref{LemM8}) yield 
\[
G'_{[a],\xi,\eta} f' = \frac{1}{2}\Big((s - u + 1)(s - u + 2)B + 2(s - u + 2)A\Big).
\]
Note that $s - u = X_g$.
\end{Lemma}

\begin{Lemma}
\label{LemM12}
Suppose $g = p = k$ and $s = \lambda'_g$, $t = 1$, $\eta = a_{g+1,1}$. For $i\in [2,\lambda'_{g+1}]$, we let $\gamma^i$ be the double path defined by $\gamma^i_{g+1,1} := 1$, 
$\gamma^i_{g+1,2} := i$. Let 
\[
A := \sumd{i\in [2,\lambda'_{g+1}]} \spi{a^{\gamma^i}}\eps_{\gamma^i}
\]
We obtain
\[
G'_{[a],\xi,\eta} f' = (X_g + 2)A,
\]
whence the map $f'$ annihilates $G'_{[a],\xi,\eta}$ modulo $X_g + 2$.

\rm
Calculations similar to those of (\ref{LemM8}) yield 
\[
G'_{[a],\xi,\eta} f' = (s - u + 2)A.
\]
\end{Lemma}

\begin{Remark}
\label{RemM12_5}
Suppose $g = p = k$. The application $M^{\lambda',-}\lraa{f'}S^\mu$ maps $\{(\ck a_\lambda)'\}^-$ to a linear combination with coefficients $\pm 1$ of standard 
$\mu$-polytabloids.

\rm
See Subsection \ref{SecSpecht} for the definition of the $\lambda$-tableau $[\ck a_\lambda]$.
\end{Remark}

We summarize.

\begin{Proposition}[{provisional version of (\ref{ThC16})}]
\label{PropM13}
Suppose $g < k$. The $\Z\Sl_n$-linear map $M^{\lambda',-}\lraa{f'} S^\mu$ factors over
\begin{center}
\begin{picture}(250,250)
\put( -10, 200){$M^{\lambda',-}$}
\put( 100, 210){\vector(1,0){80}}
\put( 120, 225){$\scm f'$}
\put( 200, 200){$S^\mu$}
\put(-130, 200){$\{ a'\}^-$}
\put(-100, 180){\vector(0,-1){130}}
\put(-130,   0){$\spi{a}$}
\put(  20, 180){\vector(0,-1){130}}
\put( -20, 110){$\scm\nu^\lambda_S$}
\put( 220, 180){\vector(0,-1){130}}
\put( 350, 200){$\spi{b}$}
\put( 380, 180){\vector(0,-1){130}}
\put( 350,   0){$\spi{b} + m' S^\mu$,}
\put( -10,   0){$S^\lambda$}
\put(  50,  10){\vector(1,0){130}}
\put( 110,  25){$\scm f$}
\put( 200,   0){$S^\mu/m'$}
\end{picture}
\end{center}
where $m'$ is an integer that divides the element
\[
\begin{array}{r}
G'_{[a],\xi,\eta} f' = (X_g + 2)\cdot\sumd{\w e\in\w E}\Big(\prodd{j\in [g+2,k]} X_j^{(2-\w e_j)}\Big) \Big(2A_{[a],\xi,\eta,\w e} + (X_g + 1)B_{[a],\xi,\eta,\w e} \\
- X_{g+1}(C_{1,[a],\xi,\eta,\w e} + C_{2,[a],\xi,\eta,\w e})\Big) \\
\end{array}
\]
of $S^\mu$ for any $\lambda$-tableau $[a]$ and any pair of subsets $\xi\tm a_g$, $\eta\tm a_{g+1}$ such that $\#\xi + \#\eta = \lambda'_p + 1$. The set $\w E$ is defined in (\ref{LemM8short}).
The elements $A_{[a],\xi,\eta,\w e}$, $B_{[a],\xi,\eta,\w e}$, $C_{1,[a],\xi,\eta,\w e}$, $C_{2,[a],\xi,\eta,\w e}$ of $S^\mu$ are defined in (\ref{LemM8short}, \ref{LemM9}, \ref{LemM10}), the 
definition depending on whether $\xi$ resp.\ $\eta$ contains $\geq 2$ elements. 

\rm
This is the conclusion we draw from (\ref{LemM8short}, \ref{LemM9}, \ref{LemM10}) concerning the modulus $m'$ and from 
(\ref{LemM2short}, \ref{LemM3}, \ref{LemM4}, \ref{LemM5}, \ref{LemM6}, \ref{LemM7})
concerning the vanishing of the remaining Garnir relations under $f'$. The assertion is {\it provisional} in that we have neither specified $m'$ nor commented on the 
possible divisibility of $M^{\lambda',-}\lraa{f'} S^\mu$ yet. 
\end{Proposition}

\begin{Proposition}[{to be included in (\ref{ThC16}), cf.\ (\ref{ThLW6}), cf.\ [K 99, 4.4.3]}]
\label{PropM14}
\Absatz
Suppose $g = k$. The $\Z\Sl_n$-linear map $M^{\lambda',-}\lraa{f'} S^\mu$ factors over 
\begin{center}
\begin{picture}(250,250)
\put( -10, 200){$M^{\lambda',-}$}
\put( 100, 210){\vector(1,0){80}}
\put( 120, 225){$\scm f'$}
\put( 200, 200){$S^\mu$}
\put(-130, 200){$\{ a'\}^-$}
\put(-100, 180){\vector(0,-1){130}}
\put(-130,   0){$\spi{a}$}
\put(  20, 180){\vector(0,-1){130}}
\put( -20, 110){$\scm\nu^\lambda_S$}
\put( 220, 180){\vector(0,-1){130}}
\put( 350, 200){$\spi{b}$}
\put( 380, 180){\vector(0,-1){130}}
\put( 350,   0){$\spi{b} + m S^\mu$,}
\put( -10,   0){$S^\lambda$}
\put(  50,  10){\vector(1,0){130}}
\put( 110,  25){$\scm f$}
\put( 200,   0){$S^\mu/m$}
\end{picture}
\end{center}
where $m := X_g + 2$ in case $X_g$ is odd and $m := (X_g + 2)/2$ in case $X_g$ is even. The resulting morphism $S^\lambda\lraa{f}S^\mu/m$ is of order $m$ in $\Hom_{\sZ\Sl_n}(S^\lambda,S^\mu/m)$. 

\rm
The assertions follow by (\ref{LemM11}, \ref{LemM12}, \ref{RemM12_5}).

We recall that in this case the double path formalism yields 
\[
\begin{array}{rcl}
M^{\lambda',-} & \lraa{f'} & S^\mu \\
\{ a'\}^-      & \lra      & \sumd{1\leq v < w \leq \lambda'_{g+1}}(-1)^{v+w}\spi{a^{v,w}},
\end{array}
\]
where we let
\[
\begin{array}{rcl}
[0,\lambda'_{g+1}]       & \lra        & [0,\mu'_{g+1}] \\
i                        & \lraa{\phi} & \# \Big( [1,i]\ohne\{v,w\}\Big) \\
\min(\phi^{-1}(\{ j\}))  & \llaa{\psi} & j \\
\end{array}
\]
in order to define the $\mu$-tableau $[a^{v,w}]$ by
\[
\begin{array}{rcll}
a^{v,w}_{j,i}              & := & a_{j,i}          & \mbox{\rm for } (j\in [1,\lambda_1]\ohne\{ g,g+1\}\mbox{\rm\ and } i\in [1,\mu'_j]) 
                                                     \mbox{\rm\ or } (j = g \mbox{\rm\ and } i\in [1,\lambda'_g]) \\
a^{v,w}_{g,\lambda'_g + 1} & := & a_{g+1,v}        & \\
a^{v,w}_{g,\lambda'_g + 2} & := & a_{g+1,w}        & \\
a^{v,w}_{g+1,i}            & := & a_{g+1,\psi(i)}  & \mbox{\rm for } i\in [1,\mu'_{g+1}]. \\
\end{array}
\]
\end{Proposition}        
\subsection{Morphism, reduced version}
\label{SubSecRedund}

\begin{footnotesize}
To achieve that the image of our morphism from $M^{\lambda',-}$ to $S^\mu$ gets is not contained in some $tS^\mu\tme S^\mu$, $t\geq 2$, we have to divide $f'$ by a factor of redundancy.
\end{footnotesize}

We assume $g < k$ throughout this section. For $g\leq p\leq q\leq k+1$, we let $\Xl[p,q]$ denote the set of {\it partial patterns,} i.e.\ the set of intermediate sets
\[
[1,2]\ti (\{ g,k+1\}\cap [p,q]) \tm \Xi\tm [1,2]\ti [p,q].
\]
For $p > q$ we let $\Xl[p,q] := \leer$. For $p\leq r\leq s\leq q$ and $\Xi\in\Xl[p,q]$, we let $\Xi_{[r,s]} := \Xi\cap ([1,2]\ti [r,s])\in\Xl[r,s]$. The {\it partial weight} of $\Xi$ 
is defined to be the map
\[
\begin{array}{rcl}
[p,q] & \lra & [0,2] \\
r     & \lra & e_{\Xi,r} = \# \Xi_{[r,r]}. \\
\end{array}
\]
Sometimes, we shall write this map as a tuple $e_\Xi = (e_{\Xi,p},\dots,e_{\Xi,q})$. An {\it ordered double path} of weight $e$ is a double path $\gamma$ that in addition satisfies the 
condition
\[
\b\gamma(p,1) < \b\gamma(p,2) \mbox{ for } p\in [g+1,k]. 
\]
Given a pattern $\Xi\in\Xl[g,k+1]$, we let $\dGamma(\Xi)$ be the set of double paths $\gamma$ of partial pattern $\Xi_\gamma = \Xi$, and we let $\v\Gamma(\Xi)$ be the set of ordered double 
paths $\gamma$ of partial pattern $\Xi_\gamma = \Xi$. We define $\Z\Sl_n$-linear maps
\[
\begin{array}{rcl}
F^\lambda & \lra & S^\mu \\
{[a]}     & \lraa{f''_\Xi}    & \sumd{\gamma\in\dGamma(\Xi)} \spi{a^\gamma}\eps_\gamma \\ 
{[a]}     & \lraa{\v f''_\Xi} & \sumd{\gamma\in\v\Gamma(\Xi)} \spi{a^\gamma}\eps_\gamma. \\ 
\end{array}
\]
We allow ourselves to write e.g.\ $\pazz{+}{-}{+}{+}$ for the partial pattern $\{ 1\ti u,\; 2\ti u,\; 2\ti (u+1)\}\in\Xl[u,u+1]$ etc.\ as long as
no confusion concerning $u$ can arise. So `$+$' stands for `coordinate contained' and `$-$' stands for `coordinate not contained'. 

Suppose given a $\lambda$-tableau $[a]$ and a double path $\gamma$. We let 
$x_\gamma := a_{j,\b\gamma(j,1)}$, $j\in [g+2,k+1]$ being minimal with $1\ti j\in\Xi_\gamma$, and 
$y_\gamma := a_{j,\b\gamma(j,2)}$, $j\in [g+2,k+1]$ being minimal with $2\ti j\in\Xi_\gamma$. We fix positions $u\neq v\in [1,\lambda'_{g+1}]$. For a $\lambda$-tableau
$[a]$ and a partial pattern $\Xi\in\Xl[g+2,k+1]$, we write $x' := a_{g+1,u}$, $y' := a_{g+1,v}$ and let
\[
\begin{array}{rcl}
F^\lambda & \lra & S^\mu \\
{[a]}     & \lraa{A^{u,v}_\Xi} & \sumd{\gamma\in\dGamma(\pazz{+}{-}{+}{-}\cup\Xi)} \spi{a^\gamma}\eps_\gamma\cdot (x_\gamma, x')\cdot (y_\gamma,y') \vspace*{1mm}\\ 
{[a]}     & \lraa{A^{u,-}_\Xi} & \sumd{\gamma\in\dGamma(\pazz{+}{-}{+}{-}\cup\Xi)} \spi{a^\gamma}\eps_\gamma\cdot (x_\gamma, x') \vspace*{1mm}\\
{[a]}     & \lraa{A^{-,v}_\Xi} & \sumd{\gamma\in\dGamma(\pazz{+}{-}{+}{-}\cup\Xi)} \spi{a^\gamma}\eps_\gamma\cdot (y_\gamma, y') \vspace*{1mm}\\
{[a]}     & \lraa{\v A^{u,v}_\Xi} & \sumd{\gamma\in\v\Gamma(\pazz{+}{-}{+}{-}\cup\Xi)} \spi{a^\gamma}\eps_\gamma\cdot (x_\gamma, x')\cdot (y_\gamma, y') \vspace*{1mm}\\
\end{array}
\]
\[
\begin{array}{rcl}
{[a]}     & \lraa{\v A^{u,-}_\Xi} & \sumd{\gamma\in\v\Gamma(\pazz{+}{-}{+}{-}\cup\Xi)} \spi{a^\gamma}\eps_\gamma\cdot (x_\gamma, x') \vspace*{1mm}\\ 
{[a]}     & \lraa{\v A^{-,v}_\Xi} & \sumd{\gamma\in\v\Gamma(\pazz{+}{-}{+}{-}\cup\Xi)} \spi{a^\gamma}\eps_\gamma\cdot (y_\gamma, y'). \\
\end{array}
\]
For instance, in this definition, it is understood that $\pazz{+}{-}{+}{-} = \{ 1\ti g,\; 2\ti g\}\tm [1,2]\ti [g,g+1]$. We note that $A^{u,v}_\Xi = A^{v,u}_\Xi$. 

For $i\in [0,1]$, let $L(i) := \{ j\in [g+1,k-1]\; |\; \lambda'_{j+1} = \lambda'_j - i \}$ be the {\it set of $i$-steps.} Let 
\[
[g,k] = L(0)\cup L(1)\cup \b L
\]
be a disjoint union. A partial pattern $\Xi\in\Xl[p,q]$, $g\leq p\leq q\leq k+1$, is called {\it bulky} in case there is some $u\in [p,q-1]\cap L(0)$ such that
\[
\Xi_{[u,u+1]} \in \{\pazz{+}{-}{+}{+},\;\pazz{+}{+}{+}{-},\;\pazz{+}{-}{-}{+},\;\pazz{-}{+}{+}{-},\;\pazz{+}{-}{-}{-},\;\pazz{-}{-}{+}{-},\;\pazz{+}{-}{+}{-}\},
\]
or in case there is some $u\in [p,q-1]\cap L(1)$ such that
\[
\Xi_{[u,u+1]} = \pazz{+}{-}{+}{-}.
\] 
Let $\Xlnb[p,q]\tm \Xl[p,q]$ be the subset of nonbulky partial patterns in $\Xl[p,q]$.

\begin{Lemma}
\label{LemC0}
Let $m'$ be a number dividing the elements
\[
\begin{array}{rl}
(X_g + 2)\cdot 2\cdot        & \sumd{\Xi\in\Xl[g+2,k+1]} (\prodd{j\in [g+2,k]} X_j^{(2 - e_{\Xi,j})}) A^{u,v}_\Xi \\
(X_g + 2)\cdot(X_g + 1)\cdot & \sumd{\Xi\in\Xl[g+2,k+1]} (\prodd{j\in [g+2,k]} X_j^{(2 - e_{\Xi,j})}) A^{u,v}_\Xi \\
(X_g + 2)\cdot X_{g+1}\cdot  & \sumd{\Xi\in\Xl[g+2,k+1]} (\prodd{j\in [g+2,k]} X_j^{(2 - e_{\Xi,j})}) (A^{u,-}_\Xi + A^{-,u}_\Xi) \\
\end{array}
\]
of $\Hom_{\sZ\Sl_n}(F^\lambda,S^\mu)$ for any choice of $u\neq v\in [1,\lambda'_{g+1}]$ resp.\ of $u\in [1,\lambda'_{g+1}]$. Then there is a factorization
\[
(M^{\lambda',-}\lraa{f'}S^\mu \lra S^\mu/m') = (M^{\lambda',-}\lraa{\nu_S^{\lambda}}S^\lambda \lra S^\mu/m').
\]

\rm
Let $[a]$ be a $\lambda$-tableau, let $\xi\tm a_g$, $\eta\tm a_{g+1}$ such that $\#\xi + \#\eta = \lambda'_p + 1$. Let $\w e$ be a map from $[1,\lambda_1]\ohne\{ g+1\}$ to $[0,2]$ that maps
$g$ and $k+1$ to $e_g = e_{k+1} = 2$, and that maps $j\in [1,\lambda_1]\ohne [g,k+1]$ to $e_j = 0$. 

{\it Case $\#\xi\geq 2$, $\#\eta\geq 2$.} Writing $a_{g+1,u} := x'$ and $a_{g+1,v} := y'$, we obtain, in the notation of (\ref{LemM8short}), including choices invoved,
\[
\begin{array}{rcl}
A_{[a],\xi,\eta,\w e}   & = & \sumd{\Xi\in\Xl[g+2,k+1],\; e_\Xi\; =\; \w e|_{[g+2,k+1]}} \Big([a] \cdot (\b\xi\ohne z,\eta\ohne\{ x,y\})\cdot (\sumd{w'\in\b\eta} (w',z))\Big) A^{u,v}_\Xi \\
B_{[a],\xi,\eta,\w e}   & = & \sumd{\Xi\in\Xl[g+2,k+1],\; e_\Xi\; =\; \w e|_{[g+2,k+1]}} \Big([a] \cdot (\b\xi\ohne z,\eta\ohne\{ x,y\})\cdot (1-\sumd{w\in\b\xi\ohne z} (w,z))\Big) A^{u,v}_\Xi\\
C_{1,[a],\xi,\eta,\w e} & = & \sumd{\Xi\in\Xl[g+2,k+1],\; e_\Xi\; =\; \w e|_{[g+2,k+1]}} \Big(- [a] \cdot (\b\xi\ohne z,\eta\ohne\{ x,y\})\cdot (x',z)\Big) A^{-,v}_\Xi \\
                        & = & \sumd{\Xi\in\Xl[g+2,k+1],\; e_\Xi\; =\; \w e|_{[g+2,k+1]}} \Big(- [a] \cdot (\b\xi\ohne z,\eta\ohne\{ x,y\})\cdot (y',z)\Big) A^{-,u}_\Xi \\
C_{2,[a],\xi,\eta,\w e} & = & \sumd{\Xi\in\Xl[g+2,k+1],\; e_\Xi\; =\; \w e|_{[g+2,k+1]}} \Big(- [a] \cdot (\b\xi\ohne z,\eta\ohne\{ x,y\})\cdot (y',z)\Big) A^{u,-}_\Xi, \\
\end{array}
\] 
where the multitransposition $(\b\xi\ohne z,\eta\ohne\{ x,y\})$ is to be read as with respect to and independent of a choice of a bijection (cf.\ Section \ref{SecGar}).

{\it Case $\#\xi = \lambda'_g$, $\#\eta = 1$.}  Writing $\{ a_{g+1,u}\} := \eta$, we obtain, in the notation of (\ref{LemM9}),
\[
\begin{array}{rcl}
A_{[a],\xi,\eta,\w e}   & = & \sumd{\Xi\in\Xl[g+2,k+1],\; e_\Xi\; =\;\w e|_{[g+2,k+1]}}\;\;\sumd{v\in [1,\lambda'_{g+1}],\; a_{g+2,v}\in\b\eta}\Big([a]\Big) A^{u,v}_\Xi \\
B_{[a],\xi,\eta,\w e}   & = & 0 \\
C_{1,[a],\xi,\eta,\w e} & = & \sumd{\Xi\in\Xl[g+2,k+1],\; e_\Xi\; =\;\w e|_{[g+2,k+1]}} \Big(-[a]\Big) A^{-,u}_\Xi \\
C_{2,[a],\xi,\eta,\w e} & = & \sumd{\Xi\in\Xl[g+2,k+1],\; e_\Xi\; =\;\w e|_{[g+2,k+1]}} \Big(-[a]\Big) A^{u,-}_\Xi. \\
\end{array}
\] 

{\it Case $\#\xi = 1$, $\#\eta = \lambda'_g = \lambda'_{g+1}$.} Writing $a_{g+1,u} := x'$ and $a_{g+1,v} := y'$, we obtain, in the notation of (\ref{LemM10}),
\[
\begin{array}{rcl}
A_{[a],\xi,\eta,\w e}   & = & 0 \\
B_{[a],\xi,\eta,\w e}   & = & \sumd{\Xi\in\Xl[g+2,k+1],\; e_\Xi\; =\; \w e|_{[g+2,k+1]}} \Big([a] \cdot (\b\xi\ohne z,\eta\ohne\{ x,y\})\cdot (1-\sumd{w\in\b\xi\ohne z} (w,z))\Big) A^{u,v}_\Xi\\
C_{1,[a],\xi,\eta,\w e} & = & \sumd{\Xi\in\Xl[g+2,k+1],\; e_\Xi\; =\; \w e|_{[g+2,k+1]}} \Big(- [a] \cdot (\b\xi\ohne z,\eta\ohne\{ x,y\})\cdot (x',z)\Big) A^{-,v}_\Xi \\
                        & = & \sumd{\Xi\in\Xl[g+2,k+1],\; e_\Xi\; =\; \w e|_{[g+2,k+1]}} \Big(- [a] \cdot (\b\xi\ohne z,\eta\ohne\{ x,y\})\cdot (y',z)\Big) A^{-,u}_\Xi \\
C_{2,[a],\xi,\eta,\w e} & = & \sumd{\Xi\in\Xl[g+2,k+1],\; e_\Xi\; =\; \w e|_{[g+2,k+1]}} \Big(- [a] \cdot (\b\xi\ohne z,\eta\ohne\{ x,y\})\cdot (y',z)\Big) A^{u,-}_\Xi, \\
\end{array}
\] 
where $(\b\xi\ohne z,\eta\ohne\{ x,y\})$ is to be read as with respect to and independent of a choice of a bijection.

The lemma now follows from (\ref{PropM13}). 
\end{Lemma}

\begin{Lemma}
\label{LemC5}
Let $\Xi\in\Xlnb[g,k+1]$. The image
\[
[\ck a_\lambda] \v f''_\Xi
\]
is a linear combination of standard $\mu$-polytabloids with coefficients $\pm 1$. For a pair of different nonbulky patterns, the corresponding pair of sets of occurring standard 
$\mu$-polytabloids is disjoint.

\rm
The summand $\pm\spi{(\ck a_\lambda)^\gamma}$ is a standard polytabloid after having ordered its columns downwards increasingly provided $\Xi_\gamma$ is not bulky.
Different ordered double paths yield different fillings of the columns, hence the coefficients remain $\pm 1$ and, moreover, the asserted disjointness holds.
\end{Lemma}

For $g \leq p\leq q\leq k$, a tuple of integers $(\theta_\Xi)_{\Xi\in\Xl[p,q]}$ is called a {\it reduced tuple of coefficients} if the following conditions hold, expressed using the 
involution introduced in (\ref{LemM0}).
\begin{itemize}
\item[$(\mbox{\rm I}_{p,q})$]  For $\Xi\in\Xl[p,q]$, the invariance $\theta_{(\Xi)\iota_{q+1}} = \theta_\Xi$ holds. Given, in addition, $r\in [p,q]$ with $e_{\Xi,r} = 2$, the invariance
                               $\theta_{(\Xi)\iota_{r+1}} = \theta_\Xi$ holds.
\item[$(\mbox{\rm II}_{p,q})$] In case $\Xi\in\Xl[p,q]$ is bulky, the coefficient $\theta_\Xi$ vanishes, in case not, it does not vanish. 
\end{itemize}
We note that the tuple $(\prodd{j\in [p,q]} X_j^{(2-e_{\Xi,j})})_{\Xi\in\Xl[g,p]}$ satisfies $(\mbox{\rm I}_{p,q})$ but in general not $(\mbox{\rm II}_{p,q})$. Let it be remarked that
e.g.\ $(\mbox{\rm II}_{p,p+1})$ implies that $\theta_{\pazz{+}{-}{-}{+}} \neq \theta_{\pazz{+}{+}{-}{-}}$ if $p\in [g+1,k-1]$ and $\lambda'_p = \lambda'_{p+1}$, whereas 
$e_{\pazz{+}{-}{-}{+}} \neq e_{\pazz{+}{+}{-}{-}}$.

\begin{Lemma}
\label{LemC4}
Let $p\in [g+1,k]$, let $(\theta_\Xi)_{\Xi\in\Xl[g,p]}$ (resp.\ $(\theta_\Xi)_{\Xi\in\Xl[g+2,p]}$) be a tuple of integers satisfying $(\mbox{\rm I}_{g,p})$ (resp.\ $(\mbox{\rm I}_{g+2,p})$), let 
$\Phi\in\Xl[p+1,k+1]$. Then
\[
\begin{array}{rcl}
\sumd{\Xi\in\Xl[g,p]} \theta_\Xi f''_{\Xi\cup\Phi}       & = & \sumd{\Xi\in\Xl[g,p]} \Big(\prod_{j\in [g+1,p]} e_{\Xi,j} !\Big)\theta_\Xi \v f''_{\Xi\cup\Phi} \\ 
\sumd{\Xi\in\Xl[g+2,p]} \theta_\Xi A^{u,v}_{\Xi\cup\Phi} & = & \sumd{\Xi\in\Xl[g+2,p]}\Big(\prod_{j\in [g+2,p]} e_{\Xi,j} !\Big) \theta_\Xi \v A^{u,v}_{\Xi\cup\Phi}\; , \\ 
\end{array}
\]
respectively, and likewise for $A^{u,-}$, $A^{-,v}$ instead of $A^{u,v}$.

\rm
This follows, separately for each partial weight $e_\Xi$, from (\ref{LemM0}, ii).
\end{Lemma}

\begin{Lemma}
\label{LemC1}
Let $p\in L(0)$, $\Psi\in\Xl[g,p-1]$, $\Phi\in\Xl[p+2,k+1]$. We obtain 
\[
\begin{array}{rrcl}
\mbox{\rm (i)}      & f''_{\Psi\cup\pazz{+}{+}{+}{+}\cup\Phi} - 2 f''_{\Psi\cup\pazz{+}{+}{+}{-}\cup\Phi} - 2 f''_{\Psi\cup\pazz{+}{-}{+}{+}\cup\Phi}
                      + 2 f''_{\Psi\cup\pazz{+}{-}{+}{-}\cup\Phi} & = & 0 \\
\mbox{\rm (ii)}     & f''_{\Psi\cup\pazz{+}{+}{+}{+}\cup\Phi} - f''_{\Psi\cup\pazz{+}{+}{+}{-}\cup\Phi} - f''_{\Psi\cup\pazz{+}{-}{+}{+}\cup\Phi} & = & 0 \\
\mbox{\rm (iii)}    & f''_{\Psi\cup\pazz{+}{+}{-}{+}\cup\Phi} - f''_{\Psi\cup\pazz{+}{-}{-}{+}\cup\Phi} - f''_{(\Psi)\iota_p\cup\pazz{-}{+}{+}{-}\cup\Phi} & = & 0 \\
\mbox{\rm (iii$'$)} & f''_{\Psi\cup\pazz{-}{+}{+}{+}\cup\Phi} - f''_{(\Psi)\iota_p\cup\pazz{+}{-}{-}{+}\cup\Phi} - f''_{\Psi\cup\pazz{-}{+}{+}{-}\cup\Phi} & = & 0 \\
\mbox{\rm (iv)}     & f''_{\Psi\cup\pazz{+}{+}{-}{-}\cup\Phi} - f''_{\Psi\cup\pazz{+}{-}{-}{-}\cup\Phi} & = & 0 \\
\mbox{\rm (iv$'$)}  & f''_{\Psi\cup\pazz{-}{-}{+}{+}\cup\Phi} - f''_{\Psi\cup\pazz{-}{-}{+}{-}\cup\Phi} & = & 0 \\
\mbox{\rm (v)}      & 2 f''_{\Psi\cup\pazz{+}{-}{+}{-}\cup\Phi} - f''_{\Psi\cup\pazz{+}{+}{+}{+}\cup\Phi} & = & 0. \\
\end{array}
\]
Equation {\rm (i)} still holds in case $p\in L(1)$. Equations {\rm (i-v)} continue to hold for $p\in L(0)\cap [g+2,k-1]$, $f''$ replaced by $A^{u,v}$, $A^{u,-}$ or $A^{-,v}$ and
$\Psi\in\Xl[g,p-1]$ replaced by $\Psi\in\Xl[g+2,p-1]$. 

\rm
Equation (i) ensues from the Garnir relation that arises, for a given double path $\gamma\in\dGamma(\Xi\cup\pazz{+}{+}{+}{+}\cup\Phi)$, from the join of the elements in column $p$ in the 
image of $\gamma$ with column $p+1$.

We deduce (ii) from the Garnir relation that arises, for a given double path $\gamma\in\dGamma(\Xi\cup\pazz{+}{+}{+}{+}\cup\Phi)$, from the join of a single element in column $p$ in the 
image of $\gamma$ with column $p+1$, using (\ref{LemM0} ii) to obtain 
$f''_{(\Psi)\iota_p\cup\pazz{+}{+}{+}{-}\cup\Phi} = f''_{(\Psi\cup\pazz{+}{+}{+}{-}\cup\Phi)\iota_{p+1}} = f''_{\Psi\cup\pazz{+}{+}{+}{-}\cup\Phi}$.

Assertion (iii) follows using the Garnir relation that arises, for a given double path $\gamma\in\dGamma(\Xi\cup\pazz{+}{+}{-}{+}\cup\Phi)$, from the join of the element in column $p$ 
in the image of $\gamma$ with column $p+1$.

Similarly, (iv) can be seen via the Garnir relation that arises, for a given double path $\gamma\in\dGamma(\Xi\cup\pazz{+}{+}{-}{-}\cup\Phi)$, to the join of the element in column $p$ in 
the image of $\gamma$ with column $p+1$.

As a consequence of (i) and (ii), we obtain (v).
\end{Lemma}

\begin{Lemma}
\label{LemC1_5}
Suppose $g+1\in L(0)$, $\Phi\in\Xl[g+3,k+1]$. We obtain
\[
\begin{array}{rrcl}
\mbox{\rm (i)}      & A^{u,v}_{\paze{+}{+}\cup\Phi} - 2 A^{u,v}_{\paze{+}{-}\cup\Phi} - 2 A^{u,v}_{\paze{-}{+}\cup\Phi}
                      + 2 A^{u,v}_{\paze{-}{-}\cup\Phi} & = & 0 \vspace*{2mm}\\
\mbox{\rm (ii)}     & A^{u,v}_{\paze{+}{+}\cup\Phi} - A^{u,v}_{\paze{+}{-}\cup\Phi} - A^{u,v}_{\paze{-}{+}\cup\Phi} & = & 0  \vspace*{2mm}\\
\mbox{\rm (iii)}    & A^{u,-}_{\paze{+}{+}\cup\Phi} - A^{u,-}_{\paze{-}{+}\cup\Phi} - A^{-,u}_{\paze{+}{-}\cup\Phi} & = & 0  \vspace*{2mm}\\
\mbox{\rm (iii$'$)} & A^{-,v}_{\paze{+}{+}\cup\Phi} - A^{v,-}_{\paze{-}{+}\cup\Phi} - A^{-,v}_{\paze{+}{-}\cup\Phi} & = & 0 \vspace*{2mm} \\
\mbox{\rm (iv)}     & A^{u,-}_{\paze{+}{-}\cup\Phi} - A^{u,-}_{\paze{-}{-}\cup\Phi} & = & 0  \vspace*{2mm}\\
\mbox{\rm (iv$'$)}  & A^{-,v}_{\paze{-}{+}\cup\Phi} - A^{-,v}_{\paze{-}{-}\cup\Phi} & = & 0  \vspace*{2mm}\\
\mbox{\rm (v)}      & 2 A^{u,v}_{\paze{-}{-}\cup\Phi} - A^{u,v}_{\paze{+}{+}\cup\Phi} & = & 0.  \\
\end{array}
\]
Equation {\rm (i)} still holds in case $g+1\in L(1)$.
\end{Lemma}

\begin{Lemma}
\label{LemC2}
Let $p\in L(0)$, let $(\theta_\Xi)_{\Xi\in\Xl[g,p]}$ be a reduced tuple of coefficients, let $\Phi\in\Xl[p+2,k+1]$. There is a reduced tuple of coefficients 
$(\h\theta_{\h\Xi})_{\h\Xi\in\Xl[g,p+1]}$ such that
\[
\sumd{\Xi\in\Xl[g,p]}\;\; \sumd{\xi\in\Xl[p+1,p+1]} \theta_\Xi \cdot  \#\xi! \cdot X_{p+1}^{(2-\#\xi)} \v f''_{\Xi\cup\xi\cup\Phi}
= \sumd{\h\Xi\in\Xl[g,p+1]} \h\theta_{\h\Xi} \v f''_{\h\Xi\cup\Phi}
\]
Fixing a partial pattern $\ck\Xi\in\Xl[g,p-1]$, $\h\theta$ can be obtained from $\theta$ letting
\[
\begin{array}{rcr}
\h\theta_{\ck\Xi\cup\pazz{+}{+}{+}{+}} & = & X_p(X_p + 1)\theta_{\ck\Xi\cup\paze{+}{+}} \\
\h\theta_{\ck\Xi\cup\pazz{-}{+}{+}{+}} & = & (X_p + 1)\theta_{\ck\Xi\cup\paze{-}{+}}\\
\h\theta_{\ck\Xi\cup\pazz{+}{+}{-}{+}} & = & (X_p + 1)\theta_{\ck\Xi\cup\paze{+}{-}} \\
\h\theta_{\ck\Xi\cup\pazz{-}{-}{+}{+}} & = & (X_p - 1)(X_p + 1)\theta_{\ck\Xi\cup\paze{-}{+}}\\
\h\theta_{\ck\Xi\cup\pazz{+}{+}{-}{-}} & = & (X_p - 1)(X_p + 1)\theta_{\ck\Xi\cup\paze{+}{-}}\\
\h\theta_{\ck\Xi\cup\pazz{-}{+}{-}{+}} & = & 2\theta_{\ck\Xi\cup\paze{-}{-}}\\
\h\theta_{\ck\Xi\cup\pazz{-}{+}{-}{-}} \;\; =\;\; \h\theta_{\ck\Xi\cup\pazz{-}{-}{-}{+}} & = & (X_p - 1)\theta_{\ck\Xi\cup\paze{-}{-}} \\
\h\theta_{\ck\Xi\cup\pazz{-}{-}{-}{-}} & = & (X_p - 1)X_p\theta_{\ck\Xi\cup\paze{-}{-}}. \\
\end{array}
\]
In case $p\in L(0)\cap [g+2,k+1]$, the assertion holds for $\v A^{u,v}$, $\v A^{u,-}$ or $\v A^{-,v}$ instead of $\v f''$, the reduced tuples of coefficients being indexed by 
$\Xi\in\Xl[g+2,p]$ on the left hand side and by $\h\Xi\in\Xl[g+2,p+1]$ on the right hand side.

\rm
Up to a rescaling, we rewrite summands of the left hand side using (\ref{LemC1}), to which the roman numbers refer. To this end, we substitute $X_{p+1} = X_p - 1$ and suppose 
$\ck\Xi\in\Xl[g,p-1]$ be given. 
\[
\begin{array}{rcl}
\theta_{\ck\Xi\cup\paze{+}{+}}\Big((X_p - 1)X_p f''_{\ck\Xi\cup\pazz{+}{-}{+}{-}\cup\Phi} \\
+ (X_p - 1) (f''_{\ck\Xi\cup\pazz{+}{+}{+}{-}\cup\Phi} + f''_{\ck\Xi\cup\pazz{+}{-}{+}{+}\cup\Phi}) 
+ f''_{\ck\Xi\cup\pazz{+}{+}{+}{+}\cup\Phi} \Big) 
& \auf{\mbox{\scr (ii, v)}}{=} & \theta_{\ck\Xi\cup\paze{+}{+}} \cdot\frac{1}{2} X_p(X_p + 1) f''_{\ck\Xi\cup\pazz{+}{+}{+}{+}\cup\Phi} \\
\theta_{\ck\Xi\cup\paze{+}{-}} \Big((X_p - 1)X_p f''_{\ck\Xi\cup\pazz{+}{-}{-}{-}\cup\Phi} + (X_p - 1)f''_{\ck\Xi\cup\pazz{+}{+}{-}{-}\cup\Phi}\Big) 
& \auf{\mbox{\scr (vi)}}{=} & \theta_{\ck\Xi\cup\paze{+}{-}} (X_p - 1)(X_p + 1) f''_{\ck\Xi\cup\pazz{+}{+}{-}{-}\cup\Phi} \\
\theta_{\ck\Xi\cup\paze{+}{-}} \Big( \frac{1}{2}(X_p - 1)(f''_{\ck\Xi\cup\pazz{+}{-}{-}{+}\cup\Phi} + f''_{(\ck\Xi)\iota_p\cup\pazz{-}{+}{+}{-}\cup\Phi}) & & \\
+ f''_{\ck\Xi\cup\pazz{+}{+}{-}{+}\cup\Phi}\Big) 
& \auf{\mbox{\scr (iii)}}{=} & \theta_{\ck\Xi\cup\paze{+}{-}} \cdot\frac{1}{2} (X_p + 1) f''_{\ck\Xi\cup\pazz{+}{+}{-}{+}\cup\Phi} \\
\end{array}
\]
For the last equation we made use of $(\mbox{\rm I}_{g,p})$ in order to have $\theta_{\ck\Xi\cup\paze{+}{-}} = \theta_{(\ck\Xi)\iota_p\cup\paze{-}{+}}$.

The rewriting of the summands having $\paze{-}{+}$ in column $p$ proceeds by symmetry. Summands of index $\Xi = \ck\Xi\cup\paze{-}{-}$ remain unchanged. We pass to $\v f''$ by inserting a 
factor $2$ where necessary. So we obtain an array of equations whose left hand sides add up to yield the left hand side of the equation we claimed. After substitution of $\theta$ by 
$\h\theta$ according to the table above, we sum up the corresponding right hand sides to obtain the right hand side of the equation we claimed.

The arguments are valid verbatim for $\v A^{u,v}$, $\v A^{u,-}$, $\v A^{-,v}$ instead of $\v f''$, except that we replace $\ck\Xi\in\Xl[g,p-1]$ by $\ck\Xi\in\Xl[g+2,p-1]$ and need to suppose 
that $p\geq g+2$.
\end{Lemma}

\begin{Lemma}
\label{LemC2_5}
Suppose $g+1\in L(0)$, let $\Phi\in\Xl[g+3,k+1]$. We obtain
\[
\begin{array}{rcl}
\sumd{\xi\in\Xl[g+2,g+2]} \#\xi!\cdot X_{g+2}^{(2-\#\xi)} \cdot \v A^{u,v}_{\xi\cup\Phi} 
& = & X_{g+1}(X_{g+1} + 1) \v A^{u,v}_{\paze{+}{+}\cup\Phi} \\ 
\sumd{\xi\in\Xl[g+2,g+2]} \#\xi!\cdot X_{g+2}^{(2-\#\xi)} \cdot (\v A^{u,-}_{\xi\cup\Phi} + \v A^{-,u}_{\xi\cup\Phi}) 
& = & (X_{g+1} + 1) (\v A^{u,-}_{\paze{+}{+}\cup\Phi} + \v A^{-,u}_{\paze{+}{+}\cup\Phi}) \\
& + & (X_{g+1} - 1)(X_{g+1} + 1) (\v A^{u,-}_{\paze{+}{-}\cup\Phi} + \v A^{-,u}_{\paze{-}{+}\cup\Phi}) \\ 
\end{array}
\]

\rm
This results from (\ref{LemC1_5}) by the argument in (\ref{LemC2}).
\end{Lemma}

\begin{Lemma}
\label{LemC3}
Let $p\in L(1)$, let $(\theta_\Xi)_{\Xi\in\Xl[g,p]}$ be a reduced tuple of coefficients, let $\Phi\in\Xl[p+2,k+1]$. There is a reduced tuple of coefficients 
$(\h\theta_{\h\Xi})_{\h\Xi\in\Xl[g,p+1]}$ such that
\[
\sumd{\Xi\in\Xl[g,p]}\;\; \sumd{\xi\in\Xl[p+1,p+1]} \theta_\Xi\cdot \#\xi !\cdot X_{p+1}^{(2-\#\xi)} \v f''_{\Xi\cup\xi\cup\Phi}
= \sumd{\h\Xi\in\Xl[g,p+1]} \h\theta_{\h\Xi} \v f''_{\h\Xi\cup\Phi} 
\]
Fixing a partial pattern $\ck\Xi\in\Xl[g,p-1]$, $\theta$ and $\h\theta$ are related by
\[
\begin{array}{rcr}
\h\theta_{\ck\Xi\cup\pazz{+}{+}{+}{+}} & = & - X_p(X_p - 3) \theta_{\ck\Xi\cup\paze{+}{+}} \\ 
\h\theta_{\ck\Xi\cup\pazz{+}{-}{+}{+}} \;\; = \;\; \h\theta_{\ck\Xi\cup\pazz{+}{+}{+}{-}} & = & X_p (X_p - 2) \theta_{\ck\Xi\cup\paze{+}{+}} \\
\h\theta_{\ck\Xi\cup\pazz{+}{+}{-}{+}} & = &  2\theta_{\ck\Xi\cup\paze{+}{-}} \\
\h\theta_{\ck\Xi\cup\pazz{+}{-}{-}{+}} \;\; = \;\;\h\theta_{\ck\Xi\cup\pazz{+}{+}{-}{-}} & = & (X_p - 2)\theta_{\ck\Xi\cup\paze{+}{-}}\\
\h\theta_{\ck\Xi\cup\pazz{+}{-}{-}{-}} & = & (X_p - 2)(X_p - 1) \theta_{\ck\Xi\cup\paze{+}{-}} \\
\h\theta_{\ck\Xi\cup\pazz{-}{+}{+}{+}} & = & 2\theta_{\ck\Xi\cup\paze{-}{+}} \vspace*{1mm}\\
\end{array}
\]
\[
\begin{array}{rcr}
\h\theta_{\ck\Xi\cup\pazz{-}{-}{+}{+}} \;\; = \;\;\h\theta_{\ck\Xi\cup\pazz{-}{+}{+}{-}} & = & (X_p - 2)\theta_{\ck\Xi\cup\paze{-}{+}} \\
\h\theta_{\ck\Xi\cup\pazz{-}{-}{+}{-}} & = & (X_p - 2)(X_p - 1)\theta_{\ck\Xi\cup\paze{-}{+}} \\
\h\theta_{\ck\Xi\cup\pazz{-}{+}{-}{+}} & = & 2\theta_{\ck\Xi\cup\paze{-}{-}} \vspace*{1mm}\\
\h\theta_{\ck\Xi\cup\pazz{-}{-}{-}{+}} \;\; = \;\;\h\theta_{\ck\Xi\cup\pazz{-}{+}{-}{-}} & = & (X_p - 2)\theta_{\ck\Xi\cup\paze{-}{-}} \\
\h\theta_{\ck\Xi\cup\pazz{-}{-}{-}{-}} & = & (X_p - 2)(X_p - 1)\theta_{\ck\Xi\cup\paze{-}{-}} \\
\end{array}
\]
In case $p\in L(1)\cap [g+2,k+1]$, the assertion holds for $\v A^{u,v}$, $\v A^{u,-}$ or $\v A^{-,v}$ instead of $f''$, the reduced tuples of coefficients being indexed by $\Xi\in\Xl[g+2,p]$ 
on the left and by $\h\Xi\in\Xl[g+2,p+1]$ on the right hand side. 

\rm
We rewrite summands of the left hand side using (\ref{LemC1} i). To this end, we substitute $X_{p+1} = X_p - 2$ and suppose $\ck\Xi\in\Xl[g,p-1]$ be given.
\[
\begin{array}{rcl}
\theta_{\ck\Xi\cup\paze{+}{+}} \Big((X_p - 2)(X_p - 1) f''_{\ck\Xi\cup\pazz{+}{-}{+}{-}\cup\Phi} \\
+ (X_p - 2)(f''_{\ck\Xi\cup\pazz{+}{+}{+}{-}\cup\Phi} + f''_{\ck\Xi\cup\pazz{+}{-}{+}{+}\cup\Phi}) + f''_{\ck\Xi\cup\pazz{+}{+}{+}{+}\cup\Phi}\Big) 
& \auf{\mbox{\scr (i)}}{=} & \theta_{\ck\Xi\cup\paze{+}{+}} \Big((X_p - 2)X_p (f''_{\ck\Xi\cup\pazz{+}{+}{+}{-}\cup\Phi} + f''_{\ck\Xi\cup\pazz{+}{-}{+}{+}\cup\Phi})\vspace*{1mm}\\
&                          &                                      - \frac{1}{2} (X_p - 3)X_p f''_{\ck\Xi\cup\pazz{+}{+}{+}{+}\cup\Phi}\Big) 
\end{array}
\]
All the other summands remain unchanged. 
\end{Lemma}

\begin{Lemma}
\label{LemC3_5}
Suppose $g+1\in L(1)$, let $\Phi\in\Xl[g+3,k+1]$. We obtain
\[
\sumd{\xi\in\Xl[g+2,g+2]} \#\xi!\cdot X_{g+2}^{(2-\#\xi)} \cdot \v A^{u,v}_{\xi\cup\Phi} 
= - X_{g+1}(X_{g+1} - 3) \v A^{u,v}_{\paze{+}{+}\cup\Phi} + X_{g+1}(X_{g+1} - 2) (\v A^{u,v}_{\paze{+}{-}\cup\Phi} + \v A^{u,v}_{\paze{-}{+}\cup\Phi}).
\]

\rm
This results from (\ref{LemC1_5}), cf.\ (\ref{LemC3}).
\end{Lemma}

{\bf Henceforth,} we denote by $\theta$ {\it the} tuple of tuples $\theta = ((\theta_\Xi)_{\Xi\in \Xl[g,r]})_{r\in [g,k]}$ of integers such that the entry 
$\theta_{\paze{+}{+}} = 1$\vspace{1mm} furnishes
the tuple in case $r = g$, such that each entry $(\theta_\Xi)_{\Xi\in \Xl[g,r]}$ is a reduced tuple of coefficients, such that for $r\in [g,k-1]$ the tuples $(\theta_\Xi)_{\Xi\in \Xl[g,r]}$ 
and $(\theta_\Xi)_{\Xi\in \Xl[g,r+1]}$ are related by the formulas in (\ref{LemC2}) (resp.\ \ref{LemC3}) in case $r\in L(0)$ (resp.\ $r\in L(1)$), and such that in the remaining case 
$r\in\b L$ the equation
\[
\theta_{\Xi\cup\xi} = \#\xi! \cdot X_{r+1}^{(2 - \#\xi)}\cdot \theta_\Xi
\]
holds `unchanged' for $\Xi\in\Xl[g,r]$, $\xi\in\Xl[r+1,r+1]$. We reindex the tuple for $r = k$ by $\theta_\Xi := \theta_{\Xi_{[g,k]}}$ for $\Xi\in\Xl[g,k+1]$. Furthermore, 
for $r,s\in [g+1,k]$, we let 
\[
R_{[r,s]} := \prodd{i\in [0,1]}\;\;\prodd{j\in L(i)\cap [r,s]} X_j^{(2 - i)},
\] 
in which a product is to be read to be equal to $1$ if no factor is present. For $r\in\Z$, we write 
\[
\dot r := \gcd(r,2).
\]

\begin{Lemma}
\label{LemC6}
Let $g+1\leq r\leq q\leq k$ such that $r - 1\in\b L$, let $\Psi\in\Xl[g,r-1]$, $\Xi\in\Xl[r,q]$ such that $\Psi\cup\Xi$ is nonbulky. Then $\theta_\Psi R_{[r,q-1]}$ divides 
$\theta_{\Psi\cup\Xi}$. If, moreover, $\Xi_{[r,r]} = \paze{+}{+}$, then $2\theta_\Xi R_{[r,q-1]}$ divides $\theta_{\Psi\cup\Xi}$.

\rm
We perform an induction on $q$ to prove the assertion that $\theta_\Psi R_{[r,q-1]}$ divides $\theta_{\Psi\cup\Xi}$ and that the quotient $\theta_{\Psi\cup\Xi}/(\theta_\Psi R_{[r,q-1]})$
is divisible by $X_q^{(2 - e_{\Xi,q})}$. The factor $2$ in the second assertion then follows from the start of the induction.

At $q = r$, we have $R_{[r,r-1]} = 1$, and the quotient $\theta_{\Psi\cup\Xi}/(\theta_\Psi R_{[r,q-1]})$ equals $\#\Xi!\cdot X_r^{(2-\#\Xi)}$ since $r-1\in\b L$.

We assume our assertion to hold for $q - 1\geq r$. Let $\ck\Xi := \Xi_{[r,q-1]}$.

{\it Case $q - 1\in \b L$.} We have $\theta_{\Psi\cup\Xi}/\theta_{\Psi\cup\ck\Xi} = \#\Xi_{[q,q]}!\cdot X_q^{(2 - e_{\Xi,q})}$ and $R_{[r,q-1]} = R_{[r,q-2]}$.

{\it Case $q - 1\in L(0)$.} We have $X_q = X_{q-1} - 1$ and $R_{[r,q-1]} = X_{q-1}(X_{q-1}+1)\cdot R_{[r,q-2]}$. By assumption, $\theta_{\Psi\cup\ck\Xi}/(R_{[r,q-2]}\theta_\Psi)$ is integral
and divisible by $X_{q-1}^{(2 - e_{\Xi,q-1})}$. We need to see that 
\[
(\theta_{\Psi\cup\Xi}/\theta_{\Psi\cup\ck\Xi})\cdot (\theta_{\Psi\cup\ck\Xi}/(R_{[r,q-2]}\theta_\Psi))
\]
is divisible by $X_{q-1}(X_{q-1} + 1)\cdot X_q^{(2 - e_{\Xi,q})}$. The left hand side factor of that expression taken from the table in (\ref{LemC2}), the assertion follows separately for 
each nonbulky subset $\Xi_{[q-1,q]}$.

{\it Case $q - 1\in L(1)$.} We have $X_q = X_{q-1} - 2$ and $R_{[r,q-1]} = X_{q-1}\cdot R_{[r,q-2]}$. By assumption, $\theta_{\Psi\cup\ck\Xi}/(R_{[r,q-2]}\theta_\Psi)$ is integral
and divisible by $X_{q-1}^{(2 - e_{\Xi,q-1})}$. We need to see that 
\[
(\theta_{\Psi\cup\Xi}/\theta_{\Psi\cup\ck\Xi})\cdot (\theta_{\Psi\cup\ck\Xi}/(R_{[r,q-2]}\theta_\Psi))
\]
is divisible by $X_{q-1}\cdot X_q^{(2 - e_{\Xi,q})}$. The left hand side factor of that expression taken from the table in (\ref{LemC3}), the assertion follows separately for each
nonbulky subset $\Xi_{[q-1,q]}$.
\end{Lemma}

\begin{Lemma}
\label{LemC7}
Let $g+1\leq r\leq q\leq k$ such that $r - 1\in\b L$ and $[r,q-1]\tm L(1)$, let $\Psi\in\Xlnb[g,r-1]$. As ideal of coefficients, we obtain
\[
\Big(\theta_{\Psi\cup\Xi} \; |\; \Xi\in\Xl[r,q] \Big) = \Big(\theta_\Psi\cdot R_{[r,q-1]}\cdot \dot X_q\Big)\tm\Z.
\]

\rm
We shall perform the following calculations.
\begin{itemize}
\item[(a)] For $\Xi\in\Xlnb[r,q]$ of partial weight $e_\Xi = (1,1,\dots,1)$, we have $\theta_{\Psi\cup\Xi} = \theta_\Psi\cdot R_{[r,q-1]}\cdot X_q$.

\item[(b)] For $\Xi\in\Xlnb[r,q]$ of partial weight $e_\Xi = (1,1,\dots,1,2)$, we have $\theta_{\Psi\cup\Xi} = \theta_\Psi\cdot R_{[r,q-1]}\cdot 2$. 

\item[(c)] Conversely, $\theta_\Psi\cdot R_{[r,q-1]}\cdot 2$ or $\theta_\Psi\cdot R_{[r,q-1]}\cdot X_q$ divides $\theta_{\Psi\cup\Xi}$ for any $\Xi\in\Xlnb[r,q]$.
\end{itemize}

Ad (a). We perform an induction on $q$, starting with $q = r$, $R_{[r,r-1]} = 1$ and $\theta_{\Psi\cup\paze{+}{-}} = \theta_{\Psi\cup\paze{-}{+}} = \theta_\Psi\cdot X_q$. Assume the 
assertion known for $\ck\Xi := \Xi_{[r,q-1]}$. We calculate
\[
\begin{array}{rcl}
\theta_{\Psi\cup\ck\Xi\cup\paze{+}{-}}
& \auf{\mbox{\scr (\ref{LemC3})}}{=} & \theta_{\Psi\cup\ck\Xi}\cdot (X_{q-1} - 2) \\
& = & \theta_\Psi\cdot R_{[r,q-2]}\cdot X_{q-1}\cdot (X_{q-1} - 2) \\
& = & \theta_\Psi\cdot R_{[r,q-1]}\cdot X_q, \\
\end{array}
\]
dito for $\xi\ohne\h\Xi = \paze{-}{+}$

Ad (b). The argument for (a) can be applied, except that we obtain a factor $2$ instead of a factor $(X_{q-1} - 2) = X_q$ in the last step.

Ad (c). We compare to the proof of (\ref{LemC6}). We need to show that there exists a step in the induction in which the quotient $\theta_{\Psi\cup\Xi}/(R_{[r,q-1]}\theta_{\Psi})$ is divisible 
by $2$ or by $X_q$. In case $e_{\Xi,r} = 2$, a factor $2$ enters at $r$ already, so that we may assume $e_{\Xi,r}\leq 1$. In case $e_{\Xi,q} \leq 1$, a factor $X_q$ enters  
at $q$, so that we may assume $e_{\Xi,q} = 2$. Since $\Xi$ is not bulky, necessarily $\pazz{-}{+}{+}{+}$ or $\pazz{+}{+}{-}{+}$ occurs as a subpattern of $\Xi$. But
then a factor $2$ enters at the induction step that adjoins the right column of that subpattern (\ref{LemC3}). 
\end{Lemma}

\begin{Lemma}
\label{LemC8}
Let $g+1\leq r\leq q\leq k$ such that $r - 1\in\b L$, $[r,q-1]\tm L(0)\cup L(1)$ and $[r,q-1]\cap L(0)\neq\leer$, let $\Psi\in\Xlnb[g,r-1]$. As ideal of coefficients, we obtain
\[
\Big(\theta_{\Psi\cup\Xi} \; |\; \Xi\in\Xl[r,q] \Big) = \Big(\theta_\Psi\cdot R_{[r,q-1]}\Big)\tm\Z.
\]

\rm
Let $s := \max(L(0)\cap [r,q-1]) + 1 \in [r+1,q]$. We shall show the following assertions.
\begin{itemize}
\item[(a)] For $t\in [s,q]$ and a partial pattern $\Xi(q,t)\in\Xlnb[r,q]$ of partial weight
\[
e_{\Xi(q,t)} = \left\{\begin{array}{ll}
               2 & \mbox{ for } j\in [s,t] \\
               1 & \mbox{ for } j\in [r,q]\ohne [s,t], \\
               \end{array}\right. 
\] 
which is a slight abuse of notation since $\Xi(q,t)$ does not only depend on $q$ and $t$, we have
\[
\theta_{\Psi\cup\Xi(q,t)} = \left\{\begin{array}{ll}
                            \pm\theta_\Psi\cdot R_{[r,q-1]}\cdot \Big(\prodd{j\in [s,t-1]} (X_j - 3)\Big)\cdot (X_{q-1} - 2) & \mbox{ for } t\in [s,q-1] \\
                            \pm\theta_\Psi\cdot R_{[r,q-1]}\cdot \Big(\prodd{j\in [s,q-1]} (X_j - 3)\Big)                    & \mbox{ for } t = q, \\
                            \end{array}\right. 
\] 
the sign $\pm$ indicating that the equation holds up to sign.

\item[(b)] In case $s\in [r+1,q-2]$, for $\Xi\in\Xlnb[r,q]$ of partial weight
\[
e_{\Xi,j} = \left\{\begin{array}{ll}
            2 & \mbox{ for } j\in \{ s, q\} \\
            1 & \mbox{ for } j\in [r,q]\ohne \{ s,q\}, \\
            \end{array}\right. 
\] 
we have
\[
\theta_{\Psi\cup\Xi} = \theta_\Psi\cdot R_{[r,q-1]}\cdot 2.
\]
\item[(c)] Assertions (a) and (b) imply the lemma.
\end{itemize}

Ad (c). The inclusion $\om$ follows from (\ref{LemC6}). In case $s = q$, the lemma follows from $\theta_{\Psi\cup\Xi(q,q)} = \theta_\Psi\cdot R_{[r,q-1]}$. In case $s = q-1$, the lemma follows 
from 
\[
\begin{array}{rcl}
\theta_{\Psi\cup\Xi(q,q-1)} & = & \pm \theta_\Psi R_{[r,q-1]} (X_{q-1} - 2) \\
\theta_{\Psi\cup\Xi(q,q)}   & = & \pm \theta_\Psi R_{[r,q-1]} (X_{q-1} - 3). \\
\end{array}
\]
In case $s\in [r+1,q-2]$, it suffices to show that
\[
\afk := \Big( 2, \prodd{j\in [s,q-1]}(X_j - 3), \Big(\prodd{j\in [s,t-1]}(X_j - 3)\Big)\cdot (X_{q-1} - 2) \; \Big|\; t\in [s,q-1] \Big) = (1).
\]
We perform an induction on $t$ from $t = q-1$ downwards to $t = s$, contending that $\prodd{j\in [s,t-1]}(X_j - 3)\in\afk$. In fact,
\[
\Big(\prodd{j\in [s,t]}(X_j - 3)\Big) - \Big(\prodd{j\in [s,t-1]}(X_j - 3)\Big)\cdot (X_{q-1} - 2) = \Big(\prodd{j\in [s,t-1]}(X_j - 3)\Big)(X_t - X_{q-1} - 1)
\]
and $X_t - X_{q-1} = 2(t-q+1)$.

Ad (a). Let $\b q \in [r,q]$. For $t\geq \b q + 1$, let $\Xi(\b q,t)\in\Xlnb[r,\b q]$ denote a partial pattern of the same partial weight as $\Xi(\b q,\b q)$. We perform an induction on 
$\b q$ running from $r$ to $q$ to prove the assertion that
\[
\theta_{\Psi\cup\Xi(\b q,t)} = \left\{
\begin{array}{ll}
\theta_{\Psi}\cdot R_{[r,\b q-1]}\cdot X_{\b q}                                                         & \mbox{ for } \b q\in [r,s-1] \\
\pm\theta_{\Psi}\cdot R_{[r,\b q-1]}\cdot \Big(\prodd{j\in [s,\b q-1]}(X_j - 3)\Big)                    & \mbox{ for } \b q\in [s,t] \\
\pm\theta_{\Psi}\cdot R_{[r,\b q-1]}\cdot \Big(\prodd{j\in [s,t-1]}(X_j - 3)\Big)\cdot (X_{\b q-1} - 2) & \mbox{ for } \b q\in [t+1,q]. \\
\end{array}
\right.
\]
To begin with, we recall our situation, viz.\ $r < r+1\leq s\leq t$ and $r\leq \b q$.

For $\b q = r$, we have $R_{[r,r-1]} = 1$ and $e_{\Xi(\b q,t),r} = 1$, whence we obtain
\[
\theta_{\Psi\cup\Xi(\b q,t)} = \theta_{\Psi}\cdot R_{[r,\b q - 1]}\cdot X_{\b q}
\]
since $r - 1\in \b L$.

Assume the assertion to be true for $r\leq \b q - 1\leq s - 2$.

Case $\b q-1\in L(0)$. We calculate
\[
\begin{array}{rcl}
\theta_{\Psi\cup\Xi(\b q,t)}
& \auf{(\ref{LemC2})}{=} & \theta_{\Psi\cup\Xi(\b q-1,t)}\cdot(X_{\b q-1} - 1)(X_{\b q-1} + 1) \\
& = & \theta_\Psi\cdot R_{[r,\b q-2]}\cdot X_{\b q-1}\cdot (X_{\b q-1} - 1)(X_{\b q-1} + 1) \\ 
& = & \theta_\Psi\cdot R_{[r,\b q-1]}\cdot X_{\b q}. \\ 
\end{array}
\]
Case $\b q-1\in L(1)$. We calculate
\[
\begin{array}{rcl}
\theta_{\Psi\cup\Xi(\b q,s)}
& \auf{(\ref{LemC3})}{=} & \theta_{\Psi\cup\Xi(\b q-1,t)}\cdot (X_{\b q-1} - 2) \\
& = & \theta_\Psi\cdot R_{[r,\b q-2]}\cdot X_{\b q-1}\cdot (X_{\b q-1} - 2) \\
& = & \theta_\Psi\cdot R_{[r,\b q-1]}\cdot X_{\b q}. \\
\end{array}
\]
Assume the assertion to be true for $\b q - 1 = s - 1$. Note that $\b q - 1\in L(0)$, so
\[
\begin{array}{rcl}
\theta_{\Psi\cup\Xi(\b q,t)}
& \auf{(\ref{LemC2})}{=} & \theta_{\Psi\cup\Xi(\b q-1,t)} \cdot (X_{\b q-1} + 1) \\
& = & \theta_\Psi\cdot R_{[r,\b q-2]}\cdot X_{\b q-1}\cdot (X_{\b q-1} + 1) \\
& = & \theta_\Psi\cdot R_{[r,\b q-1]}. \\
\end{array}
\]
Assume the assertion to be true for $s\leq \b q-1\leq t-1$. Note that $\b q - 1\in L(1)$, so
\[
\begin{array}{rcl}
\theta_{\Psi\cup\Xi(\b q,t)}
& \auf{(\ref{LemC3})}{=} & -\theta_{\Psi\cup\Xi(\b q-1,t)}\cdot X_{\b q-1}(X_{\b q-1} - 3) \\
& = & \pm\theta_\Psi\cdot R_{[r,\b q-2]}\cdot\Big(\prodd{j\in [s,\b q-2]} (X_j - 3)\Big)\cdot X_{\b q-1}(X_{\b q-1} - 3) \\
& = & \pm\theta_\Psi\cdot R_{[r,\b q-1]}\cdot\Big(\prodd{j\in [s,\b q-1]} (X_j - 3)\Big).\\
\end{array}
\]
Assume the assertion to be true for $\b q-1 = t$. Note that $\b q - 1\in L(1)$, so
\[
\begin{array}{rcl}
\theta_{\Psi\cup\Xi(\b q,t)}
& \auf{(\ref{LemC3})}{=} & \theta_{\Psi\cup\Xi(\b q-1,t)}\cdot X_{\b q-1}(X_{\b q-1} - 2) \\
& = & \pm\theta_\Psi\cdot R_{[r,\b q-2]}\cdot\Big(\prodd{j\in [s,\b q-2]} (X_j - 3)\Big)\cdot X_{\b q-1}(X_{\b q-1} - 2) \\
& = & \pm\theta_\Psi\cdot R_{[r,\b q-1]}\cdot\Big(\prodd{j\in [s,t-1]} (X_j - 3)\Big)\cdot (X_{\b q-1} - 2). \\
\end{array}
\]
Assume the assertion to be true for $t+1\leq \b q-1$. Note that $\b q - 1\in L(1)$, so
\[
\begin{array}{rcl}
\theta_{\Psi\cup\Xi(\b q,t)}
& \auf{(\ref{LemC3})}{=} & \theta_{\Psi\cup\Xi(\b q-1,t)}\cdot (X_{\b q-1} - 2) \\
& = & \pm\theta_\Psi\cdot R_{[r,\b q-2]}\cdot\Big(\prodd{j\in [s,t-1]} (X_j - 3)\Big)\cdot (X_{\b q - 2} - 2)\cdot (X_{\b q - 1} - 2) \\
& = & \pm\theta_\Psi\cdot R_{[r,\b q-1]}\cdot\Big(\prodd{j\in [s,t-1]} (X_j - 3)\Big)\cdot (X_{\b q - 1} - 2). \\
\end{array}
\]

Ad (b). We modify the last step in the proof of the formula for $\theta_{\Psi\cup\Xi(q,s)}$ in (a) according to (\ref{LemC3}), i.e.\ we replace the factor $(X_{q-1} - 2)$ by the factor $2$.
\end{Lemma}

A subinterval $[r,s]$ of $[g+1,k]$ is called {\it connected} if $[r,s-1]\tm L(0)\cup L(1)$. A maximal connected subinterval with respect to the ordering given by inclusion is called a
{\it component.} We write $[g+1,k]$ as a disjoint decomposition into components
\[
[g+1,k] = \Cup_{\kappa\in [1,K]} [p(\kappa),q(\kappa)],
\]
where $p(\kappa)\leq q(\kappa)$ for $\kappa\in [1,K]$, $p(1) = g+1$, $q(\kappa) + 1 = p(\kappa + 1)$ for $\kappa\in [1,K-1]$ and $q(K) = k$. So
$\{ q(\kappa)\; |\; \kappa\in [1,K]\} = \b L\ohne\{ g\}$.

\begin{Lemma}
\label{LemC9}
The ideal of coefficients may be factored as
\[
\Big(\theta_\Xi\; | \; \Xi\in\Xl[g,k+1]\Big) = \prod_{\kappa\in [1,K]} \tfk_\kappa,
\]
where for $\kappa\in [1,K]$ we denote
\[
\tfk_\kappa := \Big( \theta_{\Psi\cup\Xi}/\theta_\Psi \; |\; \Xi\in\Xl[p(\kappa),q(\kappa)]\Big) \tm \Z
\]
for some choice of a partial pattern $\Psi\in\Xlnb[g,p(\kappa) - 1]$, a choice which does not affect these quotients.

\rm
This factorization follows by the independence of choice just mentioned.
\end{Lemma}

\begin{Proposition}
\label{PropC10}
Let the {\rm factor of redundancy} be
\[
R := R_{[g+1,k]}\cdot\Big(\prodd{\auf{\scm\kappa\in [1,K],}{[p(\kappa),q(\kappa) - 1]\tm L(1)}} \dot X_{q(\kappa)}\Big).
\]
The ideal of coefficients is given by
\[
\Big(\theta_\Xi\; | \; \Xi\in\Xl[g,k+1]\Big) = \Big( R \Big) \tm \Z.
\]
In particular, the $\Z\Sl_n$-linear map
\[
M^{\lambda',-}\lraa{f'} S^\mu
\]
is divisible exactly by $R$, i.e.\ any matrix representing $f'$ linearly over $\Z$ has $R$ as greatest common divisor of its entries.

\rm
The assertion on the ideal of coefficients ensues from (\ref{LemC7}, \ref{LemC8}, \ref{LemC9}). We conclude by (\ref{LemC2}, \ref{LemC3}) that the linear combination
\[
f'' = \sumd{\Xi\in\Xl[g,k+1]}\Big(\prodd{j\in [g+1,k]} (e_{\Xi,j}!\cdot X_j^{(2 - e_{\Xi,j})})\Big)\cdot \v f''_\Xi: F^\lambda\lra S^\mu,
\]
cf.\ (\ref{LemC4}), can be rewritten as
\[
f'' = \sumd{\Xi\in\Xl[g,k+1]} \theta_\Xi \v f''_\Xi,
\]
whose coefficients are divisible by $R$ and vanish for $\Xi$ bulky. Applying $f''$ to the tableau $[\ck a_\lambda]$, we see by (\ref{LemC5}) that $f''$ is divisible exactly by $R$.
The assertion remains true for $f'$, being the factorization of $f''$ over the epimorphism $F^\lambda\lraa{\nu_M^\lambda} M^{\lambda',-}$.
\end{Proposition}

The modulus remains to be taken under consideration.

\begin{Lemma}
\label{LemC15}
We may rewrite
\[
\begin{array}{rcl}
2\cdot\sumd{\Xi\in\Xl[g+2,k+1]} \Big(\prodd{j\in [g+2,k]} (e_{\Xi,j}!\cdot X_j^{(2-e_{\Xi,j})})\Big) \v A^{u,v}_\Xi 
& = & \sumd{\Xi\in\Xl[g+2,k+1]} \theta_{\pazz{+}{+}{+}{+}\cup\Xi} \v A^{u,v}_\Xi \\
X_{g+1}\cdot\sumd{\Xi\in\Xl[g+2,k+1]} \Big(\prodd{j\in [g+2,k]} (e_{\Xi,j}!\cdot X_j^{(2-e_{\Xi,j})})\Big) (\v A^{u,-}_\Xi + \v A^{-,u}_\Xi) \\
= \sumd{\Xi\in\Xl[g+2,k+1]} \Big(\theta_{\pazz{+}{+}{+}{-}\cup\Xi} \v A^{u,-}_\Xi + \theta_{\pazz{+}{-}{+}{+}\cup\Xi} \v A^{-,u}_\Xi\Big)  \\
\end{array}
\]

\rm 
This ensues from (\ref{LemC2_5}, \ref{LemC2}, \ref{LemC3_5}, \ref{LemC3}).
\end{Lemma}

\begin{Lemma}
\label{LemC11}
Let $g+1\leq q\leq k$, let $\xi = \paze{+}{-}$ or $\xi = \paze{-}{+}$ in $\Xl[g+1,g+1]$. Suppose $[g+1,q-1]\tm L(1)$. As ideal of coefficients, we obtain
\[
\Big(\theta_\Xi\; |\; \Xi\in\Xl[g,q],\; \Xi_{[g+1,g+1]} = \xi\Big) = \Big( R_{[g+1,q-1]}\cdot \dot X_q\Big). 
\]
Consequently, by comparison with (\ref{PropC10}), 
\[
\Big(\theta_\Xi\; |\; \Xi\in\Xl[g,p+1],\; \Xi_{[g+1,g+1]} = \xi\Big) = \Big( R\Big). 
\]

\rm
Consider the partial weights indexing generators of the ideal of coefficients having no restriction imposed on column $g + 1$. According to assertions (a) and (b) of the proof of (\ref{LemC7}), 
they can be chosen to have column $g + 1$ equal to $\xi$.
\end{Lemma}

\begin{Lemma}
\label{LemC12}
Let $g+1\leq q\leq k$, let $\xi = \paze{+}{-}$ or $\xi = \paze{-}{+}$ in $\Xl[g+1,g+1]$. Suppose $[g+1,q-1]\tm L(0)\cup L(1)$ such that $[g+1,q-1]\cap L(0)\neq\leer$. As ideal of coefficients 
we obtain
\[
\Big(\theta_\Xi\; |\; \Xi\in\Xl[g,q],\; \Xi_{[g+1,g+1]} = \xi\Big) = \Big(R_{[g+1,q-1]}\Big). 
\]
Whence, by comparison with (\ref{PropC10}), we have
\[
\Big(\theta_\Xi\; |\; \Xi\in\Xl[g,p+1],\; \Xi_{[g+1,g+1]} = \xi\Big) = \Big( R\Big). 
\]

\rm
Consider the partial weights indexing generators of the ideal of coefficients having no restriction imposed on column $g + 1$. According to assertions (a) and (b) of the proof of (\ref{LemC8}),
they can be chosen to have column $g + 1$ equal to $\xi$.
\end{Lemma}

\begin{Lemma}
\label{LemC13}
Let $g+1\leq q\leq k$. Suppose $[g+1,q-1]\tm L(1)$. As ideal of coefficients, we obtain
\[
\Big(\theta_\Xi\; |\; \Xi\in\Xl[g,q],\; \Xi_{[g+1,g+1]} = \paze{+}{+}\Big) = \Big(2 R_{[g+1,q-1]}\Big).
\]
Therefore, by comparison with (\ref{PropC10}), 
\[
\Big(\theta_\Xi\; |\; \Xi\in\Xl[g,k+1],\; \Xi_{[g+1,g+1]} = \paze{+}{+}\Big) = \Big(2R/\dot X_{q(1)} \Big).
\]
Note that $\dot X_{q(1)} = \dot X_{g+1}$.

\rm
We claim the following assertions.
\begin{itemize}
\item[(a)] For $t\in [g+1,q]$, we let $\Xi(g,t)\in\Xlnb[g,q]$ denote a partial pattern of partial weight
\[
e_{\Xi(q,t),j} = \left\{\begin{array}{ll}
                 2 & \mbox{ for } j\in [g,t] \\
                 1 & \mbox{ for } j\in [t+1,q] \\
                 \end{array}\right.
\] 
and obtain
\[
\theta_{\Xi(q,t)} = \left\{\begin{array}{ll}
                    \pm 2 R_{[g+1,q-1]}\cdot\Big(\prodd{j\in [g+1,t-1]} (X_j - 3)\Big)\cdot (X_{q-1} - 2) & \mbox{ for } t\in [g+1,q-1] \\
                    \pm 2 R_{[g+1,q-1]}\cdot\Big(\prodd{j\in [g+1,q-1]} (X_j - 3)\Big)                    & \mbox{ for } t = q. \\
                    \end{array}\right.
\]

\item[(b)] In case $q\geq g+3$, for $\Xi\in\Xlnb[g,q]$ of partial weight $e_\Xi = (2,2,1,\dots,1,2)$, we have 
\[
\theta_\Xi = 4 R_{[g+1,q-1]}.
\]
\item[(c)] Assertions (a) and (b) imply the lemma.
\end{itemize}

(a, b, c) follow by the arguments of (\ref{LemC8}).
\end{Lemma}

\begin{Lemma}
\label{LemC14}
Let $g+1\leq q\leq k$. Suppose $[g+1,q-1]\tm L(0)\cup L(1)$, $[g+1,q-1]\cap L(0)\neq\leer$. As ideal of coefficients, we obtain
\[
\Big(\theta_\Xi\; |\; \Xi\in\Xl[g,q],\; \Xi_{[g+1,g+1]} = \paze{+}{+}\Big) = \Big( 2 R_{[g+1,q-1]}\Big).
\]
So, by comparison with (\ref{PropC10}),
\[
\Big(\theta_\Xi\; |\; \Xi\in\Xl[g,k+1],\; \Xi_{[g+1,g+1]} = \paze{+}{+}\Big) = \Big( 2R \Big).
\]

\rm
Denote 
\[
\begin{array}{rcll}
p & := & \min\Big(([g+1,q-1]\cap L(1))\cup \{ q\}\Big) & \in [g+1,q] \\
s & := & \max\Big([g+1,q-1]\cap L(0)\Big) + 1          & \in [g+2,q] \\
\end{array}
\]
so that we are in the situation $g < g + 1\leq p\leq s\leq q$, more precisely, $p = s$ or $p\leq s-2$. 

We shall show the following assertions.
\begin{itemize}
\item[(a)] For $t\in [s,q]$ and a partial pattern $\Xi(q,t)\in\Xlnb[g,q]$ of partial weight
\[
e_{\Xi(q,t),j} = 
\left\{\begin{array}{ll}
2 & \mbox{ for } j\in [g,p] \\
1 & \mbox{ for } j\in [p+1,s-1] \\
2 & \mbox{ for } j\in [s,t] \\
1 & \mbox{ for } j\in [t+1,q], \\
\end{array}\right.
\]
we obtain
\[
\theta_{\Xi(q,t)} = \left\{\begin{array}{ll}
                    \pm 2 R_{[g+1,q-1]}\cdot\Big(\prodd{j\in [s,t-1]} (X_j - 3)\Big)\cdot (X_{q-1} - 2) & \mbox{ for } t\in [s,q-1] \\
                    \pm 2 R_{[g+1,q-1]}\cdot\Big(\prodd{j\in [s,q-1]} (X_j - 3)\Big)                    & \mbox{ for } t = q. \\
                    \end{array}\right.
\]

\item[(b)] In case $s\in [g+2,q-2]$, for $\Xi\in\Xlnb[g,q]$ of partial weight
\[
e_{\Xi,j} = \left\{\begin{array}{ll}
            2 & \mbox{ for } j\in [g,p] \\
            1 & \mbox{ for } j\in [p+1,s-1] \\
            2 & \mbox{ for } j = s \\
            1 & \mbox{ for } j\in [s+1,q-1] \\
            2 & \mbox{ for } j = q, \\
            \end{array}\right. 
\]
we obtain
\[
\theta_\Xi = 4R_{[g+1,q-1]}.
\]

\item[(c)] Assertions (a) and (b) imply the lemma.
\end{itemize}

Ad (c). As in (\ref{LemC8}).

Ad (a). We perform an induction on $\b q\in [g+1,q]$ to prove
\[
\theta_{\Xi(\b q,t)} = \left\{\begin{array}{ll}
                           2 R_{[g+1,\b q - 1]}                                                            & \mbox{ for } \b q \in [g+1,p] \\
                           2 R_{[g+1,\b q - 1]}\cdot X_{\b q}                                              & \mbox{ for } \b q \in [p+1,s-1] \\
                       \pm 2 R_{[g+1,\b q - 1]}\cdot\Big(\prodd{j\in [s,\b q - 1]} (X_j - 3)\Big)          & \mbox{ for } \b q \in [s,t] \\
                       \pm 2 R_{[g+1,\b q - 1]}\cdot\Big(\prodd{j\in [s,t-1]} (X_j - 3)\Big)\cdot (X_{\b q-1} - 2) & \mbox{ for } \b q\in [t+1,q], \\
                       \end{array}\right.
\]
where $e_{\Xi(\b q,t)} = e_{\Xi(\b q,\b q)}$ for $t\geq \b q + 1$. We recall that $g < g + 1\leq p\leq s\leq t$.

We begin the induction by remarking that for $\b q = g+1$ we have $R_{[g+1,g]} = 1$, $e_{\Xi(\b q,t),g+1} = 2$, and thus $\theta_{\Xi(\b q,t)} = 2 R_{[g+1,\b q - 1]}$, as required.

{\it Case {\rm I}, $p\leq s - 2$.} 

Assume the assertion known for $\b q - 1 \in [g+1,p-1]$. Note that $\b q - 1\in L(0)$, so
\[
\begin{array}{rcl}
\theta_{\Xi(\b q,t)} 
& \auf{\mbox{\scr (\ref{LemC2})}}{=} & \theta_{\Xi(\b q - 1,t)}\cdot X_{\b q - 1}(X_{\b q - 1} + 1) \\
& = & 2 R_{[g+1,\b q - 2]} \cdot X_{\b q - 1}(X_{\b q - 1} + 1) \\
& = & 2 R_{[g+1,\b q - 1]}.\\
\end{array}
\] 
Assume the assertion known for $\b q - 1 = p$. Note that $\b q - 1\in L(1)$ since $p < s$, and that therefore $[g+1,q-1]\cap L(1)\neq\leer$, so
\[
\begin{array}{rcl}
\theta_{\Xi(\b q,t)} 
& \auf{\mbox{\scr (\ref{LemC3})}}{=} & \theta_{\Xi(\b q-1,t)}\cdot X_{\b q - 1}(X_{\b q - 1} - 2) \\
& = & 2 R_{[g+1,\b q - 2]} \cdot X_{\b q - 1}(X_{\b q - 1} - 2) \\
& = & 2 R_{[g+1,\b q - 1]} X_{\b q}. \\
\end{array}
\] 
Assume the assertion known for $\b q - 1 \in [p+1,s-2]$. In case $\b q - 1\in L(0)$, we obtain
\[
\begin{array}{rcl}
\theta_{\Xi(\b q,t)} 
& \auf{\mbox{\scr (\ref{LemC2})}}{=} & \theta_{\Xi(\b q - 1,t)}\cdot (X_{\b q-1} - 1)(X_{\b q - 1} + 1) \\
& = & 2 R_{[g+1,\b q - 2]} X_{\b q - 1}\cdot (X_{\b q - 1} - 1)(X_{\b q - 1} + 1) \\
& = & 2 R_{[g+1,\b q - 1]} X_{\b q}. \\
\end{array}
\]
In case $\b q - 1\in L(1)$, we obtain
\[
\begin{array}{rcl}
\theta_{\Xi(\b q,t)}
& \auf{\mbox{\scr (\ref{LemC3})}}{=} & \theta_{\Xi(\b q - 1,t)}\cdot (X_{\b q - 1} - 2) \\
& = & 2 R_{[g+1,\b q - 2]} X_{\b q - 1}\cdot (X_{\b q - 1} - 2) \\
& = & 2 R_{[g+1,\b q - 1]} X_{\b q}. \\
\end{array}
\] 
Assume the assertion known for $\b q - 1 = s - 1$. Note that $\b q - 1\in L(0)$, so
\[
\begin{array}{rcl}
\theta_{\Xi(\b q,t)} 
& \auf{\mbox{\scr (\ref{LemC2})}}{=} & \theta_{\Xi(\b q - 1,t)}\cdot (X_{\b q - 1} + 1) \\
& = & 2 R_{[g+1,\b q - 2]} X_{\b q - 1}\cdot (X_{\b q - 1} + 1) \\
& = & 2 R_{[g+1,\b q - 1]}. \\
\end{array}
\]
Assume the assertion known for $\b q - 1 \in [s,t-1]$. Note that $\b q - 1\in L(1)$, thus
\[
\begin{array}{rcl}
\theta_{\Xi(\b q,t)} 
& \auf{\mbox{\scr (\ref{LemC3})}}{=} & \theta_{\Xi(\b q - 1,t)}\cdot (-X_{\b q - 1}(X_{\b q - 1} - 3)) \\
& = & \pm 2 R_{[g+1,\b q - 2]} \Big(\prodd{j\in [s,\b q - 2]} (X_j - 3)\Big)\cdot X_{\b q - 1}(X_{\b q - 1} - 3)\\
& = & \pm 2 R_{[g+1,\b q - 1]} \Big(\prodd{j\in [s,\b q - 1]} (X_j - 3)\Big). \\
\end{array}
\]
Assume the assertion known for $\b q - 1 = t$. Note that $\b q - 1\in L(1)$, and consequently
\[
\begin{array}{rcl}
\theta_{\Xi(\b q,t)} 
& \auf{\mbox{\scr (\ref{LemC3})}}{=} & \theta_{\Xi(\b q - 1,t)}\cdot X_{\b q - 1}(X_{\b q - 1} - 2) \\
& = & \pm 2 R_{[g+1,\b q - 2]} \Big(\prodd{j\in [s,t-1]} (X_j - 3)\Big)\cdot X_{\b q - 1}(X_{\b q - 1} - 2)\\
& = & \pm 2 R_{[g+1,\b q - 1]} \Big(\prodd{j\in [s,t-1]} (X_j - 3)\Big)(X_{\b q - 1} - 2). \\
\end{array}
\]
Assume the assertion known for $\b q - 1 \in [t+1,q-1]$. Note that $\b q - 1\in L(1)$, whence
\[
\begin{array}{rcl}
\theta_{\Xi(\b q,t)} 
& \auf{\mbox{\scr (\ref{LemC3})}}{=} & \theta_{\Xi(\b q - 1,t)}\cdot (X_{\b q - 1} - 2) \\
& = & \pm 2 R_{[g+1,\b q - 2]} \Big(\prodd{j\in [s,t-1]} (X_j - 3)\Big)(X_{\b q - 2} - 2)\cdot (X_{\b q - 1} - 2)\\
& = & \pm 2 R_{[g+1,\b q - 1]} \Big(\prodd{j\in [s,t-1]} (X_j - 3)\Big)(X_{\b q - 1} - 2). \\
\end{array}
\]

{\it Case {\rm II}, $p = s$.} I.e.\ $L(0)\cap [g+1,q-1] = [g+1,p-1]$, $L(1)\cap [g+1,q-1] = [s,q-1]$.

The case $\b q\in [g+1,p]$ follows from the argument for $\b q - 1 \in [g+1,p-1]$ in Case I. The case $\b q \in [s+1,t]$ follows from the argument for $\b q - 1 \in [s,t-1]$ in Case I. The
case $\b q\in [t+1,q]$ follows from the arguments for $\b q - 1 = t$ and for $\b q - 1 \in [t+1,q-1]$ in Case I.

Ad (b). We modify the last step in the proof of (a) according to (\ref{LemC3}), i.e.\ we replace the factor $(X_{q-1} - 2)$ by the factor $2$. 
\end{Lemma}

\begin{Theorem}
\label{ThC16}
Suppose $g < k$. Let 
\[
m := \left\{\begin{array}{lll}
      X_g + 2                        & \mbox{ if } & [p(1),q(1)-1]\not\tm L(1) \\
     (X_g + 2)/\gcd(2, X_g, X_{g+1}) & \mbox{ if } & [p(1),q(1)-1]\tm L(1). \\
     \end{array}\right.
\]
Note that if $g+1\in\b L$, the case $\leer = [p(1),q(1)-1]\tm L(1)$ applies. The integer 
\[
X_g + 2 = (\lambda'_g - g) - (\lambda'_{k+1} - (k+1)) + 2 
\]
allows an interpretation as box shift length. The quotient $f'/R$ is integral and factors over
\begin{center}
\begin{picture}(250,250)
\put( -10, 200){$M^{\lambda',-}$}
\put( 100, 210){\vector(1,0){80}}
\put( 110, 225){$\scm f'/R$}
\put( 200, 200){$S^\mu$}
\put(  20, 180){\vector(0,-1){130}}
\put( -30, 110){$\scm\nu_S^\lambda$}
\put( 220, 180){\vector(0,-1){130}}
\put( -10,   0){$S^\lambda$}
\put(  50,  10){\vector(1,0){130}}
\put( 110,  25){$\scm f$}
\put( 200,   0){$S^\mu/m$.}
\end{picture}
\end{center}
The morphism $f$ has order $m$ as an element of $\Hom_{\sZ\Sl_n}(S^\lambda,S^\mu/m)$. Recall that $f'$ maps
\[
\begin{array}{rcl}
M^{\lambda',-} & \lraa{f'} & S^\mu \\
\{ a' \}^-     & \lra      & \sumd{e\in E} (\prodd{i\in [g+1,k]} X_i^{(2-e_i)})\sumd{\gamma\in\dGamma(e)} \spi{a^\gamma}\eps_\gamma.
\end{array}
\]
The map $f$ operates accordingly, division taking place before going modulo $m$. Alternatively, by (\ref{PropC10}) we may write
\[
\begin{array}{rcl}
S^\lambda & \lraa{f} & S^\mu/m \\
\spi{a}   & \lra     & \sumd{\Xi\in\Xl[g,k+1]} (\theta_\Xi/R) \sumd{\gamma\in\dGamma(\Xi)}\spi{a^\gamma}\eps_\gamma,
\end{array}
\]
the coefficients $\theta_\Xi/R$ being integral.

Suppose $g = k$. The formula given for $f'$ in case $g < k$ holds in case $g = k$ as well, $E$ consisting of a single element, $R = 1$ and $m = (X_g + 2)/\gcd(2, X_g)$. We have the 
factorization diagram as above, the resulting morphism $S^\lambda\lraa{f} S^\mu/m$ being of order $m$. Cf.\ (\ref{PropM14}).

\rm
Suppose $g < k$.
By (\ref{PropC10}), the ideal generated by $\Z$-linear matrix coefficients of $f'/R$ is $(1)\tm\Z$. In case $[p(1),q(1)-1]\not\tm L(1)$, the integer $m\cdot R$ 
divides each of the required expressions in (\ref{LemC0}) by (\ref{LemC15}) and by the ideals of coefficients calculated in (\ref{LemC12}, \ref{LemC14}). In case $[p(1),q(1)-1]\tm L(1)$, 
the integer $m\cdot R$ divides each of the required expressions in (\ref{LemC0}) by (\ref{LemC15}) and by the ideals of coefficients calculated in (\ref{LemC11}, \ref{LemC13}).

Suppose $g = k$. Concerning the assertion, we refer to (\ref{PropM14}) (alternatively, to \ref{ThLW6}). This case behaves as if the first component lay in $L(1)$ since $X_{g+1} = 0$.
\end{Theorem}

\begin{Remark}
\label{RemC17}
\rm
In the situation of (\ref{ThC16}), {\sc Carter} and {\sc Payne} [CP 80] have shown that 
\[
\Hom_{K\Sl_n}(K\ts_{\sZ} S^\lambda,K\ts_{\sZ} S^\mu) \neq 0,
\]
$K$ being an infinite field of characteristic $p$ such that $p$ divides $X_g + 2$ in case $p\geq 3$, respectively, such that $4$ divides $X_g + 2$ in case $p = 2$.
This part of their result is recovered by (\ref{ThC16}).
\end{Remark}

\begin{Question}
\label{QuC18}
\rm
I do not know a description of $\Hom_{\sZ\Sl_n}(S^\lambda,S^\mu/n!)$ as an abelian group. See Section (\ref{SubSecEx}) for some examples.
\end{Question}        
\subsection{Vertical examples}
\label{SubSecEx}

In the examples that follow, calculated by means of linear algebra (\footnote{I thank the computer staff at the Fakult\"at f\"ur Mathematik in Bielefeld (`root') for kind cooperation.}), we 
omit to denote the brackets that indicate polytabloids. 

\begin{Example}
\label{ExE2}
\rm
Let $\lambda = (3,3)$, $\mu = (2,2,1,1)$. Then 
\[
\Hom_{\sZ\Sl_6}(S^\lambda,S^\mu/6!) \llaiso \Hom_{\sZ\Sl_6}(S^\lambda,S^\mu/4) \iso \Z/4,
\]
and a generator of $\Hom_{\sZ\Sl_6}(S^\lambda/4,S^\mu/4)$ is given by
\begin{footnotesize}
\[
\begin{array}{rcl}
\begin{array}{lll}
1 & 3 & 5 \\
2 & 4 & 6 \\
\end{array}
& \lraa{\mbox{$\ncr f$}} &
-2\cdot
\begin{array}{ll}
1 & 5   \\
2 & 6   \\
3 &     \\
4 &     \\
\end{array} 
- \left(
\begin{array}{ll}
1 & 5   \\
2 & 4   \\
3 &     \\
6 &     \\
\end{array} 
+
\begin{array}{ll}
1 & 3   \\
2 & 5   \\
4 &     \\
6 &     \\
\end{array} 
\right) - \left( 
\begin{array}{ll}
1 & 6   \\
2 & 4   \\
5 &     \\
3 &     \\
\end{array} 
+
\begin{array}{ll}
1 & 3   \\
2 & 6   \\
5 &     \\
4 &     \\
\end{array} 
\right) - 2\cdot
\begin{array}{ll}
1 & 3   \\
2 & 4   \\
5 &     \\
6 &     \\
\end{array}. \\
\end{array}
\]
\end{footnotesize}
\end{Example}

\begin{Example}
\label{ExE3}
\rm
Let $\lambda = (3,2,2)$, $\mu = (3,1,1,1,1)$. Then 
\[
\Hom_{\sZ\Sl_7}(S^\lambda,S^\mu/7!) \llaiso \Hom_{\sZ\Sl_7}(S^\lambda,S^\mu/6) \iso \Z/2\ti\Z/3,
\]
generators of $\Hom_{\sZ\Sl_7}(S^\lambda/6,S^\mu/6)$ being given by
\begin{footnotesize}
\[
\begin{array}{rcl}
\begin{array}{lll}
1 & 4 & 7 \\
2 & 5 &   \\
3 & 6 &   \\
\end{array}
& \lraa{\mbox{$\ncr u^{\lambda',\trp}u^\mu$}} &
3\cdot\Big(\sumd{[b]\mbox{\scr\ standard $\mu$-tableau}}\spi{b}\Big) \\
\end{array}
\]
\end{footnotesize}
(cf.\ \ref{Lem11_3} below) and by
\begin{footnotesize}
\[
\begin{array}{rcl}
\begin{array}{lll}
1 & 4 & 7 \\
2 & 5 &   \\
3 & 6 &   \\
\end{array}
& \lraa{\mbox{$\ncr f$}} &
-2\cdot\left(
\begin{array}{lll}
1 & 6 & 7  \\
2 &   &    \\
3 &   &    \\
4 &   &    \\
5 &   &    \\
\end{array} 
-
\begin{array}{lll}
1 & 5 & 7   \\
2 &   &     \\
3 &   &     \\
4 &   &     \\
6 &   &     \\
\end{array} 
+ 
\begin{array}{lll}
1 & 4 & 7   \\
2 &   &     \\
3 &   &     \\
5 &   &     \\
6 &   &     \\
\end{array} 
\right). \\
\end{array}
\]
\end{footnotesize}
\end{Example}

\begin{Example}
\label{ExE3_5}
\rm
Let $\lambda = (3,3,1,1)$, $\mu = (2,2,1,1,1,1)$. Then 
\[
\Hom_{\sZ\Sl_8}(S^\lambda,S^\mu/8!) \llaiso \Hom_{\sZ\Sl_8}(S^\lambda,S^\mu/6) \iso \Z/6,
\]
and a generator of $\Hom_{\sZ\Sl_8}(S^\lambda/6,S^\mu/6)$ is given by
\begin{footnotesize}
\[
\begin{array}{rcl}
\begin{array}{lll}
1 & 5 & 7 \\
2 & 6 & 8 \\
3 &   &   \\
4 &   &   \\
\end{array}
& \lraa{\mbox{$\ncr f$}} &
-2\cdot
\begin{array}{ll}
1 & 7   \\
2 & 8   \\
3 &     \\
4 &     \\
5 &     \\
6 &     \\
\end{array} 
- \left(
\begin{array}{ll}
1 & 7   \\
2 & 6   \\
3 &     \\
4 &     \\
5 &     \\
8 &     \\
\end{array} 
+
\begin{array}{ll}
1 & 5   \\
2 & 7   \\
3 &     \\
4 &     \\
6 &     \\
8 &     \\
\end{array} 
\right) - \left( 
\begin{array}{ll}
1 & 8   \\
2 & 6   \\
3 &     \\
4 &     \\
7 &     \\
5 &     \\
\end{array} 
+
\begin{array}{ll}
1 & 5   \\
2 & 8   \\
3 &     \\
4 &     \\
7 &     \\
6 &     \\
\end{array} 
\right) - 2\cdot
\begin{array}{ll}
1 & 5   \\
2 & 6   \\
3 &     \\
4 &     \\
7 &     \\
8 &     \\
\end{array}. \\
\end{array}
\]
\end{footnotesize}
\end{Example}

\begin{Example}
\label{ExE4}
\rm
Let $\lambda = (4,4)$, $\mu = (3,3,1,1)$. Then 
\[
\Hom_{\sZ\Sl_8}(S^\lambda,S^\mu/8!) \llaiso \Hom_{\sZ\Sl_8}(S^\lambda,S^\mu/5) \iso \Z/5,
\]
and a generator of $\Hom_{\sZ\Sl_8}(S^\lambda/5,S^\mu/5)$ is given by
\begin{footnotesize}
\[
\begin{array}{l}
\begin{array}{llll}
1 & 3 & 5 & 7 \\
2 & 4 & 6 & 8 \\
\end{array}
\lraa{\mbox{$\ncr f$}} \\
- 2\cdot
\begin{array}{lll}
1 & 5 & 7  \\
2 & 6 & 8  \\
3 &        \\
4 &        \\
\end{array} 
- \left(
\begin{array}{lll}
1 & 6 & 7  \\
2 & 4 & 8  \\
5 &        \\
3 &        \\
\end{array}
+  
\begin{array}{lll}
1 & 3 & 7  \\
2 & 6 & 8  \\
5 &        \\
4 &        \\
\end{array} 
\right) - \left(
\begin{array}{lll}
1 & 5 & 7  \\
2 & 4 & 8  \\
3 &        \\
6 &        \\
\end{array} 
+
\begin{array}{lll}
1 & 3 & 7  \\
2 & 5 & 8  \\
4 &        \\
6 &        \\
\end{array} 
\right) \\
- \left(
\begin{array}{lll}
1 & 5 & 8  \\
2 & 4 & 6  \\
7 &        \\
3 &        \\
\end{array} 
+
\begin{array}{lll}
1 & 6 & 5  \\
2 & 4 & 8  \\
7 &        \\
3 &        \\
\end{array} 
+
\begin{array}{lll}
1 & 3 & 8  \\
2 & 5 & 6  \\
7 &        \\
4 &        \\
\end{array} 
+
\begin{array}{lll}
1 & 3 & 5  \\
2 & 6 & 8  \\
7 &        \\
4 &        \\
\end{array} 
\right) \\
- \left(
\begin{array}{lll}
1 & 5 & 7  \\
2 & 4 & 6  \\
3 &        \\
8 &        \\
\end{array} 
+
\begin{array}{lll}
1 & 6 & 5  \\
2 & 4 & 7  \\
3 &        \\
8 &        \\
\end{array} 
+
\begin{array}{lll}
1 & 3 & 7  \\
2 & 5 & 6  \\
4 &        \\
8 &        \\
\end{array} 
+
\begin{array}{lll}
1 & 3 & 5  \\
2 & 6 & 7  \\
4 &        \\
8 &        \\
\end{array} 
\right) \\
- 2\cdot
\begin{array}{lll}
1 & 3 & 7  \\
2 & 4 & 8  \\
5 &        \\
6 &        \\
\end{array} 
- \left(
\begin{array}{lll}
1 & 3 & 8  \\
2 & 4 & 6  \\
7 &        \\
5 &        \\
\end{array} 
+
\begin{array}{lll}
1 & 3 & 5  \\
2 & 4 & 8  \\
7 &        \\
6 &        \\
\end{array} 
\right) - \left(
\begin{array}{lll}
1 & 3 & 7  \\
2 & 4 & 6  \\
5 &        \\
8 &        \\
\end{array} 
+
\begin{array}{lll}
1 & 3 & 5  \\
2 & 4 & 7  \\
6 &        \\
8 &        \\
\end{array} 
\right) - 2\cdot
\begin{array}{lll}
1 & 3 & 5  \\
2 & 4 & 6  \\
7 &        \\
8 &        \\
\end{array}. \\
\end{array}
\]
\end{footnotesize}
Note that $R = 6$, $\b m = 1$ and
\[
\begin{array}{rcl}
\theta_{\pazd{+}{+}{+}{+}{+}{+}} & = & 12 \\
\theta_{\pazd{+}{-}{+}{+}{+}{+}} \;\; =\;\; \theta_{\pazd{+}{+}{+}{+}{-}{+}} & = & 6 \\
\theta_{\pazd{+}{-}{-}{+}{+}{+}} \;\; =\;\; \theta_{\pazd{+}{+}{+}{+}{-}{-}} & = & 6 \\
\theta_{\pazd{+}{-}{+}{+}{-}{+}} & = & 12 \\
\theta_{\pazd{+}{-}{-}{+}{-}{+}} \;\; =\;\; \theta_{\pazd{+}{-}{+}{+}{-}{-}} & = & 6 \\
\theta_{\pazd{+}{-}{-}{+}{-}{-}} & = & 12. \\
\end{array}
\]
\end{Example}

\begin{Example}
\label{ExE5}
\rm
Let $\lambda = (3,3,2)$, $\mu = (2,2,2,1,1)$. Then 
\[
\Hom_{\sZ\Sl_8}(S^\lambda,S^\mu/8!) \llaiso \Hom_{\sZ\Sl_8}(S^\lambda,S^\mu/5) \iso \Z/5,
\]
$R = 2$, and a generator of $\Hom_{\sZ\Sl_8}(S^\lambda/5,S^\mu/5)$ is given by
\begin{footnotesize}
\[
\begin{array}{rcl}
\begin{array}{lll}
1 & 4 & 7 \\
2 & 5 & 8 \\
3 & 6 &   \\
\end{array}
& \lraa{\mbox{$\ncr f$}} &
-\left(
\begin{array}{ll}
1 & 7   \\
2 & 8   \\
3 & 6   \\
4 &     \\
5 &     \\
\end{array} 
+
\begin{array}{ll}
1 & 7   \\
2 & 5   \\
3 & 8   \\
4 &     \\
6 &     \\
\end{array} 
+
\begin{array}{ll}
1 & 4   \\
2 & 7   \\
3 & 8   \\
5 &     \\
6 &     \\
\end{array} 
\right) - \left( 
\begin{array}{ll}
1 & 7   \\
2 & 5   \\
3 & 6   \\
4 &     \\
8 &     \\
\end{array} 
+
\begin{array}{ll}
1 & 4   \\
2 & 7   \\
3 & 6   \\
5 &     \\
8 &     \\
\end{array} 
+
\begin{array}{ll}
1 & 4   \\
2 & 5   \\
3 & 7   \\
6 &     \\
8 &     \\
\end{array}
\right) \\
& & - \left(
\begin{array}{ll}
1 & 8   \\
2 & 5   \\
3 & 6   \\
7 &     \\
4 &     \\
\end{array} 
+
\begin{array}{ll}
1 & 4   \\
2 & 8   \\
3 & 6   \\
7 &     \\
5 &     \\
\end{array} 
+
\begin{array}{ll}
1 & 4   \\
2 & 5   \\
3 & 8   \\
7 &     \\
6 &     \\
\end{array}
\right) - 3\cdot
\begin{array}{ll}
1 & 4   \\
2 & 5   \\
3 & 6   \\
7 &     \\
8 &     \\
\end{array}.\\
\end{array}
\]
\end{footnotesize}
\end{Example}

\begin{Example}
\label{ExE6}
\rm
Let $\lambda = (3,3,1,1)$, $\mu = (2,2,2,2)$. Then 
\[
\Hom_{\sZ\Sl_8}(S^\lambda,S^\mu/8!) \llaiso \Hom_{\sZ\Sl_8}(S^\lambda,S^\mu/3) \iso \Z/3,
\]
and a generator of $\Hom_{\sZ\Sl_8}(S^\lambda/3,S^\mu/3)$ is given by
\begin{footnotesize}
\[
\begin{array}{rcl}
\begin{array}{lll}
1 & 5 & 7 \\
2 & 6 & 8 \\
3 &   &   \\
4 &   &   \\
\end{array}
& \lraa{\mbox{$\ncr f$}} &
-
\begin{array}{ll}
1 & 5   \\
2 & 6   \\
3 & 7   \\
4 & 8   \\
\end{array}.\\
\end{array}
\]
\end{footnotesize}
\end{Example}

\begin{Example}
\label{ExE7}
\rm
Let $\lambda = (4,3)$, $\mu = (2,2,2,1)$. Then 
\[
\Hom_{\sZ\Sl_7}(S^\lambda,S^\mu/7!) \llaiso \Hom_{\sZ\Sl_7}(S^\lambda,S^\mu/4) \iso \Z/4,
\]
and a generator of $\Hom_{\sZ\Sl_7}(S^\lambda/4,S^\mu/4)$ is given by
\begin{footnotesize}
\[
\begin{array}{rcl}
\begin{array}{llll}
1 & 3 & 5 & 7 \\
2 & 4 & 6 &  \\
\end{array}
& \lraa{\mbox{$\ncr g$}} &
2\cdot
\begin{array}{ll}
1 & 5   \\
2 & 6   \\
3 & 7   \\
4 &    \\
\end{array} 
+ 
\begin{array}{ll}
1 & 5   \\
2 & 4   \\
3 & 7   \\
6 &    \\
\end{array} 
+ 
\begin{array}{ll}
1 & 3   \\
2 & 5   \\
4 & 7   \\
6 &    \\
\end{array} 
+ 
\begin{array}{ll}
1 & 6   \\
2 & 4   \\
5 & 7   \\
3 &    \\
\end{array} 
+ 
\begin{array}{ll}
1 & 3   \\
2 & 6   \\
5 & 7   \\
4 &    \\
\end{array} 
+ 2\cdot
\begin{array}{ll}
1 & 3   \\
2 & 4   \\
5 & 7   \\
6 &    \\
\end{array}. \\
\end{array}
\]
\end{footnotesize}
\end{Example}

\begin{Example}
\label{ExE8}
\rm
Let $\lambda = (4,4)$, $\mu = (2,2,2,2)$. Then 
\[
\Hom_{\sZ\Sl_8}(S^\lambda,S^\mu/8!) \llaiso \Hom_{\sZ\Sl_8}(S^\lambda,S^\mu/4)\iso\Z/4,
\]
and a generator of $\Hom_{\sZ\Sl_8}(S^\lambda/4,S^\mu/4)$ is given by
\begin{footnotesize}
\[
\begin{array}{rcl}
\begin{array}{llll}
1 & 3 & 5 & 7 \\
2 & 4 & 6 & 8 \\
\end{array}
& \lraa{\mbox{$\ncr g$}} &
2\cdot
\begin{array}{ll}
1 & 5   \\
2 & 6   \\
3 & 7   \\
4 & 8   \\
\end{array} 
+ 
\begin{array}{ll}
1 & 5   \\
2 & 4   \\
3 & 7   \\
6 & 8   \\
\end{array} 
+ 
\begin{array}{ll}
1 & 3   \\
2 & 5   \\
4 & 7   \\
6 & 8   \\
\end{array} 
+ 
\begin{array}{ll}
1 & 6   \\
2 & 4   \\
5 & 7   \\
3 & 8   \\
\end{array} 
+ 
\begin{array}{ll}
1 & 3   \\
2 & 6   \\
5 & 7   \\
4 & 8   \\
\end{array} 
+ 2\cdot
\begin{array}{ll}
1 & 3   \\
2 & 4   \\
5 & 7   \\
6 & 8   \\
\end{array}. \\
\end{array}
\]
\end{footnotesize}
\end{Example}

\begin{Question}
\label{QuE9}
\rm
I do not know a generic morphism that specializes to the morphism in (\ref{ExE7}) or to the morphism in (\ref{ExE8}).
\end{Question}          
\textheight25.1cm
\section{Semistandard morphisms}
\label{SecSemi}

\begin{footnotesize}
The following considerations are based on methods of {\sc James} [J 78, 13] and follow a hint of {\sc Ringel.} We transpose our morphism from 
(\ref{ThC16}), cf.\ (\ref{PropSST11}, \ref{CorSST17}), where transposition is given by dualization, followed by alternation
and isomorphic substitution. Moreover, we establish a connection from the morphisms from $M^{\lambda,-}$ to $S^\mu$ exhibited in Section 
\ref{SubSecMS} to those of the form $\Theta^-_\phi\nu_S^\mu$ (\ref{RemSST15}, \ref{PropSST16}), where the morphisms of type $\Theta_\phi$ occur in [J 78, 13.13].

{\it Informal remark.} The original aim to give upper bounds on the $\Hom$-groups is still out of reach (cf.\ \ref{QuC18}). It could perhaps be approached by showing that the method that has 
been sufficient to construct our morphism (\ref{ThC16}) is also necessary in some sense. The connection we establish is to be seen as the first step of this approach, namely that the 
construction of our building blocks $f''_e$ in Section \ref{SubSecMS} has been necessary in some sense, except at the prime $2$ in certain cases (\ref{ThSST6}), and that there is a 
also reason to exclude bulky patterns by means of a reduced tuple of coefficients (\ref{RemSST18}). For a particular case in which this aim could be achieved at least in the regular case, see 
(\ref{ThLW6}).
\end{footnotesize}

\subsection{Correspondoids}
\label{SubsecCorr}

\begin{footnotesize}
We shall recall the {\sc Carter-Lusztig-James} construction of the semistandard basis of 
$\Hom_{\sZ\Sl_n}(S^\lambda,M^\mu)$ [J 78, 13.13], $\lambda$ and $\mu$ partitions of $n$, in our slightly modified language convenient for the purpose of dualization (cf.\ \ref{PropSST2}). Moreover,
we shall give a recipe for the transposition (\ref{PropSST11}).
\end{footnotesize}

In this section we write maps on the right. Inverse images of subsets are written on the left. Given a $\Z\Sl_n$-lattice 
$X$, we denote its dual by $X^\ast = \Hom_{\sZ}(X,\Z)$, and, accordingly, we denote the dual of a morphism $X\lraa{g}Y$ of $\Z\Sl_n$-lattices 
by $X^\ast\llaa{g^\ast} Y^\ast$. Let $\lambda$ and $\mu$ be arbitrary partitions of $n$. A {\it correspondence} from $\mu$ to $\lambda$ is a bijection 
\[
[\mu] \lraisoa{\phi} [\lambda].
\]
The set of correspondences from $\mu$ to $\lambda$ is denoted by $\Phi^{\mu,\lambda}$. We denote by 
\[
\{\phi\} := \phi^{-1}\pi^\mu_R : [\lambda] \lra \N
\]
what is called a $\lambda$-tableau of type $\mu$ in [J 78, 13.1] and what in our context might also be called the {\it correspondoid} attached to the correspondence 
$\phi\in\Phi^{\mu,\lambda}$. The set of correspondoids from $\mu$ to $\lambda$ is denoted by $\{\Phi^{\mu,\lambda}\}$. {\sc James} uses the bijection from the set of correspondoids
$\{\Phi^{\mu,\lambda}\}$ to the set of $\mu$-tabloids that, for a fixed $\lambda$-tableau $[a]$, is given by
\[
\{\phi \} = \phi^{-1}\pi^\mu_R \;\;\lraiso\;\; [a]^{-1}\phi^{-1}\pi^\mu_R = \{ \phi [a]\}
\]
as identification. Let 
\[
\begin{array}{rcl}
R_\lambda & := & \{ \rho\in\Phi^{\lambda,\lambda}\; |\; \rho\pi^\lambda_R = \pi^\lambda_R\} \\ 
C_\lambda & := & \{ \kappa\in\Phi^{\lambda,\lambda}\; |\; \kappa\pi^\lambda_C = \pi^\lambda_C\} \\ 
\end{array}
\]
be the {\it row stabilizer} resp.\ the {\it column stabilizer} of $\lambda$. Given a correspondence $\phi\in\Phi^{\mu,\lambda}$, we write
\[
r_\phi := \# (R_\lambda\cap\phi^{-1} R_\mu \phi)
\]
Note that $r_\phi = r_{\phi^{-1}}$. Given a correspondence $\phi\in\Phi^{\mu,\lambda}$ and a $\lambda$-tableau $[a]$, conjugation with the tableau $[a]$ yields a bijection 
$R_\lambda\cap\phi^{-1} R_\mu \phi \lraiso R_{[a]} \cap R_{\phi [a]}$.

\begin{Lemma}
\label{LemSST1}
Let $\phi\in\Phi^{\mu,\lambda}$ be a correspondence. The map
\[
\begin{array}{rcl}
M^\lambda & \lraa{\Theta_\phi} & M^\mu \\
\{ a\}    & \lra               & r_\phi^{-1} \sumd{\rho\in R_{[a]}} \{ \phi [a]\}\rho = r_\phi^{-1} \sumd{\rho\in R_\lambda} \{ \phi\rho[a]\} \\
\end{array}
\]
is well defined and $\Z\Sl_n$-linear. For $\rho\in R_\mu$ and $\rho'\in R_\lambda$, we have $\Theta_{\rho\phi\rho'} = \Theta_\phi$. In particular, $\Theta_\phi$ and $\Theta_{\phi^{-1}}$
both depend only on the correspondoid $\{\phi\}$. 

\rm
First, we note that the element $\sumd{\rho\in R_{[a]}} \{ \phi [a]\}\rho\in M^\mu$ is divisible by $r_\phi = \#(R_{[a]} \cap R_{\phi [a]})$ since 
$\{ \phi [a]\}\rho = \{ \phi [a]\}\rho_0\rho$ for $\rho_0\in R_{[a]} \cap R_{\phi [a]}$. Furthermore, since $R_{[a]\sigma} = (R_{[a]})^\sigma$ for $\sigma\in\Sl_n$, the map
\[
\begin{array}{rcl}
F^\lambda & \lra & M^\mu \\
{[a]}     & \lra & r_\phi^{-1} \sumd{\rho\in R_{[a]}} \{ \phi [a]\}\rho \\
\end{array}
\]
is $\Z\Sl_n$-linear. Finally, since $\sumd{\rho\in R_{[a]}} \{ \phi [a]\}\rho = \sumd{\rho\in R_{[a]\sigma}} \{ \phi [a]\sigma\}\rho$ for $\sigma\in R_{[a]}$, this map factors over
\[
\begin{array}{rcl}
F^\lambda & \lra & M^\lambda \\
{[a]}     & \lra & \{ a\}. \\
\end{array}
\]
Suppose given $\rho\in R_\mu$ and $\rho'\in R_\lambda$. We have
\[
\begin{array}{rcl}
\{a\}\Theta_{\rho\phi\rho'}
& = & r_{\rho\phi\rho'}^{-1} \sumd{\rho''\in R_\lambda} \{ \rho\phi\rho'\rho''[a]\} \\
& = & r_\phi^{-1} \sumd{\rho''\in R_\lambda} \{ \phi\rho''[a]\} \\
& = & \{a\}\Theta_\phi. \\
\end{array}
\]
\end{Lemma}

A $\Sl_n$-invariant bilinear form on $M^\lambda$ is given by 
\[
(\{ a\},\{ b\}) = \left\{\begin{array}{ll}
                  1 & \mbox{ for } \{ a\} = \{ b\} \\
                  0 & \mbox{ for } \{ a\} \neq \{ b\}, \\
                  \end{array}\right.
\]
$\{ a\},\{ b\}$ being $\lambda$-tabloids. The resulting dualization isomorphism is denoted by
\[
\begin{array}{rcl}
M^\lambda & \lraisoa{\dell^\lambda} & M^{\lambda,\ast} \\
\{ a\}    & \lra                    & (\{ a\},-). \\
\end{array}
\]
We also write such a {\it $\lambda$-cotabloid} as $\{ a\}^\ast := (\{ a\},-)$.

\begin{Lemma}
\label{PropSST2}
Given a correspondence $\phi\in\Phi^{\mu,\lambda}$, the diagram
\begin{center}
\begin{picture}(250,250)
\put(   0, 200){$M^{\lambda,\ast}$}
\put( 180, 210){\vector(-1,0){90}}
\put( 120, 230){$\scm\Theta_\phi^\ast$}
\put( 200, 200){$M^{\mu,\ast}$}
\put(  20,  50){\vector(0,1){130}}
\put(   5, 100){$\scm\wr$}
\put(  30, 100){$\scm\dell^\lambda$}
\put( 220,  50){\vector(0,1){130}}
\put( 205, 100){$\scm\wr$}
\put( 230, 100){$\scm\dell^\mu$}
\put(   0,   0){$M^\lambda$}
\put( 180,  10){\vector(-1,0){100}}
\put( 120,  30){$\scm\Theta_{\phi^{-1}}$}
\put( 200,   0){$M^\mu$}
\end{picture}
\end{center}
commutes. 

\rm
The dual of $\Theta_\phi$ sends the $\mu$-cotabloid $\{ b\}^\ast$ to the map 
\[
\begin{array}{rcl}
M^\lambda & \lraa{\{ b\}^\ast\Theta_\phi^\ast} & \Z \\
\{ a\}    & \lra & (\{ b\},\{ a\}\Theta_\phi) \\
          & =    & r_\phi^{-1} \sumd{\rho\in R_{[a]}} (\{ b\},\{\phi [a]\}\rho) \\
          & =    & r_\phi^{-1} \sumd{\rho\in R_{[a]}}\;\;\sumd{\sigma\in R_{[b]}} \left\{\begin{array}{ll}
                                                                                  1 & \mbox{for } [b]\sigma = \phi [a]\rho \\ 
                                                                                  0 & \mbox{else}
                                                                                  \end{array}\right\} \\
          & =    & r_{\phi^{-1}}^{-1} \sumd{\sigma\in R_{[b]}}\;\;\sumd{\rho\in R_{[a]}} \left\{\begin{array}{ll}
                                                                                         1 & \mbox{for } \phi^{-1}[b]\sigma = [a]\rho \\ 
                                                                                         0 & \mbox{else}
                                                                                         \end{array}\right\} \\
          & =    & r_{\phi^{-1}}^{-1} \sumd{\sigma\in R_{[b]}} (\{ \phi^{-1}[b]\}\sigma,\{ a\}) \\
          & =    & (\{ b\}\Theta_{\phi^{-1}},\{ a\}). \\
\end{array}
\]
\end{Lemma}

We denote the inclusion that is given by the definition of the Specht lattice by $S^\lambda\lraa{\iota^\lambda} M^\lambda$.

\begin{Remark}
\label{RemSST2_5}
Let $\nu$ be still another partition of $n$, let $\phi\in\Phi^{\mu,\lambda}$, let $\psi\in\Phi^{\nu,\mu}$ and let $[a]$ be a $\lambda$-tableau. We obtain
\[
\{ a\}\Theta_\phi\Theta_\psi = \sumd{\rho\in R_{\mu}} \fracd{r_{\psi\rho\phi}}{r_\phi r_\psi} \{ a\}\Theta_{\psi\rho\phi}. 
\]

\rm
In fact,
\[
\begin{array}{rcl}
\{ a\}\Theta_\phi\Theta_\psi
& = & r_\phi^{-1}\cdot\Big(\sumd{\rho'\in R_{[a]}} \{\phi [a]\rho'\}\Big)\Theta_\psi \\
& = & r_\phi^{-1}r_\psi^{-1}\cdot\sumd{\rho'\in R_{[a]}}\;\;\sumd{\rho\in R_{\phi [a]\rho'}} \{\psi\phi [a]\rho'\rho\} \\
& = & r_\phi^{-1}r_\psi^{-1}\cdot\sumd{\rho\in R_{\phi [a]}}\;\;\sumd{\rho'\in R_{[a]}} \{(\psi\phi [a]\rho [a]^{-1}\phi^{-1}\phi)[a]\}\rho' \\
& = & r_\phi^{-1}r_\psi^{-1}\cdot\sumd{\rho\in R_{\phi [a]}} r_{\psi\phi [a]\rho [a]^{-1}\phi^{-1}\phi} \{ a\}\Theta_{\psi\phi [a]\rho [a]^{-1}\phi^{-1}\phi}  \\
& = & r_\phi^{-1}r_\psi^{-1}\cdot\sumd{\rho\in R_{\mu}} r_{\psi\rho\phi} \{ a\}\Theta_{\psi\rho\phi}.  \\
\end{array}
\]
\end{Remark}

\begin{Remark}
\label{RemSST2_6}
$X$ being a $\Z\Sl_n$-lattice, we may identify
\[
\begin{array}{rcl}
X^{\ast,-}              & \lraiso & X^{-,\ast} \\
f\ts a                  & \lra    & (x\ts b \lra xf\cdot ab) \\
(x\lra (x\ts 1)g) \ts 1 & \lla    & g.
\end{array}
\]
\end{Remark}

\begin{Proposition}[{[J 78, 6.7]}]
\label{PropSST3}
The map
\[
\begin{array}{rcl}
S^{\lambda'}              & \lraisoa{\delr^\lambda} & S^{\lambda,\ast,-} \\
\spi{a_{\lambda'}}\sigma  & \lra                    & (\{ (a_{\lambda'})'\}^{\ast}|_{S^\lambda} \ts 1)\sigma, \\
\end{array}
\]
where $\sigma\in\Sl_n$, is a well defined $\Z\Sl_n$-linear isomorphism. The diagram
\begin{center}
\begin{picture}(250,250)
\put( -50, 200){$M^{\lambda,\ast,-}$}
\put(  70, 210){\vector(1,0){110}}
\put( 110, 220){$\scm\iota^{\lambda,\ast,-}$}
\put( 200, 200){$S^{\lambda,\ast,-}$}
\put(  20,  50){\vector(0,1){130}}
\put(   5, 100){$\scm\wr$}
\put(  30, 100){$\scm\dell^{\lambda,-}$}
\put( 220,  50){\vector(0,1){130}}
\put( 205, 100){$\scm\wr$}
\put( 230, 100){$\scm\delr^\lambda$}
\put( -30,   0){$M^{\lambda,-}$}
\put(  70,  10){\vector(1,0){110}}
\put( 110,  30){$\scm\nu_S^{\lambda'}$}
\put( 200,   0){$S^{\lambda'}$}
\end{picture}
\end{center}
commutes.

\rm
We abbreviate $[a] = [(a_{\lambda'})']$ in the course of the proof and claim that the kernels of the maps 
\[
\begin{array}{rcl}
M^\lambda    & \lraa{\alpha} & S^{\lambda',-} \\
\{ a\}\sigma & \lra          & \spi{a'}\sigma\ts \eps_\sigma \\
M^\lambda    & \lraa{\beta}  & S^{\lambda,\ast} \\
\{ a\}\sigma & \lra          & (\{ a\}\sigma) \delta^\lambda\iota^{\lambda,\ast} = \{a\sigma\}^\ast|_{S^\lambda} \\
\end{array}
\]
coincide. Suppose given an element $\sumd{\sigma\in\Sl_n}x_\sigma\{ a\}\sigma\in M^\lambda$, $x_\sigma\in\Z$. To be contained in the kernel
of $\alpha$ means
\[
\begin{array}{rcl}
0 
& = & \sumd{\sigma\in\Sl_n}x_\sigma\spi{a'}\sigma\ts \eps_\sigma \\
& = & \sumd{\sigma\in\Sl_n,\; \rho\in R_{[a]}} x_\sigma\{ a'\}\rho\eps_\rho\sigma\ts \eps_\sigma \\
& = & \sumd{\sigma\in\Sl_n,\; \rho\in R_{[a]}} x_{\rho^{-1}\sigma}\{ a'\}\sigma\ts \eps_\sigma. \\
\end{array}
\]
Comparing coefficients, this is equivalent to saying that for all $\tau\in\Sl_n$ we obtain
\[
\sumd{\kappa\in C_{[a]},\; \rho\in R_{[a]}} x_{\rho^{-1}\kappa\tau}\eps_{\kappa\tau} = 0.
\]
To be contained in the kernel of $\beta$ amounts to 
\[
\begin{array}{rcl}
0 
& = & (\sumd{\sigma\in\Sl_n} x_\sigma\{ a\}\sigma,\spi{a}\tau) \\
& = & (\sumd{\sigma\in\Sl_n} x_\sigma\{ a\}\sigma\tau^{-1}\sumd{\kappa\in C_{[a]}}\kappa\eps_\kappa,\{ a\}) \\
& = & (\sumd{\kappa\in C_{[a]},\;\rho\in\Sl_n} x_{\rho\kappa^{-1}\tau}\{ a\}\rho\eps_\kappa,\{ a\}) \\
\end{array}
\] 
for all $\tau\in\Sl_n$, i.e.\ that for all $\tau\in\Sl_n$
\[
\sumd{\kappa\in C_{[a]},\; \rho\in R_{[a]}} x_{\rho\kappa^{-1}\tau}\eps_\kappa = 0.
\]
Since $\alpha$ is surjective, we may conclude that the map $S^{\lambda',-}\lraa{\delr^\lambda,-} S^{\lambda,\ast}$ is well defined and injective. To prove surjectivity, 
it remains to be shown that $\beta = \delta^\lambda\iota^{\lambda,\ast}$ is surjective, i.e.\ that $\iota^\lambda$ is a pure monomorphism. But the tuple of standard tableaux gives a 
basis of the Specht module over $\F_p$ for any prime $p$ [J 78, 8.4].
\end{Proposition}

\begin{Lemma}
\label{LemSST10}
Denoting evaluation by $X\lraa{\seva} X^{\ast\ast}$, $X$ a $\Z\Sl_n$-lattice, we may dualize to obtain
\[
\begin{array}{rclclclcl}
(M^\lambda & \lraa{\seva} & M^{\lambda,\ast\ast} & \lraa{\dell^{\lambda,\ast}}   & M^{\lambda,\ast})    & = & (M^\lambda & \lraa{\dell^\lambda}    & M^{\lambda,\ast}) \\
(S^\lambda & \lraa{\seva} & S^{\lambda,\ast\ast} & \lraa{\delr^{\lambda,-,\ast}} & S^{\lambda',-,\ast}) & = & (S^\lambda & \lraa{\delr^{\lambda'}} & S^{\lambda',-,\ast}), \\
\end{array}
\]
provided $[a_{\lambda'}] = [(a_\lambda)']$. Ignoring this condition, the second equality holds up to sign.

\rm
Suppose given $\lambda$-tableaux $[a]$ and $[b]$. Let $\sigma = [a][b]^{-1}$. We have
\[
(\{ a\},\spi{b}) = \left\{\begin{array}{ll}
                   0           & \mbox{for } \sigma\not\in R_\lambda C_\lambda \\
                   \eps_\kappa & \mbox{for } \sigma = \rho\kappa,\;\rho\in R_\lambda,\;\kappa\in C_\lambda. \\
                   \end{array}\right.
\]
Therefore,
\[
\mbox{(\#)} \hspace{2cm} (\{ a\},\spi{b}) = \eps_\sigma (\{ b'\},\spi{a'}).
\]
We abbreviate $[a] = [(a_{\lambda'})'] = [a_\lambda]$. The map $S^{\lambda'}\lraa{\delr^{\lambda}}S^{\lambda,-,\ast}$ sends $\spi{a'}\tau$ to $(\{ [a]\tau\},-)\ts\eps_\tau$, where $\tau\in\Sl_n$, 
so that the composition sends
\[
\begin{array}{rclcl}
S^\lambda & \lraa{\seva} & S^{\lambda,\ast\ast} & \lraa{\delr^{\lambda,-,\ast}} & S^{\lambda',-,\ast} \\
\spi{a}   & \lra         & \spi{a}\eva          & \lra                          & \Big(\spi{a'}\tau\ts 1 \lra \eps_\tau (\{[a]\tau\},\spi{a}) \auf{\mbox{\scr (\#)}}{=} (\{ a'\},\spi{a'}\tau)\Big).\\
\end{array}
\]
\end{Lemma}

A correspondoid $\{\phi\}\in\{\Phi^{\mu,\lambda}\}$ is called {\it semistandard} in case the following conditions hold.
\begin{itemize}
\item[(i)]   For $i\ti j,\; i\ti j'\in [\lambda]$ with $j < j'$, we have $(i\ti j)\{\phi\} \leq (i\ti j')\{\phi\}$.
\item[(ii)]  For $i\ti j,\; i'\ti j\in [\lambda]$ with $i < i'$, we have $(i\ti j)\{\phi\} < (i'\ti j)\{\phi\}$.
\end{itemize}
Let $\{\Phi^{\mu,\lambda}\}_{\mbox{\scr\rm sst}}\tm \{\Phi^{\mu,\lambda}\}$ denote the subset of {\it semistandard correspondoids.}

We attach to each correspondence $\phi\in\Phi^{\mu,\lambda}$ its {\it column distribution}
\[
\begin{array}{rclcl}
\N &\ti &\N & \lraa{|\phi|} & \N \\
j  &\ti & k & \lra          & |\phi|_{k,j} := \#\Big(\{\phi\}^{-1}(j) \cap (\N\ti\{ k\})\Big).
\end{array}
\]
The map from the correspondences to the column distribution factors over the correspondoids, $(\phi\lra |\phi|) = (\phi\lra \{\phi\}\lra |\phi|)$. For $\phi\in\Phi^{\mu,\lambda}$, $\rho\in R_\mu$ 
and $\kappa\in C_\lambda$, we have $|\rho\phi\kappa| = |\phi|$ because $\{\rho\phi\kappa\} = \{\phi\kappa\}$ and because of the bijection, $j\ti k\in \N\ti\N$,
\[
\{\phi\}^{-1}(j) \cap (\N\ti\{ k\}) \lraiso \{\phi\kappa\}^{-1}(j) \cap (\N\ti\{ k\}) 
\]
which is given by restriction of $\kappa$. Note that the column distribution may be regarded as column equivalence class in the sense of [J 78, 13.8].

A total order on the set of column distributions can be defined as follows. We order $\N\ti\N$ lexicographically via $j\ti k < j'\ti k'$ if $j < j'$ or ($j = j'$ and $k < k'$). Suppose given 
correspondences $\phi,\w\phi\in\Phi^{\mu,\lambda}$, $|\phi|\neq |\w\phi|$. Let $j\ti k$ be minimal with $|\phi|_{k,j} \neq |\w\phi|_{k,j}$. We say that $|\phi| < |\w\phi|$ in case 
$|\phi|_{k,j} < |\w\phi|_{k,j}$. 

The following arguments, needed to obtain (\ref{ThSST6}), are taken from [J 78], and we reproduce them here in a slightly adapted manner. I thank {\sc M.\ H\"arterich} for an explanation. 

\begin{Lemma}
\label{LemSST4}
For $\{\phi\}\in\{\Phi^{\mu,\lambda}\}_{\mbox{\scr\rm sst}}$ and $\rho\in R_\lambda$, we have $|\phi| \geq |\phi\rho|$.

\rm
Suppose $|\phi| \neq |\phi\rho|$ and let $j\ti k$ be minimal with $|\phi|_{k,j} \neq |\phi\rho|_{k,j}$. By induction on $j'\ti k'$, we see that 
\[
(\ast)\hspace*{1cm} \{\phi\rho\}^{-1}(j')\cap (\N\ti\{ k'\}) = \{\phi\}^{-1}(j') \cap (\N\ti\{ k'\}) 
\]
for $1\ti 1 \leq j'\ti k' < j\ti k$. In fact, assume that the two inverse images of the value $j'$ in column $k'$ differ in spite of having the same cardinality. Then there exists 
$i\ti k'\in [\lambda]$, $i\ti k'' := (i\ti k')\rho^{-1}$, such that $j' = (i\ti k'')\{\phi\}\neq (i\ti k')\{\phi\} = : j''$. Assume that $j' > j''$. But since 
$i\ti k'\in \{\phi\}^{-1}(j'')\ohne \{\phi\rho\}^{-1}(j'')$, this yields a contradiction to ($\ast$) by induction hypothesis, so that we obtain $j' < j''$. 

Let $i\ti k'^{(s)} := (i\ti k')\rho^{-s}$ for $s\geq 0$. We claim that $(i\ti k'^{(s)})\{\phi\} = j'$ for all $s\geq 1$, thus deriving a contradiction that establishes $(\ast)$. We perform an 
induction on $s$. But $(i\ti k'^{(s)})\{\phi\} = j'$ implies by semistandardness of $\{\phi\}$ that $k'^{(s)} < k'$. Thus ($\ast$) applies by induction hypothesis and shows that 
$(i\ti k'^{(s+1)})\{\phi\} = (i\ti k'^{(s)})\{\phi\rho\} = j$. 

Assume that $|\phi|_{k,j} < |\phi\rho|_{k,j}$. Then there exists $i\ti k\in [\lambda]$, $i\ti k' := (i\ti k)\rho^{-1}$, such that 
$j = (i\ti k')\{\phi\}\neq (i\ti k)\{\phi\} = : j'$. Assume that $j > j'$. But since 
$i\ti k\in \{\phi\}^{-1}(j')\ohne \{\phi\rho\}^{-1}(j')$, this yields a contradiction to ($\ast$) and shows $j < j'$. 

Let $i\ti k^{(s)} := (i\ti k)\rho^{-s}$ for $s\geq 0$. We claim that $(i\ti k^{(s)})\{\phi\} = j$ for all $s\geq 1$, thus deriving a contradiction that shows $|\phi| \geq |\phi\rho|$. We perform 
an induction on $s$. But $(i\ti k^{(s)})\{\phi\} = j$ implies by semistandardness of $\{\phi\}$ that $k^{(s)} < k$. Thus ($\ast$) shows that 
$(i\ti k^{(s+1)})\{\phi\} = (i\ti k^{(s)})\{\phi\rho\} = j$.
\end{Lemma}

\begin{Lemma}
\label{LemSST5}
Restricted to $\{\Phi^{\mu,\lambda}\}_{\mbox{\scr\rm sst}}$, the map $\{\phi\}\lra |\phi|$ becomes injective.

\rm
To determine a semistandard correspondoid $\{\phi\}$, it suffices to know the set of its values on column $k$ for each $k\geq 1$, which can be written as 
$\{ j\in \N \; |\; |\phi|_{k,j}\geq 1\}$.
\end{Lemma}

\begin{Lemma}
\label{SST5_1}
Assume $\lambda$ to be a $2$-regular partition of $n$, i.e.\ assume that $i < i'$ implies $\lambda_i > \lambda_{i'}$ for $i,i'\in [1,\lambda'_1]$.
Suppose given $\kappa,\kappa'\in C_\lambda$. Let 
\[
\begin{array}{rcl}
[\lambda] & \lraisoa{\iota} &²[\lambda] \\
i\ti k    & \lra            & i\ti(\lambda_i - k). \\
\end{array}
\]
If $\kappa\iota\kappa'\in R_\lambda$ then $\kappa = 1$, hence also $\kappa' = 1$. 

\rm
Assume that $\kappa\neq 1$, i.e.\ assume given some position $i\ti k\in [\lambda]$ such that $i\ti k \neq (i\ti k)\kappa$. First, we choose $i$ minimal such that there exists
a $k$ such that $i\ti k \neq (i\ti k)\kappa$. Having chosen our $i$, we choose $k$ maximal such that $i\ti k \neq (i\ti k)\kappa$. We write $i'\ti k := (i\ti k)\kappa$ and remark
that $i < i'$. We obtain
\[
\begin{array}{rclcrclcrclcrcl}
i    & \!\!\ti\!\! & k                          & \!\!\lraa{\kappa}\!\! & i' & \!\!\ti\!\! &  k                          & \!\!\lraa{\iota}\!\! & i' & \!\!\ti\!\! & (\lambda_{i'} - k) & 
\!\!\lraa{\kappa'}\!\! & i''' & \!\!\ti\!\! & (\lambda_{i'} - k) \\  
i    & \!\!\ti\!\! & (k+\lambda_i-\lambda_{i'}) & \!\!\lraa{\kappa}\!\! & i  & \!\!\ti\!\! &  (k+\lambda_i-\lambda_{i'}) & \!\!\lraa{\iota}\!\! & i  & \!\!\ti\!\! & (\lambda_{i'} - k) & 
\!\!\lraa{\kappa'}\!\! & i''  & \!\!\ti\!\! & (\lambda_{i'} - k), \\  
\end{array}
\]
where $i''\neq i'''$, contradicting $\kappa\iota\kappa'\in R_\lambda$.  
\end{Lemma}

\begin{Lemma}
\label{SST5_2}
Suppose given $\sigma\not\in R_\lambda\cdot C_\lambda$. There exists a transposition $\kappa_\sigma\in C_\lambda$ such that $\sigma\kappa_\sigma\sigma^{-1}\in R_\lambda$.

\rm
Assume that the conclusion does not hold, i.e.\ that $\sigma^{-1}$ sends different entries of each given column to different rows. We claim that this implies that 
$\sigma\in  R_\lambda\cdot C_\lambda$. Given a permutation $\tau$ of $[\lambda]$, we let $k_\tau$ be the maximal column position $k\geq 0$ such that $\tau$ becomes the identity when restricted to 
$[\lambda]\cap (\N\ti [1,k])$. We perform a downwards induction on $k_\sigma$, starting with $k_\sigma = \lambda_1$.

Different entries of column $k_\sigma + 1$ are mapped to different rows. Hence there is a column permutation $\kappa\in C_\lambda$ acting trivially except perhaps in column $k_\sigma + 1$ such that 
$(i\ti (k_\sigma + 1))\kappa\sigma^{-1}\pi_R^\lambda = i$ for $i\in [1,\lambda'_{k_\sigma + 1}]$. Thus there is a row permutation $\rho\in R_\lambda$ such that 
$k_{\kappa\sigma^{-1}\rho} > k_{\sigma^{-1}}$, i.e.\ $k_{\rho^{-1}\sigma\kappa^{-1}} > k_{\sigma}$. By induction, $\rho^{-1}\sigma\kappa^{-1} \in R_\lambda\cdot C_\lambda$ ensues.
\end{Lemma}

\begin{Lemma}
\label{SST5_3}
Suppose given a semistandard correspondoid $\{\phi\}\in\{\Phi^{\mu,\lambda}\}_{\mbox{\scr\rm sst}}$ and a $\lambda$-tableau $[a]$. Writing
\[
\spi{a}\Theta_\phi = \sumd{\{\psi\}\in\{\Phi^{\mu,\lambda}\}} z_\psi \{\psi [a]\} \in M^\mu,
\]
we obtain $z_{\{\phi\}} = 1$ and $z_{\{\psi\}} = 0$ for $|\psi| > |\phi|$.

\rm
Note that 
\[
\begin{array}{rcl}
\spi{a}\Theta_\phi & = & \sumd{\kappa\in C_{[a]}}\;\;\sumd{\rho\in (R_\lambda\cap\phi^{-1}R_\mu\phi)\< R_\lambda} \{ \phi\rho [a]\} \kappa\eps_\kappa \\
                   & = & \sumd{\kappa\in C_\lambda}\;\;\sumd{\rho\in (R_\lambda\cap\phi^{-1}R_\mu\phi)\< R_\lambda} \{ \phi\rho\kappa [a]\} \eps_\kappa. \\
\end{array}
\]
Since $\{\phi\}$ is semistandard, we have 
$|\phi|\auf{\mbox{\scr(\ref{LemSST4})}}{\geq} |\phi\rho| = |\phi\rho\kappa|$ for $\rho\in R_\lambda$, $\kappa\in C_\lambda$.

The equality $\{\phi\} = \{\phi\rho\kappa\}$ is equivalent to the existence of a $\rho'\in R_\mu$ such that $\kappa = \rho (\phi^{-1}\rho'\phi)$. This condition supposed 
to hold and given $i\ti k\in [\lambda]$, we obtain
\[
i\ti k \lraa{\rho} i\ti k' \mra{\phi^{-1}\rho' \phi} i'\ti k,
\] 
where $(i\ti k')\{\phi\} = (i'\ti k)\{\phi\}$. Note that the resulting map $i\lra i'$ is a bijection for each $k$. We perform an induction on $k$ to show that $i\leq i'$. Suppose $i\leq i'$ and 
hence $i = i'$ for $k_0 < k$ (supposition possibly empty). Therefore, $(i\ti k_0)\{\phi\} = (i'\ti k_0)\{\phi\} = (i\ti k_0')\{\phi\}$. Since $\{\phi\}$ is 
semistandard, this inhibits $(i\ti k)\{\phi\} > (i\ti k')\{\phi\}$, for this would mean that $\rho$ could not act on the set 
\[
\Big\{i\ti k_1\; \Big|\; k_1\in [1,\lambda_i],\; (i\ti k_1)\{\phi\} = (i\ti k')\{\phi\}\Big\} \tm \{i\}\ti [1,k-1],
\]
the inclusion given by semistandardness of $\{\phi\}$. Thus $(i\ti k)\{\phi\} \leq (i\ti k')\{\phi\} = (i'\ti k)\{\phi\}$, hence, by semistandardness of $\{\phi\}$, $i\leq i'$, forcing $i = i'$
for all $i\ti k\in [\lambda]$. We conclude that $\kappa = 1$ and that thus $\phi^{-1}\rho'\phi$ inverts $\rho$ so that $\rho$ represents the trivial coset.
\end{Lemma}

Let $m\geq 0$. $\lambda$ is called {\it $2$-singular} if there exists an $i\geq 1$ such that $\lambda_i = \lambda_{i+1} > 0$. The case of $m\geq 1$, $m\in (2)$ and $\lambda$ being $2$-singular 
shall be called the {\bf singular case.} The complementary case shall be called the {\bf regular case,} comprising in particular the case $m = 0$. 

\begin{Theorem}[{{\sc Carter}, {\sc Lusztig} [CL 74, 3.5], {\sc James} [J 78, 13.13]}]
\label{ThSST6} 
\Absatz
In the regular case, the tuple
\[
\Big(\Theta_\phi|_{S^\lambda}\ts_{\sZ}\Z/m \; \Big|\; \{\phi\}\in\{\Phi^{\mu,\lambda}\}_{\mbox{\scr\rm sst}}\Big)
\]
furnishes a $\Z/m$-linear basis of $\Hom_{\sZ\Sl_n}(S^\lambda/m,M^\mu/m)$.

In the singular case, the tuple 
\[
\Big(2\cdot\Theta_\phi|_{S^\lambda}\ts_{\sZ}\Z/m \; \Big|\; \{\phi\}\in\{\Phi^{\mu,\lambda}\}_{\mbox{\scr\rm sst}}\Big)
\]
furnishes a $\Z/(m/2)$-linear basis of $2\cdot\Hom_{\sZ\Sl_n}(S^\lambda/m,M^\mu/m)$.

\rm
{\it We claim that the respective tuple is linearly independent.} Suppose given a $\lambda$-tableau $\spi{a}$. We write down a matrix whose positions are indexed by 
$\{\Phi^{\mu,\lambda}\}_{\mbox{\scr\rm sst}} \ti \{\Phi^{\mu,\lambda}\}$, and whose entry at position $\{\phi\}\ti\{\psi\}$ is given by the multiplicity of the $\mu$-tabloid
$\{\psi[a]\}$ in $\spi{a}\Theta_\phi\in M^\mu$. It suffices to show that the $\{\Phi^{\mu,\lambda}\}_{\mbox{\scr\rm sst}} \ti \{\Phi^{\mu,\lambda}\}_{\mbox{\scr\rm sst}}$-part of
our matrix is unipotent. But this follows from (\ref{SST5_3}) in view of (\ref{LemSST5}).

{\it We claim that the respective tuple generates.} Suppose given a morphism $S^\lambda/m\lraa{\Theta}M^\mu/m$ such that $\Theta\neq 0$ in the regular case resp.\ such that $2\cdot\Theta\neq 0$
in the singular case. We fix a $\lambda$-tableau $[a]$ and write
\[
\spi{a}\Theta = \sumd{\{ \phi\}\in\{\Phi^{\mu,\lambda}\}} x_{\{\phi\}}\{\phi [a]\},
\]
where all elements of $S^\lambda$ and of $M^\mu$ which occur in this expression and in the remainder of the proof are to be read as {\it representing their residue classes modulo $m$.}
Suppose given $\{\psi\}\in\{\Phi^{\mu,\lambda}\}$ that allows for $\{\psi\tau\} = \{\psi\}$ for some $i\ti k,\; i'\ti k\in [\lambda]$, $i\neq i'$, $\tau := (i\ti k, i\ti k')$. We wish to see that 
$x_{\{\psi\}} = 0$ in the regular case and that $2\cdot x_{\{\psi\}} = 0$ in the singular case. In both cases, we may compare the coefficients of
\[
\begin{array}{rcl}
\sumd{\{ \phi\}\in\{\Phi^{\mu,\lambda}\}} x_{\{\phi\tau\}}\{\phi [a]\}
& = & \sumd{\{ \phi\}\in\{\Phi^{\mu,\lambda}\}} x_{\{\phi\}}\{\phi [a]\}(a_{k,i}, a_{k,i'}) \\
& = & \spi{a}(a_{k,i}, a_{k,i'})\Theta \\
& = & - \sumd{\{ \phi\}\in\{\Phi^{\mu,\lambda}\}} x_{\{\phi\}}\{\phi [a]\}\\
\end{array}
\]
to obtain
\[
x_{\{\psi\}} = x_{\{\psi\tau\}} = - x_{\{\psi\}}.
\]
Now suppose the partition $\lambda$ to be $2$-regular. We need to see that $x_{\{\psi\}} = 0$. Let $\iota$ be the permutation of $[\lambda]$ defined by $i\ti k\lraisoa{\iota} i\ti (\lambda_i - k)$,
$i\ti k\in [\lambda]$. We abbreviate $\iota^{[a]} := [a]^{-1}\iota [a]\in\Sl_n$ and calculate 
\[
\begin{array}{rcl}
\sumd{\kappa'\in C_{[a]}}\spi{a}\iota^{[a]}\kappa'\eps_{\kappa'}
& = & \sumd{\kappa,\kappa'\in C_{[a]}} \{a\} \kappa\iota^{[a]}\kappa'\eps_{\kappa\kappa'} \\
& = & \sumd{\kappa,\kappa'\in C_\lambda} \{\kappa\iota\kappa'[a]\} \eps_{\kappa\kappa'}, \\
\end{array}
\]
where $\kappa\in C_\lambda$ has a sign $\eps_\kappa$ as a permutation of $[\lambda]$, i.e.\ as an element of $\Sl_{[\lambda]}$. Writing the resulting expression in the form 
$\sumd{\sigma\in R_\lambda\<\Sl_{[\lambda]}} y_\sigma \{\sigma [a]\}$, $y_\sigma\in\Z$, where we choose our coset representatives $\sigma$ to lie in $C_\lambda$ whenever (uniquely) possible, we 
obtain
\[
y_\sigma = \sumd{\kappa,\kappa'\in C_\lambda,\; \kappa\iota\kappa'\in R_\lambda\sigma} \eps_{\kappa\kappa'}.
\]
In case $\sigma\in C_\lambda$, the condition on the indexing elements yields $\kappa = 1$, $\kappa' = \sigma$ and thus $y_\sigma = \eps_\sigma$ (\ref{SST5_1}). In case 
$\sigma\not\in R_\lambda\cdot C_\lambda$, we obtain, using (\ref{SST5_2}), and choosing a fixed system of representatives of $C_\lambda/\spi{\kappa_\sigma}$,
\[
\begin{array}{rcl}
y_\sigma 
& = & \sumd{\kappa,\kappa'\in C_\lambda,\; \kappa\iota\kappa'\in R_\lambda\sigma} \eps_{\kappa\kappa'} \\
& = & \sumd{\kappa\in C_\lambda,\; \kappa'\in C_\lambda/\spi{\kappa_\sigma},\; \kappa\iota\kappa'\in R_\lambda\sigma} \eps_{\kappa\kappa'} 
\;\;  + \;\;  \sumd{\kappa\in C_\lambda,\; \kappa'\in C_\lambda/\spi{\kappa_\sigma},\; \kappa\iota\kappa'\kappa_\sigma\in R_\lambda\sigma} \eps_{\kappa\kappa'\kappa_\sigma} \\
& = & 0. \\
\end{array}
\]
Hence
\[
\sumd{\kappa'\in C_{[a]}}\spi{a}\iota^{[a]}\kappa'\eps_{\kappa'} = \spi{a},
\]
and therefore
\[
\begin{array}{rcl}
\spi{a}\Theta 
& \erstgl & \sumd{\{ \phi\}\in\{\Phi^{\mu,\lambda}\}} x_{\{\phi\}}\{\phi [a]\} \\
& \zweigl & \sumd{\kappa\in C_{[a]}}\spi{a}\Theta\iota^{[a]}\kappa\eps_\kappa \\
& = & \sumd{\{\phi\}\in\{\Phi^{\mu,\lambda}\}} x_{\{\phi\}}\Big(\sumd{\kappa\in C_\lambda} \{\phi\kappa[a]\}\eps_\kappa\Big).\\
\end{array}
\]
For $\{\phi\}\in\{\Phi^{\mu,\lambda}\}$ and $\kappa\in C_\lambda$, the equality $\{\psi\} = \{\phi\kappa\}$ is equivalent to $\phi\kappa\psi^{-1}\in R_\mu$. Hence, as multiplicity of 
$\{\psi [a]\}$ in $\sumd{\kappa\in C_\lambda}\{\phi\kappa [a]\}$ we get, choosing a fixed system of representatives for $C_\lambda/\spi{\tau}$,
\[
\begin{array}{rcl}
\sumd{\kappa\in C_\lambda,\; \phi\kappa\psi^{-1}\in R_\mu} \eps_\kappa
& = & \sumd{\kappa\in C_\lambda/\spi{\tau},\; \phi\kappa\psi^{-1}\in R_\mu} \eps_\kappa \;\; +\;\; \sumd{\kappa\in C_\lambda/\spi{\tau},\; \phi\kappa\tau\psi^{-1}\in R_\mu} \eps_{\kappa\tau} \\
& = & 0.\\ 
\end{array}
\]
Therefore, a comparison of 1.\ and 2.\ yields $x_{\{\psi\}} = 0$.

Given $\kappa\in C_\lambda$, we obtain 
\[
\begin{array}{rcl}
\sumd{\{\phi\}\in\{\Phi^{\mu,\lambda}\}} x_{\{\phi\kappa^{-1}\}}\{\phi [a]\}
& = & \spi{a}\Theta [a]^{-1}\kappa [a] \\
& = & \spi{a}\Theta\eps_\kappa \\
& = & \sumd{\{\phi\}\in\{\Phi^{\mu,\lambda}\}} \eps_\kappa x_{\{\phi\}}\{\phi [a]\}, \\
\end{array}
\]
whence, up to sign, $x_{\{\phi\}}$ depends only on $|\phi|$. Let $|\phi|$ be the maximal column distribution such that $x_{\{\phi\}} \neq 0$ in the regular case, resp.\ such that
$2\cdot x_{\{\phi\}} \neq 0$ in the singular case. By our assertion on $x_{\{\psi\}}$, we may assume that $(i\ti k)\{\phi\} < (i'\ti k)\{\phi\}$ for $i\ti k,\; i'\ti k\in\lambda$, $i < i'$. 
In order to show that $\{\phi\}$ is semistandard, we assume the contrary, i.e.\ the existence of positions $i\ti k,\; i\ti (k+1)\in [\lambda]$ such that $(i\ti k)\{\phi\} > (i\ti (k+1))\{\phi\}$. 
Let $\xi := a_{k,[i,\lambda'_i]}$, let $\eta := a_{k+1,[1,i]}$. For $\sigma$ in $\Sl_{\xi\cup\eta}$ but not in $\Sl_\xi\ti\Sl_\eta$, we have $|\phi| < |\phi [a]\sigma [a]^{-1}|$ since the 
minimal value of $\{\phi\}$ changing columns under $\sigma$ cannot move to the right. Summing up over cosets, a Garnir relation gives
\[
\begin{array}{rcl}
0 
& = & \sumd{\sigma\in \Sl_\xi\ti\Sl_\eta\<\Sl_{\xi\cup\eta}}\spi{a}\Theta\sigma \\ 
& = & \sumd{\{\chi\}\in\{\Phi^{\mu,\lambda}\}}\Big(\sumd{\sigma\in\Sl_\xi\ti\Sl_\eta\<\Sl_{\xi\cup\eta}} x_{\{\chi ([a]\sigma^{-1} [a]^{-1})\}}\Big)\{\chi [a]\},\\ 
\end{array}
\]
whence $\sumd{\sigma\in\Sl_\xi\ti\Sl_\eta\<\Sl_{\xi\cup\eta}} x_{\{\phi ([a]\sigma^{-1} [a]^{-1})\}} = 0$. For each $\sigma$ in $\Sl_{\xi\cup\eta}$ but not in $\Sl_\xi\ti\Sl_\eta$, maximality
of $|\phi|$ forces $x_{\{\phi ([a]\sigma [a]^{-1})\}} = 0$ in the regular case, resp.\ $2\cdot x_{\{\phi ([a]\sigma [a]^{-1})\}} = 0$ in the singular case, contradicting $x_{\{\phi\}}\neq 0$ 
in the regular case, resp.\ $2\cdot x_{\{\phi\}}\neq 0$ in the singular case. Hence $\{\phi\}$ is semistandard, and its column distribution $|\phi|$
shall be called the leading term of $\Theta$. Given a correspondoid $\{\chi\}\in\{\Phi^{\mu,\lambda}\}$, we write $\b\Theta_\chi := \Theta_\chi|_{S^\lambda}\ts_{\sZ} \Z/m$.

We consider the regular case. By (\ref{SST5_3}), the difference $\Theta - x_{\{\phi\}}\cdot \b\Theta_\phi$ has either a strictly smaller leading term than $\Theta$ or it vanishes. Both alternatives 
allow to assume, by induction on the leading term of $\Theta$, or, respectively, directly, that $\Theta - x_{\{\phi\}}\cdot\b\Theta_\phi$ is in the linear 
span of the tuple given above in the regular case. 

We consider the singular case. By (\ref{SST5_3}), the difference $\Theta - x_{\{\phi\}}\cdot \b\Theta_\phi$ has either a strictly smaller leading term than $\Theta$ or it vanishes under 
multiplication with $2$. Both alternatives allow to assume, by induction on the leading term of $\Theta$ or, respectively, directly, that 
$2\Theta - x_{\{\phi\}}\cdot 2\b\Theta_\phi$ is in the linear span of the tuple given above in the singular case.
\end{Theorem}

\begin{Corollary}
\label{CorSST7}
In the regular case, the tuple
\[
\Big(\Theta^-_{\phi^{-1}} \nu_S^{\lambda'}\ts_{\sZ}\Z/m \; \Big|\; \{\phi\}\in\{\Phi^{\mu,\lambda}\}_{\mbox{\scr\rm sst}}\Big)
\]
furnishes a $\Z/m$-linear basis of $\Hom_{\sZ\Sl_n}(M^{\mu,-}/m,S^{\lambda'}/m)$.

In the singular case, the tuple 
\[
\Big(2\cdot\Theta^-_{\phi^{-1}} \nu_S^{\lambda'}\ts_{\sZ}\Z/m \; \Big|\; \{\phi\}\in\{\Phi^{\mu,\lambda}\}_{\mbox{\scr\rm sst}}\Big)
\]
furnishes a $\Z/(m/2)$-linear basis of $2\cdot\Hom_{\sZ\Sl_n}(M^{\mu,-}/m,S^{\lambda'}/m)$.

\rm
We write $(-)^\ast := \Hom_{\sZ}(-,\Z/m)$, neglect to denote $(-)\ts_{\sZ}\Z/m$ for maps and obtain $\Z/m$-linear isomorphisms
\[
\begin{array}{rrl}
\Hom_{\sZ\Sl_n}(S^\lambda/m,M^\mu/m) & \mraisoa{\dell^\mu\cdot (-)^\ast} & \Hom_{\sZ\Sl_n}(M^\mu/m,S^{\lambda,\ast}/m)   \\
\iota^\lambda\Theta_\phi             & \mra{\;}                          & \dell^\mu\Theta_\phi^\ast \iota^{\lambda,\ast} 
                                                                           \;\;\auf{\mbox{\scr (\ref{PropSST2})}}{=}\;\; \Theta_{\phi^{-1}}\dell^\lambda \iota^{\lambda,\ast} \\  
                                     & \cdots\mraisoa{(-)^-\cdot (\delr^\lambda)^{-1}} & \Hom_{\sZ\Sl_n}(M^{\mu,-}/m,S^{\lambda'}/m)  \\ 
                                     & \cdots\mra{\;}                                  & \Theta_{\phi^{-1}}^- \dell^{\lambda,-}  \iota^{\lambda,\ast,-}  (\delr^\lambda)^{-1} 
                                                                                         \;\;\auf{\mbox{\scr (\ref{PropSST3})}}{=}\;\; \Theta_{\phi^{-1}}^- \nu_S^{\lambda'}  \\
\end{array}
\]
which send the bases of (\ref{ThSST6}) to the tuples just described.
\end{Corollary}

\begin{Corollary}
\label{CorSST8}
In the regular case, the induced maps
\[
\begin{array}{rcl}
\Hom_{\sZ\Sl_n}(M^\lambda/m,M^\mu/m)         & \lraa{\iota^\lambda(-)}    & \Hom_{\sZ\Sl_n}(S^\lambda/m,M^\mu/m) \\
\Hom_{\sZ\Sl_n}(M^{\mu,-}/m,M^{\lambda,-}/m) & \lraa{(-)\nu_S^{\lambda'}} & \Hom_{\sZ\Sl_n}(M^{\mu,-}/m,S^{\lambda'}/m) \\
\end{array}
\]
are surjective. In the singular case, their cokernels are annihilated by multiplication by $2$. Or, slightly more precise, when passing to $2\cdot\Hom$, the maps become surjective.
\end{Corollary}

\begin{Question}
\label{QuSST9}
\rm
I do not know a basis of $\Hom_{\sF_2\Sl_n}(S^\lambda/2,M^\mu/2)$ in case $\lambda$ is $2$-singular.
\end{Question}

\begin{Proposition}
\label{PropSST11}
Suppose given a $\Z\Sl_n$-linear map of the form
\[
\begin{array}{rcl}
S^\lambda/m & \lraa{f} & S^\mu/m \\
\spi{a}     & \lra     & \sumd{\{\phi\}\in\{\Phi^{\mu',\lambda'}\}}x_{\{\phi\}}\{ a'\}^-\Theta_\phi^-\nu_S^\mu, \\
\end{array}
\]
where $x_{\{\phi\}}\in\Z/m$. Neglecting to denote $(-)\ts_{\sZ}\Z/m$ for maps and using the isomorphism $\delr$ from (\ref{PropSST3}), we obtain the {\rm transposition isomorphism}
\[
\begin{array}{rcl}
\Hom_{\sZ\Sl_n}(S^\lambda/m,S^\mu/m) & \lraisoa{\mbox{\rm\scr tra}} & \Hom_{\sZ\Sl_n}(S^{\mu'}/m,S^{\lambda'}/m) \\
f                                    & \lra                         & f^\trp := \delr^\mu f^{\ast,-} (\delr^\lambda)^{-1}. \\
\end{array}
\]
This definition can be rephrased via the formula
\[
\spi{b'} f^\trp = \sumd{\{\phi\}\in\{\Phi^{\mu',\lambda'}\}} x_{\{\phi\}}\spi{b'}\Theta_{\phi^{-1}},
\]
note that $\Theta_{\phi^{-1}}$ is applicable to $\spi{b'}\in S^{\mu'}/m\tm M^{\mu'}/m$. We have $(f^\trp)^\trp = f$.

\rm
We dualize via $(-)^\ast := \Hom_{\sZ}(-,\Z/m)$ to obtain the following diagram.
\begin{center}
\begin{picture}(550,400)
\put(  80, 350){$M^{\lambda',\ast}/m$}
\put( 430, 360){\vector(-1,0){180}}
\put( 275, 380){$\scm\sum x_{\{\phi\}}\Theta_\phi^\ast$}
\put( 450, 350){$M^{\mu',\ast}/m$}

\put( 170, 260){\vector(0,1){70}}
\put( 170, 195){\line(0,1){10}}
\put( 180, 280){$\scm\nu_S^{\lambda,\ast,-}$}
\put( 470, 195){\vector(0,1){135}}
\put( 480, 260){$\scm\nu_S^{\mu,\ast,-}$}

\put(  80, 150){$S^{\lambda,\ast,-}/m$}
\put( 430, 160){\line(-1,0){100}}
\put( 310, 160){\vector(-1,0){60}}
\put( 345, 170){$\scm f^{\ast,-}$}
\put( 450, 150){$S^{\mu,\ast,-}/m$}

\put(  50, 250){\vector(1,1){90}}
\put(  50, 285){$\scm\dell^{\lambda'}$}
\bezier{20}(85,270)(85,280)(95,280)
\bezier{20}(95,280)(105,280)(105,290)
\put( 350, 250){\vector(1,1){90}}
\put( 350, 285){$\scm\dell^{\mu'}$}
\bezier{20}(385,270)(385,280)(395,280)
\bezier{20}(395,280)(405,280)(405,290)
\put(  50,  50){\vector(1,1){90}}
\put(  55,  85){$\scm\delr^\lambda$}
\bezier{20}(85,70)(85,80)(95,80)
\bezier{20}(95,80)(105,80)(105,90)
\put( 350,  50){\vector(1,1){90}}
\put( 355,  85){$\scm\delr^\mu$}
\bezier{20}(385,70)(385,80)(395,80)
\bezier{20}(395,80)(405,80)(405,90)

\put( -40, 200){$M^{\lambda'}/m$}
\put( 280, 210){\vector(-1,0){180}}
\put( 115, 230){$\scm\sum x_{\{\phi\}}\Theta_{\phi^{-1}}$}
\put( 300, 205){$M^{\mu'}/m$}

\put(  20,  50){\vector(0,1){130}}
\put( -22, 100){$\scm\iota^{\lambda'}$}
\put( 320,  50){\vector(0,1){130}}
\put( 278, 100){$\scm\iota^{\mu'}$}

\put( -30,   0){$S^{\lambda'}/m$}
\put( 280,  10){\vector(-1,0){180}}
\put( 185,  20){$\scm f^\trp$}
\put( 300,   0){$S^{\mu'}/m$}
\end{picture}
\end{center}
The commutativity on top follows from (\ref{PropSST2}). On the left and on the right, commutativity follows by dualizing the commutativity in (\ref{PropSST3}) in the sense of (\ref{LemSST10}). 
The commutativity in the back is the dualized assumption on $f$ and the commutativity at the bottom is the definition of $f^\trp$. The claimed formula is equivalent to the commutativity in the 
front face. Moreover, (\ref{LemSST10}) allows to conclude that
\[
\begin{array}{rcl}
(f^\trp)^\trp
& = & \delr^{\lambda'}\cdot f^{\trp,\ast,-} \cdot (\delr^{\mu'})^{-1} \\
& = & \delr^{\lambda'}\cdot  (\delr^{\lambda,\ast,-})^{-1} \cdot f^{\ast\ast}\cdot\delr^{\mu,\ast,-}\cdot (\delr^{\mu'})^{-1} \\
& = & \eva \cdot f^{\ast\ast}\cdot \eva^{-1} \\
& = & f. \\
\end{array}
\]
\end{Proposition}

\begin{Question}[main technical obstacle]
\label{QuSST11_1}
\rm
Given $S^\lambda/m\lraa{f}S^\mu/m$ as in (\ref{PropSST11}), I do not know a formula for $S^{\mu'}/m\lraa{f^\trp} S^{\lambda'}/m$ in terms of a factorization over 
$M^{\mu,-}/m\lraa{\nu_S^{\mu'}}S^{\mu'}/m$ of a linear combination of maps of the form $\Theta^-_\psi\nu_S^{\lambda'}$, $\{\psi\}\in\{\Phi^{\lambda,\mu}\}$. I.e., I do not know a formula `in 
terms of polytabloids instead of tabloids'.
\end{Question}
\textheight24.7cm        
\subsection{Remarks on transposition in characteristic $2$}
\label{SubsecRemTr2}

Given an integer $a\geq 1$, we denote by $a_2 := 2^{v_2(a)}$ its $2$-part. From the more general assertion in (\ref{LemLW2}) below, but nevertheless readable independently, we take that
for $b\geq 1$ we have $v_2(\smatze{a}{i}) > 0$ for all $i\in [1,b]$ if and only if $b < a_2$.
We say that a partition $\lambda$ of $n$ is {\it $2$-convergent} if $\lambda'_{p+1} < (\lambda'_p + 1)_2$ for all $p\in [1,\lambda_1 - 1]$. For instance, if $n = 8$, the list of
$2$-convergent partitions is $(1^8)'$, $(3,1^5)'$, $(3^2,1^2)'$, $(5,1^3)'$, $(3^2,2)'$, $(7,1)'$, $(8)'$.
Let $[x] = [1\; 2\; \cdots\; n]$ denote the standard $(n)$-tableau.

\begin{Lemma}[{cf.\ [J 78, 24.4], [K 99, 4.3.35]}]
\label{LemSST11_2}
Let $\lambda$ be a partition of $n$. We have
\[
\dim\Hom_{\sF_2\Sl_n}(S^\lambda/2,S^{(n)}/2) = 
\left\{\begin{array}{ll}
1 & \mbox{ if $\lambda$ is $2$-convergent} \\
0 & \mbox{ else. } \\
\end{array}\right.
\]
In case $\lambda$ is $2$-convergent, the nonzero morphism is of the form
\[
\begin{array}{rcl}
S^\lambda/2 & \lra & S^{(n)}/2 \\
\spi{a}     & \lra & \spi{x}       \\
\end{array}
\]
By transposition (\ref{PropSST11}), we conclude that 
\[
\dim\Hom_{\sF_2\Sl_n}(S^{(n)}/2, S^\lambda/2) = 
\left\{\begin{array}{ll}
1 & \mbox{ if $\lambda'$ is $2$-convergent} \\
0 & \mbox{ else. }  \\
\end{array}\right.
\]

\rm
By $\Sl_n$-linearity, a morphism $S^\lambda/2\lra S^{(n)}/2$, if existent, is necessarily of the form just given. It remains to be seen that this map is well defined
if and only if $\lambda$ is $2$-convergent.

Welldefinedness may be rephrased as the existence of a factorization of the map
\[
\begin{array}{rcl}
F^\lambda   & \lra & S^{(n)}/2 \\
{[a]}       & \lra & \spi{x}       \\
\end{array}
\]
over $F^\lambda\lraa{\nu^\lambda} S^\lambda$. A one-step Garnir relation 
\[
G''_{[a],\xi,\eta} = \sumd{\sigma\in \Sl_\xi\ti\Sl_\eta\<\Sl_{\xi\cup\eta}} [a] \sigma\eps_\sigma,
\]
where $\xi\tm a_p$ and $\eta\tm a_{p+1}$ are subsets such that $\#\xi + \#\eta = \lambda'_p + 1$ (cf.\ Section \ref{SecSpecht}), is mapped to the element 
$\smatze{\lambda'_p + 1}{\#\eta}\cdot\spi{x}\in S^{(n)}/2$. Thus the map is well defined if and only if $\smatze{\lambda'_p + 1}{i} \con_2 0$ for all $p\in [1, \lambda_1 - 1]$ and
for all $i\in [1,\lambda'_{p+1}]$.
\end{Lemma}

\begin{Lemma}
\label{Lem11_3}
Let $\lambda$ be a partition of $n$. Let $\chi^\lambda$ be the characteristic function of the subset $C_{[a_\lambda]}R_{[a_\lambda]}\tm\Sl_n$, i.e.\ $\chi^\lambda_\sigma = 1$ for 
$\sigma\in C_{[a_\lambda]}R_{[a_\lambda]}$, and $\chi^\lambda_\sigma = 0$ for $\sigma\in\Sl_n\ohne C_{[a_\lambda]}R_{[a_\lambda]}$. Let $P^\lambda\tm\Sl_n$ be the subset of permutations 
$\sigma$ for which $[a_\lambda]\sigma$ is a standard tableau. We suppose given a nonzero $\F_2\Sl_n$-linear map (unique if existent, cf.\ \ref{LemSST11_2})
\[
\begin{array}{rcl}
S^{(n)}/2 & \lraa{u^\lambda} & S^\lambda/2 \\
\spi{x}   & \lra             & \sumd{\tau\in P^\lambda} u^\lambda_\tau\spi{a_\lambda}\tau. \\
\end{array}
\]
For each $\sigma\in P^\lambda$, we obtain
\[
1 \con_2 \sum_{\tau\in P^\lambda} u^\lambda_\tau \chi^\lambda_{\tau\sigma^{-1}}.
\]
Moreover, $(\chi^\lambda_{\tau\sigma^{-1}})_{\tau\ti\sigma\in P^\lambda\ti P^\lambda} \in\GL_{\rks\, S^\lambda}(\F_2)$. 

\rm
Note that the transpose of $u^\lambda$ maps each polytabloid to the nonzero element (\ref{LemSST11_2}).

Disregarding alternation, we obtain on the one hand, $\sigma\in\Sl_n$,
\[
\spi{x} \mra{\spi{a_\lambda'}\sigma\delr^\lambda u^{\lambda,\ast}} \sumd{\tau\in P^\lambda} u^\lambda_\tau (\{ a_\lambda\}\sigma,\spi{a_\lambda}\tau ) 
\]
and on the other hand
\[
\spi{x} \mra{\spi{a_\lambda'}\sigma u^{\lambda,\trp}\delr^{(n)}} (\{ x\},\spi{x}) = 1.
\]
It remains to be shown that $(\{ a_\lambda\}\sigma,\spi{a_\lambda}\tau ) \con_2 \chi^\lambda_{\tau\sigma^{-1}}$ for $\tau,\sigma\in P^\lambda$. But 
$(\{ a_\lambda\}\sigma,\spi{a_\lambda}\tau ) \con_2 1$ is equivalent to the
existence of $\kappa\in C_{[a_\lambda]}$ and $\rho\in R_{[a_\lambda]}$ such that $[a_\lambda]\rho = [a_\lambda]\kappa\tau\sigma^{-1}$. The invertibility of the matrix
$(\chi^\lambda_{\tau\sigma^{-1}})_{\tau\ti\sigma\in P^\lambda\ti P^\lambda}$ follows from (\ref{PropSST3}). Cf.\ [K 99, 6.2.8].
\end{Lemma}

For instance, in case $\lambda = (3,2)$, $[a_\lambda] = [\ck a_\lambda]$, and thus $P^{(3,2)} = \{ 1, (45), (23), (23)(45), (2453) \}$, this invertible matrix over $\F_2$ turns out to be
\[
(\chi^{(3,2)}_{\tau\sigma^{-1}})_{\tau\ti\sigma\in P^{(3,2)}\ti P^{(3,2)}} =
                               \enger\left[
                               \begin{array}{lllll}
                               \scm 1 & \scm 0 & \scm 0 & \scm 0 & \scm 0 \\ 
                               \scm 0 & \scm 1 & \scm 0 & \scm 0 & \scm 0 \\  
                               \scm 0 & \scm 0 & \scm 1 & \scm 0 & \scm 0 \\ 
                               \scm 0 & \scm 0 & \scm 0 & \scm 1 & \scm 0 \\  
                               \scm 1 & \scm 0 & \scm 0 & \scm 0 & \scm 1 \\  
                               \end{array}
                               \right]\weiter.
\]
In case $\lambda$ is a hook partition, i.e.\ in case $\lambda_i \leq 1$ for $i\geq 2$, the matrix $(\chi^\lambda_{\tau\sigma^{-1}})_{\tau\ti\sigma\in P^\lambda\ti P^\lambda}$ is the identity 
matrix. In particular, we recover [K 99, 4.2.11] by composition of a morphism as in (\ref{LemSST11_2}) and the transpose, as given by (\ref{Lem11_3}), of a morphism as in (\ref{LemSST11_2}).
Cf.\ (\ref{ExHE3}).     
\subsection{Paths as correspondoids}
\label{SubsecComp}

\begin{footnotesize}
In order to transpose the vertical two-box-shift morphism in (\ref{ThC16}) via (\ref{PropSST11}) to obtain a horizontal two-box-shift morphism, we translate the language of double paths to the 
language of correspondoids (\ref{RemSST15}, \ref{PropSST16}) in order to be able to use our transposition formula (\ref{PropSST11}) to obtain the provisional result 
(\ref{CorSST17}, cf.\ \ref{QuSST19}). We shall use a slightly modified and generalized variant of the setup given in Section \ref{SubSecMS}.
\end{footnotesize}

In this section maps are written on the right, with some exceptions made. Suppose given a partition $\lambda$ of $n$, an integer $d\geq 1$ and integers $1\leq g\leq k \leq \lambda_1 - 1$ such that
\[
\mu'_i := 
\left\{
\begin{array}{ll}
\lambda'_i + d & \mbox{for } i = g \\
\lambda'_i - d & \mbox{for } i = k+1 \\
\lambda'_i     & \mbox{else} \\
\end{array}
\right.
\]
defines a partition $\mu$. A {\it weight} $e$ is a map
\[
\begin{array}{rcl}
[1,\lambda_1] & \lraa{e} & [0,d] \\
j             & \lra     & e_j \\
\end{array}
\]
that maps $g$ and $k+1$ to $e_g = e_{k+1} = d$, and that maps $j$ to $e_j = 0$ in case $j\in [1,\lambda_1]\ohne [g,k+1]$. A {\it pattern} $\Xi$ of weight $e$ is a subset 
$\Xi\tm [1,d]\ti [g,k+1]$ that has  
\[
\#\Big(\Xi\cap ([1,d]\ti\{ j\})\Big) = e_j
\]
for $j\in [g,k+1]$. A {\it $d$-fold path} $\gamma$ of weight $e$ is an injection from a pattern $\Xi$ of weight $e$ to 
$[\lambda]\cup [\mu]$ of the form
\[
\begin{array}{rcl}
\Xi    & \lraa{\gamma} & [\lambda]\cup [\mu] \\
i\ti j & \lra          & \b\gamma(j,i)\ti j    \\
\end{array}
\]  
such that
\[
\b\gamma(g,i) = \lambda'_g + i   \mbox{\ \ for } i\in [1,d]. \\
\]
Sometimes, we denote its pattern by $\Xi_\gamma := \Xi$. The set of $d$-fold paths of weight $e$ is denoted by $\Gamma(e)$.

An {\it ordered $d$-fold path} of weight $e$ is a $d$-fold path $\gamma$ for which the application 
\[
\{ i\in [1,d]\; |\; i\ti j \in\Xi_\gamma \} \mra{\b\gamma(j,-)} [1,\lambda'_j] 
\]
is increasing for each $j\in [g+1,k+1]$. The set of ordered $d$-fold paths of weight $e$ is denoted by 
$\v\Gamma(e)$.

Suppose given a weight $e$. For a $d$-fold path $\gamma$ of weight $e$, we let the {\it successor permutation} be defined by  
\[
\begin{array}{rcl}
\Xi_\gamma & \lraisoa{\mbox{\scr\rm suc}_\gamma} & \Xi_\gamma \\
i\ti j     & \lra                                & \left\{\begin{array}{ll}
                                                   i\ti\min\{ j'\in [j+1,k+1] \; |\; i\ti j'\in\Xi_\gamma\} & \mbox{ for } j\in [g,k] \\
                                                   i\ti g                                                   & \mbox{ for } j = k+1. \\
                                                   \end{array}\right. \\
\end{array}
\]
Since $\gamma$ is an injection, we may define the permutation 
\[
[\lambda']\cup [\mu'] \lraisoa{\h\gamma} [\lambda']\cup [\mu']
\]
as the completion to a commutative diagram
\begin{center}
\begin{picture}(1200,300)
\put( 100, 200){$\Xi_\gamma$}
\put( 150, 210){\vector(1,0){130}}
\put( 200, 220){$\scm\gamma$}
\put( 300, 200){$[\lambda]\cup [\mu]$}
\put( 450, 210){\vector(1,0){130}}
\put( 500, 220){$\scm\tau$}
\put( 500, 190){$\scm\sim$}
\put( 380, 250){$\scm i\ti j$}
\put( 450, 260){\vector(1,0){130}}
\put( 600, 250){$\scm j\ti i$}
\put( 600, 200){$[\lambda']\cup [\mu']$}
\put( 930, 210){\vector(-1,0){150}}
\put( 950, 200){$[\lambda']\cup [\mu']\ohne (\Xi_\gamma)\gamma\tau$}

\put( 120, 180){\vector(0,-1){130}}
\put(  50, 110){$\scm\mbox{\scr\rm suc}_\gamma$}
\put( 130, 110){$\scm\wr$}
\put( 670, 180){\vector(0,-1){130}}
\put( 680, 110){$\scm\h\gamma$}
\put( 655, 110){$\scm\wr$}
\put(1050, 180){\line(0,-1){130}}
\put(1056, 180){\line(0,-1){130}}

\put( 100,   0){$\Xi_\gamma$}
\put( 150,  10){\vector(1,0){130}}
\put( 200,  20){$\scm\gamma$}
\put( 300,   0){$[\lambda]\cup [\mu]$}
\put( 450,  10){\vector(1,0){130}}
\put( 500, -10){$\scm\sim$}
\put( 500,  20){$\scm\tau$}
\put( 600,   0){$[\lambda']\cup [\mu']$}
\put( 930,  10){\vector(-1,0){150}}
\put( 950,  00){$[\lambda']\cup [\mu']\ohne (\Xi_\gamma)\gamma\tau$}
\end{picture}
\end{center}
in which the unlabeled arrows denote inclusions. Thus, roughly speaking, $\h\gamma$ is the identical prolongation of $\mb{\rm suc}_\gamma$ operating via $\gamma$.
Given a $d$-fold path $\gamma$ of weight $e$, the permutation $\sigma(\gamma)$ of $[\lambda']\cup [\mu']$ is defined as being determined by the rule
\[
\begin{array}{rcl}
[\lambda']\cup [\mu'] & \lraisoa{\sigma(\gamma)} & [\lambda']\cup [\mu'] \\
j\ti i                & \lra                     & 
\left\{
\begin{array}{rcll}
\!\! j     &\!\!\!\!\ti\!\!\!\! & i                            & \mbox{for } j\neq k+1 \\
\!\! (k+1) &\!\!\!\!\ti\!\!\!\! & \b\gamma(k+1,i - \mu'_{k+1}) & \mbox{for } j = k+1 \mbox{ and } i\in [\mu'_{k+1} + 1,\lambda'_{k+1}] \\
\end{array}
\right.
\\
\end{array}
\]
and by the requirement that its restriction to $\{ k+1\}\ti [1,\mu'_{k+1}]$ be of constant value $k+1$ in the first component and strictly increasing in the second component. 
Hence, the $[\mu'_{k+1}+1,\lambda'_{k+1}]$-part of the row $k+1$ of $[\lambda']\cup [\mu']$ is mapped under $\sigma(\gamma)$ to the image of $\gamma\tau$, whereas its $[1,\mu'_{k+1}]$-part is 
distributed, in a strictly increasing manner, over the complement of this image in that row. We define an element $\ck\gamma\in\Phi^{\mu',\lambda'}$ as the completion to a commutative diagram
\begin{center}
\begin{picture}(650,250)
\put(   0, 200){$[\mu']$}
\put(  70, 210){\vector(1,0){110}}
\put( 200, 200){$[\lambda']\cup [\mu']$}
\put( 370, 210){\vector(1,0){110}}
\put( 390, 220){$\scm\sigma(\gamma)$}
\put( 400, 190){$\scm\sim$}
\put( 500, 200){$[\lambda']\cup [\mu']$}

\put(  20, 180){\vector(0,-1){130}}
\put(  -5, 110){$\scm\ck\gamma$}
\put(  30, 110){$\scm\wr$}
\put( 520, 180){\vector(0,-1){130}}
\put( 500, 110){$\scm\wr$}
\put( 530, 110){$\scm\h\gamma$}

\put(   0,   0){$[\lambda']$}
\put(  70,  10){\vector(1,0){410}}
\put( 500,   0){$[\lambda']\cup [\mu']$,}
\end{picture}
\end{center}
the unlabeled arrows denoting inclusions. 

Given $\kappa\in C_\lambda$, we denote its transposition by $\kappa'\in R_{\lambda'}$, mapping $j\ti i\in [\lambda']$ to $(j\ti i)\kappa' := (i\ti j)\kappa\tau$. The permutation
of $[\lambda]\cup [\mu]$ that restricts to $\kappa$ on $[\lambda]$ and to the identity on $[\mu]\ohne [\lambda]$ is denoted by $\h\kappa$. Similarly, the identical prolongation of 
$\kappa'$ to a permutation of $[\lambda']\cup [\mu']$ is denoted by $\h\kappa'$. Likewise for $\mu$ instead of $\lambda$.

\begin{Lemma}
\label{LemSST12}
Let $\gamma$ be a $d$-fold path of weight $e$, let $\kappa\in C_\lambda$. We note that the composition $\gamma\h\kappa$ is again a $d$-fold path of weight $e$, of the same pattern as $\gamma$, 
and obtain
\[
\begin{array}{rcl}
(\gamma\h\kappa)\h{\,}  & = & (\h\kappa')^{-1}\h\gamma\h\kappa' \\
(\gamma\h\kappa)\ck{\,} & = & (\kappa_0')^{-1}\ck\gamma\kappa'   \\
\end{array}
\]
for some $\kappa_0\in C_\mu$.

\rm
The first equality holds since right conjugation of a cycle product in the symmetric group is performed by an application of the conjugating element to the cycle entries.
To see the second equality, we consider the diagram
\begin{center}
\begin{picture}(650,650)
\put(   0, 600){$[\mu']$}
\put(  70, 610){\vector(1,0){110}}
\put( 200, 600){$[\lambda']\cup [\mu']$}
\put( 370, 610){\vector(1,0){110}}
\put( 380, 620){$\scm\sigma(\gamma\h\kappa)$}
\put( 400, 590){$\scm\sim$}
\put( 500, 600){$[\lambda']\cup [\mu']$}

\put(  20, 580){\vector(0,-1){130}}
\put( -50, 510){$\scm\kappa_0'^{-1}$}
\put(  30, 510){$\scm\wr$}
\put( 520, 580){\vector(0,-1){130}}
\put( 505, 510){$\scm\wr$}
\put( 535, 510){$\scm\h\kappa'^{-1}$}

\put(   0, 400){$[\mu']$}
\put(  70, 410){\vector(1,0){110}}
\put( 200, 400){$[\lambda']\cup [\mu']$}
\put( 370, 410){\vector(1,0){110}}
\put( 390, 420){$\scm\sigma(\gamma)$}
\put( 400, 390){$\scm\sim$}
\put( 500, 400){$[\lambda']\cup [\mu']$}

\put(  20, 380){\vector(0,-1){130}}
\put( -20, 310){$\scm\ck\gamma$}
\put(  30, 310){$\scm\wr$}
\put( 520, 380){\vector(0,-1){130}}
\put( 505, 310){$\scm\wr$}
\put( 535, 310){$\scm\h\gamma$}

\put(   0, 200){$[\lambda']$}
\put(  70, 210){\vector(1,0){410}}
\put( 500, 200){$[\lambda']\cup [\mu']$}

\put(  20, 180){\vector(0,-1){130}}
\put( -20, 110){$\scm\kappa'$}
\put(  30, 110){$\scm\wr$}
\put( 520, 180){\vector(0,-1){130}}
\put( 505, 110){$\scm\wr$}
\put( 535, 110){$\scm\h\kappa'$}

\put(   0,   0){$[\lambda']$}
\put(  70,  10){\vector(1,0){410}}
\put( 500,   0){$[\lambda']\cup [\mu']$,}
\end{picture}
\end{center}
in which $\kappa'^{-1}_0$ is the completion by restriction and in which the unlabeled arrows denote inclusions. Note that $\sigma(\gamma\h\kappa)\h\kappa'^{-1}\sigma(\gamma)^{-1}$ restricts
identically to the row k+1 of $[\lambda']\cup [\mu']$.
\end{Lemma}

For each weight $e$, we fix a $d$-fold path $\gamma_e$ of weight $e$ and write
\[
\eps(e) := \eps_{\sigma(\gamma_e)}\eps_{[a_\mu']^{-1}\ck\gamma_e [a_\lambda']}.
\]

\begin{Lemma}
\label{LemSST12_3}
The sign $\eps(e)$ is independent of the choice of $\gamma_e$. More precisely, given a $d$-fold path $\gamma$ of weight $e$, we obtain
\[
\eps_{\ck\gamma_e^{-1}\ck\gamma} = \eps_{\sigma(\gamma_e)}\eps_{\sigma(\gamma)}.
\]

\rm
Note that $\eps_{\ck\gamma\ck\gamma_e^{-1}} = \eps_{\sigma(\gamma)\h\gamma\h\gamma_e^{-1}\sigma(\gamma_e)^{-1}}$ since $\sigma(\gamma)\h\gamma\h\gamma_e^{-1}\sigma(\gamma_e)^{-1}$
restricts to the identity on $[\lambda']\ohne [\mu']$. Moreover, we claim that $\eps_{\h\gamma}$ depends only on the weight of $\gamma$. First of all, it depends only on the pattern of $\gamma$, 
since this pattern determines the cycle type of $\h\gamma$. The sign of this cycle type in turn depends only on the cardinality $\sum_{j\in [g,k+1]} e_j$ of that pattern, since the number of 
cycles of length $\geq 2$ that occur in $\h\gamma$ equals $d$. Altogether, after a reordering we obtain $\eps_{\ck\gamma_e^{-1}\ck\gamma} = \eps_{\sigma(\gamma)}\eps_{\sigma(\gamma_e)}$.
\end{Lemma}

\begin{Lemma}
\label{LemSST12_5}
Given a $d$-fold path $\gamma$ of weight $e$ and a $\lambda$-tableau $[a]$, we obtain
\[
(\{\ck\gamma [a']\}\ts 1)\nu_S^\mu = \spi{(\ck\gamma[a'])'}\eps_{\sigma(\gamma)}\eps(e)\eps_{[a]^{-1}[a_\lambda]} \in S^\mu.
\]

\rm

We may conclude
\[
\begin{array}{rcl}
(\{\ck\gamma [a']\}\ts 1)\nu_S^\mu 
& = & \spi{(\ck\gamma [a'])'} \eps_{[a_\mu']^{-1}\ck\gamma [a']}  \\
& = & \spi{(\ck\gamma [a'])'} \eps_{[a_\mu']^{-1}\ck\gamma_e [a_\lambda']}\eps_{[a_\lambda']^{-1}[a']}\eps_{[a']^{-1} \ck\gamma_e^{-1} \ck\gamma [a']} \\
& \auf{\mb{\scr (\ref{LemSST12_3})}}{=} & \spi{(\ck\gamma [a'])'} \eps(e)\eps_{[a]^{-1}[a_\lambda]}\eps_{\sigma(\gamma)}. \\
\end{array}
\]
\end{Lemma}

\begin{Lemma}
\label{LemSST13}
Let $\gamma$ be a $d$-fold path of weight $e$. For $\xi\in\Sl_d$ and $j\in [g+1,k+1]$ we let $\xi$ operate on the columns $\leq j$ via 
\[
\begin{array}{rcl}
\Xi_\gamma & \lraa{\xi_j} & [1,d]\ti [g,k+1] \\
i\ti j'    & \lra         & \left\{\begin{array}{rcll}
                            i      & \ti & j' & \mbox{ for } j' \in \{ g\}\cup [j+1,k+1] \\
                            (i)\xi & \ti & j' & \mbox{ for } j' \in [g+1,j]. \\ 
                            \end{array}\right. \\
\end{array}
\]
We suppose in addition that $\xi_j$ restricts identically to $[1,d]\ti\{ j\}\ohne\Xi_{\gamma}$. Then the composition 
$\Big((\Xi_\gamma)\xi_j\lraisoa{\xi_j^{-1}}\Xi_\gamma\lraa{\gamma}[\lambda]\cup [\mu]\Big)$ is a $d$-fold path of weight $e$. Given a $\lambda$-tableau $[a]$, we obtain
\[
(\{(\xi_j^{-1}\gamma)\ck{\,}\; [a']\}\ts 1)\nu_S^\mu = \left\{\begin{array}{ll}
                                                           (\{\ck\gamma [a']\}\ts 1)\nu_S^\mu         & \mbox{ for } j\in [g+1,k] \\
                                                           (\{\ck\gamma [a']\}\ts 1)\nu_S^\mu\eps_\xi & \mbox{ for } j = k+1. \\
                                                           \end{array}\right.
\]

\rm
We may suppose $\xi$ to be a transposition, $\xi = (s,t)$, $s\ti j,\; t\ti j\in\Xi_j$, $s\neq t$. In case $j\in [g+1,k]$ we obtain, composing with transpositions permuting $[\mu']$, 
\[
((s,t)_j\cdot \gamma)\ck{\,} = \ck\gamma\cdot ((s\ti g)\gamma,(t\ti g)\gamma)\cdot ((s\ti j)\gamma,(t\ti j)\gamma),
\]
yielding the required equality by two column permutations applied to the resulting polytabloid.

In case $j = k+1$, we obtain
\[
((s,t)_{k+1}\cdot \gamma)\ck{\,} = \ck\gamma\cdot((s\ti g)\gamma,(t\ti g)\gamma),
\]
yielding the required equality by a single column permutation applied to the resulting polytabloid.
\end{Lemma}

Given a weight $e$, we write
\[
\begin{array}{rcl}
e!              & := & \prodd{j\in [1,\lambda_1]\ohne\{ g\}} e_j! \\
\lambda'!       & := & \prodd{j\in [1,\lambda_1]\ohne\{ g\}} \lambda'_j ! \\
(\lambda' - e)! & := & \prodd{j\in [1,\lambda_1]\ohne\{ g\}} (\lambda'_j - e_j)! \\
\end{array}
\]

\begin{Corollary}
\label{CorSST14}
Let $e$ be a weight, let $[a]$ be a $\lambda$-tableau. We obtain
\[
\sumd{\gamma\in\Gamma(e)} (\{\ck\gamma [a']\}\ts 1)\nu_S^\mu \eps_{\sigma(\gamma)} = e!\cdot\sumd{\gamma\in\v\Gamma(e)} (\{\ck\gamma [a']\}\ts 1)\nu_S^\mu \eps_{\sigma(\gamma)}.
\]

\rm
We rewrite this sum columnwise from right to left in that we use (\ref{LemSST13}) to impose an ordering condition on $\gamma$ in column $j$, starting in column $j = k+1$ and ending in column 
$j = g+1$. So we obtain a factor $e_j!$ in column $j$.
\end{Corollary}

Given a weight $e$, we let
\[
\begin{array}{rcl}
F^\lambda & \lraa{\ffk''_e} & S^\mu \\
{[a]}     & \lra          & [a]\ffk''_e := \sumd{\gamma\in\v\Gamma(e)} \spi{(\ck\gamma [a'])'}\eps_{\sigma(\gamma)}. \\
\end{array}
\]

\begin{Remark}[connection to Section \ref{SubSecMS}]
\label{RemSST15}
In case of $d = 2$, we obtain that the set of ordered double paths, defined in Section \ref{SubSecMS}, and the set of ordered $2$-fold paths, defined here, coincide. Moreover, $[a]$ being a 
$\lambda$-tableau, $e$ being a weight, we obtain
\[
[a] f''_e = - \fracd{e!}{2} [a]\ffk''_e,
\] 
the left hand side written in the notation of Section \ref{SubSecMS}, the right hand side written in the notation introduced here.

\rm
The set $\Gamma(e)$ is the disjoint union of $\dGamma(e)$ and $(1,2)_{k+1}\cdot\dGamma(e)$. We may rewrite the right hand side as
\[
\begin{array}{rcl}
- \Big(\prodd{j\in [g+1,k]} e_j!\Big)\sumd{\gamma\in\v\Gamma(e)} \spi{(\ck\gamma [a'])'}\eps_{\sigma(\gamma)}
& \auf{\mbox{\scr (\ref{CorSST14})}}{=} & - (e_{k+1}!)^{-1}\sumd{\gamma\in\Gamma(e)} \spi{(\ck\gamma [a'])'}\eps_{\sigma(\gamma)} \\
& \auf{\mbox{\scr (\ref{LemSST13})}}{=} & - \sumd{\gamma\in\dGamma(e)} \spi{(\ck\gamma [a'])'}\eps_{\sigma(\gamma)} \\
& = & \sumd{\gamma\in\dGamma(e)} \spi{(\ck\gamma [a'])'}\eps_{\gamma} \\
& = & [a] f''_e, \\
\end{array}
\]
where $\eps_\gamma = (-1)^{\b\gamma(k+1,1) + \b\gamma(k+1,2)}$ as introduced in Section \ref{SubSecMS}, and where for the last equality we translate $[a^\gamma] = (\ck\gamma [a'])'$. 
\end{Remark}

\begin{Proposition}
\label{PropSST16}
Let $e$ be a weight. We may reformulate 
\[
[a]\ffk''_e = \fracd{1}{e!(\lambda' - e)!\lambda_g'}\sumd{\gamma\in\Gamma(e)/C_\lambda} r_{\ck\gamma} \{ a'\}^- \Theta^-_{\ck\gamma}\nu_S^\mu\cdot\eps(e),
\]
where $\gamma\in\Gamma(e)/C_\lambda$ signifies that $\gamma$ runs over a set of orbit representatives of $\Gamma(e)$ under the operation of $C_\lambda$ in the sense of (\ref{LemSST12}). 
In particular, there is a factorization
\[
(F^\lambda \lraa{\ffk''_e} S^\mu) = (F^\lambda \lraa{\nu_M^\lambda} M^{\lambda,-}\lraa{\ffk'_e} S^\mu). 
\]

\rm
Given a $d$-fold path $\gamma$ of weight $e$, the stabilizer $\{\kappa\in C_\lambda \; |\; \gamma\h\kappa = \gamma\}$ of $\gamma$ in $C_\lambda$ has cardinality $(\lambda' - e)!\cdot\lambda'_g$. 
So we calculate
\[
\begin{array}{rcl}
[a]\ffk''_e
& = & \sumd{\gamma\in\v\Gamma(e)} \spi{(\ck\gamma [a'])'}\eps_{\sigma(\gamma)} \\
& \auf{\mbox{\scr (\ref{CorSST14})}}{=} & \frac{1}{e!}\sumd{\gamma\in\Gamma(e)} \spi{(\ck\gamma [a'])'}\eps_{\sigma(\gamma)} \\
& = & \frac{1}{e!(\lambda' - e)!\lambda'_g}\sumd{\gamma\in\Gamma(e)/C_\lambda} r_{\ck\gamma} r_{\ck\gamma}^{-1} 
                                 \sumd{\kappa\in C_\lambda}\spi{(\gamma\h\kappa)\ck{\,}\; [a'])'}\eps_{\sigma(\gamma)} \\
& \auf{\mbox{\scr (\ref{LemSST12_5})}}{=} & \frac{1}{e!(\lambda' - e)!\lambda'_g}\sumd{\gamma\in\Gamma(e)/C_\lambda} r_{\ck\gamma} r_{\ck\gamma}^{-1} 
                                 \sumd{\kappa\in C_\lambda}(\{(\gamma\h\kappa)\ck{\,}\; [a']\}\ts 1)\nu_S^\mu\cdot\eps(e)\eps_{[a]^{-1}[a_\lambda]} \\
\end{array}
\]
\[
\begin{array}{rcl}
& \auf{\mbox{\scr (\ref{LemSST12})}}{=} & \frac{1}{e!(\lambda' - e)!\lambda'_g}\sumd{\gamma\in\Gamma(e)/C_\lambda} r_{\ck\gamma} \Big(\Big(r_{\ck\gamma}^{-1} 
                                 \sumd{\kappa'\in R_{\lambda'}} \{\ck\gamma\kappa' [a']\}\Big)\ts 1\Big)\nu_S^\mu\cdot \eps(e)\eps_{[a]^{-1}[a_\lambda]} \\
& \auf{\mbox{\scr (\ref{LemSST1})}}{=}  & \frac{1}{e!(\lambda' - e)!\lambda'_g}\sumd{\gamma\in\Gamma(e)/C_\lambda} r_{\ck\gamma} \Big(\{a'\}\Theta_{\ck\gamma}\ts 1\Big)\nu_S^\mu\cdot 
                                 \eps(e)\eps_{[a]^{-1}[a_\lambda]}\\
& =                                     & \frac{1}{e!(\lambda' - e)!\lambda'_g}\sumd{\gamma\in\Gamma(e)/C_\lambda} r_{\ck\gamma} \{ a'\}^-\Theta^-_{\ck\gamma}\nu_S^\mu\cdot\eps(e). \\ 
\end{array}
\]
\end{Proposition}

\begin{Corollary}[to \ref{ThC16} via \ref{PropSST16}]
\label{CorSST17}
Let $d = 2$. The transpose in the sense of (\ref{PropSST11}) of the morphism given in (\ref{ThC16}) is obtained as a factorization
\begin{center}
\begin{picture}(500,250)
\put( 100, 200){$F^{\mu'}$}
\put( 170, 210){\vector(1,0){110}}
\put( 200, 220){$\scm f^{\trp,0}$}
\put( 300, 200){$S^{\lambda'}$}

\put( -50, 200){$[b']$}
\put( -30, 180){\vector(0,-1){130}}
\put(-170,   0){$\spi{b'} + mS^{\mu'}$}
\put( 120, 180){\vector(0,-1){130}}
\put( 320, 180){\vector(0,-1){130}}
\put( 450, 200){$\spi{a'}$}
\put( 470, 180){\vector(0,-1){130}}
\put( 450,   0){$\spi{a'} + mS^{\lambda'}$}

\put(  60,   0){$S^{\mu'}/m$}
\put( 200,  10){\vector(1,0){80}}
\put( 230,  20){$\scm f^\trp$}
\put( 300,   0){$S^{\lambda'}/m$,}
\end{picture}
\end{center}
where
\[
\text
[b']f^{\trp,0} := -\sumd{e\in E}\frac{1}{(\lambda'-e)!\cdot 2R\lambda'_g}\Big(\prodd{i\in [g+1,k]} X_i^{(2-e_i)}\Big)
                 \sumd{\gamma\in\Gamma(e)/C_\lambda} r_{\ck\gamma}\spi{b'}\Theta_{\ck\gamma^{-1}}\cdot\eps(e),
\]
and where $f$, $E$, $R$, $m$ and $X_i$ are defined as for (\ref{ThC16}).

\rm
Given $t, T\geq 2$, $\Z$-lattices $X$ and $Y$ and a $\Z$-linear map $X/t\lraa{u} Y/t$, we denote
\[
(X/(tT) \lraa{u\,\bullet\, T} Y/(tT)) := (X/(tT) \lra X/t \lraa{u} Y/t \lraa{T} Y/(tT)). 
\]
Let $M := R\cdot\lambda'!\cdot 2\lambda'_g$ play the role of a large enough integer and consider the morphism
\[
\begin{array}{rcl}
S^\lambda/(mM) & \lraa{f\,\bullet\, M} & S^\mu/(mM) \\
\spi{a}        & \lra                  & \lambda'!\cdot 2\lambda'_g\cdot\sumd{e\in E} \Big(\prodd{i\in [g+1,k]} X_i^{(2-e_i)}\Big)\{a'\}^- f'_e \\
               & \auf{\mb{\scr(\ref{RemSST15}, \ref{PropSST16})}}{=} & -\lambda'!\lambda'_g\cdot\sumd{e\in E} \Big(\prodd{i\in [g+1,k]} X_i^{(2-e_i)}\Big)\{a'\}^- \ffk'_e \\
               & \auf{\mb{\scr(\ref{PropSST16})}}{=} & -\sumd{e\in E} \frac{\lambda'!}{(\lambda'-e)!} \Big(\prodd{i\in [g+1,k]} X_i^{(2-e_i)}\Big)
                 \sumd{\gamma\in\Gamma(e)/C_\lambda} r_{\ck\gamma}\{ a'\}^-\Theta_{\ck\gamma}^-\nu_S^\mu\cdot \eps(e). \\
\end{array}
\]
For $\phi\in\Phi^{\mu',\lambda'}$, the coefficient $x_{\{\phi\}}$ in the sense of (\ref{PropSST11}) can be written as
\[
\textstyle
x_{\{\phi\}} = -\sumd{e\in E} \frac{\lambda'!}{(\lambda' - e)!}(\prodd{i\in [g+1,k]} X_i^{(2-e_i)})\sumd{\gamma\in\Gamma(e)/C_\lambda} r_{\ck\gamma} \cdot
\left\{\renewcommand{\arraystretch}{0.8}\begin{array}{ll}
\scm 1 & \scm\mbox{\scr\ for } \{\ck\gamma\} \; =\;  \{\phi\} \\
\scm 0 & \scm\mbox{\scr\ else } \\
\end{array}\renewcommand{\arraystretch}{1.5}\right\}
\cdot \eps(e).
\]
Thus, by (\ref{PropSST11}), the transpose of $S^\lambda/(mM) \lraa{f\,\bullet\, M} S^\mu/(mM)$ maps
\[
\begin{array}{rrl}
S^{\mu'}/(mM)  & \lraa{(f\,\bullet\, M)^\trp} & S^{\lambda'}/(mM) \\
\spi{b'}       & \lra                   & -\sumd{e\in E} \frac{\lambda'!}{(\lambda' - e)!}(\prodd{i\in [g+1,k]} X_i^{(2-e_i)})\sumd{\gamma\in\Gamma(e)/C_\lambda} 
                                          r_{\ck\gamma} \spi{b'} \Theta_{\ck\gamma^{-1}}\cdot \eps(e) \\
\end{array}
\]
Finally, there is a commutative diagram 
\begin{center}
\begin{picture}(500,250)
\put(-220, 200){$\Hom_{\sZ\Sl_n}(S^\lambda/m,S^\mu/m)$}
\put( 200, 210){\vector(1,0){80}}
\put( 220, 220){\scr tra}
\put( 230, 190){$\scm\sim$}
\put( 300, 200){$\Hom_{\sZ\Sl_n}(S^{\mu'}/m,S^{\lambda'}/m)$}

\put( -50, 180){\vector(0,-1){130}}
\put( 470, 180){\vector(0,-1){130}}

\put(-370,   0){$\Hom_{\sZ\Sl_n}(S^\lambda/(mM),S^\mu/(mM))$}
\put( 200,  10){\vector(1,0){80}}
\put( 220,  20){\scr tra}
\put( 230, -10){$\scm\sim$}
\put( 300,   0){$\Hom_{\sZ\Sl_n}(S^{\mu'}/(mM),S^{\lambda'}/(mM))$}
\end{picture}
\end{center}
with vertical injections given by multiplication by $M$, so that $f^\trp\bullet M = (f\bullet M)^\trp$.
\end{Corollary}

\begin{Remark}
\label{RemSST18}
\rm
Consider the set $\v\Gamma_0(e)$ of ordered $d$-fold paths $\gamma$ of weight $e$ that have $\b\gamma(j,i) \in [\lambda'_j-e_j,\lambda'_j]$ for $i\ti j\in\Xi_{[g+1,k+1]}$, i.e.\ that have
their positions sitting at the bottoms of the columns. We say that a correspondoid $\{\phi\}\in\{\Phi^{\mu',\lambda'}\}$ is {\it semistandard up to row permutation} if there exists 
a row permutation $\rho\in R_{\mu'}$ such that $\{\rho\phi\}\in \{\Phi^{\mu',\lambda'}\}_{\mbox{\scr\rm sst}}$ becomes semistandard. Note that $\Theta_\phi = \Theta_{\rho\phi}$ by (\ref{LemSST1}).
For $\gamma\in\v\Gamma_0(e)$, the correspondoid $\{\ck\gamma^{-1}\}$ is semistandard up to row permutation as long as $\Xi_\gamma$ is not bulky in the sense introduced at the beginning of Section 
\ref{SubSecRedund}. Cf.\ the bases in (\ref{CorSST7}), cf.\ the reduced tuple of coefficients employed in (\ref{ThC16}). 

Note that $\v\Gamma_0(e)$ is a set of representatives of $\Gamma(e)/C_\lambda$.
\end{Remark}

\begin{Question}
\label{QuSST19}
\rm 
I do not know a formula for the transpose $f^\trp$ given in (\ref{CorSST17}) in terms of $\lambda'$-polytabloids. Cf.\ (\ref{QuSST11_1}).
\end{Question}       
\textheight25.2cm

\subsection{Horizontal examples}
\label{SubsecHorEx}

We shall transpose some examples of Section \ref{SubSecEx} with respect to the convention that $[a_\lambda] = [\ck a_\lambda]$ for all partitions $\lambda$, cf.\ Section \ref{SecSpecht}. 
Again, we omit to denote the brackets that indicate polytabloids, moreover, we omit commas in cycles. If possible (cf.\ \ref{CorSST7}), we sort the resulting linear combination of polytabloids
into images under $\Theta_\phi^-\nu_S^{\lambda'}$ for various (implicitely present) correspondoids $\{\phi\}$. The formula in (\ref{CorSST17}) does not suffice to give the coefficients of
this expression (cf.\ \ref{QuSST11_1}, \ref{ExLW10}).

\begin{Example}
\label{ExHE2}
\rm
Let $\lambda = (3,3)$, $\mu = (2,2,1,1)$. The transpose of the morphism given in (\ref{ExE2}) is 

\begin{footnotesize}
\[
\begin{array}{rcl}
\mbox{\ncr $S^{\mu'}/4$} & \lraa{\mbox{\ncr $f^\trp$}} & \mbox{\ncr $S^{\lambda'}/4$} \\
\begin{array}{llll}
1 & 3 & 5 & 6 \\
2 & 4 &   &   \\
\end{array}
& \lra &
-
\begin{array}{ll}
1 & 4 \\
2 & 5 \\
3 & 6 \\
\end{array}
\cdot (1 - (34))\;\;\;
+ 2\cdot
\begin{array}{ll}
1 & 3 \\
2 & 4 \\
5 & 6 \\
\end{array}. \\
\end{array}
\]
\end{footnotesize}
\end{Example}

\begin{Example}
\label{ExHE3}
\rm
Let $\lambda = (3,2,2)$, $\mu = (3,1,1,1,1)$. The transposes of the morphisms given in (\ref{ExE3}) are 

\begin{footnotesize}
\[
\begin{array}{rcl}
\mbox{\ncr $S^{\mu'}/3$} & \lraa{\mbox{\ncr $f^\trp$}} & \mbox{\ncr $S^{\lambda'}/3$} \\
\begin{array}{lllll}
1 & 4 & 5 & 6 & 7 \\
2 &   &   &   &   \\
3 &   &   &   &   \\
\end{array}
& \lra &
-
\begin{array}{lll}
1 & 4 & 6 \\
2 & 5 & 7 \\
3 &   &   \\
\end{array}
-
\begin{array}{lll}
1 & 4 & 5 \\
2 & 6 & 7 \\
3 &   &   \\
\end{array}. \\
\end{array}
\]
\end{footnotesize}

and

\begin{footnotesize}
\[
\begin{array}{rcl}
\mbox{\ncr $S^{\mu'}/2$} & \lraa{\mbox{\ncr $u^{\mu,\trp} u^{\lambda'}$}} & \mbox{\ncr $S^{\lambda'}/2$} \\
\begin{array}{lllll}
1 & 4 & 5 & 6 & 7 \\
2 &   &   &   &   \\
3 &   &   &   &   \\
\end{array}
& \lra &
\begin{array}{lll}
1 & 2 & 6 \\
3 & 4 & 7 \\
5 &   &   \\
\end{array}
+
\begin{array}{lll}
1 & 4 & 5 \\
2 & 6 & 7 \\
3 &   &   \\
\end{array}
+
\begin{array}{lll}
1 & 3 & 5 \\
2 & 6 & 7 \\
4 &   &   \\
\end{array}
+
\begin{array}{lll}
1 & 2 & 5 \\
3 & 6 & 7 \\
4 &   &   \\
\end{array}
+
\begin{array}{lll}
1 & 3 & 4 \\
2 & 6 & 7 \\
5 &   &   \\
\end{array} \\
& & +
\begin{array}{lll}
1 & 2 & 4 \\
3 & 6 & 7 \\
5 &   &   \\
\end{array}
+
\begin{array}{lll}
1 & 2 & 3 \\
4 & 6 & 7 \\
5 &   &   \\
\end{array}
+
\begin{array}{lll}
1 & 2 & 5 \\
3 & 4 & 7 \\
6 &   &   \\
\end{array}
+
\begin{array}{lll}
1 & 3 & 4 \\
2 & 5 & 7 \\
6 &   &   \\
\end{array}
+
\begin{array}{lll}
1 & 2 & 4 \\
3 & 5 & 7 \\
6 &   &   \\
\end{array} \\
& & +
\begin{array}{lll}
1 & 2 & 3 \\
4 & 5 & 7 \\
6 &   &   \\
\end{array}
+
\begin{array}{lll}
1 & 2 & 5 \\
3 & 4 & 6 \\
7 &   &   \\
\end{array}
+
\begin{array}{lll}
1 & 3 & 4 \\
2 & 5 & 6 \\
7 &   &   \\
\end{array}
+
\begin{array}{lll}
1 & 2 & 4 \\
3 & 5 & 6 \\
7 &   &   \\
\end{array}
+
\begin{array}{lll}
1 & 2 & 3 \\
4 & 5 & 6 \\
7 &   &   \\
\end{array}. \\
\end{array}
\]
\end{footnotesize}

Let us consider the morphism $u^{\mu,\trp} u^{\lambda'}$. Using $[a_{\lambda'}] = [\ck a_{\lambda'}]$, we obtain, in the notation of (\ref{Lem11_3}),
\[
\begin{array}{l}
P^{\lambda'} = 
\{ 1, (56), (34), (34)(56), (354), (3564), (3654), (364), (37654), (3764), (234), (234)(56), \\
(2354), (23564), (23654), (2364), (237654), (23764), (24)(356), (24)(36), (24)(376)\} \tm \Sl_7. \\
\end{array}
\]
Hence, with respect to the ordering of $P^{\lambda'}$ just indicated, the matrix $(\chi^{\lambda'}_{\tau\sigma^{-1}})_{\tau\ti\sigma\in P^{\lambda'}\ti P^{\lambda'}}$ takes the form
\[
\enger\left[
\begin{array}{rrrrrrrrrrrrrrrrrrrrr}
\scm 1 &\scm 0 &\scm 0 &\scm 0 &\scm 0 &\scm 0 &\scm 0 &\scm 0 &\scm 0 &\scm 0 &\scm 0 &\scm 0 &\scm 0 &\scm 0 &\scm 0 &\scm 0 &\scm 0 &\scm 0 &\scm 0 &\scm 0 &\scm 0 \\
\scm 0 &\scm 1 &\scm 0 &\scm 0 &\scm 0 &\scm 0 &\scm 0 &\scm 0 &\scm 0 &\scm 0 &\scm 0 &\scm 0 &\scm 0 &\scm 0 &\scm 0 &\scm 0 &\scm 0 &\scm 0 &\scm 0 &\scm 0 &\scm 0 \\
\scm 0 &\scm 0 &\scm 1 &\scm 0 &\scm 0 &\scm 0 &\scm 0 &\scm 0 &\scm 0 &\scm 0 &\scm 0 &\scm 0 &\scm 0 &\scm 0 &\scm 0 &\scm 0 &\scm 0 &\scm 0 &\scm 0 &\scm 0 &\scm 0 \\
\scm 0 &\scm 0 &\scm 0 &\scm 1 &\scm 0 &\scm 0 &\scm 0 &\scm 0 &\scm 0 &\scm 0 &\scm 0 &\scm 0 &\scm 0 &\scm 0 &\scm 0 &\scm 0 &\scm 0 &\scm 0 &\scm 0 &\scm 0 &\scm 0 \\
\scm 0 &\scm 0 &\scm 0 &\scm 0 &\scm 1 &\scm 0 &\scm 0 &\scm 0 &\scm 0 &\scm 0 &\scm 0 &\scm 0 &\scm 0 &\scm 0 &\scm 0 &\scm 0 &\scm 0 &\scm 0 &\scm 0 &\scm 0 &\scm 0 \\
\scm 0 &\scm 0 &\scm 0 &\scm 0 &\scm 0 &\scm 1 &\scm 0 &\scm 0 &\scm 0 &\scm 0 &\scm 0 &\scm 0 &\scm 0 &\scm 0 &\scm 0 &\scm 0 &\scm 0 &\scm 0 &\scm 0 &\scm 0 &\scm 0 \\
\scm 0 &\scm 0 &\scm 0 &\scm 0 &\scm 0 &\scm 0 &\scm 1 &\scm 0 &\scm 0 &\scm 0 &\scm 0 &\scm 0 &\scm 0 &\scm 0 &\scm 0 &\scm 0 &\scm 0 &\scm 0 &\scm 0 &\scm 0 &\scm 0 \\
\scm 0 &\scm 0 &\scm 0 &\scm 0 &\scm 0 &\scm 0 &\scm 0 &\scm 1 &\scm 0 &\scm 0 &\scm 0 &\scm 0 &\scm 0 &\scm 0 &\scm 0 &\scm 0 &\scm 0 &\scm 0 &\scm 0 &\scm 0 &\scm 0 \\
\scm 0 &\scm 0 &\scm 0 &\scm 0 &\scm 0 &\scm 0 &\scm 0 &\scm 0 &\scm 1 &\scm 0 &\scm 0 &\scm 0 &\scm 0 &\scm 0 &\scm 0 &\scm 0 &\scm 0 &\scm 0 &\scm 0 &\scm 0 &\scm 0 \\
\scm 0 &\scm 0 &\scm 0 &\scm 0 &\scm 0 &\scm 0 &\scm 0 &\scm 0 &\scm 0 &\scm 1 &\scm 0 &\scm 0 &\scm 0 &\scm 0 &\scm 0 &\scm 0 &\scm 0 &\scm 0 &\scm 0 &\scm 0 &\scm 0 \\
\scm 0 &\scm 0 &\scm 0 &\scm 0 &\scm 0 &\scm 0 &\scm 0 &\scm 0 &\scm 0 &\scm 0 &\scm 1 &\scm 0 &\scm 0 &\scm 0 &\scm 0 &\scm 0 &\scm 0 &\scm 0 &\scm 0 &\scm 0 &\scm 0 \\
\scm 0 &\scm 0 &\scm 0 &\scm 0 &\scm 0 &\scm 0 &\scm 0 &\scm 0 &\scm 0 &\scm 0 &\scm 0 &\scm 1 &\scm 0 &\scm 0 &\scm 0 &\scm 0 &\scm 0 &\scm 0 &\scm 0 &\scm 0 &\scm 0 \\
\scm 1 &\scm 0 &\scm 0 &\scm 0 &\scm 0 &\scm 0 &\scm 0 &\scm 0 &\scm 0 &\scm 0 &\scm 0 &\scm 0 &\scm 1 &\scm 0 &\scm 0 &\scm 0 &\scm 0 &\scm 0 &\scm 0 &\scm 0 &\scm 0 \\
\scm 1 &\scm 0 &\scm 0 &\scm 0 &\scm 0 &\scm 0 &\scm 0 &\scm 0 &\scm 0 &\scm 0 &\scm 0 &\scm 0 &\scm 0 &\scm 1 &\scm 0 &\scm 0 &\scm 0 &\scm 0 &\scm 0 &\scm 0 &\scm 0 \\
\scm 0 &\scm 1 &\scm 0 &\scm 0 &\scm 0 &\scm 0 &\scm 0 &\scm 0 &\scm 0 &\scm 0 &\scm 0 &\scm 0 &\scm 0 &\scm 0 &\scm 1 &\scm 0 &\scm 0 &\scm 0 &\scm 0 &\scm 0 &\scm 0 \\
\scm 0 &\scm 1 &\scm 0 &\scm 0 &\scm 0 &\scm 0 &\scm 0 &\scm 0 &\scm 0 &\scm 0 &\scm 0 &\scm 0 &\scm 0 &\scm 0 &\scm 0 &\scm 1 &\scm 0 &\scm 0 &\scm 0 &\scm 0 &\scm 0 \\
\scm 1 &\scm 1 &\scm 0 &\scm 0 &\scm 0 &\scm 0 &\scm 0 &\scm 0 &\scm 0 &\scm 0 &\scm 0 &\scm 0 &\scm 0 &\scm 0 &\scm 0 &\scm 0 &\scm 1 &\scm 0 &\scm 0 &\scm 0 &\scm 0 \\
\scm 0 &\scm 1 &\scm 0 &\scm 0 &\scm 0 &\scm 0 &\scm 0 &\scm 0 &\scm 0 &\scm 0 &\scm 0 &\scm 0 &\scm 0 &\scm 0 &\scm 0 &\scm 0 &\scm 0 &\scm 1 &\scm 0 &\scm 0 &\scm 0 \\
\scm 0 &\scm 0 &\scm 1 &\scm 0 &\scm 1 &\scm 0 &\scm 0 &\scm 0 &\scm 0 &\scm 0 &\scm 0 &\scm 0 &\scm 0 &\scm 0 &\scm 0 &\scm 0 &\scm 0 &\scm 0 &\scm 1 &\scm 0 &\scm 0 \\
\scm 0 &\scm 0 &\scm 0 &\scm 1 &\scm 0 &\scm 0 &\scm 1 &\scm 0 &\scm 0 &\scm 0 &\scm 0 &\scm 0 &\scm 0 &\scm 0 &\scm 0 &\scm 0 &\scm 0 &\scm 0 &\scm 0 &\scm 1 &\scm 0 \\
\scm 0 &\scm 0 &\scm 0 &\scm 1 &\scm 0 &\scm 0 &\scm 0 &\scm 0 &\scm 1 &\scm 0 &\scm 1 &\scm 0 &\scm 0 &\scm 0 &\scm 0 &\scm 0 &\scm 0 &\scm 0 &\scm 0 &\scm 0 &\scm 1 \\
\end{array}
\right]\weiter. 
\]
By labelling our morphism as $u^{\mu,\trp} u^{\lambda'}$, we implicitely claimed the tuple $(u^{\lambda'}_\tau)_{\tau\in P^{\lambda'}}$ to be represented by the vector 
\[
\scm [0,1,0,1,0,1,0,1,0,1,0,1,1,1,1,1,1,1,1,1,1].
\]
But this row vector yields a row vector of constant entry $1$ by right multiplication with that matrix, which confirms this claim.

Furthermore, as remarked after (\ref{Lem11_3}), the matrix $(\chi^\mu_{\tau\sigma^{-1}})_{\tau\ti\sigma\in P^\mu\ti P^\mu}$ equals the identity matrix.
\end{Example}

\begin{Example}
\label{ExHE3_5}
\rm
Let $\lambda = (3,3,1,1)$, $\mu = (2,2,1,1,1,1)$. The transpose of the morphism given in (\ref{ExE3_5}) is 

\begin{footnotesize}
\[
\begin{array}{rcl}
\mbox{\ncr $S^{\mu'}/6$} & \lraa{\mbox{\ncr $f^\trp$}} & \mbox{\ncr $S^{\lambda'}/6$} \\
\begin{array}{llllll}
1 & 3 & 5 & 6 & 7 & 8 \\
2 & 4 &   &   &   &   \\
\end{array}
& \lra &
-
\begin{array}{llll}
1 & 4 & 7 & 8 \\
2 & 5 &   &   \\
3 & 6 &   &   \\
\end{array}
\cdot(1 - (34))\;\;\;
- 2\cdot
\begin{array}{llll}
1 & 3 & 7 & 8 \\
2 & 4 &   &   \\
5 & 6 &   &   \\
\end{array} \\
& & -
\begin{array}{llll}
1 & 4 & 6 & 8 \\
2 & 5 &   &   \\
3 & 7 &   &   \\
\end{array}
\cdot(1 - (34))\;\;\;
- 2\cdot
\begin{array}{llll}
1 & 3 & 6 & 8 \\
2 & 4 &   &   \\
5 & 7 &   &   \\
\end{array} \\
& & -
\begin{array}{llll}
1 & 4 & 5 & 8 \\
2 & 6 &   &   \\
3 & 7 &   &   \\
\end{array}
\cdot(1 - (34))\;\;\;
- 2\cdot
\begin{array}{llll}
1 & 3 & 5 & 8 \\
2 & 4 &   &   \\
6 & 7 &   &   \\
\end{array} \\
& & -
\begin{array}{llll}
1 & 4 & 6 & 7 \\
2 & 5 &   &   \\
3 & 8 &   &   \\
\end{array}
\cdot(1 - (34))\;\;\;
- 2\cdot
\begin{array}{llll}
1 & 3 & 6 & 7 \\
2 & 4 &   &   \\
5 & 8 &   &   \\
\end{array} \\
& & -
\begin{array}{llll}
1 & 4 & 5 & 7 \\
2 & 6 &   &   \\
3 & 8 &   &   \\
\end{array}
\cdot(1 - (34))\;\;\;
- 2\cdot
\begin{array}{llll}
1 & 3 & 5 & 7 \\
2 & 4 &   &   \\
6 & 8 &   &   \\
\end{array} \\
& & -
\begin{array}{llll}
1 & 4 & 5 & 6 \\
2 & 7 &   &   \\
3 & 8 &   &   \\
\end{array}
\cdot(1 - (34))\;\;\;
- 2\cdot
\begin{array}{llll}
1 & 3 & 5 & 6 \\
2 & 4 &   &   \\
7 & 8 &   &   \\
\end{array} . \\
\end{array}
\]
\end{footnotesize}
\end{Example}

\begin{Example}
\label{ExHE4}
\rm
Let $\lambda = (4,4)$, $\mu = (3,3,1,1)$. The transpose of the morphism given in (\ref{ExE3}) is \vspace*{-3mm}

\begin{footnotesize}
\[
\begin{array}{rcl}
\mbox{\ncr $S^{\mu'}/5$} & \lraa{\mbox{\ncr $f^\trp$}} & \mbox{\ncr $S^{\lambda'}/5$} \\
\begin{array}{llll}
1 & 4 & 7 & 8 \\
2 & 5 &   &   \\
3 & 6 &   &   \\
\end{array}
& \lra &
-
\begin{array}{ll}
1 & 5 \\
2 & 6 \\
3 & 7 \\
4 & 8 \\
\end{array}
\cdot(1 - (45) - (46))\;\;\;
- 3\cdot
\begin{array}{ll}
1 & 4 \\
2 & 5 \\
3 & 6 \\
7 & 8 \\
\end{array}. \\
\end{array}
\]
\end{footnotesize}
\end{Example}

\begin{Example}
\label{ExHE5}
\rm
Let $\lambda = (3,3,2)$, $\mu = (2,2,2,1,1)$. The transpose of the morphism given in (\ref{ExE3}) is \vspace*{-3mm}

\begin{footnotesize}
\[
\begin{array}{rcl}
\mbox{\ncr $S^{\mu'}/5$} & \lraa{\mbox{\ncr $f^\trp$}} & \mbox{\ncr $S^{\lambda'}/5$} \\
\begin{array}{lllll}
1 & 3 & 5 & 7 & 8 \\
2 & 4 & 6 &   &  \\
\end{array}
& \lra &
-
\begin{array}{lll}
1 & 4 & 7 \\
2 & 5 & 8 \\
3 & 6 &   \\
\end{array}
\cdot(1 - (34))\;\;\;
- 2 \cdot
\begin{array}{lll}
1 & 4 & 6 \\
2 & 5 & 8 \\
3 & 7 &   \\
\end{array}
\cdot(1 - (34))(1 - (56)) \\
& & - 
\begin{array}{lll}
1 & 3 & 6 \\
2 & 4 & 8 \\
5 & 7 &   \\
\end{array} 
\cdot(1 - (56))\;\;\;
+
\begin{array}{lll}
1 & 4 & 6 \\
2 & 5 & 7 \\
3 & 8 &   \\
\end{array}
\cdot(1 - (34))(1 - (56))\\
& & -
\begin{array}{lll}
1 & 3 & 6 \\
2 & 4 & 7 \\
5 & 8 &   \\
\end{array}
\cdot(1 - (56))\;\;\;
+ 2\cdot
\begin{array}{lll}
1 & 3 & 5 \\
2 & 4 & 6 \\
7 & 8 &   \\
\end{array}. \\
\end{array}
\]
\end{footnotesize}
\end{Example}

\begin{Example}
\label{ExHE6}
\rm
Let $\lambda = (3,3,1,1)$, $\mu = (2,2,2,2)$. The transpose of the morphism given in (\ref{ExE6}) is \vspace*{-3mm}

\begin{footnotesize}
\[
\begin{array}{rcl}
\mbox{\ncr $S^{\mu'}/3$} & \lraa{\mbox{\ncr $f^\trp$}} & \mbox{\ncr $S^{\lambda'}/3$} \\
\begin{array}{llll}
1 & 3 & 5 & 7 \\
2 & 4 & 6 & 8 \\
\end{array}
& \lra &
\begin{array}{llll}
1 & 4 & 6 & 8 \\
2 & 5 &   &   \\
3 & 7 &   &   \\
\end{array}
\cdot(1 - (34))(1 - (56))(1 - (78))\\
& & -
\begin{array}{llll}
1 & 3 & 6 & 8 \\
2 & 4 &   &   \\
5 & 7 &   &   \\
\end{array}
\cdot(1 - (56))(1 - (78)). \\
\end{array}
\]
\end{footnotesize}
\end{Example}

\begin{Remark}
\label{RemHE7}
\rm
The morphisms in (\ref{ExE7}, \ref{ExE8}) are transposed to their respective negative, $g^\trp = -g$. For instance, the morphism 
$S^{(4,3)}/4 \lraa{g} S^{(2,2,2,1)}/4$ given in (\ref{ExE7}) can be obtained as 
\[
\nu_S^\lambda \cdot g \con_4 (2\cdot\Theta^-_{\phi_1} - \Theta^-_{\phi_2} + 2\cdot\Theta^-_{\phi_3})\cdot \nu_S^\mu,
\]

\begin{footnotesize}
\[
\mbox{\ncr $\{\phi_1\} := $}
\left\{\begin{array}{ll}
1 & 3 \\
1 & 3 \\
2 & 4 \\
2 &   \\
\end{array}\right\} ,\;\;
\mbox{\ncr $\{\phi_2\} := $}
\left\{\begin{array}{ll}
1 & 2 \\
1 & 3 \\
2 & 4 \\
3 &   \\
\end{array}\right\} ,\;\;
\mbox{\ncr $\{\phi_3\} := $}
\left\{\begin{array}{ll}
1 & 2 \\
1 & 2 \\
3 & 4 \\
3 &   \\
\end{array}\right\} .\;\;
\]
\end{footnotesize}
which allows to apply (\ref{PropSST11}) directly and to compare linear combinations of tabloids.
\end{Remark}       
\section{Two columns, several boxes}
\label{SecTCSB}

\begin{footnotesize}
We shall generalize (\ref{PropM14}), treating the case of essentially two columns and two shifted boxes, to the case of essentially two columns and an arbitrary number of shifted boxes.
\end{footnotesize}

\begin{Lemma}
\label{LemLW1}
Let $x,y\geq 0$. As polynomials in the free variable $t$, we obtain
\[
\sumd{j\in [0,x]} (-1)^{x-j} \smatze{y+j}{j}\smatze{t}{x-j} = \smatze{y-t+x}{x} \in\Q[t],
\]
where 
\[
\smatze{f(t)}{z} := z!^{-1}\cdot \prodd{i\in[0,z-1]}(f(t)-i)
\]
for $f(t)\in\Q[t]$, $z\geq 0$. So, in particular, $\smatze{f(t)}{0} = 1$. 

\rm
We name the claimed equality $\Eq(x,y)$ and shall prove the following assertions.
\begin{itemize}
\item[(i)] $\Eq(0,y)$ holds for $y\geq 0$.
\item[(ii)] $\Eq(x,0)$ holds for $x\geq 0$.
\item[(iii)] $\Eq(x,y-1)$ and $\Eq(x-1,y)$ together imply $\Eq(x,y)$ for $x,y\geq 1$.
\end{itemize}

Ad (i). $\Eq(0,y)$ writes $1 = 1$.

Ad (ii). $\Eq(x,0)$ writes
\[
\sumd{j\in [0,x]} (-1)^j \smatze{t}{j} = \smatze{x-t}{x}.
\]
Proceeding by induction on $x$ and using (i) we are reduced to consider differences, i.e.\ to see that
\[
(-1)^x \smatze{t}{x} = \smatze{-(t-x+1)}{x}.
\]

Ad (iii). The right hand side of $\Eq(x,y)$ equals
\[ 
\begin{array}{rcl}
\smatze{y - t + x}{x}
& = & \smatze{y - t + (x-1)}{(x-1)} + \smatze{(y-1) - t + x}{x} \\
& \auf{\mbox{\scr by assumption}}{=} & \sumd{j\in [0,x-1]} (-1)^{x-1-j}\smatze{y+j}{j}\smatze{t}{x-1-j} + \sumd{j\in [0,x]} (-1)^{x-j}\smatze{y-1+j}{j}\smatze{t}{x-j} \\ 
& = & \sumd{j\in [1,x]} (-1)^{x-j}\smatze{y+j-1}{j-1}\smatze{t}{x-j} + \sumd{j\in [0,x]} (-1)^{x-j}\smatze{y+j-1}{j}\smatze{t}{x-j} \\ 
& = & \sumd{j\in [0,x]} (-1)^{x-j}\smatze{y+j}{j}\smatze{t}{x-j}, \\ 
\end{array}
\]
being the left hand side.
\end{Lemma}

\begin{Lemma}
\label{LemLW2}
Let $1\leq k\leq x$. We have
\[
\gcd(\smatze{x}{1},\smatze{x}{2},\dots,\smatze{x}{k}) = x\cdot\prodd{p \mbox{\rm\scr\ prime, } p | x} p^{-\min(v_p(x),\,\silog_p(k))},
\]
where 
\[
\ilog_p(k) := \max\{ i\geq 0 \;|\; p^i\leq k\} = \max v_p([1,k]).
\]

\rm
Let $p$ be a prime, let $j\in [0,x]$. For an integer $w\geq 0$, written $p$-adically as $w = \sumd{i\geq 0} w_i p^i$, $w_i\in [0,p-1]$, we denote its Quersumme by $q_p(w) := \sumd{i\geq 0} w_i$.
We obtain
\[
v_p(j!) = \sumd{i\geq 0} w_i \cdot\fracd{p^i - 1}{p-1} = \fracd{j - q_p(j)}{p-1},
\]
and thus
\[
v_p(\smatze{x}{j}) = (p-1)^{-1}(q_p(x - j) - (q_p(x) - q_p(j))).
\]
In particular, in case $j\leq p^{v_p(x)}$ we get $q_p(x - j) = (q_p(x) - 1) + (p-1)(v_p(x) - v_p(j) - 1) + p\cdot 1 - q_p(j)$, whence 
\[
v_p(\smatze{x}{j}) = v_p(x) - v_p(j).
\]
Thus the minimal such valuation for $j\in [1,k]$ takes the value claimed above.
\end{Lemma}

Let $n\geq 1$, $d\geq 1$. Assume given a partition $\lambda$ of $n$ and an integer $g\in [1,\lambda_1 - 1]$ such that
\[
\mu'_i := 
\left\{
\begin{array}{ll}
\lambda'_i + d & \mbox{for } i = g \\
\lambda'_i - d & \mbox{for } i = g+1 \\
\lambda'_i     & \mbox{else} \\
\end{array}
\right.
\] 
defines a partition $\mu$. Let $Z'$ be the set of injective maps
\[
[1,d]\lraa{\zeta'} [1,\lambda'_{g+1}].
\]
Let $Z\tm Z'$ be the subset of strictly monotone maps. Let the sign of $\zeta\in Z$ be given by $\eps_\zeta := (-1)^{\sum_{i\in [1,d]}\zeta(i)}$. 
For $\zeta'\in Z'$ there is a unique factorization $\zeta'(i) = \v\zeta'(i\sigma)$, $i\in [1,d]$, with $\sigma\in\Sl_d$ and $\v\zeta'\in Z$. Let the sign of $\zeta'\in Z'$ be given by 
$\eps_{\zeta'} := \eps_\sigma\eps_{\v\zeta'}$. Suppose given a $\lambda$-tableau $[a]$ and an injection $\zeta'\in Z'$. We let
\[
\begin{array}{rcl}
[0,\lambda'_{g+1}]       & \lra        & [0,\mu'_{g+1}] \\
i                        & \lraa{\phi} & \#\Big( [1,i]\ohne\zeta'([1,d])\Big) \\
\min(\phi^{-1}(\{ j\}))  & \llaa{\psi} & j \\
\end{array}
\]
and define the $\mu$-tableau $[a^{\zeta'}]$ by
\[
\begin{array}{rcll}
a^{\zeta'}_{j,i}              & := & a_{j,i}           & \mbox{\rm for } (j\in [1,\lambda_1]\ohne\{ g,g+1\}\mbox{\rm\ and } i\in [1,\mu'_j]) 
                                                         \mbox{\rm\ or } (j = g \mbox{\rm\ and } i\in [1,\lambda'_g]) \\
a^{\zeta'}_{g,\lambda'_g + i} & := & a_{g+1,\zeta'(i)} & \mbox{\rm for } i\in [1,d] \\
a^{\zeta'}_{g+1,i}            & := & a_{g+1,\psi(i)}   & \mbox{\rm for } i\in [1,\mu'_{g+1}]. \\
\end{array}
\]
Using this place operation, we define the $\Z\Sl_n$-linear map 
\[
\begin{array}{rcl}
F^\lambda & \lraa{f''} & S^\mu \\
{[a]}     & \lra       & \sumd{\zeta\in Z} \spi{a^\zeta}\eps_\zeta = \frac{1}{d!}\sumd{\zeta'\in Z'} \spi{a^{\zeta'}}\eps_{\zeta'}. \\
\end{array}
\]
For a $\lambda$-tableau $[a]$ and for an injection $\zeta'\in Z'$, we denote the `$[a]$-realization' of $\zeta'$ by
\[
([1,d]\lraa{\zeta'_{[a]}} [1,n]) := \Big([1,d]\lraa{\zeta'} [1,\lambda'_{g+1}] \lraiso ([1,\lambda'_{g+1}] \ti \{ g+1\})\hra [\lambda] \lraa{[a]} [1,n]\Big). 
\]

\begin{Lemma}
\label{LemLW2_5}
We have a factorization
\[
(F^\lambda\lraa{f''}S^\mu) = (F^\lambda\lraa{\nu^\lambda_M} M^{\lambda',-}\lraa{f'} S^\mu).
\]

\rm
We need to show that given $v,w\in a_{g+1}$, $v\neq w$, and a $\lambda$-tableau $[a]$, we have $([a]\cdot (v,w))f'' = -[a]f''$. Suppose given $\zeta\in Z$. In case $v,w\in\zeta_{[a]}([1,d])$ or 
in case $v,w\not\in\zeta_{[a]}([1,d])$, we obtain $\spi{a^\zeta}\eps_\zeta\cdot (v,w) = -\spi{a^\zeta}\eps_\zeta$. The maps $\zeta\in Z$ that satisfy $v\in\zeta_{[a]}([1,d])$, 
$w\not\in\zeta_{[a]}([1,d])$ furnish a subset $Z_{10}\tm Z$, the maps $\zeta\in Z$ that satisfy $v\not\in\zeta_{[a]}([1,d])$, $w\in\zeta_{[a]}([1,d])$ furnish a subset $Z_{01}\tm Z$. We have a 
bijection $\iota$ that sends $\zeta\in Z_{10}$ to the map $\iota(\zeta)\in Z_{01}$ determined by $\iota(\zeta)_{[a]}\Big([1,d]\Big) = (\zeta_{[a]}\Big([1,d]\Big))\cdot (v,w)$. We claim that 
\[
\spi{a^\zeta}\eps_\zeta\cdot (v,w) = - \spi{a^{\iota(\zeta)}}\eps_{\iota(\zeta)},
\]
thus proving the lemma. Let $v =: a_{g+1,i}$, let $w =: a_{g+1,j}$ and assume $i < j$. Comparing both sides, the entries in column $g$ of $[\mu]$ yield a sign 
$(-1)^{\#([i+1,j-1]\cap\zeta([1,d]))}$, the entries in column $g+1$ of $[\mu]$ yield a sign $(-1)^{\#([i+1,j-1]\ohne\zeta([1,d]))}$, so, altogether, the entries yield a sign $(-1)^{j-i-1}$. 
On the other hand, the quotient of the signs of the maps is $\eps_{\iota(\zeta)}/\eps_\zeta = (-1)^{j-i}$.
\end{Lemma}

\begin{Remark}
\label{RemLW2_6}
\rm
Let the correspondence $\phi\in\Phi^{\lambda',\mu'}$ be defined by
\[
\begin{array}{rcll}
[\lambda']                & \lraisoa{\phi} & [\mu']                    & \\
j\ti i                    & \lra           & j\ti i                    & \mbox{ for } j\ti i \in [\lambda']\cap [\mu'] \\
(g+1)\ti (\mu'_{g+1} + i) & \lra           & g\ti (\lambda'_{g+1} + i) & \mbox{ for } i\in [1,d] \\
\end{array}
\]
The semistandard correspondoid $\{\phi\}$ is in fact the {\it only} element in $\{\Phi^{\lambda',\mu'}\}_{\mbox{\scr\rm sst}}$ (cf.\ after \ref{LemSST10}). With respect to 
$[a_{\lambda'}] = [\ck a_{\lambda'}]$ and $[a_{\mu'}] = [\ck a_{\mu'}]$, we obtain $f' = - \Theta^-_{\phi^{-1}}\cdot\nu_S^\mu$, which reproves (\ref{LemLW2_5}) by means of (\ref{LemSST1}).

From (\ref{CorSST8}, \ref{CorSST7}) we take the following assertions. Suppose given $\w m\geq 1$. 

In case $\mu'$ is $2$-regular, the map
\[
\Hom_{\sZ\Sl_n}(M^{\lambda',-}/\w m, M^{\mu',-}/\w m) \;\;\lraa{(-)\nu_S^\mu}\;\; \Hom_{\sZ\Sl_n}(M^{\lambda',-}/\w m, S^\mu/\w m)
\]
is surjective. In particular, the group $\Hom_{\sZ\Sl_n}(S^\lambda/\w m,S^\mu/\w m)$ is cyclic since we may embed
\[
\Hom_{\sZ\Sl_n}(S^\lambda/\w m,S^\mu/\w m) \;\;\lraa{\nu_S^\lambda(-)}\;\;\Hom_{\sZ\Sl_n}(M^{\lambda',-}/\w m, S^\mu/\w m).
\]
A generator is a factorization of some multiple of $f'\ts_{\sZ}\Z/\w m$ over $\nu_S^\lambda$.

In case $\mu'$ is $2$-singular, the map
\[
2\cdot\Hom_{\sZ\Sl_n}(M^{\lambda',-}/\w m, M^{\mu',-}/\w m) \;\;\lraa{(-)\nu_S^\mu}\;\; 2\cdot\Hom_{\sZ\Sl_n}(M^{\lambda',-}/\w m, S^\mu/\w m)
\]
is surjective. In particular, the group $2\cdot\Hom_{\sZ\Sl_n}(S^\lambda/\w m,S^\mu/\w m)$ is cyclic since we may embed
\[
2\cdot\Hom_{\sZ\Sl_n}(S^\lambda/\w m,S^\mu/\w m) \;\;\lraa{\nu_S^\lambda(-)}\;\;2\cdot\Hom_{\sZ\Sl_n}(M^{\lambda',-}/\w m, S^\mu/\w m).
\]
A generator is a factorization of some multiple of $2\cdot f'\ts_{\sZ}\Z/\w m$ over $\nu_S^\lambda$.
\end{Remark}

Let $\xi\tm a_g$, $s := \#\xi$, $\b\xi := a_g\ohne\xi$, $t\in [1,\lambda'_{g+1}]$, $\eta := a_{g+1,[1,t]}$, $\b\eta := a_{g+1,[t+1,\lambda'_{g+1}]}$ be such that $s+t = \lambda'_g + 1$. Let
$u := \#\b\eta = \lambda'_{g+1} - t$.

For $x\in [0,\min(d,t)]$, we denote
\[
\begin{array}{rcl}
Z'[x]             & := & \{ \zeta'\in Z'\; |\; \zeta'(i) = i \mb{ for } i\in [1,x],\; \zeta'(i)\geq t+1 \mb{ for } i\in [x+1,d]\} \\
Z[x]              & := & Z\cap Z'[x] \\
\Lambda([a],x)    & := & \frac{1}{s!t!}\sumd{\zeta\in Z[x]} \spi{a^\zeta}\eps_\zeta\c (\xi\cup\eta) \\
                  & =  & \frac{1}{s!t!(d-x)!}\sumd{\zeta'\in Z'[x]} \spi{a^{\zeta'}}\eps_{\zeta'}\c (\xi\cup\eta) \\
\end{array}
\]
so that we can recover
\[
G'_{[a],\xi,\eta} f' = \sumd{x\in [0,\min(d,t)]} \smatze{t}{x} \Lambda([a],x).
\]
For $0\leq x\leq y\leq \min(d,t)$, there is a {\it rectification map}
\[
\begin{array}{rcl}
Z'[x]  & \lraa{[y]} & Z'[y] \\
\zeta' & \lra       & \zeta'[y] \\
\end{array}
\]
defined by $\zeta'[y](i) := i$ for $i\in [1,y]$ and $\zeta'[y](i) := \zeta'(i)$ for $i\in [y+1,d]$. In case $t\geq d$, we denote by $\zeta[d]$ the unique element in $Z[d] = Z'[d]$, mapping 
each $i\in [1,d]$ to $i\in [1,\lambda'_{g+1}]$.

\begin{Lemma}
\label{LemLW3}
Suppose $t\geq d$. Let $\b\xi = \b\xi_{-1}\cup \b\xi_2$ be a disjoint decomposition such that $\#\b\xi_{-1} = d - 1$ and such that $\#\b\xi_2 = t - d$.
Let $\eta_2 := a_{g+1,[d+1,t]}$. We obtain, in the notation of (\ref{LemG7}),
\[
G'_{[a],\xi,\eta}f' = \eps_{\zeta[d]}\cdot\sumd{i\in [0,d-1]} (-1)^i \smatze{s-u+d}{d-i}\cdot B_{[a^{\zeta[d]}]\cdot (\b\xi_2,\eta_2)_\alpha,\b\xi_{-1},\b\xi_2}(i),
\]
independent of the chosen bijection $\b\xi_2\lraisoa{\alpha}\eta_2$.

\rm
Given $d_1\in [0,d]$ and a map $\zeta\in Z[d - d_1]$, we denote
\[
\begin{array}{rcl}
\eta_1          & := & a_{g+1,[d-d_1+1,d]} \\
\b\eta_1(\zeta) & := & \zeta_{[a]}([d-d_1+1,d]). \\
\end{array}
\]
We obtain
\[
\begin{array}{rcl}
\Lambda([a],d-d_1) 
& = & \frac{1}{s!t!} \sumd{\zeta\in Z[d - d_1]} \spi{a^\zeta}\eps_\zeta\c (\xi\cup\eta) \\
& \auf{\mbox{\scr (\ref{LemG4})}}{=} & \frac{(s+d-d_1)!}{s!t!(d-1)!} \sumd{\zeta\in Z[d - d_1]} 
                                       \spi{a^\zeta}\eps_\zeta\cdot(\b\xi_2,\eta_2)\cdot(\b\eta_1(\zeta),\eta_1)\c(\b\xi\cup\b\eta_1(\zeta)) \\
& = & (-1)^{d_1}\cdot\frac{(s+d-d_1)!}{s!t!(d-1)!} \sumd{\auf{\scm\b\eta_1\tm\b\eta,}{\#\b\eta_1\; =\; d_1}} \spi{a^{\zeta[d]}}\eps_{\zeta[d]}\cdot(\b\xi_2,\eta_2)\c(\b\xi\cup\b\eta_1) \\
& \auf{\mbox{\scr (cf.\ \ref{LemG7})}}{=} & (-1)^{d_1}\cdot\frac{(s+d-d_1)!}{s!t!(d-1)!}\cdot\eps_{\zeta[d]}\cdot C_{[a^{\zeta[d]}]\cdot (\b\xi_2,\eta_2)_\alpha,\b\xi_{-1},\b\xi_2}(d_1)\\
& \auf{\mbox{\scr (\ref{LemG7})}}{=}      & (-1)^{d_1}\cdot\frac{(s+d-d_1)!(d-1)!(t - d + d_1)!}{s!t!(d-1)!}\cdot\eps_{\zeta[d]}\cdot \\
&                                         & \cdot \sumd{i\in [0,\min(d_1,d-1)]} \smatze{u-i}{d_1-i} B_{[a^{\zeta[d]}]\cdot (\b\xi_2,\eta_2)_\alpha,\b\xi_{-1},\b\xi_2}(i), \\
\end{array}
\]
whence our one-step Garnir relation is mapped to
\[
\begin{array}{rl}
  & G'_{[a],\xi,\eta} f' \\
= & \eps_{\zeta[d]}\cdot\sumd{d_1\in [0,d]}\;\; \sumd{i\in [0,\min(d_1,d-1)]}  \smatze{t}{d-d_1} (-1)^{d_1}\cdot\frac{(s+d-d_1)!(t - d + d_1)!}{s!t!}\cdot\smatze{u-i}{d_1-i} \cdot\\
  &\cdot B_{[a^{\zeta[d]}]\cdot (\b\xi_2,\eta_2)_\alpha,\b\xi_{-1},\b\xi_2}(i) \\
= & \eps_{\zeta[d]}\cdot\sumd{i\in [0,d-1]}\left(\sumd{j\in [0,d-i]} (-1)^{i+d-i-j}\smatze{s+j}{j}\smatze{u-i}{d-i-j}\right) B_{[a^{\zeta[d]}]\cdot (\b\xi_2,\eta_2)_\alpha,\b\xi_{-1},\b\xi_2}(i) \\
 \auf{\mbox{\scr (\ref{LemLW1})}}{=} & \eps_{\zeta[d]}\cdot\sumd{i\in [0,d-1]} (-1)^i \smatze{s-u+d}{d-i} B_{[a^{\zeta[d]}]\cdot (\b\xi_2,\eta_2)_\alpha,\b\xi_{-1},\b\xi_2}(i). \\
\end{array}
\]
\end{Lemma}

\begin{Lemma}
\label{RemLW3_5}
Suppose $2d - 1\leq\lambda'_{g+1}$. Let $m_0 :=  \lambda'_g - \lambda'_{g+1} + d + 1$ be the box shift length, let 
\[
m := m_0\cdot\prodd{p \mbox{\rm\scr\ prime, } p | m_0} p^{-\min(v_p(m_0),\,\silog_p(d))},
\]
In case $\mu'$ is $2$-regular, the injection 
\[
\Hom_{\sZ\Sl_n}(S^\lambda,S^\mu/m) \mra{(-)\cdot (n!/m)}\Hom_{\sZ\Sl_n}(S^\lambda,S^\mu/n!)
\]
is surjective. In case $\mu'$ is $2$-singular, the cokernel of the injection 
\[
\Hom_{\sZ\Sl_n}(S^\lambda,S^\mu/m) \mra{(-)\cdot (n!/m)} \Hom_{\sZ\Sl_n}(S^\lambda,S^\mu/n!)
\]
is annihilated by multiplication by $2$.

\rm
Consider (\ref{LemLW3}) in the case $[a] = [\ck a_\lambda]$, $t = d$, $\xi = a_{g,[1,s]}$, $\b\xi_{-1} = \b\xi$ and $\b\xi_2 = \leer$. We obtain
\[
\begin{array}{rcl}
G'_{[\ck a_\lambda],\xi,\eta} f' 
& = & \eps_{\zeta[d]}\cdot\sumd{i\in [0,d-1]} (-1)^i \smatze{s-u+d}{d-i}\cdot B_{[(\ck a_\lambda)^{\zeta[d]}],\b\xi,\leer}(i) \\
& \auf{\mb{\scr (\ref{LemG7})}}{=} & \eps_{\zeta[d]}\cdot\sumd{i\in [0,d-1]} 
      \smatze{m_0}{d-i}\cdot \sumd{\auf{\scm \b\xi_0\tm\b\xi,}{\#\b\xi_0 \; =\; i}}\;\;\sumd{\auf{\scm \phi_0\tm\b\eta,}{\#\phi_0 \; =\; i}} \spi{(\ck a_\lambda)^{\zeta[d]}} (\b\xi_0,\phi_0). \\
\end{array}
\]
Since $s = \lambda'_g + 1 - d \geq \mu'_{g+1} + 1$, the elements 
$\spi{(\ck a_\lambda)^{\zeta[d]}} (\b\xi_0,\phi_0)$ occurring in this expression are standard up to column permutation, i.e.\ up to sign. Moreover, they are pairwise different because of 
different column fillings. We note that each $i\in [0,d-1]$ indexes a nonzero summand of this expresseion, namely a nonzero linear combination of different standard polytabloids equipped 
with coefficients $\pm\smatze{m_0}{d-i}$, since we assumed $2d - 1\leq \lambda'_{g+1}$, i.e.\ since $d-1\leq\#\b\eta$ leaves space for a subset $\phi_0\tm\b\eta$ of cardinality $i$.
Therefore, by (\ref{LemLW2}), the element $G'_{[\ck a_\lambda],\xi,\eta} f'\in S^\mu$ is exactly divisible by $m$. 

Suppose given a $\Z\Sl_n$-linear map $S^\lambda \lraa{\w f} S^\mu/n!$. Let $S^\mu\lraa{\phi} S^\mu/n!$ denote the residue class morphism. 

Suppose $\mu'$ to be $2$-regular. By (\ref{RemLW2_6}), there is an integer $z$ such that $\nu_S^\lambda \w f = z\cdot f'\phi$. Thus $z\cdot G'_{[\ck a_\lambda],\xi,\eta} f'$ is divisible
by $n!$, and so $z\cdot m$ is divisible by $n!$. Hence, each element in the image of $\w f$ is contained in $(n!/m)S^\lambda/n!S^\lambda$, i.e.\ there exists a factorization
\[
(S^\lambda\lraa{\w f} S^\mu/n!) = (S^\lambda\lra S^\mu/m\lraa{n!/m}S^\mu/n!).
\]

Suppose $\mu'$ to be $2$-singular. By (\ref{RemLW2_6}), there is an integer $z$ such that $2\cdot \nu_S^\lambda \w f = 2z\cdot f'\phi$. Thus $2z\cdot G'_{[\ck a_\lambda],\xi,\eta} f'$ is divisible
by $n!$, and so $2z\cdot m$ is divisible by $n!$. Hence, each element in the image of $2 \w f$ is contained in $(n!/m)S^\lambda/n!S^\lambda$, i.e.\ there exists a factorization
\[
(S^\lambda\lraa{2 \w f} S^\mu/n!) = (S^\lambda\lra S^\mu/m\lraa{n!/m}S^\mu/n!).
\]
\end{Lemma}

\begin{Lemma}
\label{LemLW4}
Suppose $t\leq d-1$. We obtain, in the notation of (\ref{LemG8}),
\[
G'_{[a],\xi,\eta} f' = \eps_{\zeta[d]}\cdot\sumd{i\in [0,t-1]} (-1)^i \smatze{s-u+d}{t-i} B'_{[a^{\zeta[d]}],\b\xi\cup\b\eta_2,\b\eta_2,\b\eta\ohne\b\eta_2}(i). 
\]

\rm
Given $d_0\in [0,t]$ and an injection $\zeta'\in Z'[d_0]$, we denote
\[
\begin{array}{rcl}
\eta_1           & := & \zeta'_{[a]}([1,d_0]) = a_{g+1,[1,d_0]}  \\
\b\eta_1(\zeta') & := & \zeta'_{[a]}([d_0 + 1,t]) \\
\b\eta_2(\zeta') & := & \zeta'_{[a]}([t+1,d]) \\
\b\eta_2         & := & a_{g+1,[t+1,d]}. \\
\end{array}
\]
Given $\zeta'\in Z'[t]$, we get
\[
(\ast) \hspace*{2cm} \spi{a^{\zeta'}} \eps_{\zeta'}= \spi{a^{\zeta'[d]}}\sigma\eps_\sigma \eps_{\zeta'[d]}
\]
for any permutation $\sigma\in\Sl_{\b\eta}$ that maps $(a_{g+1,i})\sigma = \zeta'_{[a]}(i)$ for $i\in [t+1,d]$. Note that 
$\#\{\sigma\in\Sl_{\b\eta}\; |\; (a_{g+1,i})\sigma = \zeta'_{[a]}(i) \mbox{ for } i\in [t+1,d]\} = (u - d + t)!$. 

For $d_0\in [0,t]$, we obtain
\[
\begin{array}{cl}
  & \Lambda([a],d_0) \\
= & \frac{1}{s!t!(d-d_0)!} \sumd{\zeta'\in Z'[d_0]} \spi{a^{\zeta'}}\eps_{\zeta'}\c (\xi\cup\eta) \\
\auf{\mbox{\scr (\ref{LemG4})}}{=}      & \frac{(s+d_0)!}{s!t!(d-d_0)!(d-1)!} \sumd{\zeta'\in Z'[d_0]} \spi{a^{\zeta'}}\eps_{\zeta'}\cdot(\eta\ohne\eta_1,\b\eta_1(\zeta'))
                                          \c (\b\xi\cup\b\eta_1(\zeta')\cup\b\eta_2(\zeta')) \\
=                                       & (-1)^{t-d_0}\frac{(s+d_0)!(t - d_0)!}{s!t!(d-d_0)!(d-1)!}\sumd{\zeta'\in Z'[t]}\;\;
                                          \sumd{\auf{\scm\b\eta_1\tm\b\eta\ohne\b\eta_2(\zeta'),}{\#\b\eta_1 \; =\; t - d_0}}
                                          \spi{a^{\zeta'}}\eps_{\zeta'} \c (\b\xi\cup\b\eta_1\cup\b\eta_2(\zeta')) \\
\end{array}
\]
\[
\begin{array}{cl}
\auf{(\ast)}{=}                         & (-1)^{t-d_0}\frac{(s+d_0)!(t - d_0)!}{s!t!(d-d_0)!(d-1)!(u-d+t)!}\sumd{\auf{\scm\b\eta_1\tm\b\eta\ohne\b\eta_2,}{\#\b\eta_1 \; =\; t - d_0}}
                                          \spi{a^{\zeta[d]}}\eps_{\zeta[d]} \c (\b\xi\cup\b\eta_1\cup\b\eta_2)\c\b\eta \\
\auf{\mbox{\scr (cf.\ \ref{LemG8})}}{=} & (-1)^{t-d_0}\cdot\eps_{\zeta[d]}\cdot\frac{(s+d_0)!(t-d_0)!}{s!t!(d-d_0)!(d-1)!(u-d+t)!}\cdot  
                                          C'_{[a^{\zeta[d]}],\b\xi\cup\b\eta_2,\b\eta_2,\b\eta\ohne\b\eta_2}(t - d_0) \\
\auf{\mbox{\scr (\ref{LemG8})}}{=}      & (-1)^{t-d_0}\cdot\eps_{\zeta[d]}\cdot\frac{(s+d_0)!(t-d_0)!}{s!t!(d-d_0)!(d-1)!(u-d+t)!}\cdot \\
                                        & \cdot {\scm (u-d+t)!(d-1)!(d-d_0)!}\cdot\sumd{i\in [0,\min(t-d_0,t-1)]} \smatze{u-d+t-i}{t-d_0-i} 
                                          B'_{[a^{\zeta[d]}],\b\xi\cup\b\eta_2,\b\eta_2,\b\eta\ohne\b\eta_2}(i), \\

\end{array}
\]
whence
\[
\begin{array}{rl}
  & G'_{[a],\xi,\eta} f' \\
= & \eps_{\zeta[d]}\cdot\sumd{d_0\in [0,t]}\;\;\sumd{i\in [0,\min(t-d_0,t-1)]} \smatze{t}{d_0} (-1)^{t-d_0}\cdot\frac{(s+d_0)!(t-d_0)!}{s!t!}\cdot\\
  & \cdot\smatze{u-d+t-i}{t-d_0-i} B'_{[a^{\zeta[d]}],\b\xi\cup\b\eta_2,\b\eta_2,\b\eta\ohne\b\eta_2}(i) \\
= & \eps_{\zeta[d]}\cdot\sumd{i\in [0,t-1]}\left(\sumd{d_0\in [0,t-i]}  
    (-1)^{i+t-i-d_0}\cdot\smatze{s+d_0}{d_0}\cdot\smatze{u-d+t-i}{t-i-d_0}\right) B'_{[a^{\zeta[d]}],\b\xi\cup\b\eta_2,\b\eta_2,\b\eta\ohne\b\eta_2}(i) \\
\auf{\mbox{\scr (\ref{LemLW1})}}{=} & \eps_{\zeta[d]}\cdot\sumd{i\in [0,t-1]} (-1)^i
    \smatze{s-u+d}{t-i} B'_{[a^{\zeta[d]}],\b\xi\cup\b\eta_2,\b\eta_2,\b\eta\ohne\b\eta_2}(i). \\
\end{array}
\]
\end{Lemma}

\begin{Remark}
\label{RemLW5}
The morphism $F^\lambda\lraa{f''} S^\mu$ maps $[\ck a_\lambda]$ to a linear combination of standard $\mu$-polytabloids with coefficients $\pm 1$.
\end{Remark}

We summarize.

\begin{Theorem}
\label{ThLW6}
Let $m_0 :=  \lambda'_g - \lambda'_{g+1} + d + 1$ be the box shift length, let 
\[
m := m_0\cdot\prodd{p \mbox{\rm\scr\ prime, } p | m_0} p^{-\min(v_p(m_0),\,\silog_p(d))},
\]
where $\ilog_p(k) = \max\{i\geq 0 \; |\; p^i\leq k\}$. The $\Z\Sl_n$-linear map $M^{\lambda',-}\lraa{f'} S^\mu$ factors over 
\begin{center}
\begin{picture}(250,250)
\put( -10, 200){$M^{\lambda',-}$}
\put( 100, 210){\vector(1,0){80}}
\put( 120, 225){$\scm f'$}
\put( 200, 200){$S^\mu$}
\put(-110, 200){$\{ a'\}^-$}
\put( -80, 180){\vector(0,-1){130}}
\put(-110,   0){$\spi{a}$}
\put(  20, 180){\vector(0,-1){130}}
\put( -20, 110){$\scm\nu^\lambda_S$}
\put( 220, 180){\vector(0,-1){130}}
\put( 350, 200){$\spi{b}$}
\put( 380, 180){\vector(0,-1){130}}
\put( 350,   0){$\spi{b} + m S^\mu$,}
\put( -10,   0){$S^\lambda$}
\put(  50,  10){\vector(1,0){130}}
\put( 110,  25){$\scm f$}
\put( 200,   0){$S^\mu/m$}
\end{picture}
\end{center}
The resulting morphism $S^\lambda\lraa{f}S^\mu/m$ is of order $m$ in $\Hom_{\sZ\Sl_n}(S^\lambda,S^\mu/m)$. 

In case $\mu'$ is $2$-regular, the group $\Hom_{\sZ\Sl_n}(S^\lambda,S^\mu/m)$ is generated by $f$. 

In case $\mu'$ is $2$-singular, the group $2\cdot\Hom_{\sZ\Sl_n}(S^\lambda,S^\mu/m)$ is generated by $2f$.

In case $2d-1\leq\lambda'_{g+1}$ and $\mu'$ is $2$-regular, there is an isomorphism
\[
\Hom_{\sZ\Sl_n}(S^\lambda,S^\mu/m) \mraisoa{(-)\cdot (n!/m)}\Hom_{\sZ\Sl_n}(S^\lambda,S^\mu/n!).
\] 

In case $2d-1\leq\lambda'_{g+1}$ and $\mu'$ is $2$-singular, the cokernel of the inclusion
\[
\Hom_{\sZ\Sl_n}(S^\lambda,S^\mu/m) \mraisoa{(-)\cdot (n!/m)}\Hom_{\sZ\Sl_n}(S^\lambda,S^\mu/n!),
\]
is annihilated by multiplication with $2$. 

\rm
This follows by (\ref{LemLW3}, \ref{LemLW4}, \ref{LemLW2}, \ref{RemLW5}, \ref{RemLW2_6}, \ref{RemLW3_5}). 
\end{Theorem}

\begin{Example}
\label{ExLW6_5}
\rm
We have $\Hom_{\sZ\Sl_8}(S^{(2^3,1^2)},S^{(1^8)}/8!) = 0$. 
\end{Example}

\begin{Question}
\label{QuLW6_6}
\rm
In case $2d-1>\lambda'_{g+1}$, I do not know of a counterexample to the assertions in (\ref{ThLW6}) that presuppose $2d-1\leq\lambda'_{g+1}$.
\end{Question}

\begin{Remark}
\label{RemLW7}
\rm
In the situation of (\ref{ThLW6}), {\sc Carter} and {\sc Payne} [CP 80] have shown that 
\[
\Hom_{K\Sl_n}(K\ts_{\sZ} S^\lambda,K\ts_{\sZ} S^\mu) \neq 0,
\]
$K$ being an infinite field of characteristic $p$ such that $\ilog_p d < v_p(m_0)$. This part of their result is recovered by (\ref{ThLW6}).
\end{Remark}

\begin{Corollary}[{to (\ref{ThLW6}) via (\ref{PropSST11}, \ref{RemLW2_6})}]
\label{CorLW8}
The transpose in the sense of (\ref{PropSST11}) of the morphism given in (\ref{ThLW6}), interpreted as $S^\lambda/m\lraa{f}S^\mu/m$, is given by
\[
\begin{array}{rcl}
S^{\mu'}/m & \lraa{f^\trp} & S^{\lambda'}/m \\
\spi{b'}   & \lra          & -\spi{b'}\Theta_\phi, \\
\end{array}
\]
where the correspondoid $\{\phi\}\in\{\Phi^{\lambda',\mu'}\}$ is given by 
\[
\begin{array}{rcll}
[\mu']                    & \lraa{\{\phi\}}    & \N                        & \\
j\ti i                    & \lra               & j                         & \mbox{ for } j\ti i \in [\lambda']\cap [\mu'] \\
g\ti (\lambda'_g + i)     & \lra               & g + 1                     & \mbox{ for } i\in [1,d] \\
\end{array}
\]
All assertions on the $\Hom$-groups in (\ref{ThLW6}) have a counterpart, obtained by isomophic transport via the transposition isomorphism (\ref{PropSST11}).

\rm 
This follows by (\ref{ThLW6}, \ref{PropSST11}) using (\ref{RemLW2_6}). Concerning $\Theta_\phi$, we refer to (\ref{LemSST1}).
\end{Corollary}

\begin{Question}
\label{QuLW9}
\rm
I do not know a formula for the transpose $f^\trp$ in (\ref{CorLW8}) in terms of $\lambda'$-polytabloids. Cf.\ (\ref{QuSST11_1}, \ref{ExLW10}).
\end{Question}

\begin{Example}[{a fixed point, cf.\ [J 78, 24.4], [K 99, 4.4.1]}]
\label{ExLW10}
\rm\Absatz
Consider the case $n\geq 4$, $d = 2$, $\lambda = (2^2,1^{n-4})$ and $\mu = (1^n)$. The modulus becomes $m = n-1$ if $n-1$ is odd and $m = (n-1)/2$ if $n-1$ is even. We obtain 
\[
\begin{array}{rcl}
[\mu'] & \lraa{\{\phi\}} & \N \\
1\ti i & \lra            & \left\{\begin{array}{ll}
                           1 & \mb{for } i\in [1,n-2]\\
                           2 & \mb{for } i\in [n-1,n],\\
                           \end{array}\right.
\end{array}
\]
which allows to calculate the transpose as mapping
\[
\begin{array}{rcl}
S^{\mu'}/m & \lraa{f^\trp} & S^{\lambda'}/m \\
\spi{x}    & \lra          & -\sumd{i,j\in [1,n],\; i < j} \{ i,j\}, \\
\end{array}
\]
where $[x]$ is the unique standard $\mu'$-tableau and where we denote a $\lambda'$-tabloid by its second row. Similarly, for a $\lambda'$-polytabloid, we omit to denote the entries 
in the first row from column three onwards. I.e.\ we let a $\lambda'$-polytabloid be determined by the first two entries of the first and of the second row - which is not to be confused with a 
$(2,2)$-polytabloid. Let
\[
A := \sumd{j\in [4,n]} (j-2)\cdot \smatspzz{1}{3}{2}{j} \;\; - \sumd{i,j\in [3,n],\; i < j} \smatspzz{1}{2}{i}{j} \in S^{\lambda'}.
\]
For $k\in [3,n]$, $l\in [4,n]$, $k < l$, we obtain the following list of values of the scalar product in $M^{\lambda'}$.
\[
\begin{array}{lcclcl}
(A, \{ k,l\}) & = & + & 0                & - & 1        \\
(A, \{ 1,l\}) & = & - & (l-2)            & + & (l-3)    \\
(A, \{ 2,l\}) & = & + & (l-2)            & + & (n-l)    \\
(A, \{ 2,3\}) & = & - & (n-2)(n-1)/2 + 1 & + & (n-3) \\
(A, \{ 1,2\}) & = & + & 0                & - & (n-2)(n-3)/2 \\
(A, \{ 1,3\}) & = & + & (n-2)(n-1)/2 - 1 & + & 0. \\
\end{array}
\]
Hence, written in polytabloids, our map turns out to be
\[
\begin{array}{rcl}
S^{\mu'}/m & \lraa{f^\trp} & S^{\lambda'}/m \\
\spi{x}    & \lra          & \sumd{j\in [4,n]} (j-2)\cdot \smatspzz{1}{3}{2}{j} \;\; - \sumd{i,j\in [3,n],\; i < j} \smatspzz{1}{2}{i}{j}.\\
\end{array}
\]
\end{Example}          
\appendix
\section{Two lemmata}
\label{AppTwoLem}

\begin{footnotesize}
We present full versions of the two main lemmata of Section \ref{SubSecLinComb}. For sake of completeness of the full versions, there is a certain overlap, concerning notation and conclusion, 
with their abridged versions (\ref{LemM2short}, \ref{LemM8short}). We maintain the notation of Section \ref{SubSecLinComb}.

\begin{Lemma}[full version of (\ref{LemM2short})]
\label{LemM2}
Suppose $g < p < k$ and $s,t\geq 2$. The map $f'$ annihilates $G'_{[a],\xi,\eta}$.

\rm
We fix a map
\[
\begin{array}{rcl}
[1,\lambda_1]\ohne\{p,p+1\} & \lraa{\w e} & [0,2] \\
j                           & \lra        & \w e_j \\
\end{array}
\]
that maps $g$ and $k+1$ to $e_g = e_{k+1} = 2$, and that maps $j\in [1,\lambda_1]\ohne [g,k+1]$ to $e_j = 0$. For $\alpha,\beta\in [0,2]$, we denote by $\w e\alpha\beta$ be the prolongation 
of $\w e$ to $[1,\lambda_1]$ defined by $\w e\alpha\beta|_{[1,\lambda_1]\ohne\{p,p+1\}} := \w e$, $(\w e\alpha\beta)_p := \alpha$ and $(\w e\alpha\beta)_{p+1} := \beta$. We contend that
\[
\sumd{\alpha,\beta\in [0,2]} X_p^{(2-\alpha)} X_{p+1}^{(2-\beta)} G'_{[a],\xi,\eta} f'_{\w e\alpha\beta} = 0,
\]
from which the lemma ensues.

Given $\alpha,\beta\in [0,2]$, given a subset $\w\Xi\tm [1,2]\ti\{ p,p+1\}$ with
\[
\begin{array}{rcl}
\#(\w\Xi\cap ([1,2]\ti \{ p\}))   & = & \alpha \\
\#(\w\Xi\cap ([1,2]\ti \{ p+1\})) & = & \beta, \\
\end{array}
\]
and given prescribed inverse images $\xi^-,\eta^-\tm [1,2]$ `within $\w\Xi$', i.e.\ such that
\[
(\ast) \hspace*{3cm}
\begin{array}{rcl}
\xi^-  \ti\{ p\}  & \tm & \w\Xi\cap ([1,2]\ti\{ p\}) \\
\eta^- \ti\{p+1\} & \tm & \w\Xi\cap ([1,2]\ti\{ p+1\}), \\
\end{array}
\hspace*{3cm}
\]
we let
\[
\begin{array}{rcrcl}
\dGamma(\w e,\w\Xi,\xi^-,\eta^-) & := & \Big\{ \gamma\in\dGamma(\w e\alpha\beta) & \Big| & \Xi_\gamma \cap ([1,2]\ti\{ p,p+1\}) = \w\Xi, \\
                                 &    &                                          &       &  \gamma^{-1}(\xi) = \xi^-\ti\{ p\},\; \gamma^{-1}(\eta) = \eta^-\ti\{ p+1\} \Big\} \\
\end{array}
\]
and form the partial sum
\[
\Lambda([a],\w e,\w\Xi,\xi^-,\eta^-) := \frac{1}{s!t!} \sumd{\gamma\in\dGamma(\w e,\w\Xi,\xi^-,\eta^-)} \spi{a^\gamma}\eps_\gamma\c (\xi\cup\eta)
\]
so that we can recover
\[
(\ast\ast) \hspace*{2cm}
G'_{[a],\xi,\eta} f'_{\w e\alpha\beta} = \sumd{\w\Xi\tm [1,2]\ti\{ p,p+1\}}\;\; \sumd{\xi^-,\eta^- \tm [1,2] \mbox{\scr\ subject to } (\ast)} \Lambda([a],\w e,\w\Xi,\xi^-,\eta^-).
\hspace*{2cm}
\]
There exist elements $x,y\in\xi$, $x\neq y$, $x', y'\in\eta$, $x'\neq y'$, $z\in\b\xi$, which we choose and fix. For $\gamma\in\dGamma(\w e00)$, we let $x_\gamma := a_{j,\b\gamma(j,1)}$, where
$j\in [p+2,k+1]$ is minimal with $1\ti j\in\Xi_\gamma$, and $y_\gamma := a_{j,\b\gamma(j,2)}$, where $j\in [p+2,k+1]$ is minimal with $2\ti j\in\Xi_\gamma$. I.e.\ we pick the entries 
$x_\gamma, y_\gamma$ that `cross the columns' $p$ and $p+1$ under the operation of $\gamma$. We write 
\[
\begin{array}{rcl}
U_\gamma     & := & (s+t-2)!^{-1} \cdot\spi{a^\gamma}\eps_\gamma\cdot (\b\xi\ohne z,\eta\ohne\{ x',y'\})\cdot (x_\gamma,x,x')\cdot (y_\gamma,y,y') \c (\xi\cup\eta) \\
V_{1,\gamma} & := & (s+t-1)!^{-1} \cdot\spi{a^\gamma}\eps_\gamma\cdot (\b\xi,\eta\ohne x')\cdot (x_\gamma,x,x')\c (\xi\cup\eta) \\
V_{2,\gamma} & := & (s+t-1)!^{-1} \cdot\spi{a^\gamma}\eps_\gamma\cdot (\b\xi,\eta\ohne y')\cdot (y_\gamma,y,y')\c (\xi\cup\eta), \\
\end{array}
\]
and let
\[
\begin{array}{lcl}
A   & := & \sumd{\gamma\in \dGamma(\w e00)} U_\gamma\cdot\sumd{w'\in\b\eta} (w',z) \\
B   & := & \sumd{\gamma\in \dGamma(\w e00)} U_\gamma\cdot(1 - \sumd{w\in\b\xi\ohne z} (w,z)) \\
C_1 & := & \sumd{\gamma\in \dGamma(\w e00)} U_\gamma\cdot (z,x_\gamma) \\
C_2 & := & \sumd{\gamma\in \dGamma(\w e00)} U_\gamma\cdot (z,y_\gamma) \\
D   & := & \sumd{\gamma\in \dGamma(\w e00)} V_{1,\gamma} \cdot \sumd{w\in\b\xi} (w,y_\gamma) \\
    &  = & \sumd{\gamma\in \dGamma(\w e00)} V_{2,\gamma} \cdot \sumd{w\in\b\xi} (w,x_\gamma) \\
H   & := & \sumd{\gamma\in \dGamma(\w e00)} V_{1,\gamma} \cdot \sumd{w'\in\b\eta} (w',y_\gamma) \\
    &  = & \sumd{\gamma\in \dGamma(\w e00)} V_{2,\gamma} \cdot \sumd{w'\in\b\eta} (w',x_\gamma) \\
F_1 & := & \sumd{\gamma\in \dGamma(\w e00)} V_{1,\gamma} \\
F_2 & := & \sumd{\gamma\in \dGamma(\w e00)} V_{2,\gamma}. \\
\end{array}
\]
The equalities herein follow by the argument of (\ref{LemM0} i) and by the independence of of the respective expressions of the choice of $x,y,x',y'$.

{\bf Calculation of $\Lambda$-values.} We shall distinguish subcases and subsubcases according to the summation in $(\ast\ast)$. 

To distinguish subcases, e.g.\ $\smateckzz{+}{-}{+}{+}$ designates the subset $\w\Xi = \{ 1\ti p,\; 2\ti p,\; 2\ti (p+1)\}$; in general, the first factor of an element 
of $[1,2]\ti \{p,p+1\}$ counts rows, the second counts columns, and a plus sign $+$ denotes its appearance in $\w\Xi$. 

To distinguish subsubcases, we denote e.g.\ by $\smateckzz{\ds}{-}{+}{\ds}$ the configuration $\w\Xi = \smateckzz{+}{-}{+}{+}$, $\xi^- = \{ 1\}$ and $\eta^- = \{ 2\}$; in general, 
the sign $\ds$ in the left column means that its row number is an element of $\xi^-$, the sign $\ds$ in the right column means that its row number is an element of $\eta^-$. Concerning
the underlying set $\w\Xi$, the symbols $\ds$ and $+$ are synonymous.

{\bf Case $f'_{\w e22}$.} 

{\it Subcase $\w\Xi = \smateckzz{+}{+}{+}{+}$.}

{\it Subsubcase $(\w\Xi,\xi^-,\eta^-) = \smateckzz{\ds}{\ds}{\ds}{\ds}$.} We calculate
\[
\begin{array}{cl}
  & \Lambda([a],\w e,\pazz{+}{+}{+}{+},\{ 1,2\}, \{1,2\}) \\
= & \frac{1}{s!t!} \sumd{\gamma\in\dGamma(e,\w\Xi,\{ 1,2\}, \{1,2\})} \spi{a^\gamma}\eps_\gamma\c (\xi\cup\eta) \\
= & \frac{1}{s!t!} \sumd{\gamma\in\dGamma(\w e00)}\;\;\sumd{x_1\neq y_1\in\xi,\; x'_1\neq y'_1\in\eta} \spi{a^\gamma}\eps_\gamma \cdot(x_\gamma,x_1,x_1')\cdot(y_\gamma,y_1,y'_1)\c(\xi\cup\eta) \\
= & \frac{s(s-1)t(t-1)}{s!t!} \sumd{\gamma\in\dGamma(\w e00)} \spi{a^\gamma}\eps_\gamma\cdot(x_\gamma,x,x')\cdot(y_\gamma,y,y')\c(\xi\cup\eta) \\
\auf{\mbox{\scr(\ref{PropG6})}}{=} & \frac{s(s-1)t(t-1)\cdot s!}{s!t!(s+t-2)!} \sumd{\gamma\in\dGamma(\w e00)} \spi{a^\gamma}\eps_\gamma
                          \cdot (\b\xi\ohne z,\eta\ohne\{ x',y'\})\cdot(x_\gamma,x,x')\cdot(y_\gamma,y,y')\c(\xi\cup\eta)\c\b\xi \\
= & s(s-1)B.  \\
\end{array}
\]

{\it Subsubcase $(\w\Xi,\xi^-,\eta^-) = \smateckzz{\ds}{\ds}{\ds}{+}$.} We calculate
\[
\begin{array}{cl}
  & \Lambda([a],\w e,\pazz{+}{+}{+}{+},\{ 1,2\}, \{ 1\}) \\
= & \frac{s(s-1)t}{s!t!} \sumd{\gamma\in\dGamma(\w e00)}\;\;\sumd{y'_1\in\b\eta} \spi{a^\gamma}\eps_\gamma\cdot(x_\gamma,x,x')\cdot(y_\gamma,y,y'_1)\c (\xi\cup\eta) \\
\auf{\mbox{\scr(\ref{PropG6})}}{=} & \frac{s(s-1)t\cdot (s-1)!}{s!t!(s+t-2)!} \sumd{\gamma\in\dGamma(\w e00)}\;\;\sumd{y'_1\in\b\eta} \spi{a^\gamma}\eps_\gamma\cdot\\
                        & \cdot (x_\gamma,x,x')\cdot(y_\gamma,y,y'_1)\cdot(\b\xi\ohne z,\eta\ohne\{ x',y'\})\cdot(y'_1,y')\c (\xi\cup\eta)\c(\b\xi\cup y'_1) \\
= & - \frac{s-1}{(t-1)!(s+t-2)!} \sumd{\gamma\in\dGamma(\w e00)}\;\; \sumd{y'_1\in\b\eta} \spi{a^\gamma}\eps_\gamma\cdot\\
  & \cdot(x_\gamma,x,x')\cdot(y_\gamma,y,y'_1)\cdot(\b\xi\ohne z,\eta\ohne\{ x',y'\})\cdot(y'_1,y')\cdot (y'_1,y_\gamma)\c (\xi\cup\eta)\c(\b\xi\cup y'_1)  \\
= & - \frac{s-1}{(t-1)!(s+t-2)!} \sumd{\gamma\in\dGamma(\w e00)}\;\; \sumd{y'_1\in\b\eta}\spi{a^\gamma}\eps_\gamma\cdot
         (x_\gamma,x,x')\cdot(y_\gamma,y,y')\cdot(\b\xi\ohne z,\eta\ohne\{ x',y'\})\c (\xi\cup\eta)\c(\b\xi\cup y'_1)  \\
= & - (s-1)(uB - A). \\
\end{array}
\]

{\it Subsubcase $(\w\Xi,\xi^-,\eta^-) = \smateckzz{\ds}{+}{\ds}{\ds}$.} By symmetry, we obtain from subsubcase $\smateckzz{\ds}{\ds}{\ds}{+}$ 
\[
\Lambda([a],\w e,\pazz{+}{+}{+}{+},\{ 1,2\}, \{ 2\}) = - (s-1)(uB - A). 
\]

{\it Subsubcase $(\w\Xi,\xi^-,\eta^-) = \smateckzz{\ds}{+}{\ds}{+}$.} We calculate
\[
\begin{array}{cl}
  & \Lambda([a],\w e,\pazz{+}{+}{+}{+},\{ 1,2\}, \{\}) \\
= & \frac{s(s-1)}{s!t!} \sumd{\gamma\in\dGamma(\w e00)}\;\; \sumd{x'_1\neq y'_1\in\b\eta} \spi{a^\gamma}\eps_\gamma\cdot(x_\gamma,x,x'_1)\cdot(y_\gamma,y,y'_1)\c (\xi\cup\eta) \\
\auf{\mbox{\scr(\ref{PropG6})}}{=} & \frac{s(s-1)\cdot (s-2)!}{s!t!(s+t-2)!}\sumd{\gamma\in\dGamma(\w e00)}\;\; \sumd{x'_1\neq y'_1\in\b\eta} \spi{a^\gamma}\eps_\gamma\cdot \\
                        & \cdot(x_\gamma,x,x'_1)\cdot(y_\gamma,y,y'_1)\cdot(\b\xi\ohne z,\eta\ohne\{ x',y'\})\cdot(x'_1,x')\cdot(y'_1,y')\c (\xi\cup\eta)\c (\b\xi\cup x'_1\cup y'_1) \\
= & \frac{1}{t!(s+t-2)!}\sumd{\gamma\in\dGamma(\w e00)}\;\; \sumd{x'_1\neq y'_1\in\b\eta}\spi{a^\gamma}\eps_\gamma\cdot \\
  & \cdot(x_\gamma,x,x'_1)\cdot(y_\gamma,y,y'_1)\cdot(\b\xi\ohne z,\eta\ohne\{ x',y'\})\cdot(x'_1,x')\cdot(y'_1,y')\cdot(x'_1,x_\gamma)\cdot(y'_1,y_\gamma)\c\\
  & \c (\xi\cup\eta)\c (\b\xi\cup x'_1\cup y'_1) \\
= & \frac{1}{t!(s+t-2)!}\sumd{\gamma\in\dGamma(\w e00)}\;\; \sumd{x'_1\neq y'_1\in\b\eta}\spi{a^\gamma}\eps_\gamma\cdot
       (x_\gamma,x,x')\cdot(y_\gamma,y,y')\cdot(\b\xi\ohne z,\eta\ohne\{ x',y'\})\c (\xi\cup\eta)\c (\b\xi\cup x'_1\cup y'_1) \\
= & (u-1)(uB - 2A). \\
\end{array}
\]

{\it Subsubcase $(\w\Xi,\xi^-,\eta^-) = \smateckzz{\ds}{\ds}{+}{\ds}$.} We calculate
\[
\begin{array}{cl}
  & \Lambda([a],\w e,\pazz{+}{+}{+}{+},\{ 1\}, \{ 1,2\}) \\
= & \frac{st(t-1)}{s!t!} \sumd{\gamma\in\dGamma(\w e00)}\;\;\sumd{y_1\in\b\xi}\spi{a^\gamma}\eps_\gamma\cdot(x_\gamma,x,x')\cdot(y_\gamma,y_1,y')\c (\xi\cup\eta) \\
\auf{\mbox{\scr(\ref{PropG6})}}{=} & \frac{st(t-1)\cdot (s+1)!(t-2)!}{s!t!(s+t-1)!} \sumd{\gamma\in\dGamma(\w e00)}\;\;\sumd{y_1\in\b\xi}\spi{a^\gamma}\eps_\gamma 
                          \cdot(x_\gamma,x,x')\cdot(y_\gamma,y_1,y')\cdot(\b\xi\ohne y_1,\eta\ohne\{ x',y'\})\c (\xi\cup\eta) \\
= & \frac{s(s+1)}{(s+t-1)!} \sumd{\gamma\in\dGamma(\w e00)}\;\;\sumd{y_1\in\b\xi}\spi{a^\gamma}\eps_\gamma\cdot(\b\xi,\eta\ohne x')\cdot(x_\gamma,x,x')\cdot(y_\gamma,y_1)\c (\xi\cup\eta) \\
= & (s+1)sD. \\
\end{array}
\]

{\it Subsubcase $(\w\Xi,\xi^-,\eta^-) = \smateckzz{\ds}{\ds}{+}{+}$.} We calculate 
\[
\begin{array}{cl}
  & \Lambda([a],\w e,\pazz{+}{+}{+}{+},\{ 1\}, \{ 1\}) \\
= & \frac{st}{s!t!} \sumd{\gamma\in\dGamma(\w e00)}\;\;\sumd{y_1\in\b\xi,\; y'_1\in\b\eta} \spi{a^\gamma}\eps_\gamma\cdot(x_\gamma,x,x')\cdot(y_\gamma,y_1,y'_1)\c (\xi\cup\eta) \\
\auf{\mbox{\scr(\ref{PropG6})}}{=} & \frac{st\cdot s!(t-1)!}{s!t!(s+t-1)!} \sumd{\gamma\in\dGamma(\w e00)}\;\;\sumd{y_1\in\b\xi,\; y'_1\in\b\eta} \spi{a^\gamma}\eps_\gamma\cdot 
                        (x_\gamma,x,x')\cdot(y_\gamma,y_1,y'_1)\cdot(\b\xi\ohne y_1,\eta\ohne\{ x',y'\})\cdot(y'_1,y')\c (\xi\cup\eta) \\
= & - \frac{s}{(s+t-1)!} \sumd{\gamma\in\dGamma(\w e00)}\;\;\sumd{y_1\in\b\xi,\; y'_1\in\b\eta} \spi{a^\gamma}\eps_\gamma\cdot
           (\b\xi,\eta\ohne x')\cdot(x_\gamma,x,x')\cdot(y_\gamma,y_1,y',y'_1)\cdot(y'_1,y')\cdot(y_\gamma,y'_1)\c (\xi\cup\eta) \\
= & - \frac{su}{(s+t-1)!} \sumd{\gamma\in\dGamma(\w e00)}\;\;\sumd{y_1\in\b\xi} \spi{a^\gamma}\eps_\gamma
      \cdot(\b\xi,\eta\ohne x')\cdot(x_\gamma,x,x')\cdot(y_\gamma,y_1)\c (\xi\cup\eta) \\
= & -suD. \\
\end{array}
\]

{\it Subsubcase $(\w\Xi,\xi^-,\eta^-) = \smateckzz{\ds}{+}{+}{\ds}$.} We calculate
\[
\begin{array}{cl}
  & \Lambda([a],\w e,\pazz{+}{+}{+}{+},\{ 1\}, \{ 2\}) \\
= & \frac{st}{s!t!} \sumd{\gamma\in\dGamma(\w e00)}\;\;\sumd{y_1\in\b\xi,\; x'_1\in\b\eta} \spi{a^\gamma}\eps_\gamma\cdot(x_\gamma,x,x'_1)\cdot(y_\gamma,y_1,y')\c (\xi\cup\eta) \\
\auf{\mbox{\scr(\ref{PropG6})}}{=} & \frac{st\cdot s!(t-1)!}{s!t!(s+1-1)!} \sumd{\gamma\in\dGamma(\w e00)}\;\;\sumd{y_1\in\b\xi,\; x'_1\in\b\eta} \spi{a^\gamma}\eps_\gamma\cdot
                       (x_\gamma,x,x'_1)\cdot(y_\gamma,y_1,y')\cdot(\b\xi\ohne y_1,\eta\ohne\{ x',y'\})\cdot(x'_1,x')\c (\xi\cup\eta) \\ 
= & - \frac{s}{(s+1-1)!}\sumd{\gamma\in\dGamma(\w e00)}\;\;\sumd{y_1\in\b\xi,\; x'_1\in\b\eta} \spi{a^\gamma}\eps_\gamma\cdot
        (\b\xi,\eta\ohne x')\cdot(x_\gamma,x,x'_1)\cdot(y_\gamma,y_1)\cdot(x'_1,x')\cdot (x_\gamma,x'_1)\c (\xi\cup\eta) \\ 
= & - \frac{su}{(s+1-1)!}\sumd{\gamma\in\dGamma(\w e00)}\;\;\sumd{y_1\in\b\xi} \spi{a^\gamma}\eps_\gamma\cdot (\b\xi,\eta\ohne x')\cdot(x_\gamma,x,x')\cdot(y_\gamma,y_1)\c (\xi\cup\eta) \\ 
= & -suD. \\
\end{array}
\]

{\it Subsubcase $(\w\Xi,\xi^-,\eta^-) = \smateckzz{\ds}{+}{+}{+}$.} We calculate
\[
\begin{array}{cl}
  & \Lambda([a],\w e,\pazz{+}{+}{+}{+},\{ 1\}, \{ \}) \\
= & \frac{s}{s!t!} \sumd{\gamma\in\dGamma(\w e00)}\;\;\sumd{y_1\in\b\xi,\; x'_1\neq y'_1\in\b\eta} \spi{a^\gamma}\eps_\gamma\cdot(x_\gamma,x,x'_1)\cdot(y_\gamma,y_1,y'_1)\c (\xi\cup\eta) \\
\auf{\mbox{\scr(\ref{PropG6})}}{=} & \frac{s\cdot (s-1)!t!}{s!t!(s+t-1)!} \sumd{\gamma\in\dGamma(\w e00)}\;\;\sumd{y_1\in\b\xi,\; x'_1\neq y'_1\in\b\eta} \spi{a^\gamma}\eps_\gamma\cdot \\
                        & \cdot(x_\gamma,x,x'_1)\cdot(y_\gamma,y_1,y'_1)\cdot(\b\xi\ohne y_1,\eta\ohne\{ x',y'\})\cdot (x'_1,x')\cdot(y'_1,y')\c (\xi\cup\eta) \\
= & \frac{1}{(s+t-1)!} \sumd{\gamma\in\dGamma(\w e00)}\;\;\sumd{y_1\in\b\xi,\; x'_1\neq y'_1\in\b\eta} \spi{a^\gamma}\eps_\gamma\cdot \\
  & \cdot(\b\xi,\eta\ohne x')\cdot(x_\gamma,x,x'_1)\cdot(y_\gamma,y_1,y',y'_1)\cdot (x'_1,x')\cdot(y'_1,y')\cdot(x_\gamma,x'_1)\cdot(y_\gamma,y'_1)\c (\xi\cup\eta) \\
= & \frac{u(u-1)}{(s+t-1)!} \sumd{\gamma\in\dGamma(\w e00)}\;\;\sumd{y_1\in\b\xi} \spi{a^\gamma}\eps_\gamma\cdot(\b\xi,\eta\ohne x')\cdot(x_\gamma,x,x')\cdot(y_\gamma,y_1)\c (\xi\cup\eta) \\
= & u(u-1)D. \\
\end{array}
\]

{\it Subsubcase $(\w\Xi,\xi^-,\eta^-) = \smateckzz{+}{\ds}{\ds}{\ds}$.} By symmetry, we obtain from subsubcase $\smateckzz{\ds}{\ds}{+}{\ds}$ 
\[
\Lambda([a],\w e,\pazz{+}{+}{+}{+},\{ 2\}, \{ 1,2\}) =  (s+1)sD.\\
\]

{\it Subsubcase $(\w\Xi,\xi^-,\eta^-) = \smateckzz{+}{\ds}{\ds}{+}$.} By symmetry, we obtain from subsubcase $\smateckzz{\ds}{+}{+}{\ds}$ 
\[
\Lambda([a],\w e,\pazz{+}{+}{+}{+},\{ 2\}, \{ 1\}) = -suD. \\
\]

{\it Subsubcase $(\w\Xi,\xi^-,\eta^-) = \smateckzz{+}{+}{\ds}{\ds}$.} By symmetry, we obtain from subsubcase $\smateckzz{\ds}{\ds}{+}{+}$ 
\[
\Lambda([a],\w e,\pazz{+}{+}{+}{+},\{ 2\}, \{ 2\}) = -suD. \\
\]

{\it Subsubcase $(\w\Xi,\xi^-,\eta^-) = \smateckzz{+}{+}{\ds}{+}$.} By symmetry, we obtain from subsubcase $\smateckzz{\ds}{+}{+}{+}$ 
\[
\Lambda([a],\w e,\pazz{+}{+}{+}{+},\{ 2\}, \{\}) = u(u-1)D. \\
\]

{\it Subsubcases $\smateckzz{+}{\ds}{+}{\ds}$, $\smateckzz{+}{\ds}{+}{+}$, $\smateckzz{+}{+}{+}{\ds}$, $\smateckzz{+}{+}{+}{+}$.} The $\Lambda$-value is zero by the Garnir relations.

Summing up $\Lambda$-values over subcases and subsubcases, we obtain 
\[
G'_{[a],\xi,\eta} f'_{\w e22} = 2(s - u)A + (s - u)(s - u - 1)B + 2(s - u)(s - u + 1)D
\]

{\bf Case $f'_{\w e12}$.}

{\it Subcase $\w\Xi = \smateckzz{+}{+}{-}{+}$.}

{\it Subsubcase $(\w\Xi,\xi^-,\eta^-) = \smateckzz{\ds}{\ds}{-}{\ds}$.} We calculate
\[
\begin{array}{cl}
  & \Lambda([a],\w e,\pazz{+}{+}{-}{+},\{ 1\}, \{ 1,2\}) \\
= & \frac{st(t-1)}{s!t!} \sumd{\gamma\in\dGamma(\w e00)}\spi{a^\gamma}\eps_\gamma\cdot(x_\gamma,x,x')\cdot(y_\gamma,y')\c (\xi\cup\eta) \\
\auf{\mbox{\scr(\ref{PropG6})}}{=} & \frac{st(t-1)\cdot s!}{s!t!(s+t-2)!} \sumd{\gamma\in\dGamma(\w e00)}\spi{a^\gamma}\eps_\gamma
                          \cdot (x_\gamma,x,x')\cdot(y_\gamma,y')\cdot(\b\xi\ohne z,\eta\ohne\{ x',y'\})\c (\xi\cup\eta)\c\b\xi \\
= & - \frac{s}{(t-2)!(s+t-2)!} \sumd{\gamma\in\dGamma(\w e00)}\spi{a^\gamma}\eps_\gamma
    \cdot(x_\gamma,x,x')\cdot(y_\gamma,y,y')\cdot(\b\xi\ohne z,\eta\ohne\{ x',y'\})\c (\xi\cup\eta)\c\b\xi \\
= & -sB. \\
\end{array}
\]

{\it Subsubcase $(\w\Xi,\xi^-,\eta^-) = \smateckzz{\ds}{\ds}{-}{+}$.} We calculate
\[
\begin{array}{cl}
  & \Lambda([a],\w e,\pazz{+}{+}{-}{+},\{ 1\}, \{ 1\}) \\
= & \frac{st}{s!t!}  \sumd{\gamma\in\dGamma(\w e00)}\;\;\sumd{y'_1\in\b\eta}\spi{a^\gamma}\eps_\gamma\cdot(x_\gamma,x,x')\cdot(y_\gamma,y'_1)\c (\xi\cup\eta) \\
\auf{\mbox{\scr(\ref{PropG6})}}{=} & \frac{st\cdot s!(t-1)!}{s!t!(s+t-1)!} \sumd{\gamma\in\dGamma(\w e00)}\;\;\sumd{y'_1\in\b\eta}\spi{a^\gamma}\eps_\gamma
                          \cdot(x_\gamma,x,x')\cdot(y_\gamma,y_1')\cdot(\b\xi,\eta\ohne x')\c (\xi\cup\eta) \\
= & sH. \\
\end{array}
\]

{\it Subsubcase $(\w\Xi,\xi^-,\eta^-) = \smateckzz{\ds}{+}{-}{\ds}$.} We calculate
\[
\begin{array}{cl}
  & \Lambda([a],\w e,\pazz{+}{+}{-}{+},\{ 1\}, \{ 2\}) \\
= & \frac{st}{s!t!} \sumd{\gamma\in\dGamma(\w e00)}\;\;\sumd{x'_1\in\b\eta}\spi{a^\gamma}\eps_\gamma\cdot(x_\gamma,x,x'_1)\cdot(y_\gamma,y')\c (\xi\cup\eta) \\
\auf{\mbox{\scr(\ref{PropG6})}}{=} & \frac{st\cdot (s-1)!}{s!t!(s+t-2)!}\sumd{\gamma\in\dGamma(\w e00)}\;\;\sumd{x'_1\in\b\eta}\spi{a^\gamma}\eps_\gamma\cdot 
                          (x_\gamma,x,x'_1)\cdot(y_\gamma,y')\cdot(\b\xi\ohne z,\eta\ohne\{ x',y'\})\cdot(x'_1,x')\c (\xi\cup\eta)\c(\b\xi\cup x'_1)\\
= & \frac{1}{(t-1)!(s+t-2)!}\sumd{\gamma\in\dGamma(\w e00)}\;\;\sumd{x'_1\in\b\eta}\spi{a^\gamma}\eps_\gamma\cdot \\
  & \cdot(x_\gamma,x,x'_1)\cdot(y_\gamma,y,y')\cdot(\b\xi\ohne z,\eta\ohne\{ x',y'\})\cdot(x'_1,x')\cdot(x'_1,x_\gamma)\c (\xi\cup\eta)\c(\b\xi\cup x'_1) \\
= & \frac{1}{(t-1)!(s+t-2)!}\sumd{\gamma\in\dGamma(\w e00)}\;\;\sumd{x'_1\in\b\eta}\spi{a^\gamma}\eps_\gamma\cdot
      (x_\gamma,x,x')\cdot(y_\gamma,y,y')\cdot(\b\xi\ohne z,\eta\ohne\{ x',y'\})\c (\xi\cup\eta)\c(\b\xi\cup x'_1) \\
= & uB - A. \\
\end{array}
\]

{\it Subsubcase $(\w\Xi,\xi^-,\eta^-) = \smateckzz{\ds}{+}{-}{+}$.} We calculate
\[
\begin{array}{cl}
  & \Lambda([a],\w e,\pazz{+}{+}{-}{+},\{ 1\}, \{ \}) \\
= & \frac{s}{s!t!} \sumd{\gamma\in\dGamma(\w e00)}\;\;\sumd{x'_1\neq y'_1\in\b\eta}\spi{a^\gamma}\eps_\gamma\cdot(x_\gamma,x,x'_1)\cdot(y_\gamma,y'_1)\c (\xi\cup\eta) \\
\auf{\mbox{\scr(\ref{PropG6})}}{=} & \frac{s\cdot (s-1)!t!}{s!t!(s+t-1)!} \sumd{\gamma\in\dGamma(\w e00)}\;\;\sumd{x'_1\neq y'_1\in\b\eta}\spi{a^\gamma}\eps_\gamma\cdot
                         (x_\gamma,x,x'_1)\cdot(y_\gamma,y_1')\cdot(\b\xi,\eta\ohne x')\cdot(x'_1,x')\c (\xi\cup\eta) \\
= & - \frac{1}{(s+t-1)!} \sumd{\gamma\in\dGamma(\w e00)}\;\;\sumd{x'_1\neq y'_1\in\b\eta}\spi{a^\gamma}\eps_\gamma\cdot
        (x_\gamma,x,x'_1)\cdot(y_\gamma,y_1')\cdot(\b\xi,\eta\ohne x')\cdot(x'_1,x')\cdot(x'_1,x_\gamma)\c (\xi\cup\eta) \\
= & - \frac{u-1}{(s+t-1)!} \sumd{\gamma\in\dGamma(\w e00)}\;\;\sumd{y'_1\in\b\eta}\spi{a^\gamma}\eps_\gamma
    \cdot(x_\gamma,x,x')\cdot(y_\gamma,y_1')\cdot(\b\xi,\eta\ohne x')\c (\xi\cup\eta) \\
= & -(u-1)H.\\
\end{array}
\]

{\it Subsubcase $(\w\Xi,\xi^-,\eta^-) = \smateckzz{+}{\ds}{-}{\ds}$.} We calculate
\[
\begin{array}{cl}
  & \Lambda([a],\w e,\pazz{+}{+}{-}{+},\{ \}, \{ 1,2\}) \\
= & \frac{t(t-1)}{s!t!} \sumd{\gamma\in\dGamma(\w e00)}\;\;\sumd{x_1\in\b\xi}\spi{a^\gamma}\eps_\gamma\cdot(x_\gamma,x_1,x')\cdot(y_\gamma,y')\c (\xi\cup\eta) \\
\auf{\mbox{\scr(\ref{PropG6})}}{=} & \frac{t(t-1)\cdot (s+1)!(t-2)!}{s!t!(s+t-1)!} \sumd{\gamma\in\dGamma(\w e00)}\;\;\sumd{x_1\in\b\xi}\spi{a^\gamma}\eps_\gamma
                          \cdot(x_\gamma,x_1,x')\cdot(y_\gamma,y')\cdot(\b\xi\ohne x_1,\eta\ohne\{ x',y'\})\c (\xi\cup\eta) \\
= & - \frac{s+1}{(s+t-1)!} \sumd{\gamma\in\dGamma(\w e00)}\;\;\sumd{x_1\in\b\xi}\spi{a^\gamma}\eps_\gamma\cdot
       (\b\xi,\eta\ohne y')\cdot(x_\gamma,x_1)\cdot(y_\gamma,y,y')\c (\xi\cup\eta) \\
= & -(s+1)D.\\
\end{array}
\]

{\it Subsubcase $(\w\Xi,\xi^-,\eta^-) = \smateckzz{+}{+}{-}{\ds}$.} We calculate
\[
\begin{array}{cl}
  & \Lambda([a],\w e,\pazz{+}{+}{-}{+},\{ \}, \{ 2\}) \\
= & \frac{t}{s!t!} \sumd{\gamma\in\dGamma(\w e00)}\;\;\sumd{x_1\in\b\xi,\; x'_1\in\b\eta}\spi{a^\gamma}\eps_\gamma\cdot(x_\gamma,x_1,x'_1)\cdot(y_\gamma,y')\c (\xi\cup\eta) \\
\auf{\mbox{\scr(\ref{PropG6})}}{=} & \frac{t\cdot s!(t-1)!}{s!t!(s+t-1)!} \sumd{\gamma\in\dGamma(\w e00)}\;\;\sumd{x_1\in\b\xi,\; x'_1\in\b\eta}\spi{a^\gamma}\eps_\gamma\cdot
                        (x_\gamma,x_1,x'_1)\cdot(y_\gamma,y')\cdot(\b\xi\ohne x_1,\eta\ohne\{ x',y'\})\cdot (x'_1,x')\c (\xi\cup\eta) \\
= & \frac{1}{(s+t-1)!} \sumd{\gamma\in\dGamma(\w e00)}\;\;\sumd{x_1\in\b\xi,\; x'_1\in\b\eta}\spi{a^\gamma}\eps_\gamma\cdot
        (\b\xi,\eta\ohne y')\cdot(x_\gamma,x_1,x',x'_1)\cdot(y_\gamma,y,y')\cdot(x'_1,x')\cdot(x'_1,x_\gamma)\c (\xi\cup\eta) \\
= & \frac{u}{(s+t-1)!} \sumd{\gamma\in\dGamma(\w e00)}\;\;\sumd{x_1\in\b\xi}\spi{a^\gamma}\eps_\gamma
    \cdot(\b\xi,\eta\ohne y')\cdot(x_\gamma,x_1)\cdot(y_\gamma,y,y')\c (\xi\cup\eta) \\
= & uD. \\ 
\end{array}
\]

{\it Subsubcases $\smateckzz{+}{\ds}{-}{+}$, $\smateckzz{+}{+}{-}{+}$.} The $\Lambda$-value is zero by the Garnir relations.

{\it Subcase $\w\Xi = \smateckzz{-}{+}{+}{+}$.} Obtained, by symmetry, from subcase $\w\Xi = \smateckzz{+}{+}{-}{+}$.

Summing up $\Lambda$-values over subcases and subsubcases, we obtain 
\[
G'_{[a],\xi,\eta} f'_{\w e12} =  - 2A - 2(s - u)B - 2(s - u + 1)D + 2(s - u + 1)H.
\]

{\bf Case $f'_{\w e02}$.} 

{\it Subcase $\w\Xi = \smateckzz{-}{+}{-}{+}$.}

{\it Subsubcase $(\w\Xi,\xi^-,\eta^-) = \smateckzz{-}{\ds}{-}{\ds}$.} We calculate
\[
\begin{array}{cl}
  & \Lambda([a],\w e,\pazz{-}{+}{-}{+},\{ \}, \{ 1,2\}) \\
= & \frac{t(t-1)}{s!t!}\sumd{\gamma\in\dGamma(\w e00)}\spi{a^\gamma}\eps_\gamma\cdot(x_\gamma,x')\cdot(y_\gamma,y')\c (\xi\cup\eta) \\
\auf{\mbox{\scr(\ref{PropG6})}}{=} & \frac{t(t-1)\cdot s!}{s!t!(s+t-2)!} \sumd{\gamma\in\dGamma(\w e00)}\spi{a^\gamma}\eps_\gamma
                          \cdot(x_\gamma,x')\cdot(y_\gamma,y')\cdot(\b\xi\ohne z,\eta\ohne\{ x',y'\})\c (\xi\cup\eta)\c\b\xi \\
= & \frac{1}{(t-2)!(s+t-2)!} \sumd{\gamma\in\dGamma(\w e00)}\spi{a^\gamma}\eps_\gamma
    \cdot(x_\gamma,x,x')\cdot(y_\gamma,y,y')\cdot(\b\xi\ohne z,\eta\ohne\{ x',y'\})\c (\xi\cup\eta)\c\b\xi \\
= & B. \\
\end{array}
\]

{\it Subsubcase $(\w\Xi,\xi^-,\eta^-) = \smateckzz{-}{\ds}{-}{+}$.} We calculate
\[
\begin{array}{cl}
  & \Lambda([a],\w e,\pazz{-}{+}{-}{+},\{ \}, \{ 1\}) \\
= & \frac{t}{s!t!}\sumd{\gamma\in\dGamma(\w e00)}\;\;\sumd{y'_1\in\b\eta}\spi{a^\gamma}\eps_\gamma\cdot(x_\gamma,x')\cdot(y_\gamma,y'_1)\c (\xi\cup\eta) \\
\auf{\mbox{\scr(\ref{PropG6})}}{=} & \frac{t\cdot s!(t-1)!}{s!t!(s+t-1)!} \sumd{\gamma\in\dGamma(\w e00)}\;\;\sumd{y'_1\in\b\eta}\spi{a^\gamma}\eps_\gamma
                          \cdot(x_\gamma,x')\cdot(y_\gamma,y'_1)\cdot(\b\xi,\eta\ohne x')\c (\xi\cup\eta) \\
= & - \frac{1}{(s+t-1)!} \sumd{\gamma\in\dGamma(\w e00)}\;\;\sumd{y'_1\in\b\eta}\spi{a^\gamma}\eps_\gamma
    \cdot(x_\gamma,x,x')\cdot(y_\gamma,y'_1)\cdot(\b\xi,\eta\ohne x')\c (\xi\cup\eta) \\
= & - H.\\
\end{array}
\]

{\it Subsubcase $(\w\Xi,\xi^-,\eta^-) = \smateckzz{-}{+}{-}{\ds}$.} By symmetry, we obtain from subsubcase $\smateckzz{-}{\ds}{-}{+}$ 
\[
\Lambda([a],\w e,\pazz{-}{+}{-}{+},\{ \}, \{ 2\}) = - H.
\]

{\it Subsubcase $(\w\Xi,\xi^-,\eta^-) = \smateckzz{-}{+}{-}{+}$.} The $\Lambda$-value is zero by the Garnir relations.  

Summing up $\Lambda$-values over subcases and subsubcases, we obtain
\[
G'_{[a],\xi,\eta} f'_{\w e02} = B - 2H.
\]

{\bf Case $f'_{\w e21}$.}

{\it Subcase $\w\Xi = \smateckzz{+}{+}{+}{-}$.}

{\it Subsubcase $(\w\Xi,\xi^-,\eta^-) = \smateckzz{\ds}{\ds}{\ds}{-}$.} We calculate
\[
\begin{array}{cl}
  & \Lambda([a],\w e,\pazz{+}{+}{+}{-},\{ 1,2\}, \{ 1\}) \\
= & \frac{s(s-1)t}{s!t!}\sumd{\gamma\in\dGamma(\w e00)}\spi{a^\gamma}\eps_\gamma\cdot(x_\gamma,x,x')\cdot(y_\gamma,y)\c (\xi\cup\eta) \\
\auf{\mbox{\scr(\ref{PropG6})}}{=} & \frac{s(s-1)t\cdot (s-1)!}{s!t!(s+t-2)!} \sumd{\gamma\in\dGamma(\w e00)}\spi{a^\gamma}\eps_\gamma\cdot
                   (x_\gamma,x,x')\cdot(y_\gamma,y)\cdot(\b\xi\ohne z,\eta\ohne \{ x',y'\})\cdot (y_\gamma,y')\c (\xi\cup\eta)\c (\b\xi\cup y_\gamma) \\
= & \frac{(s-1)}{(t-1)!(s+t-2)!} \sumd{\gamma\in\dGamma(\w e00)}\spi{a^\gamma}\eps_\gamma
    \cdot(x_\gamma,x,x')\cdot(y_\gamma,y,y')\cdot(\b\xi\ohne z,\eta\ohne \{ x',y'\})\c (\xi\cup\eta)\c (\b\xi\cup y_\gamma) \\
= & (s-1)(B - C_2). \\
\end{array}
\]

{\it Subsubcase $(\w\Xi,\xi^-,\eta^-) = \smateckzz{\ds}{+}{\ds}{-}$.} We calculate
\[
\begin{array}{cl}
  & \Lambda([a],\w e,\pazz{+}{+}{+}{-},\{ 1,2\}, \{\}) \\
= & \frac{s(s-1)}{s!t!}\sumd{\gamma\in\dGamma(\w e00)}\;\;\sumd{x'_1\in\b\eta}\spi{a^\gamma}\eps_\gamma\cdot(x_\gamma,x,x'_1)\cdot(y_\gamma,y)\c (\xi\cup\eta) \rule{6cm}{0cm}\\
\end{array}
\]
\[
\begin{array}{cl}
\auf{\mbox{\scr(\ref{PropG6})}}{=} & \frac{s(s-1)\cdot (s-2)!}{s!t!(s+t-2)!}\sumd{\gamma\in\dGamma(\w e00)}\;\;\sumd{x'_1\in\b\eta}\spi{a^\gamma}\eps_\gamma\cdot \\
                        & \cdot(x_\gamma,x,x'_1)\cdot(y_\gamma,y)\cdot(\b\xi\ohne z,\eta\ohne \{ x',y'\})\cdot(x'_1,x')\cdot(y_\gamma,y')\c (\xi\cup\eta)\c (\b\xi\cup x'_1\cup y_\gamma) \\
= & - \frac{1}{t!(s+t-2)!}\sumd{\gamma\in\dGamma(\w e00)}\;\;\sumd{x'_1\in\b\eta}\spi{a^\gamma}\eps_\gamma\cdot(x_\gamma,x,x'_1)\cdot(y_\gamma,y)\cdot \\
  & \cdot(\b\xi\ohne z,\eta\ohne \{ x',y'\})\cdot(x'_1,x')\cdot(y_\gamma,y')\cdot(x'_1,x_\gamma)\c (\xi\cup\eta)\c (\b\xi\cup x'_1\cup y_\gamma) \\
= & - \frac{1}{t!(s+t-2)!}\sumd{\gamma\in\dGamma(\w e00)}\;\;\sumd{x'_1\in\b\eta}\spi{a^\gamma}\eps_\gamma\cdot
         (x_\gamma,x,x')\cdot(y_\gamma,y,y')\cdot(\b\xi\ohne z,\eta\ohne \{ x',y'\})\c (\xi\cup\eta)\c (\b\xi\cup x'_1\cup y_\gamma) \\
= & -(uB - uC_2 - A).\\
\end{array}
\]

{\it Subsubcase $(\w\Xi,\xi^-,\eta^-) = \smateckzz{\ds}{\ds}{+}{-}$.} We calculate
\[
\begin{array}{cl}
  & \Lambda([a],\w e,\pazz{+}{+}{+}{-},\{ 1\}, \{ 1\}) \\
= & \frac{st}{s!t!}\sumd{\gamma\in\dGamma(\w e00)}\;\;\sumd{y_1\in\b\xi}\spi{a^\gamma}\eps_\gamma\cdot(x_\gamma,x,x')\cdot(y_\gamma,y_1)\c (\xi\cup\eta) \\
\auf{\mbox{\scr(\ref{PropG6})}}{=} & \frac{st\cdot s!(t-1)!}{s!t!(s+t-1)!}\sumd{\gamma\in\dGamma(\w e00)}\;\;\sumd{y_1\in\b\xi}\spi{a^\gamma}\eps_\gamma
                          \cdot(x_\gamma,x,x')\cdot(y_\gamma,y_1)\cdot(\b\xi\ohne y_1,\eta\ohne \{ x',y'\})\cdot(y_\gamma,y')\c (\xi\cup\eta)\\
= & \frac{s}{(s+t-1)!}\sumd{\gamma\in\dGamma(\w e00)}\;\;\sumd{y_1\in\b\xi}\spi{a^\gamma}\eps_\gamma
    \cdot(\b\xi,\eta\ohne x') \cdot(x_\gamma,x,x')\cdot(y_\gamma,y_1)\c (\xi\cup\eta)\\
= & sD. \\
\end{array}
\]

{\it Subsubcase $(\w\Xi,\xi^-,\eta^-) = \smateckzz{\ds}{+}{+}{-}$.} We calculate
\[
\begin{array}{cl}
  & \Lambda([a],\w e,\pazz{+}{+}{+}{-},\{ 1\}, \{ \}) \\
= & \frac{s}{s!t!}\sumd{\gamma\in\dGamma(\w e00)}\;\;\sumd{y_1\in\b\xi,\; x'_1\in\b\eta}\spi{a^\gamma}\eps_\gamma\cdot(x_\gamma,x,x'_1)\cdot(y_\gamma,y_1)\c (\xi\cup\eta) \\
\auf{\mbox{\scr(\ref{PropG6})}}{=} & \frac{s\cdot (s-1)!t!}{s!t!(s+t-1)!}\sumd{\gamma\in\dGamma(\w e00)}\;\;\sumd{y_1\in\b\xi,\; x'_1\in\b\eta}\spi{a^\gamma}\eps_\gamma\cdot\\
                 & \cdot(x_\gamma,x,x'_1)\cdot(y_\gamma,y_1)\cdot(\b\xi\ohne y_1,\eta\ohne \{ x',y'\})\cdot(y_\gamma,y')\cdot(x'_1,x')\c (\xi\cup\eta)\\
= & - \frac{1}{(s+t-1)!}\sumd{\gamma\in\dGamma(\w e00)}\;\;\sumd{y_1\in\b\xi,\; x'_1\in\b\eta}\spi{a^\gamma}\eps_\gamma\cdot
         (\b\xi,\eta\ohne x')\cdot(x_\gamma,x,x'_1)\cdot(y_\gamma,y_1)\cdot(x'_1,x')\cdot(x_\gamma,x'_1)\c (\xi\cup\eta)\\
= & - \frac{u}{(s+t-1)!}\sumd{\gamma\in\dGamma(\w e00)}\;\;\sumd{y_1\in\b\xi}\spi{a^\gamma}\eps_\gamma
    \cdot(\b\xi,\eta\ohne x')\cdot(x_\gamma,x,x')\cdot(y_\gamma,y_1)\c (\xi\cup\eta)\\
= & - uD. \\
\end{array}
\]

{\it Subsubcase $(\w\Xi,\xi^-,\eta^-) = \smateckzz{+}{\ds}{\ds}{-}$.} We calculate
\[
\begin{array}{cl}
  & \Lambda([a],\w e,\pazz{+}{+}{+}{-},\{ 2\}, \{ 1\}) \\
= & \frac{st}{s!t!}\sumd{\gamma\in\dGamma(\w e00)}\;\;\sumd{x_1\in\b\xi}\spi{a^\gamma}\eps_\gamma\cdot(x_\gamma,x_1,x')\cdot(y_\gamma,y)\c (\xi\cup\eta) \\
\auf{\mbox{\scr(\ref{PropG6})}}{=} & \frac{st\cdot s!(t-1)!}{s!t!(s+t-1)!}\sumd{\gamma\in\dGamma(\w e00)}\;\;\sumd{x_1\in\b\xi}\spi{a^\gamma}\eps_\gamma
                          \cdot(x_\gamma,x_1,x')\cdot(y_\gamma,y)\cdot(\b\xi\ohne x_1,\eta\ohne\{ x',y'\})\cdot(y_\gamma,y')\c (\xi\cup\eta)\\
= & \frac{s}{(s+t-1)!}\sumd{\gamma\in\dGamma(\w e00)}\;\;\sumd{x_1\in\b\xi}\spi{a^\gamma}\eps_\gamma
    \cdot(\b\xi,\eta\ohne y')\cdot(x_\gamma,x_1)\cdot(y_\gamma,y,y')\c (\xi\cup\eta)\\
= & sD. \\
\end{array}
\]

{\it Subsubcase $(\w\Xi,\xi^-,\eta^-) = \smateckzz{+}{+}{\ds}{-}$.} We calculate
\[
\begin{array}{cl}
  & \Lambda([a],\w e,\pazz{+}{+}{+}{-},\{ 2\}, \{ \}) \\
= & \frac{s}{s!t!}\sumd{\gamma\in\dGamma(\w e00)}\;\;\sumd{x_1\in\b\xi,\; x'_1\in\b\eta}\spi{a^\gamma}\eps_\gamma\cdot(x_\gamma,x_1,x'_1)\cdot(y_\gamma,y)\c (\xi\cup\eta) \\
\auf{\mbox{\scr(\ref{PropG6})}}{=} & \frac{s\cdot (s-1)!t!}{s!t!(s+t-1)!}\sumd{\gamma\in\dGamma(\w e00)}\;\;\sumd{x_1\in\b\xi,\; x'_1\in\b\eta}\spi{a^\gamma}\eps_\gamma\cdot\\
                      & \cdot(x_\gamma,x_1,x'_1)\cdot(y_\gamma,y)\cdot(\b\xi\ohne x_1,\eta\ohne\{ x',y'\})\cdot(y_\gamma,y')\cdot(x'_1,x')\c (\xi\cup\eta)\\
= & - \frac{1}{(s+t-1)!}\sumd{\gamma\in\dGamma(\w e00)}\;\;\sumd{x_1\in\b\xi,\; x'_1\in\b\eta}\spi{a^\gamma}\eps_\gamma\cdot
        (\b\xi,\eta\ohne y')\cdot(x_\gamma,x_1,x',x'_1)\cdot(y_\gamma,y,y')\cdot(x'_1,x')\cdot(x_\gamma,x'_1)\c (\xi\cup\eta)\\
= & - \frac{u}{(s+t-1)!}\sumd{\gamma\in\dGamma(\w e00)}\;\;\sumd{x_1\in\b\xi}\spi{a^\gamma}\eps_\gamma
    \cdot(\b\xi,\eta\ohne y')\cdot(x_\gamma,x_1)\cdot(y_\gamma,y,y')\c (\xi\cup\eta)\\
= & - uD.\\
\end{array}
\]

{\it Subsubcases $\smateckzz{+}{\ds}{+}{-}$, $\smateckzz{+}{+}{+}{-}$.} The $\Lambda$-value is zero by the Garnir relations.

{\it Subcase $\w\Xi = \smateckzz{+}{-}{+}{+}$.} Obtained, by symmetry, from subcase $\w\Xi = \smateckzz{+}{+}{+}{-}$.

Summing up $\Lambda$-values over subcases and subsubcases, we obtain 
\[
G'_{[a],\xi,\eta} f'_{\w e21} = 2A + 2(s - u - 1)B - (s - u - 1)(C_1 + C_2) + 4(s - u)D.
\]

{\bf Case $f'_{\w e11}$.}

{\it Subcase $\w\Xi = \smateckzz{+}{+}{-}{-}$.}

{\it Subsubcase $(\w\Xi,\xi^-,\eta^-) = \smateckzz{\ds}{\ds}{-}{-}$.} We calculate
\[
\begin{array}{cl}
  & \Lambda([a],\w e,\pazz{+}{+}{-}{-},\{ 1\}, \{ 1\}) \\
= & \frac{st}{s!t!}\sumd{\gamma\in\dGamma(\w e00)}\spi{a^\gamma}\eps_\gamma\cdot(x_\gamma,x,x')\c (\xi\cup\eta) \\
\auf{\mbox{\scr(\ref{PropG6})}}{=} & \frac{st\cdot s!(t-1)!}{s!t!(s+t-1)!}\sumd{\gamma\in\dGamma(\w e00)}\spi{a^\gamma}\eps_\gamma
                          \cdot(x_\gamma,x,x')\cdot(\b\xi,\eta\ohne x')\c (\xi\cup\eta)\\
= & s F_1. \\
\end{array}
\]

{\it Subsubcase $(\w\Xi,\xi^-,\eta^-) = \smateckzz{\ds}{+}{-}{-}$.} We calculate
\[
\begin{array}{cl}
  & \Lambda([a],\w e,\pazz{+}{+}{-}{-},\{ 1\}, \{ \}) \\
= & \frac{s}{s!t!}\sumd{\gamma\in\dGamma(\w e00)}\;\;\sumd{x'_1\in\b\eta}\spi{a^\gamma}\eps_\gamma\cdot(x_\gamma,x,x'_1)\c (\xi\cup\eta) \\
\auf{\mbox{\scr(\ref{PropG6})}}{=} & \frac{s\cdot (s-1)!t!}{s!t!(s+t-1)!}\sumd{\gamma\in\dGamma(\w e00)}\;\;\sumd{x'_1\in\b\eta}\spi{a^\gamma}\eps_\gamma
                          \cdot(x_\gamma,x,x'_1)\cdot(\b\xi,\eta\ohne x')\cdot(x'_1,x')\c (\xi\cup\eta)\\
= & - \frac{1}{(s+t-1)!}\sumd{\gamma\in\dGamma(\w e00)}\;\;\sumd{x'_1\in\b\eta}\spi{a^\gamma}\eps_\gamma
    \cdot(x_\gamma,x,x'_1)\cdot(\b\xi,\eta\ohne x')\cdot(x'_1,x')\cdot(x'_1,x_\gamma)\c (\xi\cup\eta)\\
= & - \frac{u}{(s+t-1)!}\sumd{\gamma\in\dGamma(\w e00)}\spi{a^\gamma}\eps_\gamma
    \cdot(\b\xi,\eta\ohne x')\cdot(x_\gamma,x,x')\c (\xi\cup\eta)\\
= & - uF_1. \\
\end{array}
\]

{\it Subsubcases $\smateckzz{+}{\ds}{-}{-}$, $\smateckzz{+}{+}{-}{-}$.} The $\Lambda$-value is zero by the Garnir relations.

{\it Subcase $\w\Xi = \smateckzz{-}{-}{+}{+}$.} Obtained, by symmetry, from subcase $\w\Xi = \smateckzz{+}{+}{-}{-}$.

{\it Subcase $\w\Xi = \smateckzz{+}{-}{-}{+}$.}

{\it Subsubcase $(\w\Xi,\xi^-,\eta^-) = \smateckzz{\ds}{-}{-}{\ds}$.} We calculate
\[
\begin{array}{cl}
  & \Lambda([a],\w e,\pazz{+}{-}{-}{+},\{ 1\}, \{ 2\}) \\
= & \frac{st}{s!t!}\sumd{\gamma\in\dGamma(\w e00)}\spi{a^\gamma}\eps_\gamma\cdot(x_\gamma,x)\cdot(y_\gamma,y')\c (\xi\cup\eta) \\
\auf{\mbox{\scr(\ref{PropG6})}}{=} & \frac{st\cdot (s-1)!}{s!t!(s+t-2)!}\sumd{\gamma\in\dGamma(\w e00)}\spi{a^\gamma}\eps_\gamma\cdot
                     (x_\gamma,x)\cdot(y_\gamma,y')\cdot(\b\xi\ohne z,\eta\ohne \{ x',y'\})\cdot(x_\gamma,x')\c (\xi\cup\eta)\c(\b\xi\cup x_\gamma)\\
= & - \frac{1}{(t-1)!(s+t-2)!}\sumd{\gamma\in\dGamma(\w e00)}\spi{a^\gamma}\eps_\gamma\cdot
           (x_\gamma,x,x')\cdot(y_\gamma,y,y')\cdot(\b\xi\ohne z,\eta\ohne \{ x',y'\})\c (\xi\cup\eta)\c(\b\xi\cup x_\gamma)\\
= & -(B - C_1). \\
\end{array}
\]

{\it Subsubcase $(\w\Xi,\xi^-,\eta^-) = \smateckzz{+}{-}{-}{\ds}$.} We calculate
\[
\begin{array}{cl}
  & \Lambda([a],\w e,\pazz{+}{-}{-}{+},\{ \}, \{ 2\}) \\
= & \frac{t}{s!t!}\sumd{\gamma\in\dGamma(\w e00)}\;\;\sumd{x_1\in\b\xi}\spi{a^\gamma}\eps_\gamma\cdot(x_\gamma,x_1)\cdot(y_\gamma,y')\c (\xi\cup\eta) \\
\auf{\mbox{\scr(\ref{PropG6})}}{=} & \frac{t\cdot s!(t-1)!}{s!t!(s+t-1)!}\sumd{\gamma\in\dGamma(\w e00)}\;\;\sumd{x_1\in\b\xi}\spi{a^\gamma}\eps_\gamma
                          \cdot(x_\gamma,x_1)\cdot(y_\gamma,y')\cdot(\b\xi\ohne x_1,\eta\ohne \{ x',y'\})\cdot(x_\gamma,x')\c (\xi\cup\eta)\\
= & - \frac{1}{(s+t-1)!}\sumd{\gamma\in\dGamma(\w e00)}\;\;\sumd{x_1\in\b\xi}\spi{a^\gamma}\eps_\gamma
    \cdot(\b\xi,\eta\ohne y')\cdot(x_\gamma,x_1)\cdot(y_\gamma,y,y')\c (\xi\cup\eta)\\
= & - D.\\
\end{array}
\]

{\it Subsubcase $(\w\Xi,\xi^-,\eta^-) = \smateckzz{\ds}{-}{-}{+}$.} We calculate
\[
\begin{array}{cl}
  & \Lambda([a],\w e,\pazz{+}{-}{-}{+},\{ 1\}, \{ \}) \\
= & \frac{s}{s!t!}\sumd{\gamma\in\dGamma(\w e00)}\;\;\sumd{y'_1\in\b\eta}\spi{a^\gamma}\eps_\gamma\cdot(x_\gamma,x)\cdot(y_\gamma,y'_1)\c (\xi\cup\eta) \\
\auf{\mbox{\scr(\ref{PropG6})}}{=} & \frac{s\cdot (s-1)!t!}{s!t!(s+t-1)!}\sumd{\gamma\in\dGamma(\w e00)}\;\;\sumd{y'_1\in\b\eta}\spi{a^\gamma}\eps_\gamma
                          \cdot(x_\gamma,x)\cdot(y_\gamma,y'_1)\cdot(\b\xi,\eta\ohne x')\cdot(x_\gamma,x')\c (\xi\cup\eta)\\
= & \frac{1}{(s+t-1)!}\sumd{\gamma\in\dGamma(\w e00)}\;\;\sumd{y'_1\in\b\eta}\spi{a^\gamma}\eps_\gamma
    \cdot(x_\gamma,x,x')\cdot(y_\gamma,y'_1)\cdot(\b\xi,\eta\ohne x')\c (\xi\cup\eta)\\
= & H. \\
\end{array}
\]

{\it Subsubcase $(\w\Xi,\xi^-,\eta^-) = \smateckzz{+}{-}{-}{+}$.} The $\Lambda$-value is zero by the Garnir relations.

{\it Subcase $\w\Xi = \smateckzz{-}{+}{+}{-}$.} Obtained, by symmetry, from subcase $\w\Xi = \smateckzz{+}{-}{-}{+}$.

Summing up $\Lambda$-values over subcases and subsubcases, we obtain 
\[
G'_{[a],\xi,\eta} f'_{\w e11} = - 2B + (C_1 + C_2) - 2D + 2H + (s - u)(F_1 + F_2).
\]

{\bf Case $f'_{\w e01}$.}

{\it Subcase $\w\Xi = \smateckzz{-}{+}{-}{-}$.}

{\it Subsubcase $(\w\Xi,\xi^-,\eta^-) = \smateckzz{-}{\ds}{-}{-}$.} We calculate
\[
\begin{array}{cl}
  & \Lambda([a],\w e,\pazz{-}{+}{-}{-},\{ \}, \{ 1\}) \\
= & \frac{t}{s!t!}\sumd{\gamma\in\dGamma(\w e00)}\spi{a^\gamma}\eps_\gamma\cdot(x_\gamma,x')\c (\xi\cup\eta) \\
\end{array}
\]
\[
\begin{array}{cl}
\auf{\mbox{\scr(\ref{PropG6})}}{=} & \frac{t\cdot s!(t-1)!}{s!t!(s+t-1)!}\sumd{\gamma\in\dGamma(\w e00)}\spi{a^\gamma}\eps_\gamma
                          \cdot(x_\gamma,x')\cdot(\b\xi,\eta\ohne x')\c (\xi\cup\eta)\\
= & - \frac{1}{(s+t-1)!}\sumd{\gamma\in\dGamma(\w e00)}\spi{a^\gamma}\eps_\gamma
    \cdot(x_\gamma,x,x')\cdot(\b\xi,\eta\ohne x')\c (\xi\cup\eta)\\
= & - F_1. \\
\end{array}
\]

{\it Subsubcase $(\w\Xi,\xi^-,\eta^-) = \smateckzz{-}{+}{-}{-}$.} The $\Lambda$-value is zero by the Garnir relations.

{\it Subcase $\w\Xi = \smateckzz{-}{-}{-}{+}$.} Obtained, by symmetry, from subcase $\w\Xi = \smateckzz{-}{+}{-}{-}$.

Summing up $\Lambda$-values over subcases and subsubcases, we obtain 
\[
G'_{[a],\xi,\eta} f'_{\w e01} = -(F_1 + F_2).
\]

{\bf Case $f'_{\w e20}$.} 

{\it Subcase $\w\Xi = \smateckzz{+}{-}{+}{-}$.}

{\it Subsubcase $(\w\Xi,\xi^-,\eta^-) = \smateckzz{\ds}{-}{\ds}{-}$.} We calculate
\[
\begin{array}{cl}
  & \Lambda([a],\w e,\pazz{+}{-}{+}{-},\{ 1,2\}, \{ \}) \\
= & \frac{s(s-1)}{s!t!}\sumd{\gamma\in\dGamma(\w e00)}\spi{a^\gamma}\eps_\gamma\cdot(x_\gamma,x)\cdot(y_\gamma,y)\c (\xi\cup\eta) \\
\auf{\mbox{\scr(\ref{PropG6})}}{=} & \frac{s(s-1)(s-2)!}{s!t!(s+t-2)!}\sumd{\gamma\in\dGamma(\w e00)}\spi{a^\gamma}\eps_\gamma\cdot\\
                        & \cdot(x_\gamma,x)\cdot(y_\gamma,y)\cdot(\b\xi\ohne z,\eta\ohne\{ x',y'\})\cdot(x_\gamma,x')\cdot(y_\gamma,y')\c (\xi\cup\eta)\c (\b\xi\cup x_\gamma\cup y_\gamma)\\
= & \frac{1}{t!(s+t-2)!}\sumd{\gamma\in\dGamma(\w e00)}\spi{a^\gamma}\eps_\gamma\cdot
     (x_\gamma,x,x')\cdot(y_\gamma,y,y')\cdot(\b\xi\ohne z,\eta\ohne\{ x',y'\})\c (\xi\cup\eta)\c (\b\xi\cup x_\gamma\cup y_\gamma)\\
= & B - C_1 - C_2.\\
\end{array}
\]

{\it Subsubcase $(\w\Xi,\xi^-,\eta^-) = \smateckzz{\ds}{-}{+}{-}$.} We calculate
\[
\begin{array}{cl}
  & \Lambda([a],\w e,\pazz{+}{-}{+}{-},\{ 1\}, \{ \}) \\
= & \frac{s}{s!t!}\sumd{\gamma\in\dGamma(\w e00)}\;\;\sumd{y_1\in\b\xi}\spi{a^\gamma}\eps_\gamma\cdot(x_\gamma,x)\cdot(y_\gamma,y_1)\c (\xi\cup\eta) \\
\auf{\mbox{\scr(\ref{PropG6})}}{=} & \frac{s(s-1)!t!}{s!t!(s+t-1)!}\sumd{\gamma\in\dGamma(\w e00)}\;\;\sumd{y_1\in\b\xi}\spi{a^\gamma}\eps_\gamma\cdot
                     (x_\gamma,x)\cdot(y_\gamma,y_1)\cdot(\b\xi\ohne y_1,\eta\ohne\{ x',y'\})\cdot(x_\gamma,x')\cdot(y_\gamma,y')\c (\xi\cup\eta)\\
= & \frac{1}{(s+t-1)!}\sumd{\gamma\in\dGamma(\w e00)}\;\;\sumd{y_1\in\b\xi}\spi{a^\gamma}\eps_\gamma
    \cdot(\b\xi,\eta\ohne x')\cdot(x_\gamma,x,x')\cdot(y_\gamma,y_1)\c (\xi\cup\eta)\\
= & D.\\
\end{array}
\]

{\it Subsubcase $(\w\Xi,\xi^-,\eta^-) = \smateckzz{+}{-}{\ds}{-}$.} By symmetry, we obtain from subsubcase $\smateckzz{\ds}{-}{+}{-}$ 
\[
\Lambda([a],\w e,\pazz{+}{-}{+}{-},\{ 2\}, \{ \}) = D.
\]

{\it Subsubcase $(\w\Xi,\xi^-,\eta^-) = \smateckzz{+}{-}{+}{-}$.} The $\Lambda$-value is zero by the Garnir relations.

Summing up $\Lambda$-values over subcases and subsubcases, we obtain 
\[
G'_{[a],\xi,\eta} f'_{\w e20} = B - (C_1 + C_2) + 2D.
\]

{\bf Case $f'_{\w e10}$.}

{\it Subcase $\w\Xi = \smateckzz{+}{-}{-}{-}$.}

{\it Subsubcase $(\w\Xi,\xi^-,\eta^-) = \smateckzz{\ds}{-}{-}{-}$.} We calculate
\[
\begin{array}{cl}
  & \Lambda([a],\w e,\pazz{+}{-}{-}{-},\{ 1\}, \{ \}) \\
= & \frac{s}{s!t!}\sumd{\gamma\in\dGamma(\w e00)}\spi{a^\gamma}\eps_\gamma\cdot(x_\gamma,x)\c (\xi\cup\eta) \\
\auf{\mbox{\scr(\ref{PropG6})}}{=} & \frac{s(s-1)!t!}{s!t!(s+t-1)!}\sumd{\gamma\in\dGamma(\w e00)}\spi{a^\gamma}\eps_\gamma
                          \cdot(x_\gamma,x)\cdot(\b\xi,\eta\ohne x')\cdot(x_\gamma,x')\c (\xi\cup\eta)\\
= & \frac{1}{(s+t-1)!}\sumd{\gamma\in\dGamma(\w e00)}\spi{a^\gamma}\eps_\gamma
    \cdot(x_\gamma,x,x')\cdot(\b\xi,\eta\ohne x')\c (\xi\cup\eta)\\
= & F_1.\\
\end{array}
\]

{\it Subsubcase $(\w\Xi,\xi^-,\eta^-) = \smateckzz{+}{-}{-}{-}$.} The $\Lambda$-value is zero by the Garnir relations.

{\it Subcase $\w\Xi = \smateckzz{-}{-}{+}{-}$.} Obtained, by symmetry, from subcase $\w\Xi = \smateckzz{+}{-}{-}{-}$.

Summing up $\Lambda$-values over subcases and subsubcases, we obtain 
\[
G'_{[a],\xi,\eta} f'_{\w e10} = (F_1 + F_2).
\]

{\bf Case $f'_{\w e00}$.} The Garnir relations yield 
\[
G'_{[a],\xi,\eta} f'_{\w e00} = 0.
\]

We note that $s - u = X_p - X_{p+1}$ and evaluate the linear combination
\[
\begin{array}{l}
\sumd{\alpha,\beta\in [0,2]} X_p^{(2-\alpha)} X_{p+1}^{(2-\beta)} G'_{[a],\xi,\eta} f'_{\w e\alpha\beta} \\
\begin{array}{clll}
= & 1            & \cdot 1                    & \cdot\Big( 2(X_p - X_{p+1})A + (X_p - X_{p+1})(X_p - X_{p+1} - 1)B \\
  &              &                            & \;\; + 2(X_p - X_{p+1})(X_p - X_{p+1} + 1)D\Big) \\
+ & X_p          & \cdot 1                    & \cdot\Big( - 2A - 2(X_p - X_{p+1})B - 2(X_p - X_{p+1} + 1)D \\
  &              &                            & \;\; + 2(X_p - X_{p+1} + 1)H\Big) \\
+ & X_p(X_p + 1) & \cdot 1                    & \cdot\Big( B - 2H\Big) \\
+ & 1            & \cdot X_{p+1}              & \cdot\Big( 2A + 2(X_p - X_{p+1} - 1)B  \\
  &              &                            & \;\; - (X_p - X_{p+1} - 1)(C_1 + C_2) + 4(X_p - X_{p+1})D\Big) \\
+ & X_p          & \cdot X_{p+1}              & \cdot\Big( - 2B + (C_1 + C_2) - 2D + 2H \\
  &              &                            & \;\; + (X_p - X_{p+1})(F_1 + F_2)\Big)\\
+ & X_p(X_p + 1) & \cdot X_{p+1}              & \cdot\Big( -(F_1 + F_2)\Big)\\
+ & 1            & \cdot X_{p+1}(X_{p+1} + 1) & \cdot\Big( B - (C_1 + C_2) + 2D\Big) \\
+ & X_p          & \cdot X_{p+1}(X_{p+1} + 1) & \cdot\Big( (F_1 + F_2)\Big)\\
+ & X_p(X_p + 1) & \cdot X_{p+1}(X_{p+1} + 1) & \cdot\Big( 0 \Big) \\
= & 0.\\
\end{array} \\
\end{array}
\]
\end{Lemma}

\pagebreak[4]

\begin{Lemma}[full version of (\ref{LemM8short})]
\label{LemM8}
Suppose $g = p < k$ and $s,t\geq 2$. There exist elements $x', y'\in\eta$, $x'\neq y'$, $z\in\b\xi$, which we choose and fix. For $\gamma\in\dGamma(\w e 0)$, we let 
$x_\gamma := a_{j,\b\gamma(j,1)}$, where $j\in [g+2,k+1]$ is minimal with $1\ti j\in\Xi_\gamma$, and $y_\gamma := a_{j,\b\gamma(j,2)}$, where $j\in [g+2,k+1]$ is minimal with 
$2\ti j\in\Xi_\gamma$. I.e.\ we pick the entries $x_\gamma, y_\gamma$ that `cross the column' $g+1$ under the operation of $\gamma$. The set of maps
\[
\begin{array}{rcl}
[1,\lambda_1]\ohne\{ g+1\} & \lraa{\w e} & [0,2] \\
j                          & \lra        & \w e_j \\ 
\end{array}
\]
that send $g$ and $k+1$ to $\w e_g = \w e_{k+1} = 2$, and that map $j\in [1,\lambda_1]\ohne [g,k+1]$ to $e_j = 0$, is denoted by $\w E$. For $\w e\in\w E$ and $\beta\in [0,2]$, we denote by 
$\w e\beta$ the prolongation of $\w e$ to $[g,k+1]$ by $(\w e\beta)_{g+1} = \beta$. For $\gamma\in\dGamma(\w e0)$, we write 
\[
U_\gamma := \spi{a^\gamma}\eps_\gamma\cdot(\b\xi\ohne z,\eta\ohne\{ x',y'\})\cdot(x_\gamma,x')\cdot(y_\gamma,y') \\
\]
and let
\[
\begin{array}{lcl}
A_{[a],\xi,\eta,\w e}   & := & \sumd{\gamma\in\dGamma(\w e0)} U_\gamma\cdot\sumd{w'\in\b\eta} (w',z) \\
B_{[a],\xi,\eta,\w e}   & := & \sumd{\gamma\in\dGamma(\w e0)} U_\gamma\cdot(1 - \sumd{w\in\b\xi\ohne z} (w,z)) \\
C_{1,[a],\xi,\eta,\w e} & := & \sumd{\gamma\in\dGamma(\w e0)} U_\gamma\cdot (z,x_\gamma) \\
C_{2,[a],\xi,\eta,\w e} & := & \sumd{\gamma\in\dGamma(\w e0)} U_\gamma\cdot (z,y_\gamma). \\
\end{array}
\]
We obtain
\[
\begin{array}{r}
G'_{[a],\xi,\eta} f' = (X_g + 2)\sumd{\w e\in\w E}\Big(\prodd{j\in [g+2,k]} X_j^{(2-\w e_j)}\Big) \Big( 2 A_{[a],\xi,\eta,\w e} + (X_g + 1)B_{[a],\xi,\eta,\w e} \\
- X_{g+1}(C_{1,[a],\xi,\eta,\w e} + C_{2,[a],\xi,\eta,\w e})\Big).\\
\end{array}
\]

\rm
Having fixed a map $\w e\in\w E$, we need to evaluate the expression
\[
\sumd{\beta\in [0,2]} X_{g+1}^{(2-\beta)} G'_{[a],\xi,\eta} f'_{\w e\beta}.
\]
Given $\beta\in [0,2]$, given a subset $\w\Xi\tm [1,2]\ti \{ g+1\}$ of cardinality $\#\w\Xi = \beta$ and given a subset $\eta^-\tm [1,2]$ such that
\[
(\ast) \hspace*{3cm}
\eta^- \ti (g+1) \tm\w\Xi, 
\hspace*{3cm}
\]
we let
\[
\dGamma(\w e,\w\Xi,\eta^-) := \Big\{ \gamma\in\dGamma(\w e\beta) \; \Big|\; \Xi_\gamma \cap ([1,2]\ti\{ g+1\}) = \w\Xi,\; \gamma^{-1}(\eta) = \eta^-\ti\{ g+1\}\Big\} 
\]
and form the partial sum
\[
\Lambda([a],\w e,\w\Xi,\eta^-) := \frac{1}{s!t!} \sumd{\gamma\in\dGamma(\w e,\w\Xi,\eta^-)} \spi{a^\gamma}\eps_\gamma\c (\xi\cup\eta)
\]
so that we can recover
\[
G'_{[a],\xi,\eta} f'_{\w e\beta} = \sumd{\w\Xi\tm [1,2]\ti\{ g+1\}}\;\; \sumd{\eta^- \tm [1,2] \mbox{\scr\ subject to } (\ast)} \Lambda([a],\w e,\w\Xi,\eta^-).
\]

{\bf Calculation of $\Lambda$-values.} To distinguish subcases and subsubcases, we adapt the according notation of (\ref{LemM2}).

{\bf Case $f'_{\w e2}$.}

{\it Subcase $\w\Xi = \smateckze{+}{+}$.}

{\it Subsubcase $(\w\Xi,\eta^-) = \smateckze{\ds}{\ds}$.} We calculate
\[
\begin{array}{cl}
  & \Lambda([a],\w e,\paze{+}{+},\{ 1,2\}) \\
= & \frac{t(t-1)}{s!t!} \sumd{\gamma\in\dGamma(\w e0)}\spi{a^\gamma}\eps_\gamma\cdot(x_\gamma,x')\cdot(y_\gamma,y')\c (\xi\cup\eta) \\
\auf{\mbox{\scr(\ref{LemG4})}}{=} & \frac{t(t-1)(s+2)!}{s!t!} \sumd{\gamma\in\dGamma(\w e0)}\spi{a^\gamma}\eps_\gamma
                         \cdot (x_\gamma,x')\cdot(y_\gamma,y')\cdot(\b\xi\ohne z,\eta\ohne\{ x',y'\})\c \b\xi \\
= & (s+2)(s+1)B_{[a],\xi,\eta,\w e}. \\
\end{array}
\]

{\it Subsubcase $(\w\Xi,\eta^-) = \smateckze{\ds}{+}$.} We calculate
\[
\begin{array}{cl}
  & \Lambda([a],\w e,\paze{+}{+},\{ 1\}) \\
= & \frac{t}{s!t!} \sumd{\gamma\in\dGamma(\w e0)}\;\;\sumd{y'_1\in\b\eta}\spi{a^\gamma}\eps_\gamma\cdot(x_\gamma,x')\cdot(y_\gamma,y'_1)\c (\xi\cup\eta) \\
\auf{\mbox{\scr(\ref{LemG4})}}{=} & \frac{t(s+1)!}{s!t!} \sumd{\gamma\in\dGamma(\w e0)}\;\;\sumd{y'_1\in\b\eta}\spi{a^\gamma}\eps_\gamma\cdot
                           (x_\gamma,x')\cdot(y_\gamma,y'_1)\cdot(\b\xi\ohne z,\eta\ohne\{ x',y'\})\cdot(y'_1,y')\c(\b\xi\cup y'_1) \\
= & - \frac{(s+1)}{(t-1)!} \sumd{\gamma\in\dGamma(\w e0)}\;\;\sumd{y'_1\in\b\eta}\spi{a^\gamma}\eps_\gamma\cdot 
         (x_\gamma,x')\cdot(y_\gamma,y'_1)\cdot(\b\xi\ohne z,\eta\ohne\{ x',y'\})\cdot(y'_1,y')\cdot(y'_1,y_\gamma)\c(\b\xi\cup y'_1) \\
= & - \frac{(s+1)}{(t-1)!} \sumd{\gamma\in\dGamma(\w e0)}\;\;\sumd{y'_1\in\b\eta}\spi{a^\gamma}\eps_\gamma
    \cdot (x_\gamma,x')\cdot(\b\xi\ohne z,\eta\ohne\{ x',y'\})\cdot(y_\gamma,y')\c(\b\xi\cup y'_1) \\
= & -(s+1)(uB_{[a],\xi,\eta,\w e} - A_{[a],\xi,\eta,\w e}). \\
\end{array}
\]

{\it Subsubcase $(\w\Xi,\eta^-) = \smateckze{+}{\ds}$.} By symmetry, we obtain from subsubcase $\smateckze{\ds}{+}$
\[
\Lambda([a],\w e,\paze{+}{+},\{ 2\}) = -(s+1)(uB_{[a],\xi,\eta,\w e} - A_{[a],\xi,\eta,\w e}).
\]

{\it Subsubcase $(\w\Xi,\eta^-) = \smateckze{+}{+}$.} We calculate
\[
\begin{array}{cl}
  & \Lambda([a],\w e,\paze{+}{+},\{\}) \\
= & \frac{1}{s!t!} \sumd{\gamma\in\dGamma(\w e0)}\;\;\sumd{x'_1\neq y'_1\in\b\eta}\spi{a^\gamma}\eps_\gamma\cdot(x_\gamma,x'_1)\cdot(y_\gamma,y'_1)\c (\xi\cup\eta) \\
\auf{\mbox{\scr(\ref{LemG4})}}{=} & \frac{s!}{s!t!} \sumd{\gamma\in\dGamma(\w e0)}\;\;\sumd{x'_1\neq y'_1\in\b\eta}\spi{a^\gamma}\eps_\gamma\cdot
                              (x_\gamma,x'_1)\cdot(y_\gamma,y'_1)\cdot(\b\xi\ohne z,\eta\ohne\{ x',y'\})\cdot(x'_1,x')\cdot(y'_1,y')\c(\b\xi\cup x'_1\cup y'_1) \\
= & \frac{1}{t!} \sumd{\gamma\in\dGamma(\w e0)}\;\;\sumd{x'_1\neq y'_1\in\b\eta}\spi{a^\gamma}\eps_\gamma\cdot\\
  & \cdot (x_\gamma,x'_1)\cdot(y_\gamma,y'_1)\cdot(\b\xi\ohne z,\eta\ohne\{ x',y'\})\cdot(x'_1,x')\cdot(y'_1,y')\cdot(x'_1,x_\gamma)\cdot(y'_1,y_\gamma)\c(\b\xi\cup x'_1\cup y'_1) \\
= & \frac{1}{t!} \sumd{\gamma\in\dGamma(\w e0)}\;\;\sumd{x'_1\neq y'_1\in\b\eta}\spi{a^\gamma}\eps_\gamma
    \cdot(\b\xi\ohne z,\eta\ohne\{ x',y'\})\cdot(x_\gamma,x')\cdot(y_\gamma,y')\c(\b\xi\cup x'_1\cup y'_1) \\
= & (u-1)uB_{[a],\xi,\eta,\w e} - 2(u-1)A_{[a],\xi,\eta,\w e}.\\
\end{array}
\]

{\bf Case $f'_{\w e1}$.}

{\it Subcase $\w\Xi = \smateckze{+}{-}$.}

{\it Subsubcase $(\w\Xi,\eta^-) = \smateckze{\ds}{-}$.} We calculate
\[
\begin{array}{cl}
  & \Lambda([a],\w e,\paze{+}{-},\{ 1\}) \\
= & \frac{t}{s!t!} \sumd{\gamma\in\dGamma(\w e0)}\spi{a^\gamma}\eps_\gamma\cdot(x_\gamma,x')\c (\xi\cup\eta) \\
\auf{\mbox{\scr(\ref{LemG4})}}{=} & \frac{t(s+1)!}{s!t!} \sumd{\gamma\in\dGamma(\w e0)}\spi{a^\gamma}\eps_\gamma
                         \cdot (x_\gamma,x')\cdot(\b\xi\ohne z,\eta\ohne\{ x',y'\})\cdot(y_\gamma,y')\c(\b\xi\cup y_\gamma) \\
= & (s+1)(B_{[a],\xi,\eta,\w e} - C_{2,[a],\xi,\eta,\w e}).\\
\end{array}
\]

{\it Subsubcase $(\w\Xi,\eta^-) = \smateckze{+}{-}$.} We calculate
\[
\begin{array}{cl}
  & \Lambda([a],\w e,\paze{+}{-},\{ \}) \\
= & \frac{1}{s!t!} \sumd{\gamma\in\dGamma(\w e0)}\;\;\sumd{x'_1\in\b\eta}\spi{a^\gamma}\eps_\gamma\cdot(x_\gamma,x'_1)\c (\xi\cup\eta) \\
\auf{\mbox{\scr(\ref{LemG4})}}{=} & \frac{s!}{s!t!} \sumd{\gamma\in\dGamma(\w e0)}\;\;\sumd{x'_1\in\b\eta}\spi{a^\gamma}\eps_\gamma
                         \cdot (x_\gamma,x'_1)\cdot(\b\xi\ohne z,\eta\ohne\{ x',y'\})\cdot(x'_1,x')\cdot(y_\gamma,y')\c(\b\xi\cup x'_1\cup y_\gamma) \\
= & - \frac{1}{t!} \sumd{\gamma\in\dGamma(\w e0)}\;\;\sumd{x'_1\in\b\eta}\spi{a^\gamma}\eps_\gamma\cdot
      (x_\gamma,x'_1)\cdot(\b\xi\ohne z,\eta\ohne\{ x',y'\})\cdot(x'_1,x')\cdot(y_\gamma,y')\cdot(x_\gamma,x'_1)\c(\b\xi\cup x'_1\cup y_\gamma) \\
= & - \frac{1}{t!} \sumd{\gamma\in\dGamma(\w e0)}\;\;\sumd{x'_1\in\b\eta}\spi{a^\gamma}\eps_\gamma
    \cdot(\b\xi\ohne z,\eta\ohne\{ x',y'\})\cdot(x_\gamma,x')\cdot(y_\gamma,y')\c(\b\xi\cup x'_1\cup y_\gamma) \\
= & -(uB_{[a],\xi,\eta,\w e} - uC_{2,[a],\xi,\eta,\w e} - A_{[a],\xi,\eta,\w e}). \\
\end{array}
\]

{\it Subcase $\w\Xi = \smateckze{-}{+}$.} Obtained, by symmetry, from subcase $\w\Xi = \smateckze{+}{-}$.

{\bf Case $f'_{\w e0}$.}

{\it Subcase $\w\Xi = \smateckze{-}{-}$.}

{\it Subsubcase $(\w\Xi,\eta^-) = \smateckze{-}{-}$.} We calculate
\[
\begin{array}{cl}
  & \Lambda([a],\w e,\paze{-}{-},\{ \}) \\
= & \frac{1}{s!t!} \sumd{\gamma\in\dGamma(\w e0)}\spi{a^\gamma}\eps_\gamma\c (\xi\cup\eta) \\
\auf{\mbox{\scr(\ref{LemG4})}}{=} & \frac{s!}{s!t!} \sumd{\gamma\in\dGamma(\w e0)}\spi{a^\gamma}\eps_\gamma
                         \cdot(\b\xi\ohne z,\eta\ohne\{ x',y'\})\cdot(x_\gamma,x')\cdot(y_\gamma,y')\c(\b\xi\cup x_\gamma\cup y_\gamma) \\
= & B_{[a],\xi,\eta,\w e} - C_{1,[a],\xi,\eta,\w e} - C_{2,[a],\xi,\eta,\w e}. \\
\end{array}
\]
We note that $s - u = X_g - X_{g+1}$ and evaluate the linear combination
\[
\begin{array}{l}
\sumd{\beta\in [0,2]} X_{g+1}^{(2 - \beta)} G'_{[a],\xi,\eta} f'_{\w e\beta} \\
\begin{array}{cll}
= & 1                    & \cdot\Big( 2(X_g - X_{g+1} + 2)A_{[a],\xi,\eta,\w e} \\
  &                      & + (X_g - X_{g+1} + 2)(X_g - X_{g+1} + 1)B_{[a],\xi,\eta,\w e} \Big)  \\
+ & X_{g+1}              & \cdot\Big( 2A_{[a],\xi,\eta,\w e} + 2(X_g - X_{g+1} + 1)B_{[a],\xi,\eta,\w e} \\
  &                      & - (X_g - X_{g+1} + 1)(C_{1,[a],\xi,\eta,\w e} + C_{2,[a],\xi,\eta,\w e}) \Big) \\
+ & X_{g+1}(X_{g+1} + 1) & \cdot\Big( B_{[a],\xi,\eta,\w e} - (C_{1,[a],\xi,\eta,\w e} + C_{2,[a],\xi,\eta,\w e}) \Big)  \\
\end{array} \\
= 2(X_g + 2)A_{[a],\xi,\eta,\w e} + (X_g + 2)(X_g + 1)B_{[a],\xi,\eta,\w e} - X_{g+1}(X_g + 2)(C_{1,[a],\xi,\eta,\w e} + C_{2,[a],\xi,\eta,\w e}). \\
\end{array}
\]
\end{Lemma}
\end{footnotesize}      
\section{References}
\label{SecRef}

\begin{footnotesize}
{\sc Carter, R.W., Lusztig, G.}\upl
\begin{itemize}
\item[{[CL 74]}] {\it On the Modular Representations of the General Linear and Symmetric Groups}, Math.\ Z.\ 136, p.\ 139-242, 1974.
\end{itemize}
   
{\sc Carter, R.W., Payne, M.T.J.}\upl
\begin{itemize}
\item[{[CP 80]}] {\it On homomorphisms between Weyl modules and Specht modules}, Math.\ Proc.\ Camb.\ Phil.\ Soc.\ 87, p.\ 419-425, 1980. 
\end{itemize}

{\sc Garnir, H.}\upl
\begin{itemize}
\item[{[G 50]}] {\it Th\'eorie de la r\'epresentation lin\'eaire des groupes sym\'etriques}, Th\`ese, M\'em. Soc.\ Roy.\ Sc.\ Li\`ege, (4), 10, 1950. 
\end{itemize}

{\sc James, G.D.}\upl
\begin{itemize}
\item[{[J 78]}] {\it The Representation Theory of the Symmetric Groups}, SLN 682, 1978. 
\end{itemize}

{\sc K\"unzer, M.}\upl
\begin{itemize}
\item[{[K 99]}] {\it Ties for the $\Z\Sl_n$,} thesis, http://www.mathematik.uni-bielefeld.de/$\scm\sim$kuenzer, Bielefeld, 1999.
\end{itemize}

{\sc Specht, W.}\upl
\begin{itemize}
\item[{[Sp 35]}] {\it Die irreduziblen Darstellungen der Symmetrischen Gruppe,} Math.\ Z.\ 39, p.\ 696-711, 1935.
\end{itemize}

\vspace*{8cm}

\begin{flushright}
Invariant address:

Matthias K\"unzer\\
Fakult\"at f\"ur Mathematik\\
Universit\"at Bielefeld\\
Postfach 100131\\
D-33501 Bielefeld\\
kuenzer@mathematik.uni-bielefeld.de\\
\end{flushright}
\end{footnotesize}
\end{document}